\documentclass{svmono}

\usepackage{makeidx}     
\usepackage{graphicx}    
\usepackage{multicol}    

\makeindex             
\usepackage{amscd}
\usepackage{amsmath}
\usepackage{amssymb}
\usepackage{latexsym}
\usepackage{theorem}

\usepackage[all]{xy}

\usepackage{mykomma}
\usepackage{eepic}

\newenvironment{pf}
{\noindent {\it Proof.}}
{\hfill $\Box$}

\begin{document}

\author{Markus Rosellen}
\title{A Course in Vertex Algebra}
\maketitle

\frontmatter

\preface

This book presents vertex algebra from an algebraic perspective.
Thus we emphasize algebraic concepts and methods.
Definitions and constructions that come from physics
play a less important role and are always motivated algebraically.

My intention was to write a book very much like other textbooks in algebra,
in particular textbooks about noncommutative rings and Lie algebras.
My hope is that through this book more mathematicians will
appreciate, use, and study vertex algebra.

A more detailed introduction is given in chapter \ref{C:intro}.
References to the literature are collected in
Appendix \ref{C:notes}.

\bigskip

{\bf Acknowledgements.}\:
First of all, 
I would like to thank Sergei Merkulov for the opportunity to work
at his department and thus making it possible for me to write this book.

Since this book grew out of seminars in Bonn, Mainz, and Stockholm,
I like to thank the organizers and the participants for their
contributions.
In Bonn were
Vicente Cort\'es,
Daniel Huybrechts,
Werner Nahm,
Holger Eberle, 
Daniel Roggenkamp, 
Stephan Stolz,
Peter Teichner,
Bertrand Toen, and
Katrin Wendland.
In Mainz were
Theo de Jong,
Manfred Lehn,
Duco van Straten,
Christian van Enckevort, and
Christian Sevenheck.
And in Stockholm were
Julius Borcea,
Torsten Ekedahl,
Clas L\"ofwall, and
Sergei Merkulov.

I am also very grateful to
Toshiyuki Abe,
Bojko Bakalov,
Keith Hubbard,
Haisheng Li,
Yuri Manin,
Atsushi Matsuo,
Kiyokazu Nagatomo,
Dmitri Orlov,
Markus Spitzweck, and
Hiroshi Yamauchi
for discussions and support.

\newpage

{\bf Terminology.}\:
Our terminology concerning superalgebra is explained in 
appendix \ref{C:supalg}.
Here we only mention the following conventions.

We always work with super objects without making
this explicit in our terminology,
e.g. we work with super vector spaces
and just call them vector spaces.
Morphisms of vector spaces need not preserve the supergrading; 
in other words, they need not be even.

An {\bf even} vector space is 
\index{even vector space}
a vector space $E=E\even\oplus E\odd$ with $E\odd=0$.
An {\bf algebra} is 
\index{algebra}
a vector space $V$ together with 
an even linear map $V\otimes V\to V, a\otimes b\mapsto ab$.
{\bf Commutative algebras} are 
\index{commutative algebra}
always assumed to be associative and unital.

Supersigns are written as powers of 
$$
\zeta
\; :=\; 
-1
$$
in order to distinguish them from other signs. 
In the exponent of $\zeta$ we omit to write the parity $\ta$ of 
an element $a$,
instead of $\zeta^{\ta\tb}$ we just write $\zeta^{ab}$.

We sometimes omit supersigns to make the text more readable.
Supersigns are usually omitted in proofs.
If we do not want to specify a supersign we just write $\pm_s$.

\bigskip

{\bf Standard Notation.}\:
We denote by $\N, \Z, \Q, \R, \C$
the monoid of non-negative integers,
the ring of integers, and 
the fields of rational, real, and complex numbers.  
For $\F=\Z, \Q$, or $\R$, define $\F_{\geq}:=\set{a\in\F\mid a\geq 0}$.
The semigroups $\F_{\leq}, \F_>, \F_<$ are defined in the same way.

Let $a$ be an element of an associative unital algebra and $n\in\K$. 
The 
\index{divided power}
divided powers and the 
\index{binomial coefficient}
binomial coefficients are
$$
a^{(n)}:=
\begin{cases}
a^n/n! \\
\;\; 0 
\end{cases}
\qquad
\binom{a}{n}:=
\begin{cases}
(n!)\inv\prod_{i=0}^{n-1}(a-i) &\quad\text{if $n\in\N$,}\\
\qquad\quad 0 &\quad\text{otherwise}.
\end{cases}
$$
We denote by $e^a=\exp a=\sum a^{(n)}$ the exponential.

For a vector space $E$, 
we denote by $E[\la]$ and $E\pau{\la}$ the spaces of polynomials and 
power series in the variable $\la$ with coefficients in $E$.

We denote by $\fgl(E)$ the general linear Lie algebra.

\bigskip

{\bf Our Notation.}\:
We work over a fixed ground field $\K$ of characteristic $0$.

We write $A\cong B$ if two objects $A$ and $B$ are isomorphic. 
We often write $A=B$ if $A$ and $B$ are canonically isomorphic.

If $S$ is a set with an equivalence relation and $a\in E$ then
$[a]$ denotes the equivalence class of $a$.

If $R$ is a ring then $R\uptimes$ denotes the group of units of $R$.

If $V$ is an algebra with multiplication $a\otimes b\mapsto a\subal b$
and $E, F\subset V$ are subspaces then we define
$E\subal F:=\rspan\set{a\subal b\mid a\in E, b\in F}$.

We sometimes consider $\rho\Z$-gradations.
In this case $\rho\in\Q\uptimes$ such that $\rho\inv\in\N$.

\tableofcontents

\mainmatter

\chapter{Introduction}
\label{C:intro}

\noindent
{\bf The New Approach.}\:
This book presents a new approach to vertex algebras.
The new approach says that a vertex algebra is an associative algebra
such that the underlying Lie algebra is a vertex Lie algebra.
In particular, vertex algebras can be formulated in terms of a single
multiplication and they behave like associative algebras
with respect to it.

This statement is not just an analogy but is based on a theorem
of Bakalov and Kac.
They observed that the data of a vertex algebra are equivalent to
an algebra and a $\la$-bracket and proved that the usual vertex algebra axioms 
are equivalent to five identities for the algebra and the $\la$-bracket.
These five identities correspond to the five identities in the 
following definition of an associative algebra.

We first recall the notion of a {\it pre-Lie algebra}.
This is a nonassociative algebra $V$ such that left multiplication 
$\rho: V\to\End(V), a\mapsto a\cdot$, is an algebra morphism
with respect to the commutator $[a,b]\subast:=ab-ba$ on $V$ and $\End(V)$. 
Any associative algebra is a pre-Lie algebra because an algebra is associative 
iff $\rho(ab)=(\rho a)(\rho b)$ and this implies 
$\rho[a,b]\subast=[\rho a,\rho b]\subast$. 
Moreover, if $V$ is a pre-Lie algebra then $[\, ,]\subast$ is a Lie bracket
because $\fgl(V)$ is a Lie algebra and one can make $\rho$ into a monomorphism
by adjoining an identity $1$ to $V$.

Now an equivalent but redundant definition of an associative algebra is:
An associative algebra is a pre-Lie algebra together with 
a Lie bracket $[\, ,]$ such that $[\, ,]\subast=[\, ,]$ and 
$[a,\; ]$ is a derivation: $[a,bc]=[a,b]c+b[a,c]$.

In the case of vertex algebras, the Lie bracket is replaced by
a Lie $\la$-bracket and the redundancy in the axioms disappears.
This is explained next.

\medskip

\noindent
{\bf Vertex Lie Algebras.}\:
Vertex algebras live in the category of $\K[T]$-modules, 
i.e.~vector spaces endowed with an operator $T$.
Algebras in this category are {\it differential algebras},
i.e.~$\K[T]$-modules with a multiplication such that $T$ is a derivation.

A {\it $\la$-bracket} on a $\K[T]$-module $R$ is a linear map 
$a\otimes b\mapsto [a\subla b]$, from $R\otimes R$ to the space $R[\la]$ of
polynomials with coefficients in $R$, such that 
$$
T[a\subla b]
\; =\;
[Ta\subla b]+[a\subla Tb]
\qquad\text{and}\qquad
[Ta\subla b]
\; =\;
-\la[a\subla b].
$$
A {\it Lie $\la$-bracket} is a $\la$-bracket that satisfies
{\it conformal skew-symmetry} $[a\subla b]=-[b_{-\la-T}a]$
and the {\it conformal Jacobi identity}
$$
[a\subla[b\submu c]]
\; =\;
[[a\subla b]_{\mu+\la}c]
\; +\; 
[b\submu[a\subla c]].
$$
Here we use that $c_T:=\sum T^i(c_i)$ for $c\subla=\sum c_i\la^i\in R[\la]$.
A {\it vertex Lie algebra} is a $\K[T]$-module with a Lie $\la$-bracket.

One can show that any vertex Lie algebra has 
a natural Lie bracket $[\, ,]\subslie$ defined by the formal integral 
$[a,b]\subslie:=\int_{-T}^0 d\la [a\subla b]$.

\medskip

\noindent
{\bf Associative Vertex Algebras.}\:
We define a {\it vertex algebra} to be a differential algebra with 
a $\la$-bracket.
A {\it conformal derivation} of a vertex algebra $V$
is a map $d\subla: V\to V[\la]$ that satisfies 
$d\subla[a\submu b]=[(d\subla a)_{\mu+\la}b]+[a\submu(d\subla b)]$ and 
the {\it Wick formula}:
$$
d\subla(ab)
\; =\;
(d\subla a)b
\; +\;
a(d\subla b)
\; +\;
\int_0\upla d\mu\, [(d\subla a)\submu b].
$$
A vertex algebra is {\it associative} if the underlying algebra is a
pre-Lie algebra, the $\la$-bracket is a Lie $\la$-bracket, 
the commutator $[\, ,]\subast$ is equal to $[\, ,]\subslie$, and 
$[a\subla\;]$ is a conformal derivation for any $a$.

Thus an associative vertex algebra is a differential pre-Lie algebra with
a Lie $\la$-bracket such that $[\, ,]\subast=[\, ,]\subslie$ and 
$[a\subla\;]$ is a conformal derivation.

A vertex algebra is {\it unital} if the underlying algebra is unital,
i.e.~if it has an identity element $1$.
Bakalov and Kac proved the following result.

\bigskip

{\bf Theorem.}\: {\it
The notion of an associative unital vertex algebra is equivalent to
the usual definition of a vertex algebra due to Borcherds.
}

\bigskip

Reversing the historical development, we take the above definition of 
an associative vertex algebra as the starting point of this book.

\medskip

\noindent
{\bf The Original Approach.}\:
The usual formulation of vertex algebras is in terms of a linear map
$$
Y:\;\: V\;\to\;\End(V)\pau{z\uppm},\quad 
a\;\mapsto\; a(z)\; =\;\sum_{t\in\Z}\: a_t\, z^{-t-1}.
$$
The property $(Ta)(z)=\del_z a(z)$ implies that $a_{-1-t}b=(T^{(t)}a)_{-1}b$
for $t\geq 0$ where $T^{(t)}:=T^t/t!$.
The relation between the two formulations is given by $ab=a_{-1}b$ and 
$[a\subla b]=\sum_{t\geq 0}a_t b\,\la^{(t)}$.

The proof of the theorem goes roughly as follows.
The map $Y$ defines a vertex algebra in Borcherds' sense iff
$Y$ satisfies skew-symmetry and the commutator formula.
Skew-symmetry is $a(z)b=e^{zT}b(-z)a$ or equivalently
$$
a_r b
\; =\;
\sum_{i\geq 0}\: (-1)^{r+i+1}\, T^{(i)}(b_{r+i}a).
$$ 
It is easy to see that skew-symmetry for $r\geq 0$ is 
conformal skew-symmetry and for $r=-1$ it is $[\, ,]\subast=[\, ,]\subslie$.
Moreover, skew-symmetry for $r\ne 0$ implies skew-symmetry for $r-1$.
Thus skew-symmetry is equivalent to two identities:
conformal skew-symmetry and $[\, ,]\subast=[\, ,]\subslie$.

Similarly, the commutator formula $[a_t,b_s]=\sum\binom{t}{i}(a_i b)_{t+s-i}$
can be decomposed into three identities that express that the product and
the $\la$-bracket are compatible with themselves and with each other.
For $t, s\geq 0$ it is the conformal Jacobi identity, 
for $t\geq 0, s=-1$ it is the Wick formula, and 
for $t=s=-1$ it is the pre-Lie identity if $[\, ,]\subast=[\, ,]\subslie$.
Conformal skew-symmetry implies that the commutator formula holds for $a, b$ 
and $t, s$ iff it holds for $b, a$ and $s, t$.
Finally, the commutator formula for $t, s\ne 0$ implies 
the commutator formula for $t-1, s$ and $t, s-1$.
Thus the theorem follows.

\medskip

\noindent
{\bf Advantages of the New Approach.}\:
The new formulation of vertex algebras in terms of pre-Lie algebras
makes it possible to study vertex algebras as a class of associative algebras
and to present the subject like in an ordinary algebra textbook. 
The original approach, using infinite families of multiplications,
distributions, the Jacobi identity, and locality,
made it hard to get used to vertex algebras and to apply standard ideas
from algebra to them.

Various notions in algebra appear now more naturally for vertex algebras:
identity element $1$ and idempotent; derivation; 
commuting elements, centre, and commutative vertex algebra; 
vertex Poisson algebra; 
gradation and differential and invariant filtration.
We define conformal derivations of a vertex algebra and prove that
they form an unbounded vertex Lie algebra.

Various results and proofs in algebra now have more natural analogues 
for vertex algebras.
For example, results about generating subspaces, the PBW-property,
and $C_1(V)$ and results about generating subspaces without repeats 
and $C_2(V)$, see section \ref{S:filt span sets}.
These results are proven using filtrations.
Another example is the statement that if an associative algebra is
endowed with a compatible bracket such that $[a,b]=1$ for some $a, b$
then $[\, ,]$ is equal to the commutator.
The vertex analogue of this statement is the following theorem
of Bakalov and Kac: 
If a vertex algebra $V$
satisfies the associativity formula and $a_t b=1$ for some 
$a, b\in V, t\geq 0$ then $V$ is associative,
see sections \ref{SS:assoc alg bracket}--\ref{SS:eq commut}.
The block decomposition into indecomposable subalgebras 
also exists for vertex algebras, 
see section \ref{SS:indecomp va}.

Various constructions of ordinary algebras from vertex algebras 
are now more natural.
An easy example is the commutative subalgebra $V_0$ of an $\N$-graded
vertex algebra that is the starting point for the vertex algebras of 
chiral differential operators. 
Another example is the Gerstenhaber algebra $H_Q(V)$ of a 
topological vertex algebra, see section \ref{SS:topl va}.
The associated graded vertex algebra $\rgr V$ of a filtered vertex algebra $V$
is often a commutative algebra.

We emphasize universal properties and adjoint functors.
In particular, we construct associative vertex algebras in terms of
generators and relations.
We systematically use the convenient formalism of $\la$-products due to Kac.

\medskip

\noindent
{\bf Further Highlights.}\:
In chapter \ref{C:vlie} on vertex Lie algebras we present 
Chambert-Loir's beautiful proof that locality implies the 
weak commutator formula and Roitman's result that, roughly speaking, 
the functors $\fg\to\fg\pau{z\uppm}$ and $R\mapsto\fg(R)$ are adjoint.
This last result implies that vertex Lie algebras form a full subcategory
of the category of local Lie algebras. 
Instead of locality, we primarily use the weak commutator formula to motivate
the notions of a local Lie algebra and a vertex Lie algebra. 

In chapter \ref{C:va} on associative vertex algebras we prove that 
a vertex algebra of fields is associative iff it is local.
More precisely, we show that the pre-Lie identity, 
the conformal Jacobi identity, and the Wick formula are satisfied for any
fields.
Conformal skew-symmetry and $[\, ,]\subast=[\, ,]\subslie$ for fields
are equivalent to locality. 
We present Matsuo and Nagatomo's proof of the Wick formula that
exhibits the symmetry of the Jacobi identity.
We also give a new construction of the $t$-th products of fields
and use it to prove that skew-symmetry for fields is equivalent to locality
and to prove Dong's lemma. 

We prove that associative vertex algebras form a full subcategory
of the category of local associative algebras. 

In chapter \ref{C:resul} on basic results we give a unified
presentation of the thirteen vertex algebra identities and the
implications between them. 
Most results are derived from the fundamental recursion, $\bbS_3$-symmetry,
and the second recursion. 
In particular, we give a new proof of the above theorem of Bakalov and Kac.

In chapter \ref{C:env} on the enveloping vertex algebra $U(R)$,
we explain three constructions of $U(R)$: as a quotient of
the free vertex algebra generated by $R$, as a quotient of the
tensor vertex algebra of $R$, and as a Verma module $V(\fg)$ of the
Borcherds Lie algebra $\fg$ of $R$. 
For an arbitrary local Lie algebra $\fg$, we define the notion of
a vertex algebra over $\fg$ and prove that $V(\fg)$ is the universal
vertex algebra over $\fg$.
We prove the Poincar\'e-Birkhoff-Witt theorem for $U(R)$.
It states that the associated graded vertex algebra $\rgr U(R)$ with respect to
the invariant filtration generated by $R$ is isomorphic 
to the symmetric algebra $S\dual R$.

In chapter \ref{C:reprva} on representation theory
we prove the basic properties of the Zhu algebra $A(V)$ using
an isomorphism $A(V)\to(\hV\subast)_0/K$ where $\hV$ is the affinization of $V$.

\chapter{Vertex Lie Algebras}
\label{C:vlie}

The purpose of this chapter is to discuss examples of vertex Lie algebras and
to explain the relation between vertex Lie algebras and local Lie
algebras.
We emphasize that the relation between Lie brackets and Lie $\la$-brackets
is provided by the weak commutator formula.

In section \ref{S:loclie} 
we show that many important infinite-dimensional Lie algebras 
have the structure of a local Lie algebra.

In sections \ref{S:loc weak comm f}--\ref{S:vlie to loclie}
we show that vertex Lie algebras form a full subcategory of the category
of local Lie algebras.
We first prove that the weak commutator formula is equivalent to locality.
Then we show that the distributions $c^i(z)$ appearing in the 
weak commutator formula define the structure of a vertex Lie algebra on 
the space $R(\fg)$ of distributions of a local Lie algebra $\fg$.
Finally, we construct a fully faithful functor $R\mapsto\fg(R)$ 
from vertex Lie algebras to local Lie algebras such that $R(\fg(R))=R$ and
the $c^i(z)$s of the weak commutator formula for $\fg(R)$ are given by
the $\la$-bracket of $R$.

In section \ref{S:ex vlie}
we show that all the examples of local Lie algebras from 
section \ref{S:loclie} are in fact vertex Lie algebras and 
discuss free vertex Lie algebras and 
vertex Lie algebras of conformal operators and conformal derivations.

Note that our {\bf conventions} regarding terminology and notation
are listed at the end of the preface.

\section{Examples of Local Lie Algebras}
\label{S:loclie}

The most important infinite-dimensional Lie algebras are
the affine and (super) Virasoro Lie algebras.
An important auxiliary role is played by certain infinite-dimensional 
Heisenberg and Clifford Lie algebras.
These Lie algebras $\fg$ have the additional structure of a set 
of families $(a_t)$, where $a_t\in\fg$ and $t\in\Z$, such that 
the $a_t$ span $\fg$ and satisfy the weak commutator formula.
Such Lie algebras are called local Lie algebras.

Two families $(a_t), (b_t)$ satisfy the weak commutator formula 
if there exist families $(d^1_t), \dots, (d^m_t)$
such that $[a_t,b_s]$ is equal to a linear combination of $d^j_{t+s}$ 
with polynomial coefficients in $t, s$.
Note that this condition is particularly natural if $\fg$ is a $\Z$-graded
Lie algebra and $a_t, b_t, d^j_t\in\fg_t$.

In section \ref{SS:loclie}
we define graded local Lie algebras. 
In the remaining sections \ref{SS:loop affi lie}--\ref{SS:super vir} 
we consider examples.

\subsection{Weak Commutator Formula and Local Lie Algebras}
\label{SS:loclie}

We define local Lie algebras in terms of the weak commutator formula.
We also define graded local Lie algebras since 
it makes it easier to describe examples like the Virasoro algebra.

\bigskip

Let $\fg$ be a Lie algebra and $a_t, b_t\in\fg$ for $t\in\Z$. 
We say that the families $(a_t), (b_t)$ satisfy the 
{\bf weak commutator formula} if there exist 
\index{weak commutator formula}
families $(d^j_t)$ in $\fg$ and polynomials $p_j(t,s)\in\K[t,s]$ such that
$$
[a_t,b_s]
\; =\;
\sum_{j=1}^m\:
p_j(t,s)\, d^j_{t+s}.
$$

Everything we do takes place in the category of $\K[T]$-modules.
Here $\K[T]$ is the polynomial ring in one even variable.
Thus a $\K[T]$-module is the same thing as a vector space with an 
even operator $T$.
The category of $\K[T]$-modules is a monoidal category
with $T$ acting on $V\otimes W$ as $T\otimes 1+1\otimes T$.

A {\bf differential algebra} is a $\K[T]$-module $V$ with 
\index{differential algebra}
a morphism $V\otimes V\to V$.
In other words, it is an ordinary algebra with a derivation $T$. 

Let $V$ be a differential algebra.
A {\bf derivation} of $V$ is a $\K[T]$-module
\index{derivation!of a differential algebra}
morphism $d$ such that $d(ab)=(da)b+\zeta^{da}a(db)$.
A {\bf differential} of $V$ is 
\index{differential!of a differential algebra}
an odd derivation $d$ such that $d^2=0$.
Note that $T$ is in general {\it not} a differential.

A {\bf left ideal} of $V$ is
\index{ideal}
a $\K[T]$-submodule $I$ such that $VI\subset I$.
Right ideals and (two-sided) ideals are defined in the same way.
If $I$ is an ideal then $V/I$ is 
\index{quotient}
the {\bf quotient} differential algebra.

\bigskip

{\bf Definition.}\:
A {\bf local Lie algebra} is 
\index{local!Lie algebra}
a differential Lie algebra $\fg$ together with 
a set $F\subfg$ of families $a=(a_t)$ in $\fg$ such that
any $a, b\in F\subfg$ satisfy the weak commutator formula,
$\fg=\rspan\set{a_t\mid a\in F\subfg, t\in\Z}$, and $Ta_t=-ta_{t-1}$.

\bigskip

Morphisms of local Lie algebras are defined in section \ref{SS:vlie}.
Since a local Lie algebra is spanned by the elements $a_t$, 
the operator $T$ is uniquely determined by $Ta_t=-ta_{t-1}$.

A {\bf gradation} of a $\K[T]$-module $V$ is a $\K$-gradation 
$V=\bigoplus V_h$ such that $TV_h\subset V_{h+1}$ for $h\in\K$.
This is equivalent to giving an even diagonalizable operator $H$ 
such that $[H,T]=T$.
The correspondence is given by $H|_{V_h}=h$.
We 
\index{Hamiltonian}
call $H$ the {\bf Hamiltonian} and $h_a:=h$ the {\bf weight} of $a\in V_h$.
A subset $S\subset V$ is {\bf homogeneous} if $S\subset\bigoplus (S\cap V_h)$.
A subspace $E\subset V$ is homogeneous iff $HE\subset E$,
see Proposition \ref{S:supalg}.

A {\bf gradation} of a differential algebra $V$ is 
\index{graded!differential algebra}
a gradation of the $\K[T]$-module such that $V_h\cdot V_k\subset V_{h+k}$.

If $\fg$ is a $\K$-graded Lie algebra
then we denote families in $\fg$ by $(a_{(t)})$ where $t\in\Z$.
Let $(\fg^{\Z})_h$ be the space of families $(a_{(t)})$ such that 
$a_{(t)}\in\fg_{h-t-1}$.
If $(a_{(t)})\in(\fg^{\Z})_h$, define $a_n:=a_{(n+h-1)}$ for $n\in \Z-h$.
Then $a_n\in\fg_{-n}$.

\bigskip

{\bf Definition.}\:
A {\bf graded local Lie algebra} is 
\index{graded!local Lie algebra}
a local Lie algebra $\fg$ together with 
a differential Lie algebra gradation $\fg=\bigoplus\fg_h$ 
such that $F\subfg$ is a homogeneous subset of $\bigoplus(\fg^{\Z})_h$.

\bigskip

We usually specify a graded local Lie algebra by giving a 
differential Lie algebra $\fg$, a set of families $a=(a_n)_{n\in\Z-h_a}$, 
and numbers $h_a$.
Then $F\subfg:=\set{(a_{t+1-h_a})_{t\in\Z}\mid a}$ and 
$\fg_n:=\rspan\set{a_{-n}\mid a}$.
The data $\fg, a, h_a$ define a graded local Lie algebra iff
$\fg=\rspan\set{a_n}$, $Ta_n=-(n+h_a-1)a_{n-1}$, and for any 
$a, b$ there exist $e^j_n\in\fg_{-n}$ and $q_j\in\K[n,m]$ such that
$[a_n,b_m]=\sum q_j(n,m)e^j_{n+m}$.

\subsection{Affine Lie Algebras}
\label{SS:loop affi lie}

We show that loop and affine Lie algebras are graded local Lie algebras.

\bigskip

Let $\fg$ be a Lie algebra.
If $C$ is an even commutative differential algebra then $\fg\otimes C$ is 
a differential Lie algebra with $T=1\otimes T$ and bracket
$$
[a\otimes f,b\otimes g]
\; :=\;
[a,b]\otimes fg.
$$
The functor $\fg\mapsto\fg\otimes C$ is a base change functor 
from Lie algebras over $\K$ to Lie algebras over $C$.

The {\bf loop Lie algebra} of $\fg$ is
\index{loop Lie algebra}
$\tfg:=\fg\otimes\K[x\uppm]$ where $\K[x\uppm]$ is endowed with 
the derivation $T=-\del_x$.
Define $a_n:=a\otimes x^n$ for $a\in\fg$ and $n\in\Z$.
Then $[a_n,b_m]=[a,b]_{n+m}$ and
$Ta_n=-\del_x(a\otimes x^n)=-n a_{n-1}$.
Thus $\tfg$ with the families $(a_n), a\in\fg, h=1$, 
is a graded local Lie algebra.

A bilinear form 
\index{invariant bilinear form}
on $\fg$ is {\bf invariant} if $([a,b],c)=(a,[b,c])$.
For example, if $\fg$ is even then the {\bf Killing form} $(a,b):=\Tr [a,][b,]$
is 
\index{Killing form}
an invariant symmetric bilinear form.
Moreover, if $\K$ is algebraically closed and $\fg$ is simple
then any invariant bilinear form is a multiple of the Killing form.
This follows from Schur's lemma and the fact that a bilinear form
is invariant iff $\fg\to\fg\dual, a\mapsto (a,\;)$, is a $\fg$-module morphism.

Let $\fg$ be a Lie algebra with an invariant symmetric bilinear form.
Then $\ep:\rLa^2\tfg\to\K, a_n\wedge b_m\mapsto n(a,b)\de_{n+m}$
is a Chevalley-Eilenberg 2-cocycle because, omitting supersigns, we have 
\begin{align}
\notag
&\ep([a_n,b_m],c_k)+\ep([c_k,a_n],b_m)+\ep([b_m,c_k],a_n)
\\
\notag
=\;
&\big(
(n+m)([a,b],c)
+
(k+n) ([c,a],b)
+
(m+k) ([b,c],a)
\big)
\de_{n+m+k}
\\
\notag
=\;
&2(n+m+k)\, ([a,b],c)\, \de_{n+m+k}
\; =\;
0.
\end{align}
 
The {\bf affine Lie algebra} $\hfg$ is
\index{affine!Lie algebra}
the central extension of $\tfg$ given by the 2-cocycle $\ep$.
Thus $\hfg=\tfg\oplus\K\hk$ such that $\hk$ is central and
$$
[a_n,b_m]
\; =\;
[a,b]_{n+m}
\; +\;
n(a,b)\de_{n+m}\hk.
$$
The operator $T: a_n+\la\hk\mapsto -n a_{n-1}$ is a derivation of $\hfg$ 
because $T$ is a derivation of $\tfg$ and the projection of 
$[Ta_n,b_m]+[a_n,Tb_m]$ onto $\K\hk$ is equal to
$(-n(n-1)-mn)(a,b)\de_{n+m-1}\hk=0$.
Moreover, $T(\de_n \hk)=0=-(n-1)\de_{n-1}\hk$.
Thus $\hfg$ with the families $(a_n), a\in\fg, h=1$,
and $(\de_n \hk), h=0$, is a graded local Lie algebra.

Let $\K=\C$ and $\fg$ be an even finite-dimensional simple Lie algebra.
By Cartan's criterion the Killing form is non-degenerate.
Hence the non-degenerate invariant symmetric bilinear forms are the non-zero
multiples of the Killing form. 
We state without proof the following result.

\bigskip

{\bf Proposition.}\: {\it
Let $\K=\C$ and $\fg$ be an even finite-dimensional simple Lie algebra.
Then $\hfg$ is the universal central extension of $\tfg$.
\hfill $\square$
}

\subsection{Heisenberg and Clifford Lie Algebras}
\label{SS:heis cliff lie}

We show that Heisenberg and Clifford Lie algebras are graded 
local Lie algebras and explain their relation to Weyl and Clifford algebras.

\bigskip

A {\bf Heisenberg Lie algebra} is a Lie algebra $\fh$ such
\index{Heisenberg!Lie algebra}
that $[\fh,\fh]$ is even and one-dimensional and equal to the centre of $\fh$.

Let $\fh$ be a Heisenberg Lie algebra.
Choose an isomorphism between $[\fh,\fh]$ and $\K$.
Then the bracket of $\fh$ induces 
a non-degenerate skew-symmetric bilinear form on $W:=\fh/[\fh,\fh]$.
In particular, if $\fh=\fh\even$ and $\dim\fh<\infty$ then 
$W$ is a symplectic vector space and $\dim\fh\in 3+2\N$.

Conversely, let $W$ be a vector space with 
a skew-symmetric bilinear form $\om$.
In other words, $W$ is given by a pair of even vector spaces $W\even, W\odd$ 
such that $W\even$ is endowed with a skew-symmetric and $W\odd$ with a 
symmetric bilinear form.
Then $\cH(W):=W\oplus\K\hk$ with $[a,b]:=\om(a,b)\hk$ and $[a,\hk]:=0$
is a Lie algebra.
Moreover, $\cH(W)$ is a Heisenberg Lie algebra iff 
$W\ne 0$ and $\om$ is non-degenerate.
One can view $\cH(W)$ as the central extension of the abelian Lie algebra $W$
given by the $2$-cocycle $\rLa^2 W\to\K, a\wedge b\mapsto\om(a,b)$.

If $W$ is a symplectic vector space and $2n=\dim W<\infty$ then 
the {\bf Weyl algebra} $A_n$ is the enveloping algebra of $\cH(W)$ modulo 
the relation $\hk=1$.
It is isomorphic to the algebra of differential operators of 
$\K[x_1, \dots, x_n]$.

If $B$ is an even vector space with a symmetric bilinear form then
the {\bf Clifford algebra} is 
\index{Clifford algebra}
the associative unital algebra generated by $B$ with relations $a^2=(a,a)$.
These relations are equivalent to $ab+ba=2(a,b)$. 
Thus the Clifford algebra is the enveloping algebra of $\cH(\rPi B)$ 
modulo the relation $\hk=2$.
Here $\rPi$ is the parity-change functor.

\medskip

There are two ways to construct $W$ from a vector space $B$ with 
a symmetric bilinear form in such a way that $\cH(W)$ is a 
local Lie algebra.

The {\bf Heisenberg Lie algebra} 
\index{Heisenberg!Lie algebra}
of $B$ is $\cH(B\otimes\K[x\uppm])$ where
$\om(a_n,b_m):=n(a,b)\de_{n+m}$ for $a_n:=a\otimes x^n$.
This is just the affine Lie algebra $\hB$ where $B$ is viewed 
as an abelian Lie algebra with an invariant symmetric bilinear form.
In particular, it is a graded local Lie algebra.
Since $[a_0,\hB]=0$, the Heisenberg Lie algebra of $B$ is {\it never}
a Heisenberg Lie algebra.
But $\K\hk\oplus\bigoplus_{n\ne 0}B\otimes\K x^n$ is a Heisenberg Lie algebra
if $B$ is non-degenerate.

The {\bf Clifford Lie algebra} is 
\index{Clifford Lie algebra}
$C(B):=\cH(B\otimes\K[x\uppm]\theta)$ where $\theta$ is an odd variable and
$\om(a_n,b_m):=(a,b)\de_{n+m}$ for $a_n:=a\otimes x^{n-1/2}\theta$ and 
$n\in\Z+1/2$.
The operator $T:a_n\mapsto -(n-1/2) a_{n-1}, \hk\mapsto 0$,
is a derivation of $C(B)$ since
$[Ta_n,b_m]+[a_n,Tb_m]=-(n+m-1)(a,b)\de_{n+m-1}\hk=0$.
Thus $C(B)$ with the families $(a_n), a\in V, h=1/2$, and $(\de_n\hk), h=0$, 
is a graded local Lie algebra.

If $W$ is a symplectic vector space then $\widehat{\rPi W}$ is
\index{symplectic fermion Lie algebra}
called a {\bf symplectic fermion Lie algebra} and $C(\rPi W)$ is
\index{bosonic ghost Lie algebra}
called a {\bf bosonic ghost Lie algebra}.

\subsection{Witt Algebra}
\label{SS:witt alg}

We show that the Witt algebra is a simple graded local Lie algebra.

\bigskip

Let $\Der(C)$ denote the space of derivations of a commutative algebra $C$.
This is a $C$-module and a Lie subalgebra of $\fgl(C)$.
If $\del, \del'\in\Der(C)$ such that $[\del,\del']=0$ and $f, g\in C$ then 
$$
[f\del,g\del']
\; =\;
f\del(g)\del'
\; -\; 
\zeta^{(f+\del)(g+\del')}g\del'(f)\del.
$$

The {\bf Witt algebra} $\Witt$ is
\index{Witt algebra} 
the Lie algebra of derivations of $\K[z\uppm]$.
It is a free $\K[z\uppm]$-module of rank $1$ with basis $\del_z$.
Define $\ell_n:=-z^{n+1}\del_z$.
Then $[\ell_n,\ell_m]=(n-m)\ell_{n+m}$. 
In particular, $[\ell_{-1},\ell_m]=-(m+1)\ell_{m-1}$.
Since $[\ell_{-1},\, ]$ is a derivation of $\Witt$, we see
that $\Witt$ with the family $(\ell_n), h=2$, is a graded local Lie algebra.

We now show that $\Witt$ is a simple Lie algebra.
Since $[\ell_0,I]\subset I$, any ideal $I\subset\Witt$ is $\Z$-graded by
Proposition \ref{S:supalg}.
If $\ell_n\in I$ then $\ell_m\in I$ for any $m\ne 2n$ since
$[\ell_n,\ell_k]=(n-k)\ell_{n+k}$.
Thus if $n=0$ then $I=\Witt$.
If $n\ne 0$ then also $I=\Witt$ since 
$\ell_{4n}\in I$ and hence $\ell_{2n}\in I$.

\subsection{Virasoro Algebra}
\label{SS:vir alg}

We show that the universal central extension of the Witt algebra 
is a graded local Lie algebra.

\bigskip

It is well-known that a Lie algebra $\fg$ has a universal central extension
\index{perfect}
iff $\fg$ is {\bf perfect}: $\fg=[\fg,\fg]$.
Any simple Lie algebra is perfect, in particular $\Witt$ is perfect.
We define the {\bf Virasoro algebra} $\Vir$ to be 
\index{Virasoro algebra}
the universal central extension of $\Witt$.

\bigskip

{\bf Lemma.}\: {\it
Let $\fh$ be a trivial $\Witt$-module.

\smallskip

\iti\:
Any Chevalley-Eilenberg 2-cocycle $\ep:\rLa^2\,\Witt\to\fh$ is 
cohomologous to a cocycle of the form 
$\ell_n\wedge\ell_m\mapsto c_n\de_{n+m}$ for some $c_n\in\fh$.

\smallskip

\itii\:
Let $\ep: \ell_n\wedge\ell_m\mapsto c_n\de_{n+m}$ be a 2-cochain.
Then $\ep$ is a cocycle iff $c_n=n^3a+nb$ for some $a, b\in\fh$.
Moreover, $\ep$ is a coboundary iff $a=0$.
}

\bigskip

\begin{pf}
\iti\:
Define $c_{n,m}:=\ep(\ell_n\wedge\ell_m)$.
The cocycle condition is  
$\ep([\ell_n,\ell_m],\ell_k)+\ep([\ell_k,\ell_n],\ell_m)+
\ep([\ell_m,\ell_k],\ell_n)=0$.
Since $c_{n,m}=-c_{m,n}$ and $[\ell_0,\ell_n]=-n\ell_n$ 
the cocycle condition for $n=0$ yields $(m-k)c_{0,m+k}=-(m+k)c_{m,k}$.

The coboundary of a 1-cochain $f:\Witt\to\fh$ is 
$df: a\wedge b\mapsto -f([a,b])$.
Define $f$ by $\ell_n\mapsto -c_{0,n}/n$ if $n\ne 0$ and $\ell_0\mapsto 0$.
Then $(df)(\ell_n\wedge\ell_m)=-(n-m)f(\ell_{n+m})=
(n-m)c_{0,n+m}/(n+m)=-c_{n,m}$ if $n+m\ne 0$ and $(df)(\ell_n\wedge\ell_m)=0$
otherwise.
Thus $\ep+df:\ell_n\wedge\ell_m\mapsto \de_{n+m}c_{n,-n}$.
 
\smallskip

\itii\:
The cocycle condition for $\ep$ is  
$$
\de_{n+m+k}
\big(
(n-m)c_{n+m}
\; +\;
(k-n)c_{k+n}
\; +\;
(m-k)c_{m+k}
\big)
=0.
$$

\smallskip

`$\Rightarrow$'\:
We have to show that $(n-m)c_{n+m}-(2n+m)c_{-m}+(2m+n)c_{-n}=0$.
This is clear for $c_n=na$. 
For $c_n=n^3b$ this follows from
$(n-m)(n+m)^3=(n^2-m^2)(n^2+m^2+2nm)=n^4-m^4+2n^3m-2nm^3$.

\smallskip

`$\Leftarrow$'\:
For any $h\in\fh$,
the 2-cochain $\ell_n\wedge\ell_m\mapsto n\de_{n+m}h$ is the coboundary of
$f:\ell_n\mapsto -\de_n h/2$ because
$-f([\ell_n,\ell_m])=(n-m)\de_{n+m}h/2=n\de_{n+m}h$.
Thus $\ell_n\wedge\ell_m\mapsto d_n\de_{n+m}$ is also a cocycle
where $d_n:=c_n-nc_1$. 

We have $d_1=0$. 
Since $\ep$ is a cochain we have $c_{-n}=-c_n$ and hence also $d_{-n}=-d_n$.
Thus $d_{-1}=d_0=0$.
The cocycle condition for $n=1$ and $k=-1-m$ is 
$(1-m)d_{1+m}+(-m-2)d_{-m}=0$ using that $d_{-1}=0$.
This implies that $d_{m+1}=d_m(m+2)/(m-1)$ for any $m\geq 2$.
Since $d_{-m}=-d_m$ we get $d_m=\binom{m+1}{3}d_2$ for any $m\in\Z$.

Suppose that $n^3\de_{n+m}a=-(n-m)f(\ell_{n+m})$ for some $f$.
For $n=-m$ we get $n^3a=-2nf(\ell_0)$ for any $n$. 
Contradiction.
\end{pf}

\bigskip

{\bf Remark.}\: {\it
Let $\fg$ be a perfect Lie algebra and $\fe, \fe'$ two central extensions 
of $\fg$.
Then there exists at most one morphism $\fe\to\fe'$ over $\fg$.
}

\bigskip

\begin{pf}
Let $\phi, \psi: \fe\to\fe'$ be two morphisms over $\fg$ and $\de:=\phi-\psi$.
Then $\de\fe\subset\ker(\fe'\to\fg)$.
Thus $\phi[a,b]=[\psi a+\de a,\psi b+\de b]=\psi[a,b]$.
Since $\fe$ is perfect we get $\phi=\psi$.
\end{pf}

\bigskip

{\bf Proposition.}\: {\it
The universal central extension of $\Witt$ is given by the 2-cocycle 
$\ell_n\wedge\ell_m\mapsto n^3\de_{n+m}$.
}

\bigskip

\begin{pf}
Let $\Vir'$ be the central extension of $\Witt$
given by the 2-cocycle $\ell_n\wedge\ell_m\mapsto n^3\de_{n+m}$
and $L_n, \hc, n\in\Z$, be a basis of $\Vir'$ 
such that $\hc$ is central and $[L_n,L_m]=(n-m)L_{n+m}+n^3\de_{n+m}\hc$. 
Since $[L_1,L_{-1}]=2L_0+\hc$ and $[L_2,L_{-2}]=4L_0+8\hc$
we see that $\hc\in[\Vir',\Vir']$. 
Because $\Witt$ is perfect it follows that $\Vir'$ is perfect, too. 

Let $\fg$ be a central extension of $\Witt$ by the trivial $\Witt$-module 
$\fh$. 
The Lemma shows that $H^2(\Witt,\fh)$ is isomorphic to $\fh$ via 
$[\ep_a]\mapsto a$ where $\ep_a: \ell_n\wedge\ell_m\mapsto n^3\de_{n+m}a$.
Thus the relation between $H^2(\Witt,\fh)$ and central extensions of $\Witt$
by $\fh$ directly implies that there exists a morphism $\Vir'\to\fg$ over
$\Witt$.
This morphism is unique because of the Remark.
\end{pf}

\bigskip

The Proposition and the Lemma show that $\Vir$ has a basis 
$L_n, \hc, n\in\Z$, such that $\hc$ is central and 
$$
[L_n,L_m]
\; =\;
(n-m)L_{n+m}
\; +\;
(n^3-n)\de_{n+m}\hc/12.
$$
Since $[L_{-1},L_n]=-(n+1)L_{n-1}$ we see that $\Vir$ with the families
$(L_n), h=2$, and $(\de_n\hc), h=0$, is a graded local Lie algebra.

\subsection{Super Witt, Contact, and Special Lie Algebras}
\label{SS:super witt}

We show that the super Witt algebras $W_N$ are graded local Lie algebras.
We state that the same is true for the contact and special Lie algebras 
$K_N$ and $S'_N$ and we describe $K_1$ explicitly.

\bigskip

\subsubsection{Super Witt Algebras.}
Let $\La(N):=\K[\xi_1, \dots, \xi_N]$ be the polynomial ring in $N$ odd
variables. 
We consider $\La(N)$ as a graded ring with $h_{\xi_i}:=1/2$.
Define $\La(1,N):=\K[z\uppm]\otimes\La(N)$ where $z$ is even. 
The $N$=$N$ {\bf super Witt algebra} $W_N$ is the Lie algebra of 
derivations of $\La(1,N)$. 
In particular, $W_0=\Witt$.

The $\La(1,N)$-module $W_N$ is free of rank $N+1$ with basis
$\del_z, \del_i:=\del_{\xi_i}$.
Let $a\in\La(N)$. 
Define $a_n^0:=-z^{n+1-h_a} a\del_z$ for $n\in\Z+h_a$ and 
$a_n^i:=z^{n+(1/2)-h_a}a\del_i$ for $n\in\Z+h_a-1/2$ and $i=1, \dots, N$.
Then 
\begin{align}
\notag
[a_n^0,b_m^0]
\; &=\;
(n-m-h_a+h_b)(ab)_{n+m}^0,
\\
\notag
[a_n^0,b_m^j]
\; &=\;
-(m+(1/2)-h_b)(ab)_{n+m}^j-\zeta^{a+b}((\del_j a)b)^0_{n+m},
\\
\notag
[a_n^i,b_m^j]
\; &=\;
(a\del_i b)_{n+m}^j+\zeta^a((\del_j a)b)^i_{n+m}.
\end{align}
In particular, $[1_{-1}^0,a_n^0]=-(n+1-h_a) a_{n-1}^0$ and 
$[1_{-1}^0,a_n^i]=-(n+(1/2)-h_a)a_{n-1}^i$.
Thus $W_N$ with the families $(a_n^0), h=2-h_a$, and 
$(a_n^i), h=(3/2)-h_a$ for $a\in\La(N), i=1, \dots, N$,
is a graded local Lie algebra.
One can show that $W_N$ is simple.

It is clear that $W_{N-1}\subset W_N$ is a subalgebra.
In particular, we have $\Witt\subset W_N$ and $\ell_n=1_n^0$.
A calculation shows that for any $\la\in\K^N$ the map
$\ell_n\mapsto\ell_n^{\la}:=1_n^0-(n+1)\sum_i \la_i (\xi_i)^i_n$
is a monomorphism $\Witt\to W_N$.
As will be clear later, this fact also follows from Proposition 
\ref{SS:chodos thorn}.

\subsubsection{Contact Lie Algebras.}
Let $\Om(1,N)$ be the $\La(1,N)$-module of K\"ahler differentials with
universal differential $d: \La(1,N)\to\Om(1,N)$.
Then $df=(dz)\del_zf+\sum(d\xi_i)\del_i f$.
The Lie algebra $W_N$ acts on $\Om(1,N)$ by derivations such that
$d: \La(1,N)\to\Om(1,N)$ is an odd $W_N$-module morphism. 
The action of $W_N$ on $\Om(1,N)$ is called 
\index{Lie derivative}
the {\bf Lie derivative}.

Define $\om:=dz-\sum\xi_i d\xi_i$.
The $N$=$N$ {\bf contact Lie algebra} is
\index{contact Lie algebra}
$$
K_N
\; :=\;
\set{\del\in W_N\mid \del\om=p\om \text{ for some } p\in\La(1,N)}.
$$
This is a subalgebra because $\del\om=p\om$ and $\del'\om=p'\om$ implies
$[\del,\del']\om=(\del p'+(-1)^{p'\del}p'p-
(-1)^{\del\del'}\del'p-(-1)^{\del\del'+p\del'}pp')\om
=(\del p'-(-1)^{\del\del'}\del'p)\om$.
We have $K_0=W_0=\Witt$ since $(p\del_z)dz=(\del_z p)dz$ for any 
$p\in\K[z\uppm]$.
One can show that $K_2\cong W_1$.

One can show that $K_N$ is a graded local Lie algebra, 
$F_{K_N}$ has $2^N$ elements, and $K_N$ is simple iff $N\ne 4$.
Moreover, $K'_4:=[K_4,K_4]$ is a simple graded local Lie algebra.
We have $\ell_n^{\la}\in K_N$ for any $n$ iff $\la_i=1/2$ for any $i$.

One can show that a basis of $K_1$ is given by the derivations $\ell_n$ and
$g_m$ where $\ell_n:=\ell_n^{1/2}=-z^{n+1}\del_z-(1/2)(n+1)z^n\xi\del_{\xi}$ 
and $g_m:=z^{m+1/2}(\xi\del_z+\del\subxi)$ for $n\in\Z$ and $m\in\Z+1/2$.
We have
$$
[\ell_n,g_m]
\; =\;
((n/2)-m)g_{n+m},
\qquad
[g_n,g_m]
\; =\;
2\ell_{n+m}.
$$
Since $[\ell_{-1},g_m]=-(m+1/2)g_{m-1}$ we see that $K_1$ 
with the families $(\ell_n), h=2$, and $(g_n), h=3/2$,
is a graded local Lie algebra.

\subsubsection{Special Lie Algebras.}
The {\bf divergence} of $X=\sum X_i\del_i\in W_N$ is
$$
\rdiv X
\; :=\;
\sum_{i=0}^N\:
\zeta^{X_i\del_i}\del_i X_i
$$
where $\del_0:=\del_z$.
The $N$=$N$ {\bf special} Lie algebra is 
$$
S_N
\; :=\;
\set{X\in W_N\mid \rdiv X=0}.
$$
This is a subalgebra because $\rdiv[X,Y]=X\rdiv Y-\zeta^{XY}Y\rdiv X$.

One can show that $S_N=S'_N\oplus\K\xi_1\dots\xi_N\del_z$ where 
$S'_N:=[S_N,S_N]$.
Moreover, $S'_N$ is a graded local Lie algebra, 
$F_{S'_N}$ has $N2^N$ elements, and $S'_N$ is simple iff $N\geq 2$.

\subsection{Super Virasoro Algebras}
\label{SS:super vir}

We give an overview of the super Virasoro algebras.
They are universal central extensions of the simple Lie algebras 
$K_N, N\leq 3, S_2'$, and $K_4'$ and are also graded local Lie algebras.

\bigskip

Let $0\leq N\leq 3$.
The $N$=$N$ {\bf super Virasoro algebra} $\NVir$ is the 
universal central extension of the contact Lie algebra $K_N$.
Thus $\zeroVir=\Vir$.
One can show that the centre of $\NVir$ is one-dimensional and
$\NVir$ is a graded local Lie algebra.
For $N=0$ these results were proven in section \ref{SS:vir alg}.

The $N$=1 super Virasoro algebra $\oneVir$ is also called the 
{\bf Neveu-Schwarz algebra}.
It has a basis $L_n, G_m, \hc$ for $n\in\Z, m\in\Z+1/2$
such that $L_n, \hc$ span a Virasoro algebra and $G_m$ are odd with
$$
[L_n,G_m]
=
((n/2)-m)G_{n+m},
\quad
[G_n,G_m]
=
2L_{n+m}
+
(4n^2-1)\de_{n+m}\hc/12.
$$
Thus $\oneVir$ with the families $(L_n), (G_m), (\de_n\hc)$ and
$h=2, 3/2, 0$ is a graded local Lie algebra.
We have $T=[L_{-1},\;]$.

The $N$=2 super Virasoro algebra $\twoVir$ has a basis 
$L_n, G^+_m, G^-_m, J_n, \hc$ for $n\in\Z, m\in\Z+1/2$ such that
$L_n, \hc$ span a Virasoro algebra, $J_n, \hc$ span a Heisenberg Lie algebra
with $[J_n,J_m]=n\de_{n+m}\hc/3$ and $[L_n,J_m]=-mJ_{n+m}$, and 
$G_m\suppm$ are odd with 
$$
[L_n,G_m\suppm]
=
((n/2)-m)G_{n+m}\suppm,
\qquad
[J_n,G_m\suppm]
=
\pm G_{n+m}\suppm,
$$
$$
[G_n^+,G_m^-]
=
2L_{n+m}
+
(n-m)J_{n+m}
+
(4n^2-1)\de_{n+m}\hc/12,
$$
$[G_n^+,G_m^+]=[G_n^-,G_m^-]=0$.
Thus $\twoVir$ with the families $(L_n), (G_m\suppm), (J_n), (\de_n\hc)$ and
$h=2, 3/2, 1, 0$ is a graded local Lie algebra.

The $N$=4 {\bf super Virasoro algebra} $\fouVir$ is the 
universal central extension of $S'_2$.
The {\bf large} $N$=4 {\bf super Virasoro algebra} $\cA$
is the universal central extension of $K'_4$.
One can show that $\fouVir$ and $\cA$ are graded local Lie algebras.

The centre of $\fouVir$ and of the universal central extension of 
$W_N$ for $N\leq 2$ are one-dimensional.
The centre of $\cA$ is two-dimensional.
The Lie algebras $W_N$ for $N\geq 3$, $K_N$ for $N\geq 5$, and 
$S'_N$ for $N\geq 3$ have no non-trivial central extensions.

It seems that $\Vir, \oneVir, \twoVir, \fouVir$, and $\cA$
are the only physically interesting Lie algebras that contain
the Virasoro algebra.

\section{Weak Commutator Formula and Locality}
\label{S:loc weak comm f}

In section \ref{SS:distr} 
we write the weak commutator formula in terms of distributions.
In sections \ref{SS:diff op}--\ref{SS:loc diff op}
we prove that the weak commutator formula is equivalent to locality
using the correspondence between distributions and morphisms defined
on test functions.

\subsection{Weak Commutator Formula and Distributions}
\label{SS:distr}

We rewrite the weak commutator formula in terms of distributions.

\bigskip

Let $\fg$ be a Lie algebra and $a_t, b_t\in\fg$ for any $t\in\Z$. 
Since $\K[t,s]=\K[t,t+s]$ and since the polynomials $\binom{t}{i}$ 
for $i\geq 0$ form a basis of $\K[t]$,
we see that $(a_t), (b_t)$ satisfy the weak commutator formula iff
there exist families $(c^i_t)$ in $\fg$ such that
$$
[a_t,b_s]
\; =\;
\sum_{i=0}^n\:
\binom{t}{i}\, c^i_{t+s-i}.
$$
This is 
\index{weak commutator formula}
the {\bf weak commutator formula}.

We identify a family $(a_t)$ of vectors of a vector space $E$ 
with the formal sum 
$$
a(z)
\; =\;
\sum_{t\in\Z}\: a_t z^{-t-1}
$$ 
where $z$ is an even formal variable. 
We 
\index{distribution}
call $a(z)$ a {\bf distribution}.

The space $E\pau{z\uppm}$ of $E$-valued distributions
is a module over the ring $\K[z\uppm]$ of 
Laurent polynomials, $z^s a(z):=\sum a_{t+s}z^{-t-1}$.
The module $E\pau{z\uppm}$ has a derivation $\del_z$ defined by 
$\del_z a(z):=-\sum t\, a_{t-1}z^{-t-1}$.
We identify $a\in E$ with the constant distribution 
$a=\sum\de_{t,-1}a z^{-t-1}$.

A linear map $\ka:E\otimes F\to G$ induces 
\index{operator product} 
a morphism 
\begin{align}
\notag
E\pau{z\uppm}\otimes F\pau{w\uppm}&\; \to\; \qquad G\pau{z\uppm,w\uppm},
\\
\notag
a(z)\otimes b(w)\quad &\; \mapsto\; 
\sum_{t,s}\:\ka(a_t\otimes b_s)z^{-t-1}w^{-s-1},
\end{align}
called the {\bf operator product}. 
The operator product is in general not well-defined for $z=w$.
The singularity at $z=w$ is studied using the operator product expansion,
see section \ref{SS:ope sing}.

The {\bf delta distribution} is 
\index{delta distribution}
$$
\de(z,w)
\; :=\;
\sum_{t\in\Z}\:
w^t\, z^{-t-1}
\; \in\;
\K[w\uppm]\pau{z\uppm}.
$$
Then the weak commutator formula, written in terms of distributions, is 
$$
[a(z),b(w)]
\; =\;
\sum_{i=0}^n\:
c^i(w)\, \del_w^{(i)}\de(z,w)
$$
where $\del_z^{(i)}:=\del_z^i/i!$.
    
If $\fg$ is a differential algebra then we denote by $\fg\pau{z\uppm}_T$ 
the space of $a(z)\in\fg\pau{z\uppm}$ such that 
$T a(z)=\del_z a(z)$, in other words $Ta_t=-t a_{t-1}$.
Distributions in $\fg\pau{z\uppm}_T$ are 
\index{translation covariant}
called {\bf translation covariant}.
The third axiom in the definition of a local Lie algebra is that
$F\subfg\subset\fg\pau{z\uppm}_T$.

If $R$ is a $\K[T]$-module then we consider $\fgl(R)$ as a 
differential algebra with $T=[T,\;]$.
Thus $a(z)\in\fgl(R)\pau{z\uppm}_T$ iff $[T,a(z)]=\del_z a(z)$.

If $\fg$ is a graded local Lie algebra then
$a(z)=\sum_{n\in\Z-h_a}a_n z^{-n-h_a}$.

\subsection{Distribution Kernels and Differential Operators}
\label{SS:diff op}

We show that distributions of the form $\sum_i c^i(w)\del_w^{(i)} \de(z,w)$ 
are distribution kernels of differential operators and give a formula
for $c^i(z)$.

\bigskip

For $a(z)\in E\pau{z\uppm}$, define the linear map $\al:\K[z\uppm]\to E$ by
$$
\al:\:
p(z)
\;\mapsto\;
\res_z\, p(z)a(z)
$$
where $\res_z a(z):=a_0$ is the {\bf residue}.
Thus $\al: z^t\mapsto a_t$ and we get a natural isomorphism $a(z)\mapsto\al$ 
from $E\pau{z\uppm}$ to the space of linear maps $\K[z\uppm]\to E$.
This is a $\K[z\uppm]$- and a $\K[\del_z]$-module isomorphism
where $p\al:=\al(p\cdot\:)$ and $\del_z\al:=-\al\circ\del_z$.
The second claim amounts to the {\bf integration-by-parts} formula
\index{integration-by-parts formula}
$$
-\res_z \del_z p(z)a(z)
\; =\;
\res_z p(z)\del_z a(z)
$$
that follows from $\res_z \del_z a(z)=0$ and the product formula.
We call $a(z)$ the {\bf distribution kernel} of  $\al$.
The delta distribution $\de(z,w)$ is the distribution kernel of 
the linear map $\K[z\uppm]\to \K[w\uppm], p(z)\mapsto p(w)$.

In particular, there is an isomorphism $c(z,w)\mapsto\al$ from
$E\pau{z\uppm,w\uppm}$ to the space of linear maps 
$\K[z\uppm]\to E\pau{w\uppm}$.
We identify such linear maps with $\K[w\uppm]$-module morphisms 
$\K[z\uppm,w\uppm]\to E\pau{w\uppm}$.
The morphism $\al$ is given by 
$$
\al:\:
p(z,w)
\;\mapsto\;
\res_z\, p(z,w)c(z,w).
$$

A {\bf differential operator} of order $n$ is 
\index{differential operator}
a morphism $\al:\K[z\uppm,w\uppm]\to E\pau{w\uppm}$ 
of the form
$$
\al:\: 
p(z,w)
\;\mapsto\; 
\sum_{i=0}^n\: c^i(w)\,\del_z^{(i)} p(z,w)|_{z=w}
$$
for some $c^i(z)\in E\pau{z\uppm}, c^n(z)\ne 0$.
The distribution kernel of $\al$ is 
$$
c(z,w)
\; =\;
\sum_{i=0}^n\: 
c^i(w)\, \del_w^{(i)}\de(z,w)
$$
because for any $p(z)\in\K[z\uppm]$ we have
$$
\res_z p(z)c^i(w)\del_w^{(i)}\de(z,w)
\; =\; 
c^i(w)\del_w^{(i)}\res_z p(z)\de(z,w)
\; =\; 
c^i(w)\del_w^{(i)}p(w).
$$
The coefficient $c^i(z)$ is given by 
$c^i(w)=\al((z-w)^i)=\res_z (z-w)^i c(z,w)$.
In particular, we have:

\bigskip

{\bf Proposition.}\: {\it
The distributions $c^i(z)$ of the weak commutator formula are
$$
c^i(w)
\; =\;
\res_z (z-w)^i [a(z),b(w)].
$$
\hfill $\square$
}

\subsection{Locality of Differential Operators}
\label{SS:delta dis diff op}

We prove that the distribution kernel of a differential operator is local.

\bigskip

A distribution $c(z,w)\in E\pau{z\uppm,w\uppm}$ is {\bf local} of 
\index{local!distribution}
order $\leq n$ if
$$
(z-w)^n c(z,w) 
\; =\; 
0.
$$

\bigskip

{\bf Proposition.}\: {\it
For any $i\geq 0$, we have
$$
(z-w)\del_w^{(i)}\de(z,w) 
\; =\; 
\del_w^{(i-1)}\de(z,w).
$$
In particular, the distribution kernel of a differential operator of order $n$
is local of order $n+1$.
}

\bigskip

\begin{pf}
For any $p(z)\in\K[z\uppm]$, we have
\begin{align}
\notag
\res_z p(z)(z-w)\del_w^{(i)}\de(z,w)
\; &=\; 
\del_z^{(i)}(p(z)(z-w))|_{z=w}
\; =\; 
\del_w^{(i-1)}p(w)
\\
\notag
&=\;
\res_z p(z)\del_z^{(i-1)}\de(z,w).
\end{align}
\end{pf}

\bigskip

Denote by $E\lau{z}$ the $\K[z\uppm]$-submodule of $E\pau{z\uppm}$ 
generated by the subspace $E\pau{z}$ of power series.
The module $E\lau{z_1, \dots, z_r}$ is defined in the same way.
Its elements are 
\index{Laurent series}
called {\bf Laurent series}.

\bigskip

{\bf Remark.}\: {\it
Let $c(z,w)$ be in $E\pau{w\uppm}\lau{z}$ or
$E\pau{w\uppm}\lau{z\inv}$ or $E\pau{z\uppm}\lau{w}$ or 
$E\pau{z\uppm}\lau{w\inv}$.
If $c(z,w)$ is local then $c(z,w)=0$.
}

\bigskip

\begin{pf}
Let $c(z,w)\in E\pau{w\uppm}\lau{z}$.
If $c_{p,q}\ne 0$ and $c_{t,s}=0$ for any $t>p$ and $s\in\Z$ then
$(z-w)c(z,w)=\sum_{t,s}(c_{t+1,s}-c_{t,s+1})z^{-t-1}w^{-s-1}$
and $c_{p+1,q-1}-c_{p,q}=-c_{p,q}\ne 0$.
The other cases are proven in the same way.
\end{pf}

\subsection{Locality of Distributions}
\label{SS:loc diff op}

We prove that a distribution is local iff it is the distribution kernel of 
a differential operator.
This result corresponds to the fact from functional analysis that
a distribution is supported on the diagonal iff it is the distribution kernel
of a differential operator.

\bigskip

{\bf Lemma.} (Taylor's Formula)\: {\it
Let $c(z,w)\in E\lau{z,w}$ and $n\geq 0$. 
Then $c^i(w)=\del_z^{(i)}c(z,w)|_{z=w}$ are 
\index{Taylor's formula}
the unique distributions in $E\lau{w}$ such that
$$
c(z,w)
\; =\;
\sum_{i=0}^{n-1}\:
c^i(w)\, (z-w)^i
\; +\;
(z-w)^n\, r(z,w)
$$
for some $r(z,w)\in E\lau{z,w}$.
If $c(z,w)\in E[z\uppm,w\uppm]$ then $r(z,w)\in E[z\uppm,w\uppm]$.
If $c^i(z)=0$ for any $i\geq 0$ then $c(z,w)=0$.
}

\bigskip

\begin{pf}
We do induction on $n$.
For $n=1$ we have to show that if $d(z,w)\in E\lau{z,w}$ vanishes for $z=w$ 
then $d(z,w)=(z-w)r(z,w)$ for some $r(z,w)\in E\lau{z,w}$. 
We may assume that $d(z,w)\in E\pau{z,w}$ and show $r(z,w)\in E\pau{z,w}$. 
By projecting onto a basis of $E$ we may assume that $E=\K$.
Write $d(z,w)=\sum_i d_i(z,w)$ with $d_i(z,w)$ homogeneous of degree $i$. 
Then $d_i(w,w)=0$ for any $i$. 
Thus $d_i(z,w)$ is divisible by $z-w$ and hence so is $d(z,w)$.

Using induction and applying the induction beginning $n=1$ to $r(z,w)$
we see that there exist $c^i(z)\in E\lau{z}$ and $r(z,w)\in E\lau{z,w}$ 
such that the above identity for $c(z,w)$ is satisfied.
Applying $\del_z^{(i)}$ to both sides of this identity
and setting $z=w$ we get $\del_z^{(i)}c(z,w)|_{z=w}=c^i(w)$.

The second claim follows from the arguments we already gave.
The last claim follows from the fact that polynomials like $d_i(z,w)$ 
have only a finite order of vanishing at $z=w$.
\end{pf}

\bigskip

{\bf Proposition.\;}{\it
A distribution in $E\pau{z\uppm,w\uppm}$ is local
iff it is the distribution kernel of a differential operator.
}

\bigskip

\begin{pf}
Suppose that $c(z,w)$ is local of order $\leq n$.
Let $\al:\K[z\uppm,w\uppm]\to E\pau{w\uppm}$ be the module morphism 
with distribution kernel $c(z,w)$.
By Taylor's formula for any $p(z,w)\in\K[z\uppm,w\uppm]$
there exists $r(z,w)\in\K[z\uppm,w\uppm]$ such that
\begin{align}
\notag
&\al(p(z,w))
\\
\notag
=\;
&\res_z\Bigg( 
\sum_{i=0}^{n-1}\:
\del_z^{(i)}p(z,w)|_{z=w}\,(z-w)^i
\; +\;
(z-w)^n r(z,w)\Bigg) \: c(z,w)
\\
\notag
=\;
&\sum_{i=0}^{n-1}\;\;
\del_z^{(i)}p(z,w)|_{z=w}\;\;
\res_z (z-w)^i c(z,w).
\end{align}
Thus $\al$ is a differential operator.
The converse is Proposition \ref{SS:delta dis diff op}.
\end{pf}

\bigskip

Let $\fg$ be an algebra.
Then $a(z), b(z)\in\fg\pau{z\uppm}$ are {\bf local} if $[a(z),b(w)]$ is 
\index{local!pair of distributions}
local. 
The Proposition implies the following result.

\bigskip

{\bf Corollary.}\: {\it
Let $\fg$ be an algebra.
Then $a(z), b(z)\in\fg\pau{z\uppm}$ are local iff  
they satisfy the weak commutator formula.
\hfill $\square$
}

\bigskip

A subset $S\subset\fg\pau{z\uppm}$ is {\bf local} if 
\index{local!set of distributions}
any pair of elements of $S$ is local.
Thus a local Lie algebra is a differential Lie algebra $\fg$
together with a local subset $F\subfg\subset\fg\pau{z\uppm}_T$ such that
$\fg=\rspan\set{a_t\mid a(z)\in F\subfg, t\in\Z}$.

\section{Vertex Lie Algebras of Distributions}
\label{S:loclie to vlie}

In sections \ref{SS:confa}--\ref{SS:weak comm f ex}
we show that $\fg\pau{z\uppm}$ is an unbounded conformal algebra with 
$i$-th products given by the distributions $c^i(z)$ of 
the weak commutator formula.
In sections \ref{SS:conf jac}--\ref{SS:vlie}
we prove that if $\fg$ is a local Lie algebra then
$F\subfg$ generates a vertex Lie subalgebra of $\fg\pau{z\uppm}$.

\subsection{Conformal Algebras}
\label{SS:confa}

We show that unbounded conformal algebras are equivalent to
$\N$-fold algebras with a translation operator.

\bigskip

An {\bf unbounded $\la$-product} on a $\K[T]$-module $R$ is 
\index{unbounded!lambda-product@$\la$-product}
an even linear map 
$R\otimes R\to R\pau{\la}, a\otimes b\mapsto a\subla b$, such that
$$
T(a\subla b)
\; =\;
(Ta)\subla b
\; +\;
a\subla(Tb)
\qquad\text{and}\qquad
(Ta)\subla b
\; =\;
-\la\, a\subla b.
$$
One obtains an equivalent definition if one replaces 
either one of the above identities by
$$
a\subla Tb
\; =\;
(T+\la)(a\subla b).
$$

An {\bf unbounded conformal algebra} is 
\index{unbounded!conformal algebra}
a $\K[T]$-module with an unbounded $\la$-product.
A {\it morphism} of unbounded conformal algebras is a $\K[T]$-module morphism 
$\phi$ such that $\phi(a\subla b)=(\phi a)\subla(\phi b)$.

Elements $a, b$ of an unbounded conformal algebra $R$ are {\bf weakly local} 
\index{weakly local}
if $a\subla b\in R[\la]$.
If any $a, b\in R$ are weakly local
then $R$ is called a {\bf conformal algebra}
\index{conformal algebra}
and $a\subla b$ is called 
\index{lambda-product@$\la$-product}
a {\bf $\la$-product}.
We interchangeably denote 
\index{unbounded!lambda-bracket@$\la$-bracket}
$\la$-products by $[a\subla b]$ and call them
\index{lambda-bracket@$\la$-bracket}
{\bf $\la$-brackets}. 

What we have just defined are {\it left} $\la$-products.
Right $\la$-products satisfy $a\subla Tb=-\la\, a\subla b$ instead of
$(Ta)\subla b=-\la\, a\subla b$.
If $a\subla b$ is a left $\la$-product then 
$a\otimes b\mapsto\paraab b\subla a$ is a right $\la$-product.
More interesting is the {\bf opposite} $\la$-product
\index{opposite!lambdaproduct@$\la$-product}
$a\otimes b\mapsto\paraab b_{-\la-T}a$.
Here $c_T:=\sum T^i(c_i)$ for $c\subla=\sum c_i\la^i\in R[\la]$.
The opposite $\la$-product is a left $\la$-product.

A {\bf derivation} of an unbounded conformal algebra $R$ is 
\index{derivation!of a conformal algebra}
a $\K[T]$-module morphism $d$ such that 
$d(a\subla b)=(da)\subla b+\zeta^{da}a\subla(db)$. 

A {\bf left ideal} of $R$ is
\index{ideal}
a $\K[T]$-submodule $I$ such that $R\subla I\subset I\pau{\la}$.
Right ideals and (two-sided) ideals are defined in the same way.
If $I$ is an ideal then $R/I$ is 
\index{quotient}
the {\bf quotient} conformal algebra.

The {\bf $t$-th products} $a_t b$ of $R$ are 
\index{tth product@$t$-th product}
defined by
$$
[a\subla b]
\; =\;
\sum_{t\geq 0}\: a_t b\, \la^{(t)}.
$$
We shall always use Greek letters for $\la$-products $a\subla b$ 
and Latin letters for $t$-th products $a_t b$ in order to distinguish them.

Let $S$ be an even set.
An {\bf $S$-fold algebra} is 
\index{Sfold algebra@$S$-fold algebra}
a vector space $V$ together with a family of multiplications indexed by $S$:
$$
V\otimes V\;\to\; V, 
\quad 
a\otimes b\;\mapsto\; a_s b,
\quad s\in S.
$$ 
We 
\index{sth product@$s$-th product}
call $a_s b$ the {\bf $s$-th product} and denote by
$a_s$ the operator $b\mapsto a_s b$.
We always use the {\bf left-operator notation}
\index{left!operator notation}
$a_t b_s:=a_t\circ b_s$ so that $a_t b_s c=a_t(b_s c)$.

A {\it morphism} of $S$-fold algebras is an even linear map $\phi$ 
such that $\phi(a_s b)=(\phi a)_s (\phi b)$.
A {\bf derivation} of an $S$-fold algebra is an operator $d$ such that 
$d(a_s b)=(da)_s b+\zeta^{da}\, a_s (db)$.

A {\bf translation operator} of an $\N$-fold algebra 
\index{translation operator}
is an even derivation $T$ such that $(Ta)_t b=-t\, a_{t-1}b$,
where $t\, a_{t-1}b:=0$ if $t=0$.
An $\N$-fold algebra with a translation operator 
is the same thing as an unbounded conformal algebra because
$-\la\sum a_t\la^{(t)}=-\sum t\,a_{t-1}\la^{(t)}$.

If $R$ is a {\it graded} $\K[T]$-module, we change notation and
denote the $t$-th product of a $\la$-bracket by $a_{(t)}b$ so that 
$a\subla b=\sum a_{(t)} b\,\la^{(t)}$.
A {\bf gradation} of a conformal algebra $R$ is 
\index{graded!conformal algebra}
a $\K[T]$-module gradation $R=\bigoplus_{h\in\K}R_h$ such that 
$(R_h)_{(t)} R_k\subset R_{h+k-t-1}$.
For $a\in R_h$ and $n\in\N-h+1$, 
we define $a_n:=a_{(n+h-1)}$.
Then $[H,a_n]=-na_n$.

A conformal algebra of {\bf CFT-type} is 
\index{CFT-type!conformal algebra of}
a $\rho\N$-graded conformal algebra $R$, for some $\rho\in\Q_>$,
such that $\dim R_0=1$ and $TR_0=0$.

\subsection{Unbounded Conformal Algebras of Distributions}
\label{SS:alg of distr}

We show that $\fg\pau{z\uppm}$ is an unbounded conformal algebra.
The $i$-th products of $\fg\pau{z\uppm}$ are the distributions $c^i(z)$ of 
the weak commutator formula.
We also consider the subalgebras $\fg\pau{z\uppm}_T, \fg\pau{z\uppm}_H$,
and $\Endv(E)$.

\bigskip

If $\fg$ is an algebra then $\fg\pau{z\uppm}$ is an unbounded conformal algebra
with $T=\del_z$ and 
$$
[a(w)\subla b(w)]
\; :=\;
\res_z e^{(z-w)\la}[a(z),b(w)]
$$
because $\del_w\res_z\: e^{(z-w)\la}[a(z),b(w)]=-\la[a(w)\subla b(w)]
+[a(w)\subla \del_w b(w)]$ and 
$\res_z e^{(z-w)\la}[\del_z a(z),b(w)]=-\la[a(w)\subla b(w)]$
by the integration-by-parts formula.

The $i$-th product of $\fg\pau{z\uppm}$ is the coefficient $c^i(z)$ of 
the weak commutator formula: 
$$
a(w)_i b(w)
\; =\;
\res_z (z-w)^i[a(z),b(w)].
$$
In particular, $a(z)_0 b(z)=[a_0,b(z)]$.
If $a(z), b(z)$ are local then they are weakly local.

If $\fg$ is a differential algebra then $\fg\pau{z\uppm}_T$ consists of $a(z)$ 
such that $Ta(z)\linebreak[0]=\del_z a(z)$, see section \ref{SS:distr}.
This is an unbounded conformal subalgebra because $T$ and $\del_z$
are derivations of $\fg\pau{z\uppm}$.

If $\fg$ is a $\K$-graded Lie algebra
then we denote by $a_{(t)}$ the coefficients of $a(z)\in\fg\pau{z\uppm}$ 
so that $a(z)=\sum a_{(t)}z^{-t-1}$.

For $h\in\K$, let $\fg\pau{z\uppm}_h$ be the space of $a(z)\in\fg\pau{z\uppm}$
such that $a_{(t)}\in\fg_{h-t-1}$. 
This is equivalent to $H a(z)=h a(z)+z\del_z a(z)$.
There is an inclusion $\fg_h\to\fg\pau{z\uppm}_h, a\mapsto\sum\de_t\, a z^t$.

The space $\fg\pau{z\uppm}_h$ corresponds to $(\fg^{\Z})_h$ via
the isomorphism $\fg\pau{z\uppm}\to\fg^{\Z}, a(z)\mapsto(a_{(t)})$,
see section \ref{SS:loclie}.
As in section \ref{SS:loclie}, we define $a_n:=a_{(n+h-1)}$ for 
$a(z)\in\fg\pau{z\uppm}_h$ and $n\in\Z-h$.
Then $a_n\in\fg_{-n}$ and $a(z)=\sum a_n z^{-n-h}$.

The space $\fg\pau{z\uppm}_H:=\bigoplus\fg\pau{z\uppm}_h$
is a graded unbounded conformal algebra because if $a(z), b(z)$
are homogeneous then $(\del_z a(z))_{(t)}=-ta_{(t-1)}\in\fg_{h_a+1-t-1}$ and
$$
a(w)_i b(w)
\; =\;
\sum (-1)^j\binom{i}{j}\: [a_{(i-j)},b_{(t+j)}]\, w^{-t-1}
$$
with $[a_{(i-j)},b_{(t+j)}]\in\fg_{h_a+h_b-i-1-t-1}$.
We have $(\del_z a(z))_n=(\del_z a(z))_{(n+h_a)}=-(n+h_a)a_n$.
The weak commutator formula for the $a_n$ is
$$
[a_n,b_m]
\; =\;
\sum_{i=0}^N\: \binom{n+h_a-1}{i}\: c^i_{n+m}
$$
since $(n+h_a-1)+(m+h_b-1)-i=n+m+(h_a+h_b-i-1)-1$.

Let $E$ be a vector space.
A distribution $a(z)\in\End(E)\pau{z\uppm}$ is
\index{field}
a {\bf field} on $E$ if $a(z)b\in E\lau{z}$ for any $b\in E$.
In section \ref{SS:fields} we discuss this notion in detail.
Here we only remark that the space $\Endv(E)$ of fields is 
\index{EndvE@$\Endv(E)$}
an unbounded conformal subalgebra of 
$\fgl(E)\pau{z\uppm}$ since $a(w)_t b(w)=\sum\binom{t}{i}[a_i,b(w)](-w)^{t-i}$.
This example is important because $\Endv(E)$ is in fact an
unbounded vertex algebra, see section \ref{SS:fields}.

\subsection{Examples of $\la$-Brackets of Distributions}
\label{SS:weak comm f ex}

We calculate the $\la$-brackets of any $a(z), b(z)\in F\subfg$ 
for affine, Clifford, and $N$=0, 1, 2 Virasoro Lie algebras.
We always have $a(z)_i b(z)\in\K[\del_z]F\subfg$.

\bigskip

Let $\fg$ be a $\K$-graded Lie algebra and 
$a(z), b(z)\in\fg\pau{z\uppm}_H$ be homogeneous.
The ``constant'' weak commutator formula $[a_n,b_m]=c_{n+m}$
is equivalent to $a(z)\subla b(z)=c(z)$.

In particular, the weak commutator formula $[a_n,b_m]=[a,b]_{n+m}$
for a loop Lie algebra $\tfg$ is equivalent to
$a(z)\subla b(z)=[a,b](z)$
and the weak commutator formula $[a_n,b_m]=(a,b)\de_{n+m}\hk$
for a Clifford Lie algebra $C(B)$ is equivalent to  
$a(z)\subla b(z)=(a,b)\hk$.
Here we identify $\hk\in C(B)_0$ with 
$\hk=\sum\de_n\hk z^{-n}$.

The weak commutator formula $[a_n,b_m]=[a,b]_{n+m}+n(a,b)\de_{n+m}\hk$
for an affine Lie algebra $\hfg$ is equivalent to 
$$
a(z)\subla b(z)
\; =\;
[a,b](z)
\; +\; 
(a,b)\hk\la.
$$
In particular, for a Heisenberg Lie algebra $\hB$ we have
$a(z)\subla b(z)=(a,b)\hk\la$.

Suppose that $a(z)\subla b(z)=\del_z c(z)+\al c(z)\la$ for some $\al\in\K$.
Thus $\al c(z)=a(z)_1 b(z)$ and hence $h_c=h_a+h_b-2$.
The weak commutator formula in this case is 
$[a_n,b_m]=(-(n+m+(h_a+h_b-2))+\al(n+h_a-1))c_{n+m}$.
If $\al=(h_a+h_b-2)(h_a-1)\inv$ then 
$[a_n,b_m]=((h_b-1)(h_a-1)\inv n-m)c_{n+m}$.
In particular, if $h_a=2$ then $[a_n,b_m]=((h_b-1)n-m)b_{n+m}$ is equivalent to
$$
a(z)\subla b(z)
\; =\;
\del_z b(z)
\; +\; 
h_b b(z)\la.
$$

The last remark shows that the weak commutator formula 
$[\ell_n,\ell_m]=(n-m)\ell_{n+m}$ 
for the Witt algebra is equivalent to
$\ell(z)\subla\ell(z)=\del_z\ell(z)+2\ell(z)\la$.
The weak commutator formula 
$[L_n,L_m]=(n-m)L_{n+m}+(n^3-n)\de_{n+m}\linebreak[0]\hc/12$
for the Virasoro algebra is equivalent to
$$
L(z)\subla L(z)
\; =\;
\del_z L(z)
\; +\; 
2L(z)\la
\; +\; 
(\hc/2)\la^{(3)}
$$
since $\binom{n+1}{3}(\hc/2)=(n^3-n)\hc/12$.

The $\la$-brackets for the $N$=1 and $N$=2 super Virasoro algebra 
are as follows.
The $\la$-brackets with $L(z)$ are $L(z)\subla a(z)=\del_z a(z)+h_a a(z)\la$ 
for $a=G,G\suppm, J$ since $[L_n,a_m]=((h_a-1)n-m)a_{n+m}$.
The $\la$-brackets with $J(z)$ are $J(z)\subla J(z)=(\hc/3)\la$
and $J(z)\subla G\suppm(z)=\pm G\suppm(z)$ since 
$[J_n,J_m]=n\de_{n+m}\hc/3$ and $[J_n,G\suppm_m]=\pm G\suppm_{n+m}$.

The weak commutator formula $[G_n,G_m]=2L_{n+m}+(4n^2-1)\de_{n+m}\hc/12$
is equivalent to $G(z)\subla G(z)=2L(z)+(2\hc/3)\la^{(2)}$ because
$\binom{n+1/2}{2}(2\hc/3)=(n^2-1/4)\hc/3$.
Finally, the general remark above in the case that $h_a=h_b=3/2$ and
$c(z)=J(z)$ shows that 
$[G_n^+,G_m^-]=2L_{n+m}+(n-m)J_{n+m}+(4n^2-1)\de_{n+m}\hc/12$
is equivalent to
$$
G^+(z)\subla G^-(z)
\; =\;
2L(z)+\del_z J(z)
\; +\;
2J(z)\la
\; +\;
(2\hc/3)\la^{(2)}.
$$

\subsection{Conformal Jacobi Identity}
\label{SS:conf jac}

We prove that the Leibniz identity for an algebra $\fg$ is equivalent to 
the conformal Jacobi identity for $\fg\pau{z\uppm}$.

\bigskip

The {\bf Leibniz identity} for an algebra $\fg$ is 
\index{Leibniz!identity}
$$
[[a,b],c]
\; =\;
[a,[b,c]]
\; -\; 
\paraab\,[b,[a,c]].
$$
It is satisfied iff left multiplication $\fg\to\fgl(\fg), a\mapsto [a,\, ]$, 
is an algebra morphism iff $[a,\, ]$ is a derivation.

What we have defined is the {\it left} Leibniz identity.
The right Leibniz identity is $[a,[b,c]]=[[a,b],c]-\zeta^{bc}[[a,c],b]$.
If $[\, ,]$ is skew-symmetric then they are both equivalent to 
the Jacobi identity $[a,[b,c]]+[b,[c,a]]+[c,[a,b]]=0$.

The {\bf conformal Jacobi identity} for 
\index{conformal!Jacobi identity}
an unbounded conformal algebra is
$$
[[a\subla b]\submu c]
\; =\;
[a\subla[b_{\mu-\la}c]]
\; -\; 
\paraab\, [b_{\mu-\la}[a\subla c]].
$$
Define
$$
[a,b]
\; :=\;
[a_0 b].
$$ 
Setting $\la=\mu=0$ in the conformal Jacobi identity, we get:

\bigskip

{\bf Remark.}\: {\it
If an unbounded conformal algebra satisfies the conformal Jacobi identity 
then $[\, ,]$ satisfies the Leibniz identity.
\hfill $\square$
}

\bigskip

The indices $\la, \mu$, and $\mu-\la$ in the conformal Jacobi identity
are determined by the fact that 
all three terms of the identity transform in the same way
with respect to $T$: replacing $a$ by $Ta$, $b$ by $Tb$, or $c$ by $Tc$ 
is equivalent to multiplication with $-\la$, $\la-\mu$, and $T+\mu$, resp.
This fact also yields the following result.

\bigskip

{\bf Lemma.}\: {\it
Let $R$ be an unbounded conformal algebra and $S\subset R$ a subset.
If any $a, b, c\in S$ satisfy the conformal Jacobi identity then 
so do any $a, b, c\in\K[T]S$.
${}_{}{}_{}{}_{}$\hfill $\square$
}

\bigskip

{\bf Proposition.}\: {\it
An algebra $\fg$ satisfies the Leibniz identity iff
$\fg\pau{z\uppm}$ satisfies the conformal Jacobi identity.
}

\bigskip

\begin{pf}
`$\Rightarrow$'\:
This follows from
\begin{align}
\notag
&\res_{z,w}e^{(z-w)\la} e^{(w-y)\mu}[[a(z),b(w)],c(y)]
\\
\notag
=\;
&\res_{z,w}e^{(z-y)\la} e^{(w-y)(\mu-\la)}
([a(z),[b(w),c(y)]]-[b(w),[a(z),c(y)]]).
\end{align}

\smallskip

`$\Leftarrow$'\:
The map $\fg\to\fg\pau{z\uppm}, a\mapsto az\inv$, is an algebra monomorphism 
since $\res_z [az\inv,bw\inv]=[a,b]w\inv$.
By the Remark $\fg\pau{z\uppm}$ satisfies the Leibniz identity. 
Thus the claim follows.
\end{pf}

\bigskip

The conformal Jacobi identity is equivalent to 
$[a_t,b_s]c=\sum_i\binom{t}{i}(a_i b)_{t+s-i} c$ for any $t, s\geq 0$
because replacing $\mu$ by $\mu+\la$ in the conformal Jacobi identity we obtain
$[[a\subla b]_{\mu+\la}c]\linebreak[0]
=[a\subla[b\submu c]]-\paraab[b\submu[a\subla c]]$
and 
\begin{align}
\notag
&[[a\subla b]_{\mu+\la}c]
\; =\;
\sum_{t,s,i}(a_i b)_s c\: \la^{(i)}\mu^{(s-t)}\la^{(t)}
\\
\notag
&\quad
\; =\;\sum_{t,s,i}
\binom{i+t}{i}(a_i b)_s c\: \la^{(i+t)}\mu^{(s-t)}
\; =\;
\sum_{t,s,i}
\binom{t}{i}(a_i b)_{t+s-i} c\: \la^{(t)}\mu^{(s)}.
\end{align}
In section \ref{SS:bracket vlie} we give a reformulation of the above 
formula in terms of local Lie algebras.
Proposition \ref{SS:j id conf j id} states other identities 
that are equivalent to the conformal Jacobi identity.

\subsection{The Delta Distribution}
\label{SS:delta dis}

We prove three basic properties of the delta distribution
that are used in section \ref{SS:conf skew sym}.

\bigskip

For a vector space $E$, define a subspace of $E\pau{z\uppm,w\uppm}$:
$$
E[(z/w)\uppm]\pau{w\uppm}
\; :=\;
\big\{
\, a(z,w)\:\big| \: \forall\: n\in\Z\;\, 
\sum_{m\in\Z}a_{n-m,m}x^m\in E[x\uppm]\:
\big\}.
$$
If $a(z,w)\in E[(z/w)\uppm]\pau{w\uppm}$ then
$a(w,w):=\sum_{n,m\in\Z} a_{n,m-n-1} w^{-m-1}$ is well-defined. 
For $a(z)\in E\pau{z\uppm}$, define 
$$
a(z+w)
\; :=\;
e^{w\del_z}a(z)
\;\in\;
E\pau{z\uppm}\pau{w}.
$$
Note that $a(z+w)=a(w+z)$ iff $a(z)\in E\pau{z}$.

\bigskip

{\bf Proposition.}\: {\it
\iti\: 
Let $a(z,w)\in E[(z/w)\uppm]\pau{w\uppm}$. 
Then $\de(z,w)a(z,w)$ is well-defined and 
$$
\de(z,w)a(z,w)
\; =\;
\de(z,w)a(w,w).
$$
In particular, $\res_z \de(z,w)a(z,w)=a(w,w)$.

\smallskip

\itii\:
$\del_w \de(z,w)\, =\, -\del_z \de(z,w)$ and $\de(z-x,w)=\de(z,w+x)$.

\smallskip

\itiii\: 
$\de(z,w)=\de(w,z)$.
}

\bigskip

\begin{pf}
\iti\:
The product $\de(z,w)a(z,w)$ is well-defined because
$$
\de(z,w)a(z,w)
\; =\;
\sum_{n,m,k\in\Z}\: a_{n-k-1,m+k}\,z^{-n-1}\,w^{-m-1}.
$$
We have $p(z)\de(z,w)=p(w)\de(z,w)$ for any $p(z)\in\K[z\uppm]$ since
$$
\res_z q(z)p(z)\de(z,w)
\; =\;
q(w)p(w)
\; =\;
\res_z q(z)p(w)\de(z,w).
$$
Thus $\de(z,w)a(z,w)=\de(z,w)a(w,w)$ follows from a continuity argument
since $E[z\uppm]$ is dense is $E\pau{z\uppm}$.

\smallskip

\itii\:
The integration-by-parts formula implies
$$
\res_z p(z)\del_w \de(z,w)
\, =\, 
\del_w p(w)
\, =\, 
\res_z (\del_z p(z))\de(z,w)
\, =\, 
-\res_z p(z)\del_z\de(z,w).
$$
This proves the first identity. 
The second identity follows from the first:
$$
\de(z-x,w)\; =\; e^{-x\del_z}\de(z,w)\; =\; e^{x\del_w}\de(z,w)\; =\; 
\de(z,w+x).
$$

\smallskip

\itiii\:
This follows from the fact that $n\mapsto -n-1$ is an involution of $\Z$.
\end{pf}

\subsection{Conformal Skew-Symmetry}
\label{SS:conf skew sym}

We prove that local distributions with values in a skew-symmetric algebra
satisfy conformal skew-symmetry.

\bigskip

{\bf Conformal skew-symmetry} for  
\index{conformal!skew-symmetry}
an unbounded conformal algebra is
$$
[a\subla b]
\; =\;
-\,\paraab\, [b_{-\la-T}a].
$$
Here $b, a$ are assumed to be weakly local.
Note that $\paraab[b_{-\la-T}a]$ is the opposite $\la$-bracket, 
see section \ref{SS:confa}.
Conformal skew-symmetry is equivalent to
$$
a_r b
\; =\;
\paraab\sum_{i\geq 0} (-1)^{r+1+i}\,T^{(i)}(b_{r+i}a)
$$ 
since
$[b_{-\la-T}a]=\sum_{r, i}(-1)^r T^{(i)}(b_r a)\la^{(r-i)}=
\sum_{r, i}(-1)^{r+i}T^{(i)}(b_{r+i}a)\la^{(r)}$.

\smallskip

Conformal skew-symmetry implies that the identities 
$[(Ta)\subla b]=-\la[a\subla b]$ and $[a\subla Tb]=(T+\la)[a\subla b]$ 
are equivalent.
Setting $\la=0$ in these two identities and in conformal skew-symmetry, 
we obtain:

\bigskip

{\bf Remark.}\: {\it
Let $R$ be a conformal algebra.
Then $TR$ is a two-sided ideal with respect to $[\,,]$.
If $R$ satisfies conformal skew-symmetry 
then $[\,,]$ on $R/TR$ is skew-symmetric.
\hfill $\square$
}

\bigskip

Since $[a\subla b]$ and $a\otimes b\mapsto\paraab[b_{-\la-T}a]$ are 
$\la$-brackets, we obtain:

\bigskip

{\bf Lemma.}\: {\it
Let $R$ be a conformal algebra and $S\subset R$ a subset.
If any $a, b\in S$ satisfy conformal skew-symmetry then 
so do any $a, b\in\K[T]S$.
\hfill $\square$
}

\bigskip

{\bf Proposition.}\: {\it
Let $\fg$ be a skew-symmetric algebra.
If $a(z), b(z)\in\fg\pau{z\uppm}$ are local
then they satisfy conformal skew-symmetry.
}

\bigskip

\begin{pf}
Corollary \ref{SS:loc diff op} and 
Proposition \ref{SS:delta dis}\,\itii, \itiii\ imply
$$
[b(z),a(w)]
\; =\;
-[a(w),b(z)]
\; =\;
-\sum_i\: (-1)^i \, a(z)_i b(z)\del_w^{(i)}\de(z,w).
$$
By Proposition \ref{SS:delta dis diff op} we have
$e^{(z-w)\la}\del_w^{(i)}\de(z,w)=\sum_j \la^{(j)}\del_w^{(i-j)}\de(z,w)=
(\la+\del_w)^{(i)}\de(z,w)$.
Thus we get
$$
[b(w)\subla a(w)]
\; =\;
-\sum_i\: (-1)^i\,\res_z a(z)_i b(z)(\la+\del_w)^{(i)}\de(z,w).
$$
By Proposition \ref{SS:delta dis}\,\iti\ we obtain
$$
[b(w)_{-\la-\del_w}a(w)]
\; =\;
-\sum_i \res_z a(z)_i b(z)\de(z,w)\,\la^{(i)}
\; =\;
-[a(w)\subla b(w)].
$$
\end{pf}

\subsection{Vertex Lie Algebras}
\label{SS:vlie}

We prove that $F\subfg$ generates a vertex Lie subalgebra 
$R(\fg)\subset\fg\pau{z\uppm}$.
In section \ref{SS:vlie into loclie} we describe $R(\fg)$ explicitly 
for some examples.

\bigskip

{\bf Definition.}\: 
A {\bf vertex Lie algebra} is 
\index{vertex Lie algebra}
a conformal algebra that satisfies the conformal Jacobi identity and 
conformal skew-symmetry.

\bigskip

The $\la$-bracket of a vertex Lie algebra is 
\index{Lie!lambdabracket@$\la$-bracket}
called a {\bf Lie $\la$-bracket}.
The zero $\la$-bracket $[\subla]=0$ is a Lie $\la$-bracket. 
A vertex Lie algebra is 
\index{abelian vertex Lie algebra}
{\bf abelian} if $[\subla]=0$.

Remarks \ref{SS:conf jac} and \ref{SS:conf skew sym} yield:

\bigskip

{\bf Remark.}\: {\it
If $R$ is a vertex Lie algebra then $R/TR$ is a Lie algebra.
\hfill $\square$
}

\bigskip

For example, if $R=\hfg$ is an affine vertex Lie algebra then 
$R/TR=\fg\oplus\K\hk$, see section \ref{SS:aff vlie}.

From Propositions \ref{SS:conf jac} and \ref{SS:conf skew sym}
follows that if $\fg$ is a Lie algebra and 
$R$ a local conformal subalgebra of $\fg\pau{z\uppm}$
then $R$ is a vertex Lie algebra.

Conversely, a vertex Lie subalgebra of $\fg\pau{z\uppm}$ need {\it not} be 
local.
For example, $\fg\subset\fg\pau{z\uppm}$ is an abelian vertex Lie subalgebra 
but $\fg$ is not local if $\fg$ is non-abelian.

\bigskip

{\bf Dong's Lemma.}\: {\it
Let $\fg$ be 
\index{Dong's lemma}
a Lie algebra.
If $a(z), b(z), c(z)\in\fg\pau{z\uppm}$ are pairwise local
then $a(z)_i b(z)$ and $c(z)$ are local for any $i\geq 0$.
}

\bigskip

\begin{pf}
Suppose that $a(z), b(z), c(z)$ are pairwise local of order $\leq n$. 
Then
\begin{align}
\notag
&(w-x)^{3n}[[a(z),b(w)],c(x)]
\\
\notag
=\;
&\sum_k\binom{2n}{k} (w-x)^n\, (w-z)^k\, (z-x)^{2n-k}[[a(z),b(w)],c(x)]
\; =\;
0
\end{align}
because the summand for $k\geq n$ is $0$ since $a(z), b(z)$ are local and
the summand for $k\leq n$ is $0$ since $a(z), c(z)$ and $b(z), c(z)$ are
local and we may apply the Leibniz identity.
Applying $\res_z(z-w)^i\cdot$ the claim follows.
\end{pf}

\bigskip

If $R$ is an unbounded conformal algebra then we denote by $\bS\subset R$ 
the unbounded conformal subalgebra generated by a subset $S\subset R$.

If $S\subset\fg\pau{z\uppm}$ is a local subset then $\bS$ is local 
because of Dong's lemma and because locality of $a(z), b(z)$ implies 
locality of $\del_z a(z), b(z)$.
Thus $\bS$ is a vertex Lie algebra.
In particular, if $\fg$ is a local Lie algebra then
$R(\fg):=\bF_{\fg}$ is a vertex Lie algebra.

A {\it morphism} of local Lie algebras is a differential algebra morphism 
$\phi:\fg\to\fg'$ such that $\phi F\subfg\subset\bF_{\fg'}$.
It suffices to require that $\phi$ is an algebra morphism
since $\bF_{\fg'}\subset\fg'\pau{z\uppm}_T$ and hence 
$\phi(Ta_t)=-t\phi(a_{t-1})=T\phi(a_t)$.

We thus obtain a functor $\fg\mapsto R(\fg)$ from local Lie algebras to 
vertex Lie algebras.

If $\fg$ is a graded local Lie algebra then, by definition,
$F\subfg\subset\bigoplus(F\subfg\cap\fg\pau{z\uppm}_h)$.
This implies that $R(\fg)\subset\fg\pau{z\uppm}_H$ is 
a graded vertex Lie algebra.

\section{The Borcherds Lie Algebra}
\label{S:vlie to loclie}

In sections \ref{SS:vlie to loclie}--\ref{SS:loclie of vlie}
we use the weak commutator formula to 
construct a local Lie algebra $\fg(R)$ from a vertex Lie algebra $R$.
We prove that the functor $R\mapsto \fg(R)$ is fully faithful and 
left adjoint to $\fg\mapsto R(\fg)$.

In section \ref{SS:vlie into loclie}
we show that many local Lie algebras from section \ref{S:loclie}
lie in the image of the functor $R\mapsto \fg(R)$.
In section \ref{SS:bracket vlie} 
we prove that any vertex Lie algebra has a natural Lie bracket.

\subsection{The Borcherds Lie Algebra}
\label{SS:vlie to loclie}

We construct a functor $R\mapsto\fg(R)$ from vertex Lie algebras
to local Lie algebras and show that it is left adjoint to $\fg\mapsto R(\fg)$.
The bracket of $\fg(R)$ is determined by the $\la$-bracket of $R$
via the weak commutator formula.
The algebra $\fg(R)$ is
\index{Borcherds Lie algebra} 
the {\bf Borcherds Lie algebra} of $R$.

\bigskip

Let $R$ be a vertex Lie algebra.
There exists a Lie algebra $\fg(R)$ and a morphism $Y: R\to\fg(R)\pau{z\uppm}$
of unbounded conformal algebras such that $Y(R)$ is local and
the following universal property holds:
for any Lie algebra $\fg$ and any morphism $\psi: R\to\fg\pau{z\uppm}$ 
such that $\psi(R)$ is local, there exists a unique algebra morphism 
$\phi: \fg(R)\to\fg$ such that $\psi=\phi\circ Y$.
The pair $(\fg(R),Y)$ is unique up to a unique isomorphism.

Roughly speaking, the universal property says that the functor 
$R\mapsto\fg(R)$ is left adjoint to $\fg\mapsto\fg\pau{z\uppm}$.

We construct $\fg(R)$ as follows.
Let $\fg(R)$ be the Lie algebra generated by the vector space 
$R[x\uppm]$ with relations given by $(Ta)_t=-t a_{t-1}$ and
$$
[a_t,b_s]
\; =\;
\sum_{i\geq 0}\:\binom{t}{i}\: (a_i b)_{t+s-i}
$$
for any $a, b\in R$ and $t, s\in\Z$ where $a_t:=ax^t$.
We also denote by $a_t$ the image of $a_t$ in $\fg(R)$.
Define $Y: R\to\fg(R)\pau{z\uppm}$ by $a\mapsto a(z):=\sum a_t z^{-t-1}$.

By Corollary \ref{SS:loc diff op} a linear map 
$\psi: R\to\fg\pau{z\uppm}$ is a morphism with $\psi(R)$ local iff
$(\psi Ta)_t=-t (\psi a)_{t-1}$ and 
$[(\psi a)_t,(\psi b)_s]=\sum\binom{t}{i}(\psi(a_i b))_{t+s-i}$.
Thus $Y$ is a morphism, $Y(R)$ is local, and 
$\phi: \fg(R)\to\fg, a_t\mapsto (\psi a)_t$, is the unique 
algebra morphism such that $\psi=\phi\circ Y$.

The map $a_t\mapsto (Ta)_t$ induces a derivation $T$ of $\fg(R)$ since
it preserves the relations.
Thus $\fg(R)$ together with $Y(R)$ is a local Lie algebra.
We have $R(\fg(R))=Y(R)$.

\bigskip

{\bf Proposition.}\: {\it
The functor $R\mapsto\fg(R)$ from vertex Lie algebras to local Lie algebras
is left adjoint to $\fg\mapsto R(\fg)$.
The morphism $Y: R\to R(\fg(R))$ is the unit of adjunction.
}

\bigskip

\begin{pf}
Let $R$ be a vertex Lie algebra and $\fg$ a local Lie algebra.
By definition, a local Lie algebra morphism $\fg(R)\to\fg$ is 
an algebra morphism $\phi:\fg(R)\to\fg$ such that $\phi Y(R)\subset R(\fg)$. 
By the universal property of $\fg(R)$,
the map $\phi\mapsto\phi\circ Y$ is a bijection onto the set of
morphisms $\psi: R\to\fg\pau{z\uppm}$ such that 
$\psi(R)$ is local and $\psi(R)=\phi Y(R)\subset R(\fg)$.
Since $R(\fg)$ is local, 
these $\psi$s are just the vertex Lie algebra morphisms $R\to R(\fg)$.
\end{pf}

\bigskip

A gradation of a vertex Lie algebra $R$ induces a gradation of the
local Lie algebra $\fg(R)$ such that $Y: R\to\fg(R)\pau{z\uppm}_H$ 
is a morphism of graded unbounded conformal algebras.
Thus $a_{(t)}\in\fg(R)_{h_a-t-1}$.

\subsection{Modules and the Commutator Formula}
\label{SS:module vlie}

We prove that $R$ and $\fg(R)$ have the same module categories and that
a $\Z$-fold module $M$ is a module iff $M$ satisfies the commutator formula.

\bigskip

{\bf Definition.}\:
A {\bf module} over a conformal algebra $R$ is 
\index{module}
a vector space $M$ together with a morphism $Y_M: R\to\fgl(M)\pau{z\uppm}$ 
such that $Y_M(R)$ is local.

\bigskip

We often just write $a(z)$ instead of $Y_M(a)$.

There is also another notion of module over a conformal algebra,
called conformal module, see section \ref{SS:conf op}.
It is given by a map $R\to\fgl(M)\pau{\mu}$.

The above definition is motivated by two results.
First, bounded modules over a vertex Lie algebra $R$ are the same as
modules over the enveloping vertex algebra $U(R)$, 
see Proposition \ref{SS:vertex envelope}.
Second, the universal property of $\fg(R)$ applied to $\fg=\fgl(M)$
yields the following fact:

\bigskip

{\bf Remark.}\: {\it
Let $R$ be a vertex Lie algebra.
Then the categories of $R$-modules and 
of modules over the Lie algebra $\fg(R)$ are isomorphic.
Moreover, $Y_M=\rho\circ Y$, where $\rho: \fg(R)\to\fgl(M)$.
\hfill $\square$
}

\bigskip

For example, from the adjoint representation $\fg(R)\to\fgl(\fg(R))$
we obtain that $\fg(R)$ is an $R$-module.
The Remark and the existence of the counit of adjunction 
$\fg(R(\fg))\to\fg$ imply that any module over a local Lie algebra $\fg$
is an $R(\fg)$-module.

A {\bf $\Z$-fold module} over a vector space $R$ is 
\index{Zfold@$\Z$-fold!module}
a vector space $M$ together with an even linear map 
$Y_M: R\to\End(M)\pau{z\uppm}, a\mapsto a(z)=\sum a_t z^{-t-1}$.
Equivalently, it is a vector space $M$ together with a family
of even linear maps $V\otimes M\to M, a\otimes b\mapsto a_t b$, 
indexed by $t\in\Z$.
In particular, for $M=V$ this is equivalent to the notion of 
a $\Z$-fold algebra, see section \ref{SS:confa}.

Let $M$ be a $\Z$-fold $R$-module.
We call $a\in R$ and $b\in M$ {\bf weakly local} if $a(z)b\in M\lau{z}$.
We call $M$ {\bf bounded} 
\index{bounded!zfold module@$\Z$-fold module}
if any $a\in R$ and $b\in M$ are weakly local.
In other words, $Y_M(R)\subset\Endv(M)$, see section \ref{SS:alg of distr}.
A special case of this definition are the notions of 
a bounded $R$-module and of a bounded $\Z$-fold algebra. 

A module $M$ over a local Lie algebra $\fg$ is 
\index{bounded!module}
{\bf bounded} if $\rho F\subfg\subset\Endv(M)$ where $\rho: \fg\to\fgl(M)$.
Since $Y_M(R)=\rho Y(R)=\rho F_{\fg(R)}$,
a module $M$ over a vertex Lie algebra $R$ 
is bounded iff the $\fg(R)$-module $M$ is bounded.

Let $M$ be a $\Z$-fold module over a conformal algebra $R$.
The {\bf commutator formula} is
\index{commutator formula}
$$
[a(z),b(w)]
\; =\;
\sum_{i\geq 0}\:
(a_i b)(w)\; \del_w^{(i)}\de(z,w).
$$
This is an identity in $\fgl(M)\pau{z\uppm,w\uppm}$. 
It is equivalent to 
$$
[a_t,b_s]
\; =\;
\sum_{i\geq 0}\:
\binom{t}{i}\,
(a_i b)_{t+s-i}.
$$
We call this the commutator formula for {\bf indices} $t, s$.
It is well-defined for an $\N$-fold algebra $R$ if $t\geq 0$ or 
if $a, b$ are weakly local.

The commutator formula is equivalent to the weak commutator formula
together with the identities $(a_r b)(z)=a(z)_r b(z)$ for $r\geq 0$.
Thus Corollary \ref{SS:loc diff op} yields the following result.

\bigskip

{\bf Proposition.}\: {\it
Let $R$ be a conformal algebra and $M$ a $\Z$-fold $R$-module.
Then $M$ is an $R$-module iff $Y_M(Ta)=\del_z Y_M(a)$ and
$M$ satisfies the commutator formula.
\hfill $\square$
}

\subsection{Base Change and Affinization}
\label{SS:affin vlie}

We discuss the base change functor $R\mapsto C\otimes R$ and obtain
as a special case the affinization $\tR=\K[x\uppm]\otimes R$ 
of a vertex Lie algebra $R$.

\bigskip

Let $C=(C,\del)$ be an even commutative differential algebra.
Recall that there is a functor $\fg\mapsto C\otimes \fg$ from Lie algebras to 
Lie algebras over $C$.

A conformal algebra {\bf over} $C$ is 
\index{conformal algebra!over C@over $C$}
a conformal algebra $R$ with a $C$-module structure such that
$T(fa)=(\del f)a+fTa$ and 
$$
[fa\subla b]
\; =\; 
(e^{\del\del\subla}f)[a\subla b],
\qquad
[a\subla fb]
\; =\; 
f[a\subla b]
$$
for any $a, b\in R$ and $f\in C$. 

We note that if we require instead that $[fa\subla b]=f[a\subla b]$ 
then $[\subla]=0$ for any $C$ such that $1\in\del(C)$.
The factor $e^{\del\del\subla}f$ ensures that both sides of the above
identity have the same transformation property with respect to $T$.

If $R$ is a conformal algebra then $C\otimes R$ is a conformal algebra
over $C$ with $T=\del\otimes 1+1\otimes T$ and
$$
[fa\subla gb]
\; :=\; 
(e^{\del\del\subla}f)g[a\subla b]
$$
because
\begin{align}
\notag
[(\del+T)fa\subla gb]
\; =\; 
(e^{\del\del\subla}\del f)g[a\subla b]-(e^{\del\del\subla}f)g\la[a\subla b]
\; =\; 
-\la[fa\subla gb]
\end{align}
and 
\begin{align}
\notag
[fa\subla(\del+T)gb]
\; =\; 
&
(e^{\del\del\subla}f)\del g[a\subla b]+(e^{\del\del\subla}f)g(T+\la)[a\subla b]
\\
\notag
\; =\; 
&(\del+T+\la)[fa\subla gb].
\end{align}
Both times we used that $[e^{\del\del\subla},\la\cdot]=e^{\del\del\subla}\del$
which follows from $[\del\del\subla,\la\cdot]=\del$ since
$[a^{(i)},b]=a^{(i-1)}[a,b]$ if $[a,[a,b]]=0$. 

The {\bf base change} functor $R\mapsto C\otimes R$ from 
\index{base change}
conformal algebras to conformal algebras over $C$ is left adjoint to 
the forgetful functor.
We do not prove this easy fact here. 

If $R$ satisfies conformal skew-symmetry then so does $C\otimes R$ because
$$
-[gb_{-\la-\del-T}fa]
\; =\; 
-e^{\del\del\subla}((e^{\del\del_{-\la}}g)f)[b_{-\la-T}a]
\; =\; 
[fa\subla gb].
$$
If $R$ satisfies the conformal Jacobi identity then so does $C\otimes R$ 
because
\begin{align}
\notag
[[fa\subla gb]\submu hc]
\; =\;
&(e^{\del\del\submu}((e^{\del\del\subla}f)g))h\,
[[a\subla b]\submu c]
\\
\notag
\; =\;
&(e^{\del(\del\submu+\del\subla)}f)(e^{\del\del\submu}g)h\,
([a\subla[b_{\mu-\la}c]]-[b_{\mu-\la}[a\subla c]])
\\
\notag
=\;
&[fa\subla[gb_{\mu-\la}hc]]-[gb_{\mu-\la}[fa\subla hc]]
\end{align}
where we use that $e^{\del(\del\submu+\del\subla)}p(\la,\mu-\la)=
e^{\del(-\del_{\la'}+\del\subla+\del_{\la'})}p(\la,\mu-\la')|_{\la=\la'}=
e^{\del\del\subla}p(\la,\mu-\la')|_{\la=\la'}$.

We note that if $\del=0$ and $R$ is a vertex Lie algebra 
then $\fg(C\otimes R)=C\otimes\fg(R)$.
The {\bf affinization} of a conformal algebra $R$ is 
\index{affinization}
$\tR:=\K[x\uppm]\otimes R$ with $\del=\del_x$.

\subsection{The Embedding}
\label{SS:loclie of vlie}

We prove that the functor $R\mapsto\fg(R)$ is fully faithful.

\bigskip

If $R$ is a conformal algebra then $\fg_R:=\tR /T\tR$ is an algebra
by Remark \ref{SS:conf skew sym}.
Define $a_t:=\io(x^t a)$ where $\io:\tR\to\fg_R$ is the quotient map.
We have $(Ta)_t=-t a_{t-1}$ since 
$T(x^t a)=tx^{t-1}a+x^tTa$.
Since 
$$
(e^{\del_x\del\subla}f)g[a\subla b]
\; =\;
\sum_{t,i\geq 0}
(\del_x^{(i)}f)g\: a_{t+i}b\: \la^{(t)},
$$
the bracket of $\fg_R$ is
$$
[a_t,b_s]
\; =\;
\io\sum_{i\geq 0}(\del_x^{(i)}x^t)x^s\: a_i b
\; =\;
\sum_{i\geq 0}\binom{t}{i}(a_i b)_{t+s-i}.
$$
It follows that $\fg_R$ is the algebra generated by the vector space 
$R[x\uppm]$ with relations $(Ta)_t=-t a_{t-1}$ and
$[a_t,b_s]=\sum\binom{t}{i}(a_i b)_{t+s-i}$.
Thus there exists a unique 
algebra epimorphism $\io:\fg_R\to\fg(R)$ such that $a_t\mapsto a_t$.

If $R$ is a vertex Lie algebra then $\fg_R$ is 
a Lie algebra by Remark \ref{SS:vlie}.
In this case $\io$ is an isomorphism.

\bigskip

{\bf Remark.}\: {\it
If $R$ is a conformal algebra then $R\to\fg_R, a\mapsto a_{-1}$, is injective. 
}

\bigskip

\begin{pf}
The map $p:\tR\to R, x^t a\mapsto T^{(-1-t)}(a)$, induces a map 
$p:\fg_R\to R$ because 
$pT(x^t a)=p(tx^{t-1}a+x^tTa)=tT^{(-t)}(a)+T^{(-1-t)}(Ta)=0$.
From $p(a_{-1})=a$ follows that $a\mapsto a_{-1}$ is injective. 
\end{pf}

\bigskip

{\bf Proposition.}\: {\it
The functor $R\mapsto\fg(R)$ from vertex Lie algebras to local Lie algebras
is fully faithful.
}

\bigskip

\begin{pf}
The Remark implies that $Y: R\to R(\fg(R))$ is an isomorphism.
By Proposition \ref{SS:vlie to loclie}
the map $Y$ is the unit of adjunction
of the pair of functors $R\mapsto\fg(R)$ and $\fg\mapsto R(\fg)$.
It is a general fact that if the unit of adjunction is an isomorphism
then the left adjoint functor is fully faithful.
\end{pf}

\subsection{The Image}
\label{SS:vlie into loclie}

We show that many of the local Lie algebras from section \ref{S:loclie}
are vertex Lie algebras, i.e.~they are contained in the image of the 
embedding $R\mapsto\fg(R)$.

\bigskip

A local Lie algebra $\fg$ is {\bf regular} 
\index{regular local Lie algebra}
if $\fg$ is isomorphic to $\fg(R)$ for some vertex Lie algebra $R$.
In other words, the category of regular local Lie algebras is 
the essential image of the fully faithful functor $R\mapsto\fg(R)$ and is thus 
equivalent to the category of vertex Lie algebras.
We identify the category of vertex Lie algebras with its essential image in 
the category of local Lie algebras.

The counit of adjunction $\fg(R(\fg))\to\fg$ is given by 
$a(z)_t\mapsto a_t$ where 
$a(z)\in R(\fg)\subset\fg\pau{z\uppm}$. 
This is an epimorphism.
Hence $\fg\mapsto R(\fg)$ is faithful.
A local Lie algebra $\fg$ is regular iff the counit of 
adjunction $\fg(R(\fg))\to\fg$ is an isomorphism.

Let $\fg$ be a local Lie algebra and $K:=F\subfg\cap\fg$ the set 
of constant distributions in $F\subfg$.
Since $[\hk,a(z)]$ is local we get $[\hk,a(z)]=0$
for any $\hk\in K$ and $a(z)\in F\subfg$.
Thus $K$ is contained in the centre of $\fg$.

\bigskip

{\bf Proposition.}\: {\it
Let $\fg$ be a local Lie algebra such that $\bF\subfg=\K[\del_z]F\subfg$.
If $\hk\in K:=F\subfg\cap\fg$ and $a_t$ for 
$a(z)\in F\subfg\setminus K, t\in\Z$ form a basis of $\fg$ then 
$\fg$ is regular and the map
$$
\K^{\oplus K}\oplus\K[\del_z]^{\oplus F\subfg\setminus K}\; \to\; 
\bF\subfg,
\quad
p\mapsto \sum_{a\in F\subfg}\: p_a a,
$$
is an isomorphism. 
}

\bigskip

\begin{pf}
We have $-t\hk_{t-1}=(\del_z \hk)_t=0\in\fg(R(\fg))$ and hence
$\hk_t=0\in\fg(R(\fg))$ for $t\ne -1$.
From $(\del_z a(z))_t=-t a(z)_t$ thus follows that 
$\hk_{-1}$ for $\hk\in K$ and $a(z)_t$ for $a(z)\in F\subfg\setminus K, 
t\in\Z$ form a spanning set of $\fg(R(\fg))$. 
Hence the counit of adjunction $\fg(R(\fg))\to\fg, a(z)_t\mapsto a_t$,
is an isomorphism since it maps a spanning set to a basis.

The map $p\mapsto\sum p_a a$ is obviously surjective.
It is injective since its composition with the map $\bF\subfg\to\fg, 
a(z)\mapsto a_{-1}$, is injective.
\end{pf}

\bigskip

Note that $\bF\subfg=\K[\del_z]F\subfg$ iff $a(z)_t b(z)\in\K[\del_z]F\subfg$ 
for any $a(z), b(z)\in F\subfg$.

By section \ref{SS:weak comm f ex} the Proposition applies to 
the local Lie algebras $\tfg, \hfg, \linebreak[0] C(B), 
\Witt, \Vir, \oneVir$, and $\twoVir$.

\subsection{The Lie Bracket of a Vertex Lie Algebra}
\label{SS:bracket vlie}

We prove that any vertex Lie algebra has a natural Lie bracket 
$[\, ,]\subslie$.

\bigskip

For a vector space $E$, define the formal integral
$$
\int\subal\upbe d\la :\;
E\pau{\la}\to E\pau{\al,\be},
\quad
\sum a_t\la^{(t)}\mapsto\sum a_t(\be^{(t+1)}-\al^{(t+1)}).
$$
A calculation shows that the integral satisfies the {\bf substitution formula} 
\index{substitution formula}
$$
\int\subal\upbe d\la\, p(\la)
\; =\; 
\int_{(\al-\ka)/c}^{(\be-\ka)/c} d\la\; c\, p(c\la+\ka),
$$
for $c\in\K\uptimes$, and the {\bf Fubini formula}
\index{Fubini formula}
$$
\int\subal\upbe d\la\int\subal\upla d\mu\:p(\la,\mu)
\; =\; 
\int\subal\upbe d\mu\int\submu\upbe d\la\:p(\la,\mu).
$$

For a conformal algebra $R$, define 
$$
[a,b]\subslie
\; :=\;
\int_{-T}^0 d\la\: [a\subla b]
\; =\; 
\sum_{i\geq 0}
(-1)^i\, T^{(i+1)}(a_i b).
$$  
Then $R$ with $[\, ,]\subslie$ is a differential algebra, that we 
\index{Rlie@$R\subsmlie$}
denote by $R\subsmlie$.
The algebra $R\subsmlie$ should not be confused with the conformal algebra
$V\subslie$ that is defined for any vertex algebra $V$, 
see section \ref{SS:va}.

For a local Lie algebra $\fg$, define 
$\fg_+:=\rspan\set{a_t\mid a(z)\in F\subfg, t\geq 0}$ and
$\fg_-:=\rspan\set{a_t\mid a(z)\in F\subfg, t<0}$.
These are $\K[T]$-submodules.

If $R$ is a vertex Lie algebra then $\fg(R)_+$ and $\fg(R)_-$ are 
subalgebras of $\fg(R)$ and $\fg(R)=\fg(R)_-\oplus\fg(R)_+$
since $\fg(R)=\tR/T\tR$ and 
$R\otimes\K[x]$ and $R\otimes x\inv\K[x\inv]$ are $\K[T]$-submodules of $\tR$.

The formula proved at the end of section \ref{SS:conf jac} shows that 
there is an algebra morphism $\fg(R)_+\to\fgl(R), a_t\mapsto a_t$.

\bigskip

{\bf Proposition.}\: {\it
Let $R$ be a vertex Lie algebra.
The map $R\subsmlie\to\fg(R)_-, a\mapsto a_{-1}$, is 
an isomorphism of differential algebras. 
In particular, $[\, ,]\subslie$ is a Lie bracket.
}

\bigskip

\begin{pf}
The map is injective by Remark \ref{SS:loclie of vlie}.
It is surjective since $a_{-1-t}=(T^{(t)}a)_{-1}$ for $t\geq 0$.
It is an algebra morphism because
$$
([a,b]\subslie)_{-1}
=
\sum_i (-1)^i\, T^{(i+1)}(a_i b)_{-1}
=
\sum_i \binom{-1}{i}(a_i b)_{-2-i}
=
[a_{-1},b_{-1}].
$$
\end{pf}

\bigskip

The Proposition and Remark \ref{SS:module vlie} imply that
any $R$-module $M$ is an $R\subsmlie$-module since we have a morphism
$R\subsmlie\to\fg(R)\to\fgl(M), a\mapsto a_{-1}$.

\section{Examples of Vertex Lie Algebras}
\label{S:ex vlie}

In sections \ref{SS:unb confa power series}--\ref{SS:conf der}
we discuss vertex Lie algebras of power series, conformal operators, and 
conformal derivations.
In section \ref{SS:free vlie} we construct free vertex Lie algebras.

In section \ref{SS:finite confa} 
we explain the general structure of finite vertex Lie algebras.
Section \ref{SS:aff vlie}
is about affine and Clifford vertex Lie algebras.

In sections \ref{SS:witt vlie}--\ref{SS:n1 n2 top vir}
we discuss examples of conformal vertex Lie algebras:
$\Witt, \Vir$, the GKO coset and the Chodos-Thorn construction, 
$\oneVir, \twoVir$, and the topological Virasoro algebra.
Moreover, we discuss conformal and primary vectors, the Griess algebra, 
and $U(1)$-currents.

\subsection{Unbounded Vertex Lie Algebras of Power Series}
\label{SS:unb confa power series}

We prove that $\fg\pau{\mu}$ is an unbounded vertex Lie algebra 
iff $\fg$ is a Lie algebra.

\bigskip

Let $\fg$ be an algebra.
Then $\fg\pau{\mu}$ is an unbounded conformal algebra
with $T=-\mu\cdot$ and 
$$
[(a\submu)\subla b\submu]
\; :=\;
[a\subla,b_{\mu-\la}]
$$ 
since 
$[(-\mu a\submu)\subla b\submu]=-\la[a\subla,b_{\mu-\la}]$ and
$[(a\submu)\subla(-\mu b\submu)]=(\la-\mu)[a\subla,b_{\mu-\la}]$.

There is a morphism of unbounded conformal algebras
$$
\fg\pau{z\uppm}\;\to\;\fg\pau{\mu},
\quad
a(z)\;\mapsto\; a\submu\; :=\; \res_z e^{z\mu}a(z)\; =\; \sum a_t\, \mu^{(t)}
$$
because the integration-by-parts formula implies 
$\del_z a(z)\mapsto -\mu a\submu$ and 
$$
\res_w e^{w\mu}\res_z e^{\la(z-w)}[a(z),b(w)]
\; =\;
[\res_z e^{\la z}a(z),\res_w e^{(\mu-\la)w}b(w)].
$$
Thus there is a short exact sequence
$0\to\fg\pau{z}\to\fg\pau{z\uppm}\to\fg\pau{\mu}\to 0$.

We often identify $z\inv\fg\pau{z\inv}$ with $\fg\pau{\mu}$
using the $\K[\del_z]$-module isomorphism 
$z\inv\fg\pau{z\inv}\to\fg\pau{\mu}, a(z)\mapsto a\submu$.
The map $a(z)\mapsto a\submu$ can be viewed as a Fourier transformation.
We have $\del_z\leftrightarrow-\mu\cdot$ and $z\cdot\leftrightarrow\del\submu$.

In general, conformal skew-symmetry is only well-defined for 
weakly local elements.
Yet conformal skew-symmetry for $\fg\pau{\mu}$ is always well-defined because
$$
[(b\submu)_{-\la-T} a\submu]
\; =\;
e^{T\del\subla}[(b\submu)_{-\la} a\submu]
\; =\;
e^{-\mu\del\subla}[(b\submu)_{-\la} a\submu]
$$
and $\sum p_n\mu^{(n)}$ is well-defined for any $p_n\in\fg\pau{\la,\mu}$.
Here we used the fact that $a(z+w)=e^{w\del_z}a(z)$ for $a(z)\in\fg\pau{z}$
since $e^{w\del_z}z^n=(z+w)^n$ for $n\geq 0$.

An {\bf unbounded vertex Lie algebra} is an unbounded conformal algebra 
\index{unbounded!vertex Lie algebra}
for which the conformal Jacobi identity and conformal skew-symmetry are 
well-defined and satisfied.
``Well-defined'' can be made rigorous by requiring that the infinite sum
$b_{-\la-T}a=e^{T\del\subla}b_{-\la}a$ converges in some topology.

\bigskip

{\bf Proposition.}\: {\it
An algebra $\fg$ is a Lie algebra iff 
$\fg\pau{\mu}$ is an unbounded vertex Lie algebra.
More precisely: 

\smallskip

\iti\:
$\fg$ is skew-symmetric iff $\fg\pau{\mu}$ satisfies conformal skew-symmetry.

\smallskip

\itii\:
$\fg$ satisfies the Leibniz identity 
iff $\fg\pau{\mu}$ satisfies the conformal Jacobi identity.
}

\bigskip

\begin{pf}
`$\Rightarrow$'\:
\iti\:
This follows from
$$
[(b\submu)_{-\la-T} a\submu]
\; =\;
e^{T\del\subla}[(b\submu)_{-\la} a\submu]
\; =\;
e^{-\mu\del\subla}[b_{-\la}, a_{\mu+\la}]
\; =\;
[b_{\mu-\la}, a\subla].
$$

\smallskip

\itii\:
This follows from $[[(a\subka)\subla b\subka]\submu c\subka]=
[[a\subla,b_{\ka-\la}]\submu c\subka]=[[a\subla,b_{\mu-\la}],c_{\ka-\mu}]$,
$$
[(a\subka)\subla[(b\subka)_{\mu-\la}c\subka]]
\; =\;
[(a\subka)\subla[b_{\mu-\la},c_{\ka-\mu+\la}]]
\; =\;
[a\subla,[b_{\mu-\la},c_{\ka-\mu}]],
$$
and $[(b\subka)_{\mu-\la}[(a\subka)\subla c\subka]]=
[(b\subka)_{\mu-\la}[a\subla,c_{\ka-\la}]]=
[b_{\mu-\la},[a\subla,c_{\ka-\mu}]]$.

\smallskip

`$\Leftarrow$'\:
This follows from Remarks \ref{SS:conf jac} and \ref{SS:conf skew sym} 
and the fact that there exists 
an algebra isomorphism $\fg\to\fg[\mu]/T\fg[\mu]$.
\end{pf}

\subsection{Conformal Operators}
\label{SS:conf op}

We discuss the general linear unbounded vertex Lie algebra $\fgl\subsc(R)$.

\bigskip

Let $\fg$ be a differential Lie algebra. 
The {\bf translation covariant} $\fg$-valued power series 
\index{translation covariant}
form an unbounded vertex Lie algebra
$\fg\pau{\mu}_T:=\fg\pau{\mu}\cap\fg\pau{z\uppm}_T$.
By definition, $a\submu\in\fg\pau{\mu}_T$ iff $Ta\submu=-\mu a\submu$,
see sections \ref{SS:distr} and \ref{SS:unb confa power series}.
Note that $\fg\pau{\mu}_T\cap\fg[\mu]=0$.
The projection $a(z)\mapsto a\submu$ maps $\fg\pau{z\uppm}_T$ to 
$\fg\pau{\mu}_T$.

If $R$ is a $\K[T]$-module
then $\fgl(R)$ is a differential Lie algebra with $T=[T,\; ]$.
To give an unbounded $\la$-product on $R$ is equivalent to giving 
a $\K[T]$-module morphism $R\to\fgl(R)\pau{\mu}_T, a\mapsto a\submu$. 
The confomal Jacobi identity holds for $R$ iff $a\mapsto a\submu$ is 
a morphism of unbounded conformal algebras.

Recall that the space $\Endv(R)$ of fields consists of 
$a(z)\in\End(R)\pau{z\uppm}$ such that $a(z)b\in R\lau{z}$ for any $b\in R$, 
see section \ref{SS:alg of distr}.

The {\bf general linear} unbounded 
\index{general linear vertex Lie algebra}
vertex Lie algebra is 
\index{gaaaglcr@$\fgl\subsc(R)$}
$$
\fgl\subsc(R)
\; :=\;
\fgl(R)\pau{\mu}_T\,\cap\,\Endv(R).
$$
Thus $a\submu\in\fgl\subsc(R)$ iff $[T,a\submu]=-\mu a\submu$ and 
$a\submu b\in R[\mu]$ for any $b\in R$.
The elements of $\fgl\subsc(R)$ are 
\index{conformal operator}
called {\bf conformal operators}.

To give a $\la$-product on $R$ is equivalent to giving a 
$\K[T]$-module morphism $R\to \fgl\subsc(R)$. 
A {\bf conformal module} over a conformal algebra $R$ is 
\index{conformal!module}
a $\K[T]$-module $M$ with a morphism $R\to\fgl\subsc(M)$.

We remark that $\fgl\subsc(R)$ is in general {\it not} a 
finitely generated $\K[T]$-module, even if $R$ is.
For example, if $R=E[T]$ is a free $\K[T]$-module then 
$\fgl\subsc(R)\to\Hom(E,R[\mu]), a\submu\mapsto\al=a\submu|_E$,
is a $\K[T]$-module isomorphism where $T\al=-\mu\al$.
Of course, $\Hom(E,R[\mu])$ is isomorphic to $\Hom(E,R)[T]$.

\bigskip

{\bf Proposition.}\: {\it
\iti\: 
If $R$ is a finitely generated $\K[T]$-module then 
$\fgl\subsc(R)$ is a vertex Lie algebra.

\smallskip

\itii\:
The $\K[T]$-modules $\fg\pau{\mu}$ and $\fgl\subsc(R)$ are torsionfree.
In particular, if $R$ is an unbounded conformal algebra and 
$c\in R$ is torsion then $c\submu=0$.

\smallskip

\itiii\:
If $R$ is a $\K[T]$-module, $a\submu\in\fgl\subsc(R)$, and $c$ is torsion 
then $a\submu c=0$.
}

\bigskip

\begin{pf}
\iti\:
We have to show that $(a\submu)\subla b\submu\in\fgl\subsc(R)[\la]$
for any $a\submu, b\submu\in\fgl\subsc(R)$.
This follows from $((a\submu)\subla b\submu)c\in R[\mu,\la]$ and 
$((a\submu)\subla b\submu)Tc=(T+\mu)((a\submu)\subla b\submu)c$ 
for any $c\in R$.

\smallskip

\itii\:
This is clear.

\smallskip

\itiii\:
Let $p(T)\in\K[T], p\ne 0$, such that $p(T)c=0$.
Considering the leading coefficient of $p(T+\mu)(a\submu c)=a\submu p(T)c=0$
we get $a\submu c=0$.
\end{pf}

\subsection{Conformal Derivations}
\label{SS:conf der}

We discuss the unbounded vertex Lie algebra $\Der\subsc(R)$ of conformal
derivations.

\bigskip

A {\bf conformal derivation} of a conformal algebra $R$ is 
\index{conformal derivation!of a conformal algebra}
a conformal operator $d\submu$ such that
$$
d\submu (a\subla b)
\; =\;
(d\submu a)_{\la+\mu}b
\; +\;
\zeta^{da}\, a\subla(d\submu b).
$$
This is equivalent to $[d\submu, a\subla]=(d\submu a)_{\la+\mu}$. 

The conformal Jacobi identity holds iff 
$[a\submu\;]$ is a conformal derivation for any $a$.
Conformal derivations of the form $[a\submu\;]$ are 
\index{conformal derivation!inner}
called {\bf inner}.

Setting $\mu=0$ in the above identity and in $[T,d\submu]=-\mu d\submu$,
we obtain:

\bigskip

{\bf Remark.}\: {\it
Let $R$ be a conformal algebra and
$d\submu$ a conformal derivation of $R$. 
Then $d_0$ is a derivation of $R$.
\hfill $\square$
}

\bigskip

The space $\Der\subsc(R)$ of conformal derivations is 
an unbounded vertex Lie subalgebra of $\fgl\subsc(R)$ 
because $T=-\mu\cdot$ and
\begin{align}
\notag
[[(d\submu)\subla d'\submu], a\subnu]
\; =\; 
&[[d\subla, d'_{\mu-\la}], a\subnu]
\; =\; 
[d\subla,[d'_{\mu-\la},a\subnu]]-[d'_{\mu-\la},[d\subla,a\subnu]]
\\
\notag
\; =\; 
&(d\subla(d'_{\mu-\la}a))_{\nu+\mu}-(d'_{\mu-\la}(d\subla a))_{\nu+\mu}
\; =\; 
([(d\submu)\subla d'\submu]a)_{\nu+\mu}.
\end{align}

\subsection{Free Vertex Lie Algebras}
\label{SS:free vlie}

We construct the free vertex Lie algebra $R(S,o)$ generated by 
\index{free vertex Lie algebra}
a set $S$ with locality function $o$.
In section \ref{SS:free vlie rev}
we give a second construction of $R(S,o)$.

\bigskip

A {\bf locality function} on a set $S$ is 
\index{locality function}
a map $o:S^2\to\Z\cup\set{-\infty}$.
A {\it morphism} of sets with locality function
is an even map $\phi:S\to S'$ such that $o(\phi a,\phi b)\leq o(a,b)$.
A locality function $o$ is {\bf non-negative}
\index{locality function!non-negative}
if $o(S^2)\subset\N$.

For elements $a, b$ of an unbounded conformal algebra,
define $o_+(a,b)$ to be the least $n\geq 0$ such that $a_m b=0$ for $m\geq n$.
If $R$ is a conformal algebra then $o_+:R^2\to\N$ is 
\index{locality function!of a conformal algebra}
the {\bf locality function} of $R$.

Let $\fg$ be an algebra and $a(z), b(z)\in\fg\pau{z\uppm}$.
Define $N_+(a(z),b(z))$ to be the least $n\geq 0$ such that
$a(z), b(z)$ are local of order $\leq n$.
It is clear that $o_+(a(z),b(z))\leq N_+(a(z),b(z))$. 
Proposition \ref{SS:delta dis diff op} implies that
if $a(z), b(z)$ are local then $N_+(a(z),b(z))=o_+(a(z),b(z))$.
In particular, if $R$ is a vertex Lie algebra then $o_+(a,b)=N_+(Ya,Yb)$.

\bigskip

{\bf Proposition.}\: {\it
The functor $R\mapsto (R,o_+)$ from vertex Lie algebras to sets with 
non-negative locality function has a left adjoint $(S,o)\mapsto R(S,o)$. 
}

\bigskip

\begin{pf}
Let $F$ be the free $\N$-fold algebra with derivation generated by the set $S$.
A basis of $F$ is given by $B=\bigsqcup B_n$ where $B_1:=\set{T^k a\mid a\in S,
k\geq 0}$ and 
$B_n:=\set{a_i b\mid a\in B_m, b\in B_{n-m}, 1\leq m<n, i\geq 0}$
for $n\geq 2$.

Define a map $o: B^2\to\N$ by induction on $n$: $o(T^k a,T^l b):=o(a,b)+k+l$
for $a, b\in S$, $o(T^k a,b_i c):=3\max(o(T^k a,b),o(T^k a,c),o(b,c))$
for $a\in S, b_i c\in B_n$, and $o(a_i b,c):=3\max(o(a,b),o(a,c),o(b,c))$
for $a_i b\in B_n, c\in B$.
Let $F'$ be the quotient of $F$ by the relations $a_i b=0$ and 
$(Ta)_t b=-t a_{t-1}b$ for $a, b\in B, i\geq o(a,b)$, and $t\geq 0$.
Then $F'$ is a conformal algebra.
Let $R(S,o)$ be the quotient of $F'$ by the conformal Jacobi identity and
conformal skew-symmetry.

Let $R$ be a vertex Lie algebra with a morphism $(S,o)\to (R,o_+)$.
Since $F$ is free, 
there exists a unique morphism $F\to R'$ compatible with $S$.

We have $o_+(a,b)=N_+(Ya,Yb)$ for $a, b\in R$.
The proof of Dong's lemma shows that
$N_+(a_i b,c)\leq 3\max(N_+(a,b),N_+(a,c),N_+(b,c))$ 
for any $a, b, c\in Y(R)$.
Moreover, $N_+(\del_z a,b)\leq N_+(a,b)+1$.
Thus $F\to R$ induces a morphism $F'\to R$.
Since $R$ is a vertex Lie algebra, the map $F'\to R$ induces a morphism
$R(S,o)\to R$.
\end{pf}

\bigskip

Since the composition of left adjoint functors is left adjoint,
the Proposition implies that the functor $\fg\mapsto (R(\fg),o_+)$ 
from local Lie algebras to sets with non-negative locality function has 
a left adjoint $(S,o)\mapsto\fg(R(S,o))$. 

Let $\cR\subset R(S,o)$ be a subset and $I$ the ideal generated by $\cR$.
Then $R(S,o)/I$ is the vertex Lie algebra 
\index{relations!for a vertex Lie algebra}
generated by $(S,o)$ with {\bf relations} $\cR$.

Suppose that $S$ is a $\K$-graded set, $S=\bigsqcup S_h$.
In the following we use the algebras $F, F', R(S,o)$ defined in
the proof of the Proposition. 
The gradation of $S$ induces a gradation of the free $\N$-fold algebra $F$.
Thus $F'$ is also graded.
Since the conformal Jacobi identity and conformal skew-symmetry are 
homogeneous identities, 
it follows that $R(S,o)$ is a graded vertex Lie algebra.

\subsection{Finite Vertex Lie Algebras and Cosets}
\label{SS:finite confa}

We explain the general structure of finite conformal algebras and
present finite vertex Lie algebras in terms of generators and relations.

\bigskip

A conformal
\index{finite conformal algebra}
algebra is {\bf finite} if it is finitely generated as a $\K[T]$-module.

The polynomial ring $\K[T]$ is a principal ideal domain.
A basic theorem of algebra says that any finitely generated module over a 
principal ideal domain is the direct sum of the torsion submodule
and a free module of finite rank. 
We first consider the torsion submodule of a conformal algebra.

Elements $a, b$ of a 
\index{commuting elements}
conformal algebra $R$ {\bf commute} if $[a\subla b]=[b\subla a]=0$.
Note that $a$ need {\it not} commute with $a$, even if $R$ is a vertex
Lie algebra.

The {\bf centralizer} or {\bf coset} of 
\index{centralizer}
a subset $S\subset R$ is
\index{coset}
$$
C_R(S)
\; :=\;
\set{a\in R\mid a\;\text{commutes with any}\;b\in S}.
$$
The {\bf centre} of $R$ 
\index{centre!of a conformal algebra}
is $C(R):=C_R(R)$.
The elements of the centre are 
\index{central element!of a conformal algebra}
called {\bf central}.
If $R$ is a vertex Lie algebra then the coset $C_R(S)$ is a 
vertex Lie subalgebra and $C_R(S)=C_R(\bS)$.
Proposition \ref{SS:conf op}\,\itiii\ implies:

\bigskip

{\bf Proposition.}\: {\it
The torsion submodule of a conformal algebra is contained in the centre.
\hfill $\square$
}

\bigskip

The Proposition implies that any finite-dimensional vertex Lie algebra is 
abelian.

Let $R=F\oplus K$ be a $\K[T]$-module, 
where $F$ is a free and $K$ a torsion $\K[T]$-module.
Let $E\subset F$ be a subspace such that $F=E[\mu]$ and $T|_F=-\mu\cdot$.

Let $p\subla: E\otimes E\to R[\la]$ be an even linear map.
It induces a $\K[\mu,\ka]$-module morphism 
$p\subla: E[\mu]\otimes E[\ka]\to R[\mu,\ka,\la]$.
The $\K[T]$-module $R$ is a conformal algebra with $\la$-product
$$
(a\submu+c)\subla(b\submu+d)
\; :=\;
p\subla(a\subla,b_{\mu-\la}),
$$
for $a\submu, b\submu\in E[\mu]$ and $c, d\in K$, since
$(-\mu a\submu)\subla b\submu=-\la p\subla(a\subla,b_{\mu-\la})$ and
$(a\submu)\subla(-\mu b\submu)\linebreak[0]
=(-\mu+\la)p\subla(a\subla,b_{\mu-\la})$.
The Proposition implies that this $\la$-product is the
unique extension of $p\subla$ from $E$ to $R$.

Lemmata \ref{SS:conf jac} and \ref{SS:conf skew sym} show that 
if conformal skew-symmetry and the conformal Jacobi identity are 
satisfied for the elements of $E$ then $R$ is a vertex Lie algebra.
The verification of the conformal Jacobi identity can be simplified
by using $\bbS_3$-symmetry, see Proposition \ref{SS:s3 symm conf J}. 
We often define vertex Lie algebras in terms of $E, K$, and the map
$E\otimes E\to R[\la]$.

Assume that $E=\bigoplus E_h$ is $\K$-graded.
The $\K[T]$-module $R$ has a unique $\K$-gradation such that 
$E_h\subset R_h$ and $K\subset R_0$.
Of course, $R$ is graded as a conformal algebra iff
$(E_h)_{(t)}(E_k)\subset R_{h+k-t-1}$.

\bigskip

{\bf Remark.}\: {\it
Let $R=E[\mu]\oplus K$ be a vertex Lie algebra such that 
$T|_{E[\mu]}=-\mu\cdot$ and $T|_K=0$.
Denote by $\io: R\to E[T]\oplus K$ the natural $\K[T]$-module isomorphism.
Then $R$ is the vertex Lie algebra generated by the vector space $E\oplus K$
with relations $a\subla b=\io(a\subla b)$ for $a, b\in E$ and 
$T\hk=0$ for $\hk\in K$. 
}

\bigskip

\begin{pf}
Let $\phi$ be an even linear map from $E\oplus K$ to a 
vertex Lie algebra $R'$ such that $(\phi a)\subla(\phi b)=\phi\io(a\subla b)$
and $T\phi\hk=0$.
Then $\phi$ can be extended uniquely to a $\K[T]$-module morphism 
$\phi: R\to R'$.
This is a vertex Lie algebra morphism since 
$\phi(a\subla b)=(\phi a)\subla(\phi b)$ for any $a, b\in E\oplus K$.
\end{pf}

\subsection{Affine and Clifford Vertex Lie Algebras}
\label{SS:aff vlie}

We show that any vertex Lie algebra $R$ of CFT-type 
contains an image of the affine vertex Lie algebra $\hR_1$.

\bigskip

Let $\fg$ be a Lie algebra.
The {\bf loop vertex Lie algebra} of 
\index{loop!vertex Lie algebra}
$\fg$ is $\tfg:=\fg[\mu]$ with $a\subla b:=[a,b]$ for $a, b\in\fg$.
This is a graded vertex Lie algebra with $\tfg_1=\fg$.
It is a subalgebra of $\fg\pau{\mu}$.

Endow $\fg$ with an invariant, symmetric bilinear form.
The {\bf affine vertex Lie algebra} of 
\index{affine!vertex Lie algebra}
$\fg$ is the central extension $\hfg=\tfg\oplus\K\hk$ defined by $T\hk=0$ and
$$
a\subla b
\; :=\;
[a,b]
\; +\;
(a,b)\hk\,\la
$$
for $a, b\in\fg$.
This is a graded vertex Lie algebra since, omitting supersigns, 
the conformal Jacobi identity 
for $a, b, c\in\fg$ is
$$
[[a,b],c]+([a,b],c)\hk\mu
=
[a,[b,c]]+(a,[b,c])\hk\la
-
[b,[a,c]]-(b,[a,c])\hk(\mu-\la)
$$
and conformal skew-symmetry is $[a,b]+(a,b)\hk\la=-[b,a]+(b,a)\hk\la$.

Let $R$ be a vertex Lie algebra of CFT-type, see section \ref{SS:confa}.
Let $\hk\in R_0, \hk\ne 0$.
Then $R_1$ with $[a,b]:=a_0 b$ is an algebra with a bilinear form defined by 
$(a,b)\hk=a_1 b$.

\bigskip

{\bf Proposition.}\: {\it
Let $R$ be a vertex Lie algebra of CFT-type.
Then $R_1$ is a Lie algebra with an invariant, symmetric bilinear form.
Moreover, there is a natural epimorphism $\hR_1\to\K[T]R_1\oplus R_0 $
of graded vertex Lie algebras.
}

\bigskip

\begin{pf}
We have $a\subla b=[a,b]+(a,b)\hk\la$ for $a, b\in R_1$.
The above proof that $\hfg$ is a vertex Lie algebra
also shows that $R_1$ is a Lie algebra and 
$(\, ,)$ is invariant and symmetric.
The second claim is clear.
\end{pf}

\bigskip

Let $B$ be a vector space with a symmetric bilinear form.
We consider $B$ as an abelian Lie algebra.
If $B=B\even$ then $\hB$ is 
\index{Heisenberg!vertex Lie algebra}
called a {\bf Heisenberg vertex Lie algebra} and if $B=B\odd$ then $\hB$ is 
\index{symplectic fermion!vertex Lie algebra}
called a {\bf symplectic fermion vertex Lie algebra}.

The {\bf Clifford vertex Lie algebra} $C(B):=(\rPi B)[\mu]\oplus\K\hk$ is 
\index{Clifford!vertex Lie algebra}
defined by $T\hk=0$ and 
$$
a\subla b
\; :=\;
(a,b)\hk
$$
for $a, b\in\rPi B$.
Here $\rPi$ is the parity-change functor.
The Clifford vertex Lie algebra is graded with $C(B)_{1/2}=\rPi B$.

By Proposition \ref{SS:vlie into loclie} and section \ref{SS:weak comm f ex}
we have $R(\tfg)=\tfg, R(\hfg)=\hfg$, and $R(C(B))=C(B)$.
The following result is proven in the same way as the Proposition.

\bigskip

{\bf Remark.}\: {\it
Let $R$ be a vertex Lie algebra of CFT-type.
Then $\rPi R_{1/2}$ is a vector space with a symmetric bilinear form
defined by $(a,b)\hk=a_{1/2}b$.
Moreover, there is a natural epimorphism 
$C(\rPi R_{1/2})\to\K[T]R_{1/2}\oplus R_0$ of graded vertex Lie algebras.
\hfill $\square$
}

\subsection{Witt Vertex Lie Algebra}
\label{SS:witt vlie}

We prove that $\Witt$ is the unique non-abelian vertex Lie algebra of rank $1$.

\bigskip

The {\bf Witt vertex Lie algebra} $\Witt:=\K[T]L$ is defined by
$$
L\subla L
\; :=\;
(T+2\la)L.
$$
By Proposition \ref{SS:vlie into loclie} and section \ref{SS:weak comm f ex}
we have $R(\Witt)=\Witt$.
In particular, $\Witt$ is indeed a vertex Lie algebra.
It is graded, $h_L=2$.

\bigskip

{\bf Proposition.}\: {\it
Let $R$ be an even vertex Lie algebra such that $R$ is a free $\K[T]$-module
of rank $1$.
Then $R$ is abelian or $R\cong\Witt$.
}

\bigskip

\begin{pf}
We may identify $R$ with $\K[T]$.
Let $p(T,\la):=1\subla 1\in\K[T,\la]$ and 
$p(T,\la)=\sum_{i=0}^n p_i(\la)T^i$.

The conformal Jacobi identity for $a=b=c=1$ is
$$
p(-\mu,\la)p(T,\mu)
=
p(T+\la,\mu-\la)p(T,\la)
-
p(T+\mu-\la,\la)p(T,\mu-\la).
$$
Suppose that $p_n(\la)\ne 0$ and $n\geq 2$.
Then the coefficient of $T^{2n-1}$ of these three terms are $0$,
$$
p_{n-1}(\mu-\la)p_n(\la)+
n p_n(\mu-\la)\la p_n(\la)+
p_n(\mu-\la)p_{n-1}(\la),
$$
and $p_{n-1}(\la)p_n(\mu-\la)+n p_n(\la)(\mu-\la)p_n(\mu-\la)+
p_n(\la)p_{n-1}(\mu-\la)$.
We obtain the contradiction $0=n p_n(\la)p_n(\mu-\la)(2\la-\mu)$.
Thus $1\subla 1=p_0(\la)+Tp_1(\la)$.

The conformal Jacobi identity for $\la=\mu-\la$ yields
$p(-2\la,\la)p(T,2\la)=0$.
Thus $p(-2\la,\la)=0$ and $1\subla 1=(T+2\la)p_1(\la)$.

Conformal skew-symmetry implies $(T+2\la)p_1(\la)=(T+2\la)p_1(-\la-T)$.
This shows that $p:=p_1(\la)\in\K$.
If $p=0$ then $1\subla 1=0$ and $R$ is abelian.
If $p\ne 0$ then $p\inv{}\subla p\inv=(T+2\la)p\inv$.
\end{pf}

\bigskip

A vertex Lie algebra $R$ is {\bf simple} if $R$ is not abelian and 
\index{simple!vertex Lie algebra}
the only ideals of $R$ are $0$ and $R$.

\bigskip

{\bf Remark.}\: {\it
\iti\:
A Lie algebra $\fg$ is simple iff $\tfg$ is simple.

\smallskip

\itii\:
The Witt vertex Lie algebra is simple.
}

\bigskip

\begin{pf}
\iti\:
If $I\subset\fg$ is an ideal then $I[\mu]\subset\fg[\mu]=\tfg$ is an ideal.
Thus if $\tfg$ is simple then $\fg$ is simple.
Conversely, let $I\subset\tfg$ be a non-zero ideal. 
Then $I\cap\fg$ is an ideal of $\fg$. 
If $a\submu\in I$ and $b\in\fg$ then 
$[(a\submu)\subla b]=\sum[a_n,b]\la^{(n)}\in I[\la]$.
Thus $I\cap\fg\ne 0$. 
We get $I\cap\fg=\fg$ and $I=\tfg$. 

\smallskip

\itii\:
Let $a=\sum_{i=0}^n k_i T^i L\in\Witt$ with $k_i\in\K$ and $k_n\ne 0$.
Then $a\subla L=\sum_i k_i(-\la)^i(T+2\la)L$.
The coefficient of $a\subla L$ of degree $n+1$ is a non-zero
scalar multiple of $L$. 
Thus the claim follows.
\end{pf}

\bigskip

We state without proof the following deep result.

\bigskip

{\bf Theorem.}\: {\it
Let $R$ be an even finite simple vertex Lie algebra. 
Then $R$ is isomorphic to a loop vertex Lie algebra or 
the Witt vertex Lie algebra.
\hfill $\square$
}

\subsection{Virasoro and Conformal Vertex Lie Algebras}
\label{SS:vir vlie}

We discuss Virasoro and conformal vectors of vertex Lie algebras.

\bigskip

The {\bf Virasoro} vertex Lie algebra $\Vir=\Witt\oplus\K\hc$ is 
\index{Virasoro!vertex Lie algebra}
the central extension of $\Witt$ defined by $T\hc=0$ and
$$
L\subla L
\; :=\;
(T+2\la)L
\; +\;
(\hc/2)\la^{(3)}.
$$
By Proposition \ref{SS:vlie into loclie} and section \ref{SS:weak comm f ex}
we have $R(\Vir)=\Vir$.
It is graded, $h_L=2$.

A {\bf Virasoro vector} of a vertex Lie algebra $R$ 
is a vector $L$ such that there exists a morphism 
$\Vir\to R$ with $L\mapsto L$.

Section \ref{SS:finite confa} provides a description of $\Vir$
in terms of generators and relations.
It follows that an even vector $L$ is a Virasoro vector iff 
$L\subla L=(T+2\la)L+L_{(3)}L\la^{(3)}$ and $\hc_L:=2L_{(3)}L\in\ker T$.

A {\bf dilatation operator} of $R$ is 
\index{dilatation operator}
an even diagonalizable operator $H$ such that $[H,T]=T$ and 
$H(a_t b)=(Ha)_t b+a_t(Hb)-(t+1)a_t b$.
To give a gradation of $R$ is equivalent to giving a dilatation operator.

A {\bf conformal vector} of a vertex Lie algebra $R$
is a Virasoro vector such that $L_{(0)}=T$ and $L_{(1)}$ is a 
dilatation operator.
In this case $R$ is graded, $L\in R_2$, and hence $T=L_{-1}, H=L_0$.

A {\bf conformal vertex Lie algebra} is a vertex Lie algebra together
with a conformal vector. 
A {\bf conformal vector} of a {\it graded} vertex Lie algebra 
is a Virasoro vector such that $L_{-1}=T$ and $L_0=H$.

\bigskip

{\bf Remark.}\: {\it
Let $R$ be a conformal vertex Lie algebra with 
\index{centre}
centre $C$.
Then $C=\ker T$ and $C\subset V_0$.
}

\bigskip

\begin{pf}
By Proposition \ref{SS:finite confa} we have $\ker T\subset C$.
We have $C\subset\ker T$ since $T=L_{(0)}$.
We have $C\subset V_0$ since $H=L_{(1)}$.
\end{pf}

\bigskip

{\bf Proposition.}\: {\it
Let $R$ be a vertex Lie algebra and $S\subset R$ a subspace that generates $R$.
If $L_{(0)}|_S=T$ and $L_{(1)}|_S$ is diagonalizable
then $L_{(0)}=T$ and $L_{(1)}$ is a dilatation operator.
}

\bigskip

\begin{pf}
By Remark \ref{SS:conf der} the operator $L_{(0)}$ is a derivation of $R$.
Thus $\ker(T-L_{(0)})$ is a conformal subalgebra of $R$.
Since $S\subset\ker(T-L_{(0)})$, we get $L_{(0)}=T$.

The conformal Jacobi identity implies
$$
[L_{(1)},a_{(t)}]b
\; =\;
(L_{(0)}a)_{(t+1)}b+(L_{(1)}a)_{(t)}b
\; =\; 
(L_{(1)}a)_{(t)}b-(t+1)a_{(t)}b.
$$
Moreover, $[L_{(1)},T]=-[T,L_{(1)}]=L_{(0)}=T$.
This shows that the span of the eigenvectors of $L_{(1)}$ is 
a conformal subalgebra.
Hence $L_{(1)}$ is diagonalizable.
The above two identities also show that $L_{(1)}$ is a dilatation operator.
\end{pf}

\subsection{Griess Algebra}
\label{SS:griess}

We show that the Virasoro vectors are the idempotents of the Griess algebra.

\bigskip

Let $R$ be a vertex Lie algebra of CFT-type and $\hk\in R_0, \hk\ne 0$.
Then $R_2$ with $ab:=a_0 b$ is an algebra with a bilinear form defined by 
$(a,b)\hk=a_2 b$.
The algebra $R_2$ is called 
\index{Griess algebra}
the {\bf Griess algebra}.

\bigskip

{\bf Remark.}\: {\it
Let $R$ be a vertex Lie algebra of CFT-type such that $R_1=0$.
Then $ab$ is commutative and $(\, ,)$ is invariant and symmetric.
If $L$ is a conformal vector then $L/2$ is an identity of $R_2$.
}

\bigskip

\begin{pf}
We have $a\subla b=a_{-1}b+ab\la+(a,b)\hk\la^{(3)}$ and
$-b_{-\la-T}a=-b_{-1}a+T(ba)+ba\la+(b,a)\hk\la^{(3)}$.
Thus $ab$ is commutative and $(\, ,)$ invariant. 
Replacing $\la, \mu$ by $\la+\ka, \mu+\ka$, 
the conformal Jacobi identity becomes 
$[[a_{\la+\ka}b]_{\mu+\ka}c]=[a_{\ka+\la}[b_{\mu-\la}c]-
[b_{\mu-\la}[a_{\ka+\la}c]$.
The coefficient of $\la\mu\ka^{(2)}$ is $(ab,c)\hk=(a,bc)\hk$.
The last claim follows from $L_0|_{R_2}=2$.
\end{pf}

\bigskip

{\bf Proposition.}\: {\it
Let $R$ be a vertex Lie algebra of CFT-type.

\smallskip

\iti\:
An even vector $L\in R_2$ is a Virasoro vector iff $L_0 L=2L$;
in other words, iff $L/2$ is an idempotent of the Griess algebra $R_2$.

\smallskip

\itii\:
Let $L\in R_2$ be an even vector and $S\subset R$ a homogeneous subset 
that generates $R$.
If $L_{-1}|_S=T$ and $L_0|_S=H$ then $L$ is a conformal vector.
}

\bigskip

\begin{pf}
\iti\:
Since $R$ is of CFT-type we have $L_2 L\in\ker T$ and $L_n L=0$ for $n\geq 3$.
Conformal skew-symmetry implies $L_1 L=-L_1 L$ and $L_{-1}L=-L_{-1}L+T(2L)$.

\smallskip

\itii\:
By Proposition \ref{SS:vir vlie} we have $L_{-1}=T$ and $L_0=H$.
Since $L_0 L=HL=2L$, the claim follows from \iti.
\end{pf}

\subsection{Primary Vectors and GKO Coset Construction}
\label{SS:gko vector}

We prove that 
\index{coset!construction}
the coset of a conformal vertex Lie subalgebra is conformal.
This 
\index{GKO construction}
is 
\index{Goddard-Kent-Olive construction}
the {\bf Goddard-Kent-Olive} coset construction.

\bigskip

A homogeneous vector $a$ of 
\index{primary vector}
a conformal vertex Lie algebra $R$ is {\bf primary} if $L\subla a=(T+H\la)a$.
In other words, $L_n a=0$ for $n\geq 1$.
A 
\index{quasi-primary vector}
homogeneous vector $a$ is {\bf quasi-primary} if $L_1 a=0$.

For example, $L$ is quasi-primary and primary iff $\hc_L=0$.
For $a$ primary, $Ta$ is primary iff $h_a=0$.
Central elements are primary.
Many conformal vertex Lie algebras are generated by primary vectors.

A conformal vertex Lie {\bf subalgebra} of $R$ is 
\index{conformal!vertex Lie subalgebra}
a graded vertex Lie subalgebra $R'$ together with a quasi-primary 
conformal vector $L'\in R'$.

\bigskip

{\bf Remark.}\: {\it 
\iti\:
Let $S$ be a subset of a vertex Lie algebra $R$.
If $L\in S$ such that $L_{(0)}|_S=T$ then $C_R(S)=\ker L_{(0)}$.

\smallskip

\itii\:
Let $R$ be a conformal vertex Lie algebra of CFT-type.
Then a Virasoro vector $L'\in R_2$ 
is quasi-primary iff $L'$ commutes with $L-L'$. 
}

\bigskip

\begin{pf}
\iti\:
Since $L\in S$ we have $C_R(S)\subset\ker L_{(0)}$.
Conversely, let $a\in\ker L_{(0)}$ and $b\in S$. 
Then $L_{(0)}(a\subla b)=(L_{(0)}a)\subla b+a\subla(L_{(0)}b)=a\subla (Tb)=
(T+\la)a\subla b$.
This implies $a\subla b=0$.

\smallskip

\itii\:
Let $L^c:=L-L'$.
If $L', L^c$ commute then $L_1 L'=L^c_1 L'+L'_1 L'=0$.
Conversely, suppose that $L'$ is quasi-primary.
Since $L_{-1}L'=TL'=L'_{-1}L'$ and $L_0 L'=2L'=L'_0 L'$ 
and $L_1 L'=0=L'_1 L'$ we have $L^c_n L'=0$ for $n=-1,0,1$.
Conformal skew-symmetry and the fact that $R$ is of CFT-type imply 
$L'_{-1}L^c=0$.
Proposition \ref{SS:griess}\,\itii\ 
shows that $L'$ is a conformal vector of $\overline{\set{L'}}$.
Hence $L', L^c$ commute by \iti.
\end{pf}

\bigskip

{\bf Proposition.}\: {\it
Let $R$ be a conformal vertex Lie algebra of CFT-type and 
$R'$ a conformal vertex Lie subalgebra. 
Then $C_R(R')$ is a conformal vertex Lie subalgebra with 
conformal vector $L-L'$ and $\hc_{L-L'}=\hc_L-\hc_{L'}$.
}

\bigskip

\begin{pf}
The Remark shows that $R^c:=C_R(R')=\ker L'_{-1}$ and that 
$L'$ and $L^c:=L-L'$ commute.
Thus $L^c\in R^c$ and $R^c$ is a graded vertex Lie subalgebra.
We have $L^c_{-1}=L_{-1}=T$ and $L^c_0=L_0=H$ on $R^c$.
Moreover, $L^c=L-L'\in R_2$.
Thus $L^c$ is a conformal vector of $R^c$ by 
Proposition \ref{SS:griess}\,\itii.
The vector $L^c$ is quasi-primary because $L^c, L'$ commute.
We have $\hc_{L^c}=\hc_L-\hc_{L'}$ because $L_2 L=L'_2 L'+L^c_2 L^c$.
\end{pf}

\subsection{$U(1)$-Currents and Chodos-Thorn Construction}
\label{SS:chodos thorn}

We show how conformal vectors can be modified using a $U(1)$-current.

\bigskip

A {\bf $U(1)$-vector} of a vertex Lie algebra $R$ is 
\index{U1-vector@$U(1)$-vector}
an even vector $J$ such that $J\subla J=\hk_J\la$ for some $\hk_J\in\ker T$.
Proposition \ref{SS:aff vlie} implies that if $R$ is of CFT-type 
then any even vector $J\in R_1$ is a $U(1)$-vector.

Let $J$ be a $U(1)$-vector.
If $a\in R$ is an eigenvector of $J_0$ then the eigenvalue $q_a$ is 
called 
\index{Jcharge@$J$-charge}
the {\bf $J$-charge} of $a$.
Since $J_0$ is a derivation, 
the charge of $a_i b$ is $q_a+q_b$ and the charge of $Ta$ is $q_a$.
Hence if $J_0$ is diagonalizable on a subspace $S$ then
$J_0$ is diagonalizable on $\bS$.
A vector $a$ 
\index{Jprimary@$J$-primary}
is {\bf $J$-primary} if $J\subla a=q_a a$.
In this case $Ta$ is $J$-primary 
\index{primary!wrt a $U(1)$-current}
iff $q_a=0$.

A {\bf $U(1)$-current} of a conformal vertex Lie algebra $R$ is 
\index{U1-current@$U(1)$-current}
a primary $U(1)$-vector $J\in R_1$ such that $J_0$ is diagonalizable.
If $J$ is a $U(1)$-current then $L\subla J=(T+\la)J$.
Hence $J\subla L=-L_{-\la-T}J=\la J$ and $L$ has charge $0$.
Of course, $J$ has also charge $0$.

\bigskip

{\bf Proposition.}\: {\it
Let $J\in R_1$ be a primary $U(1)$-vector of 
a conformal vertex Lie algebra $R$.
Then $L':=L+TJ$ is a Virasoro vector with $\hc_{L'}=\hc_L-12\hk_J$.
If $J$ is a $U(1)$-current then $L'$ is a conformal vector with
$L'_0=L_0-J_0$ and 
any vector, that is both $L$- and $J$-primary, is $L'$-primary.
}

\bigskip

\begin{pf}
We have $L\subla TJ=(T+\la)L\subla J=T^2 J+2T\la J+\la^2 J$.
From $J\subla L=-L_{-\la-T}J=\la J$ we get $(TJ)\subla L=-\la^2 J$.
Finally, $(TJ)\subla TJ=-\hk_J\la^3$.
Thus we obtain
\begin{align}
\notag
L'\subla L'
\; &=\;
TL+2\la L+\hc_L\la^{(3)}/2+T^2 J+2TJ\la+\la^2 J-\la^2 J-\hk_J\la^3
\\
\notag
\; &=\;
T(L+TJ)
\; +\;
2\la(L+TJ)
\; +\;
(\hc_L-12\hk_J)\la^{(3)}/2.
\end{align}

Suppose that $J$ is a $U(1)$-current.
Since $[L_0,J_0]=[H,J_0]=0$, 
the operators $L_0$ and $J_0$ have a simultaneous eigenspace decomposition. 
From $(TJ)\subla=-\la J\subla$ and Proposition \ref{SS:vir vlie} follows that
$L'$ is a conformal vector with $L'_0=L_0-J_0$.
The last claim is clear.
\end{pf}

\subsection{$N$=1, $N$=2, and Topological Virasoro Algebra}
\label{SS:n1 n2 top vir}

We apply the Chodos-Thorn construction to $\twoVir$ and obtain
the topological Virasoro vertex Lie algebra.

\bigskip

The $N$=1 {\bf super Virasoro} vertex Lie algebra 
\index{N1 super Virasoro@$N$=1 super Virasoro!vertex Lie algebra}
is the conformal vertex Lie algebra $\oneVir=\K[T]L\oplus\K[T]G\oplus\K\hc$ 
with conformal vector $L$ and $\hc_L=\hc$ and an odd, primary vector $G$
of weight $3/2$ such that $G\subla G=2L+(2\hc/3)\la^{(2)}$.

The $N$=2 {\bf super Virasoro} vertex Lie algebra 
\index{N2 super Virasoro@$N$=2 super Virasoro!vertex Lie algebra}
is the conformal vertex Lie algebra $\twoVir=
\K[T](\K L\oplus\K G^+\oplus\K G^-\oplus\K J)\oplus\K\hc$
with conformal vector $L$ and $\hc_L=\hc$,
a $U(1)$-current $J$ with $\hk_J=\hc/3$, and 
two odd, primary, $J$-primary vectors $G\suppm$ of weight $3/2$ and 
charge $\pm 1$ with $G^+{}\subla G^-=2L+(T+2\la)J+(2\hc/3)\la^{(2)}$
and $G\suppm{}\subla G\suppm=0$.

By Proposition \ref{SS:vlie into loclie} and section \ref{SS:weak comm f ex}
we have $R(\oneVir)=\oneVir$ and $R(\twoVir)$=$\twoVir$.
In particular, $\oneVir$ and $\twoVir$ are indeed vertex Lie algebras.

The conformal vertex Lie algebra 2-$\Vir$ has an involution
given by $G\suppm\mapsto G^{\mp}$ and $J\mapsto -J$.
This follows from 
$G^-{}\subla G^+=G^+{}_{-\la-T}G^-=2L+(-T-2\la)J+(2\hc/3)\la^{(2)}$.
This involution is 
\index{mirror involution}
the {\bf mirror involution}.

Proposition \ref{SS:chodos thorn} shows that $L+TJ/2$ is a conformal vector
of $\twoVir$.
The {\bf topological Virasoro} vertex Lie algebra $\tVir$ is 
\index{topological Virasoro algebra}
the vertex Lie algebra $\twoVir$ endowed with the conformal vector $L+TJ/2$.

Proposition \ref{SS:chodos thorn} yields the following 
description of the topological Virasoro algebra.
We have $\tVir=\K[T](\K L\oplus\K Q\oplus\K G\oplus\K J)\oplus\K\hd$
where $L$ is the conformal vector and $\hc_L=0$, 
$J$ is a $U(1)$-vector with $\hk_J=\hd$ and $L\subla J=(T+\la)J-\hd\la^{(2)}$,
$Q, G$ are odd, primary, $J$-primary vectors of weight $1, 2$ and 
charge $1, -1$, and $Q\subla G=L+J\la+\hd\la^{(2)}, Q\subla Q=G\subla G=0$.

A vertex Lie algebra isomorphism $\twoVir\to\tVir$
is given by $L\mapsto L-TJ/2, G^+\mapsto 2Q, G^-\mapsto G, J\mapsto J$, and
$\hc\mapsto 3\hd$.

The non-trivial commutators of $\fg(\tVir)$ are
\begin{align}
\notag
[Q_n,G_m]
\; &=\;
L_{n+m}
+
nJ_{n+m}
+(n^2-n)\de_{n+m}\hd/2,
\\
\notag
[L_n,J_m]
\; &=\;
-mJ_{n+m}
-(n^2+n)\de_{n+m}\hd/2.
\end{align}

\section{Supplements}
\label{S:vlie suppl}

In sections \ref{SS:s3 symm conf J}--\ref{SS:frob alg vlie}
we present further examples of vertex Lie algebras:
semidirect products, the $\la$-commutator of an
associative conformal algebra, and vertex Lie algebras constructed from
Frobenius algebras.
In section \ref{SS:free vlie rev}
we give a second construction of free vertex Lie algebras.

In section \ref{SS:1trunc conf} and \ref{SS:1trunc vlie}
we discuss the functor $R\mapsto (R_0, R_1)$ from 
$\N$-graded vertex Lie algebras to 1-truncated vertex Lie algebras 
and construct a left adjoint for it.

\subsection{$\bbS_3$-Symmetry of the Conformal Jacobi Identity}
\label{SS:s3 symm conf J}

We prove that conformal skew-symmetry implies an $\bbS_3$-symmetry of the
conformal Jacobi identity that permutes the elements $a, b, c$.

\bigskip

{\bf Proposition.}\: {\it
Let $R$ be a conformal algebra that satisfies conformal skew-symmetry.
Then the conformal Jacobi identity holds for $a, b, c$ iff it holds
for any permutation of $a, b, c$.
}

\bigskip

\begin{pf}
The conformal Jacobi identity is 
$[[a\subla b]\submu c]=[a\subla[b_{\mu-\la}c]]-[b_{\mu-\la}[a\subla c]]$.
It holds for $a, b, c$ iff it holds for $b, a, c$ since 
$$
[[b_{-\la-T}a]\submu c]
\; =\;
[e^{T\del\subla}[b_{-\la}a]\submu c]
\; =\;
e^{-\mu\del\subla}[[b_{-\la}a]\submu c]
\; =\;
[[b_{\mu-\la}a]\submu c].
$$
It holds for $a, b, c$ iff it holds for $a, c, b$ since 
$[[a\subla b]\submu c]=-e^{T\del\submu}[c_{-\mu}[a\subla b]]$,
$$
[a\subla[b_{\mu-\la}c]]
\; =\;
-[a\subla e^{(T-\la)\del\submu}[c_{-\mu}b]]
\; =\;
-e^{T\del\submu}[a\subla[c_{-\mu}b]],
$$
and $[b_{\mu-\la}[a\subla c]]=-e^{T\del\submu}[[a\subla c]_{-\mu+\la}b]$.
\end{pf}

\subsection{Products and Semidirect Products}
\label{SS:semidir prod}

We discuss products and semidirect products of vertex Lie algebras.

\bigskip

The {\bf product} $\prod R_i$ of a family of vertex Lie algebras $R_i$ is  
\index{product of vertex Lie algebras}
an unbounded vertex Lie algebra with $\la$-bracket
$[(a^i)\subla(b^i)]:=([a^i{}\subla b^i])$.
If the family is finite then $\prod R_i$ is a vertex Lie algebra.

Let $R$ and $I$ be vertex Lie algebras together with a morphism 
$R\to\Der\subsc(I), a\mapsto a\subla$.
The {\bf semidirect product} $R\ltimes I$ is 
\index{semidirect product}
the conformal algebra $R\oplus I$ with $\la$-bracket
$$
[(a+e)\subla(b+f)]
\; :=\;
[a\subla b]
\: + \:
a\subla f
\: -\:
\zeta^{eb}\, b_{-\la-T} e
\: +\:
[e\subla f]
$$
for $a, b\in R$ and $e, f\in I$.

\bigskip

{\bf Proposition.}\: {\it
The semidirect product $R\ltimes I$ is a vertex Lie algebra.
}

\bigskip

\begin{pf}
It is clear that conformal skew-symmetry is satisfied.
The conformal Jacobi identity holds for $a, b, c\in R$ and for $a, b, c\in I$
since $R$ and $I$ are vertex Lie algebras.
It holds for $a, b\in R, c\in I$ since $R\to\fgl\subsc(I)$ is a morphism.
It holds for $a\in R, b, c\in I$ since $a\subla\in\Der\subsc(I)$.
Proposition \ref{SS:s3 symm conf J} implies that it holds for any 
$a, b, c\in R\ltimes I$.
\end{pf}

\bigskip

A short exact sequence 
$0\to I\to R^e\overset{p}{\to}R\to 0$ of vertex Lie algebras
{\bf splits} if there exists a section of $p$, that is,  
a morphism $\io: R\to R^e$ such that $p\circ\io=\id_R$. 

The inclusion $I\to R\ltimes I$ and the projection $R\ltimes I\to R$
yield a split short exact sequence $0\to I\to R\ltimes I\to R\to 0$.
Conversely, any split short exact sequence $0\to I\to R^e\to R\to 0$ is 
isomorphic to such a sequence where $R\to\Der\subsc(I)$ is given
by the composition of a section $\io: R\to R^e$ and the adjoint
action $R^e\to\Der\subsc(I)$.

We have $R\ltimes I=R\times I$ iff the map $R\to\Der\subsc(I)$ is zero
iff the extension $0\to I\to R\ltimes I\to R\to 0$ is central: 
$I\subset C(R\ltimes I)$.

Here is a concrete example of a semidirect product.
Let $R=E[\mu]\oplus K$ be a $\K$-graded $\K[T]$-module such that 
$T|_{E[\mu]}=-\mu\cdot$, $E\subset R$ is a $\K$-graded subspace, and
$K\subset R_0$ is torsion. 
Define a $\K[T]$-module morphism $\Witt\to\fgl\subsc(R)$ by
$L\subla a:=(T+H\la)a$ for $a\in E$. 
A direct calculation shows that $R$ thus becomes a $\Witt$-module.

In particular, loop and affine vertex Lie algebras $\tfg$ and $\hfg$
with $E=\fg$ and 
Clifford vertex Lie algebras $C(B)$ with $E=\rPi B$ are $\Vir$-modules.
In these three cases we have $\Witt\to\Der\subsc(R)$ 
as direct calculations show.
Thus semidirect products such as $\Witt\ltimes\tfg$ and
$\Vir\ltimes C(B)$ are defined.

For example, if $\fg$ is a one-dimensional abelian Lie algebra 
then $\Witt\ltimes\tfg$ is isomorphic to the vertex Lie algebra
$\set{p\del_x+q\mid p, q\in\K[x\uppm]}$
of first order differential operators on
$\K[x\uppm]$ because $[-x^{n+1}\del_x,x^m]=-mx^{n+m}$.
This is 
\index{Atiyah vertex Lie algebra}
the {\bf Atiyah vertex Lie algebra}.

\subsection{Associative and Commutative Conformal Algebras}
\label{SS:assoc confa}

We show that the Jacobi identity, skew-symmetry, associativity, and
commutativity for the algebra $\fg_R$ are equivalent to the 
corresponding ``conformal'' identities for the conformal algebra $R$.

\bigskip

In section \ref{SS:loclie of vlie}
we showed that if $R$ is a conformal algebra then $\fg_R=\tR /T\tR$
is an algebra with $[a_t,b_s]=\sum\binom{t}{i}(a_i b)_{t+s-i}$.
From Corollary \ref{SS:loc diff op} and Remark \ref{SS:loclie of vlie}
follows that the linear map 
$$
Y:\: 
R\;\to\;\fg_R\pau{z\uppm}, 
\quad
a\;\mapsto\; a(z)\; =\; \sum a_t z^{-t-1},
$$
is a monomorphism of unbounded conformal algebras and $Y(R)$ is local.

\bigskip

{\bf Proposition.}\: {\it
\iti\:
A conformal algebra $R$ satisfies the conformal Jacobi identity iff
$\fg_R$ satisfies the Leibniz identity.

\smallskip

\itii\:
A conformal algebra $R$ satisfies conformal skew-symmetry iff
$\fg_R$ is skew-symmetric.
}

\bigskip

\begin{pf}
\iti\:
This follows from Remark \ref{SS:conf jac}, Proposition \ref{SS:conf jac},
and the fact that $Y:R\to\fg_R\pau{z\uppm}$ is a monomorphism.

\smallskip

\itii\:
This follows from Remark \ref{SS:conf skew sym}, 
Proposition \ref{SS:conf skew sym}, and the fact that 
$Y:R\to\fg_R\pau{z\uppm}$ is a monomorphism and $Y(R)$ is local.
\end{pf}

\bigskip

{\bf Conformal associativity} for an unbounded conformal algebra $R$ is
\index{conformal!associativity}
$$
(a\subla b)\submu c
\; =\;
a\subla(b_{\mu-\la}c).
$$
{\bf Conformal commutativity} for weakly local $b, a\in R$ is
\index{conformal!commutativity}
$$
a\subla b
\; =\;
\paraab\, b_{-\la-T}a.
$$

A 
\index{associative!conformal algebra}
conformal algebra $R$ is {\bf associative} or {\bf commutative}
if $R$ satisfies 
\index{commutative!conformal algebra}
conformal associativity or conformal commutativity, resp.
Vertex Lie algebras are 
\index{Lie conformal algebra}
also called {\bf Lie conformal algebras}.

The following statements complement 
Remarks \ref{SS:conf jac} and \ref{SS:conf skew sym},
Propositions \ref{SS:conf jac}, \ref{SS:conf skew sym}, and 
\ref{SS:unb confa power series}, and the Proposition above.

If $R$ is associative then $ab:=a_0 b$ is associative.
If $R$ is commutative then $ab$ is commutative on $R/TR$.
An algebra $A$ is associative iff $A\pau{z\uppm}$ is associative iff 
$A\pau{\mu}$ is associative.
The multiplication of an algebra $C$ is commutative 
iff $C\pau{\mu}$ is commutative.
In this case local $C$-valued distributions satisfy conformal commutativity.
Dong's lemma also holds for associative algebras.

Calculations as in section \ref{SS:affin vlie} show that
if $R$ is commutative or associative then so is $\tR$.
It follows that $R$ is commutative or  associative 
iff the multiplication of $\fg_R$ is commutative or associative, resp.
As in section \ref{SS:loclie of vlie}
one shows that $\fg_R$ satisfies a universal property depending
on the class of algebras and conformal algebras one looks at.

\subsection{$\la$-Commutator}
\label{SS:lambda commut}

In chapter \ref{C:va} we will see that any associative vertex algebra has 
an underlying vertex Lie algebra.
Here we prove that also associative conformal algebras have
an underlying vertex Lie algebra.

\bigskip

Let $R$ be an unbounded conformal algebra such that 
conformal commutativity is well-defined.
The {\bf $\la$-commutator} is 
$$
[a\subla b]
\; :=\;
a\subla b\; -\; \paraab\, b_{-\la-T}a.
$$
It is clear that the $\la$-commutator is a $\la$-bracket and 
satisfies conformal skew-symmetry.

If $A$ is an associative algebra then the commutator $ab-\paraab ba$
is a Lie bracket. 
Let $A'$ denote this Lie algebra.
The identity in the proof of Proposition \ref{SS:unb confa power series}\,\iti\
shows that the $\la$-commutator of $A\pau{\mu}$ is equal to 
the unbounded $\la$-bracket of $A'\pau{\mu}$.
In particular, the unbounded $\la$-bracket of $\fgl(R)\pau{\mu}$
is the $\la$-commutator of the associative conformal algebra 
$\End(R)\pau{\mu}$.

\bigskip

{\bf Proposition.}\: {\it
Let $R$ be an associative conformal algebra. 
Then the $\la$-commutator of $R$ is a Lie $\la$-bracket.
}

\bigskip

\begin{pf}
Since $R$ is associative, $\fg_R$ is associative.
Thus $\fg_R$ endowed with the commutator is a Lie algebra.
Let $R'$ be $R$ endowed with the $\la$-commutator.
As in section \ref{SS:affin vlie} we see that the $\la$-product of 
the affinization $\tR'$ is
$$
(e^{\del\del\subla}f)g\,(a\subla b-b_{-\la-T}a)
\; =\;
(fa)\subla(gb)-(gb)_{-\la-\del-T}(fa).
$$
Setting $\la=0$, this shows that the algebra $\fg_{R'}$ is equal to
the Lie algebra $\fg_R$.
Proposition \ref{SS:assoc confa} implies that $R'$ is a vertex Lie algebra.
\end{pf}

\subsection{Frobenius Algebras and Vertex Lie Algebras}
\label{SS:frob alg vlie}

We construct a vertex Lie algebra from a Frobenius algebra.

\bigskip

A {\bf Frobenius algebra} is
\index{Frobenius algebra}
a commutative algebra $C$ together with an {\it invariant} symmetric
bilinear form: $(ab,c)=(a,bc)$.

Our definition differs from the standard one in two respects.
First, Frobenius algebras in general need not be commutative.
Second, the bilinear form is usually assumed to be non-degenerate.
Its existence imposes then a strong restriction on the algebra $C$.

Let $C$ be a commutative algebra.
Consider $C$ as a differential algebra with $T=0$.
By section \ref{SS:affin vlie}
the tensor product $C[T]=\Witt\otimes C$ is a vertex Lie algebra
with $\la$-bracket $[a\subla b]=(T+2\la)ab$ for $a, b\in C$.

If $C$ is a Frobenius algebra then $C[T]$ has a central extension
$R=C[T]\oplus\K\hc$ defined by $T\hc=0$ and
$$
[a\subla b]
\; =\;
(T+2\la)ab
\; +\;
(a,b)\hc\,\la^{(3)}
$$
for $a, b\in C$.
A direct calculation shows that $R$ is indeed a vertex Lie algebra.
It is of CFT-type with $R_1=0$ and $R_2=C$.
Moreover, the Griess algebra of $R$ is equal to the Frobenius algebra $C$.

Taking $C=\K$ with bilinear form defined by $(1,1):=1/2$, 
we get $R=\Vir$ with $1\mapsto L$ and $\hc\mapsto\hc$.

\subsection{Free Vertex Lie Algebras Revisited}
\label{SS:free vlie rev}

We construct the free vertex Lie algebra $R(S,o)$ 
\index{free vertex Lie algebra}
by describing the Lie algebra $\fg(R(S,o))$ explicitly 
in terms of generators and relations.

\bigskip

{\bf Proposition.}\: {\it
The functor $R\mapsto (R,o_+)$ from vertex Lie algebras to sets with 
non-negative locality function has a left adjoint $(S,o)\mapsto R(S,o)$. 
}

\bigskip

\begin{pf}
Let $\fg$ be the Lie algebra generated by $S\times\Z$ with relations
$$
\sum_{i\geq 0}\:
(-1)^i \binom{o(a,b)}{i}\; [(a,t+i),(b,s+o(a,b)-i)]
\; =\; 0
$$
for any $a, b\in S$ and $t, s\in\Z$.
Denote by $a_t$ the image of $(a,t)$ in $\fg$ and
define $F:=\set{a(z)\mid a\in S}$.
Then any $a(z), b(z)$ are local of order $\leq o(a,b)$, $R:=\bF$ is 
a vertex Lie algebra, and $(S,o)\to (R,o_+), a\mapsto a(z)$, is a morphism.

The Lie algebra $\fg$ satisfies the following universal property.
Let $\fg'$ be a local Lie algebra. 
Then a morphism $\psi: (S,o)\to(R(\fg'),o_+)$ yields a unique 
Lie algebra morphism $\fg\to\fg'$ such that $a(z)\mapsto \psi a$ 
for any $a\in S$.
Conversely, a Lie algebra morphism $\phi:\fg\to\fg'$ such that 
$\phi F\subset\bF_{\fg'}$ yields a morphism 
$(S,o)\to (R(\fg'),o_+), a\mapsto\phi a(z)$.

The composition $(S,o)\to (R,o_+)\to (R(\fg(R)),o_+), a\mapsto Y(a(z))$, 
yields a morphism $\phi:\fg\to\fg(R)$ such that $a_t\mapsto a(z)_t$.
On the other hand, by the universal property of $\fg(R)$ the inclusion 
$R\subset\fg\pau{z\uppm}$ induces a morphism $\phi':\fg(R)\to\fg$ such that 
$a(z)_t\mapsto a_t$ for any $a(z)\in R$.

In general, if $R$ is a vertex Lie algebra and $F\subset R$ a subset 
such that $\bF=R$ then the algebra $\fg(R)$ is generated by 
$\set{a_t\mid a\in F, t\in\Z}$ because $(a_r b)_t=(a(z)_r b(z))_t=
\sum_i (-1)^i\binom{r}{i}[a_{r-i},b_{t+i}]$ and $(Ta)_t=-ta_{t-1}$ 
for $a, b\in R$ and $r\geq 0, t\in\Z$.

It follows that $\phi, \phi'$ are inverse to each other because
$\phi, \phi'$ are inverse to each other on the sets of generators 
$\set{a_t\mid a\in S, t\in\Z}\subset\fg$ and 
$\set{a(z)_t\mid a\in S, t\in\Z}\subset\fg(R)$. 
It is clear that $\phi, \phi'$ respect $R$ and $F_{\fg(R)}$.
Since $\fg(R)$ is a local Lie algebra, so is $(\fg,R)$.

The universal property of $\fg$ implies that $(S,o)\mapsto(\fg,R)$
is left adjoint to $\fg'\mapsto(R(\fg'),o_+)$.
Since $R'\mapsto\fg(R')$ is fully faithful and the counit of adjunction
$\phi':\fg(R)\to\fg$ is an isomorphism, the claim follows.
\end{pf}

\subsection{Skew-Symmetric 1-Truncated Conformal Algebras}
\label{SS:1trunc conf}

Note that if $R$ is an $\N$-graded conformal algebra then 
$R_0\oplus R_1$ is an $\N$-fold subalgebra of $R$.
We describe the algebraic structure on $R_0$ and $R_1$ and 
construct a universal conformal algebra from it.

We denote the $t$-th product by $a_t b$, not by $a_{(t)}b$.

\bigskip

Let $\cO$ be a vector space.
Denote the elements of $\cO$ by $f, g, h$. 

A {\bf skew-symmetric 1-truncated conformal algebra} over $\cO$ is 
\index{skew-symmetric 1-truncated conformal algebra}
an algebra $\fa, x\otimes y\mapsto [x,y]$, 
with an even symmetric bilinear map 
$\fa\otimes\fa\to\cO, x\otimes y\mapsto (x,y)$, and even linear maps 
$\fa\otimes\cO\to\cO, x\otimes f\mapsto xf$, and $d: \cO\to\fa$
such that $[x,y]+[y,x]=d(x,y)$, $(df)g=[df,x]=0$, and $(x,df)=xf$.

Here and in the following we denote the elements of $\fa$ by $x, y, z$. 
The map $\fa\to\End(\cO), x\mapsto x\cdot$, is 
\index{anchor}
the {\bf anchor} of $\fa$.
Since $[x,df]+[df,x]=d(x,df)=d(xf)$, the axiom $[df,x]=0$ is equivalent to 
$[x,df]=d(xf)$.

For example, skew-symmetric 1-truncated conformal algebras $\fa$ over $\K$ with
$d=0$ and $x\cdot=0$ for any $x$ are exactly skew-symmetric algebras $\fa$ 
with an even symmetric bilinear form.

Let $R$ be an $\N$-graded conformal algebra that satisfies 
conformal skew-symmetry.
Then $f\subla g=0, x\subla f=x_0 f, f\subla x=-x_0 f$, and 
$x\subla y=x_0 y+x_1 y\,\la$.
It is easy to see that $R_1$ is a skew-symmetric 1-truncated conformal algebra
over $R_0$ with $[x,y]:=x_0 y, (x,y):=x_1 y, xf:=x_0 f$, and $d:=T$.

It is clear how to define morphisms $(\cO,\fa)\to (\cO',\fa')$ 
of skew-symmetric 1-truncated conformal algebras.
We write $\fa=(\cO,\fa)$.

\bigskip

{\bf Proposition.}\: {\it
The functor $R\mapsto R_1$, from $\N$-graded conformal algebras that satisfy 
conformal skew-symmetry to skew-symmetric 1-truncated conformal algebras,
has a left adjoint $\fa\mapsto R(\fa)$.
}

\bigskip

\begin{pf}
We construct the left adjoint $\fa\mapsto R(\fa)$ as follows.
Let $R(\fa):=\cO\oplus\fa[T]$ be the $\K[T]$-module defined by 
$T:=d$ on $\cO$ and $T=T\cdot$ on $\fa[T]$.
It is $\N$-graded with $R_0=\cO$ and $R_1=\fa$.

Define $f\subla g:=0, x\subla f:=xf, f\subla x:=-xf$, and 
$x\subla y:=[x,y]+(x,y)\la$.
Then it is clear that conformal skew-symmetry holds for any 
$a, b\in\cO\oplus\fa$.

We have $(Tf)\subla g=-\la f\subla g$ since $(df)g=0$.
We have $(Tf)\subla x=-\la f\subla x$ since $[df,x]=0$ and 
$(df,x)=(x,df)=xf=-f_0 x$.
Conformal skew-symmetry implies that also $f\subla Tg=(T+\la)f\subla g$
and $x\subla Tf=(T+\la)x\subla f$.
Thus we can uniquely extend the $\la$-product on $\cO\oplus\fa$ 
to a $\la$-product on $R(\fa)$.
By Lemma \ref{SS:conf skew sym} the conformal algebra $R(\fa)$ satisfies
conformal skew-symmetry since any $a, b\in\cO\oplus\fa$ do.

It is clear that $\fa\mapsto R(\fa)$ is left adjoint to $R\mapsto R_1$.
\end{pf}

\subsection{1-Truncated Vertex Lie Algebras}
\label{SS:1trunc vlie}

We continue the investigation of section \ref{SS:1trunc conf} 
and consider now the consequences of the conformal Jacobi identity.

\bigskip

A {\bf 1-truncated vertex Lie algebra} is 
\index{onetruncated vertex Lie algebra@1-truncated vertex Lie algebra}
a skew-symmetric 1-truncated conformal algebra $\fa$ such that 
$R(\fa)$ is a vertex Lie algebra.

A {\bf Leibniz algebra} is 
\index{Leibniz!algebra}
an algebra $\fa$ that satisfies the Leibniz identity
$[[a,b],c]=[a,[b,c]]-\paraab\,[b,[a,c]]$, see section \ref{SS:conf jac}.
Equivalently, left multiplication $\fa\to\fgl(\fa), a\mapsto [a,\, ]$, 
is an algebra morphism.

A {\bf module} over a Leibniz algebra $\fa$ is 
\index{module!over a Leibniz algebra}
a vector space $M$ with an algebra morphism $\fa\to\fgl(M), a\mapsto a\cdot$.

If $M$ and $N$ are $\fa$-modules then $M\otimes N$ is an $\fa$-module with
$a(m\otimes n)=am\otimes n+\zeta^{am}\, m\otimes an$.

\bigskip

{\bf Proposition.}\: {\it
Let $\fa$ be a skew-symmetric 1-truncated conformal algebra.
Then $\fa$ is a 1-truncated vertex Lie algebra iff $\fa$ is a Leibniz algebra,
$\cO$ with $x\mapsto x\cdot$ is an $\fa$-module, and
$\fa\otimes\fa\to\cO$ is an $\fa$-module morphism.
}

\bigskip

\begin{pf}
By Lemma \ref{SS:conf jac} it suffices to consider the conformal Jacobi
identity for $a, b, c\in\cO\oplus\fa$.
By section \ref{SS:conf jac} the conformal Jacobi identity is equivalent to
$[a_t,b_s]c=\sum\binom{t}{i}(a_i b)_{t+s-i}c$ for $t, s\geq 0$.
Because of the gradation we only need to consider $t, s$ with $t+s\leq 1$.

For $t+s=1$, we get $a_0(b_1 c)=(a_0 b)_1 c+b_1(a_0 c)$ and
$b_1(a_0 c)=a_0(b_1 c)+(b_0 a)_1 c+(b_1 a)_0 c$.
Conformal skew-symmetry yields $(b_0 a)_1 c=-(a_0 b)_1 c+(T(b_1 a))_1 c=
-(a_0 b)_1 c-(b_1 a)_0 c$.
Thus the two identities for $t+s=1$ are equivalent. 
This also follows from Lemma \ref{SS:s3 symm}\,\itii.

Because of the gradation the identity for $t+s=1$ is only non-trivial
if $a, b, c\in\fa$. 
In this case it is just the condition that $\fa\otimes\fa\to\cO$
is a module morphism.

The identity for $t=s=0$ is only non-trivial if at least two elements
of $a, b, c$ are in $\fa$.
If $a, b, c\in\fa$ then it is the Leibniz identity for $\fa$.
If one element is in $\cO$ then we may assume that $c\in\cO$
because of Proposition \ref{SS:s3 symm conf J}.
Then it is the statement that $x\mapsto x\cdot$ is an algebra morphism.
\end{pf}

\bigskip

Thus a 1-truncated vertex Lie algebra over $\cO$ is a Leibniz algebra $\fa$
with a symmetric $\fa$-module morphism $\fa\otimes\fa\to\cO$, 
an $\fa$-module structure on $\cO$, and an even $\fa$-module morphism
$d: \cO\to\fa$ such that $[x,y]+[y,x]=d(x,y)$, $(df)g=0$, and $(x,df)=xf$.

It is clear that if $R$ is a vertex Lie algebra then $R_1$ is a 
1-truncated vertex Lie algebra over $R_0$. 
Moreover, the functor $\fa\mapsto R(\fa)$ from 1-truncated vertex Lie algebras 
to vertex Lie algebras is left adjoint to $R\mapsto R_1$.

For example, 1-truncated vertex Lie algebras $\fa$ over $\K$ with
$d=0$ and $x\cdot=0$ for any $x$ are exactly Lie algebras $\fa$ with 
an invariant symmetric bilinear form.
In this case the vertex Lie algebra $R(\fa)$ is the affine vertex Lie
algebra $\hat{\fa}=\fa[T]\oplus\K$, see section \ref{SS:aff vlie}.

\chapter{Associative Vertex Algebras}
\label{C:va}

Sections \ref{S:va def} and \ref{S:state f corres}
are about vertex algebras and their reformulation in terms of the
state-field correspondence.
We discuss associative, commutative, and Poisson vertex algebras.
We unite the algebra and the conformal algebra structure of a 
vertex algebra into one family of multiplications indexed by $\Z$.
In this way various vertex algebra identities are united as well.

In sections \ref{S:vas of distr} and \ref{S:fz alg}
we show that associative vertex algebras form a full subcategory of 
the category of local associative algebras.
We first prove that the normal ordered product and the distributions $c^i(z)$
define the structure of an associative vertex algebra on 
the space $V(\cA)$ of distributions of a local associative algebra $\cA$.
Then we construct a fully faithful functor $V\mapsto\cA(V)$ 
from associative vertex algebras to local associative algebras such that 
$V(\cA(V))=V$.

\section{Associative, Commutative, and Poisson Vertex Algebras}
\label{S:va def}

In sections \ref{SS:prelie}--\ref{SS:conf deriv va}
we discuss pre-Lie algebras, vertex algebras, the Wick formula, and 
conformal derivations of vertex algebras.
In sections \ref{SS:assoc va def}--\ref{SS:poisson va}
we show that if $V$ is an associative vertex algebra or a
vertex Poisson algebra then $V/C_2(V)$ is a Poisson algebra and 
we prove that commutative vertex algebras and 
commutative differential algebras are equivalent notions.

\subsection{Pre-Lie Algebras}
\label{SS:prelie}

We explain in which sense pre-Lie algebras generalize associative algebras
and prove that the commutator of a pre-Lie algebra is a Lie bracket.

\bigskip

Let $V$ be an algebra. 
We always use the {\bf left-operator notation} that is 
\index{left-operator notation}
defined inductively by $a_1\dots a_r:=a_1(a_2\dots a_r)$
for $a_i\in V$.
One must always be aware that in general 
$(a_1\dots a_r)\cdot a\ne a_1\dots a_ra$.

The {\bf commutator} of $V$ is
\index{commutator}
$$
[a,b]\subast
\; :=\;
ab
\, -\,
\paraab\, ba.
$$  
We often just write $[\,,]$ instead of $[\,,]\subast$.

\bigskip

{\bf Definition.}\: 
A {\bf pre-Lie algebra} is 
\index{pre-Lie algebra}
an algebra $V$ such that left multiplication $V\to\End(V), a\mapsto a\cdot$, 
is an algebra morphism with respect to the commutator on $V$ and $\End(V)$:\;
$[a,b]c=[a\cdot,b\cdot]c$.

\bigskip

Note that an algebra $V$ is associative iff $V\to\End(V), a\mapsto a\cdot$,
is an algebra morphism and a skew-symmetric algebra $\fg$ is a Lie algebra iff 
$\fg\to\fgl(\fg), a\mapsto [a,\;]$, is an algebra morphism.
Since $\rho(ab)=(\rho a)(\rho b)$ implies $\rho[a,b]=[\rho a,\rho b]$,
any associative algebra is a pre-Lie algebra.

What we have defined are {\it left} pre-Lie algebras.
Right pre-Lie algebras satisfy $a[b,c]=a[\cdot b,\cdot c]$.

The {\bf associator} of an algebra is 
\index{associator}
$\fa(a,b,c):=(ab)c-abc$.

\bigskip

{\bf Proposition.}\: {\it
\iti\:
An algebra is a pre-Lie algebra iff its associator is symmetric in the 
first two arguments: \; $\fa(a,b,c)=\paraab\fa(b,a,c)$. 

\smallskip

\itii\:
An algebra $V$ is associative iff $V$ is a pre-Lie algebra and
$[a,\;]$ is a derivation of $V$ for any $a$. 
}

\bigskip

\begin{pf}
\iti\:
We have $(ab)c-(ba)c=abc-bac$ iff $(ab)c-abc=(ba)c-bac$.

\smallskip

\itii\:
We already remarked that any associative algebra is a pre-Lie algebra.
This also follows from \iti.

The map $[a,\;]$ is a derivation iff $abc-(bc)a=[a,b]c+bac-bca$. 
Thus any two of the following three properties imply the third one:
$[a,\;]$ is a derivation; $[a\cdot,b\cdot]c=[a,b]c$; and $(bc)a=bca$.
This implies the claim.
\end{pf}

\bigskip

{\bf Remark.}\: {\it
Let $V$ be a pre-Lie algebra.

\smallskip

\iti\:
The algebra $V\oplus\K 1$, obtained by adjoining an identity $1$ to $V$,
is a pre-Lie algebra as well.

\smallskip

\itii\:
The commutator of $V$ is a Lie bracket.

\smallskip

\itiii\:
If $e, f$, and $ef$ are central then $(ab)e=abe$ and $(ef)a=efa$.

\smallskip

\itiv\:
If the multiplication of $V$ is commutative then $V$ is associative.
}

\bigskip

\begin{pf}
\iti\:
This follows from part \iti\ of the Proposition and from 
$\fa(1,a,b)=\fa(a,1,b)=\fa(a,b,1)=0$.

\smallskip

\itii\:
By \iti\ the map 
$V\oplus\K 1\to\fgl(V\oplus\K 1),\linebreak[1] a\mapsto a\cdot$, is a
monomorphism.
Since $\fgl(V\oplus\K)$ is a Lie algebra, so is $V$.

\smallskip

\itiii\:
We have $(ab)e=eab=aeb=abe$ and $(ef)a=afe=(af)e=efa$.

\smallskip

\itiv\:
This follows from \itiii\ and also from part \itii\ of the Proposition.
\end{pf}

\subsection{Vertex Algebras}
\label{SS:va}

We define vertex algebras and some standard notions for them.

\bigskip

A {\bf vertex algebra} is 
\index{vertex algebra}
a differential algebra with a $\la$-bracket.
In other words, it is an algebra $V$ with an even derivation $T$ and 
an even linear map $V\otimes V\to V[\la], a\otimes b\mapsto [a\subla b]$, 
such that $[T,[a\subla\;]]=[(Ta)\subla\;]=-\la[a\subla\;]$, 
see section \ref{SS:confa}.

An {\bf unbounded vertex algebra} is 
\index{unbounded!vertex algebra}
a differential algebra $V$ with an unbounded $\la$-bracket
$V\otimes V\to V\pau{\la}$.

If $V$ is a vertex algebra then $V\subast$ denotes the differential algebra 
and $V\subslie$ the conformal algebra underlying $V$.
A {\it morphism} $V\to W$ of vertex algebras 
is a morphism $V\subast\to W\subast$ that is also a morphism
$V\subslie\to W\subslie$.

An {\bf identity} $1$ of a vertex algebra $V$ is 
\index{identity element}
an identity of $V\subast$.
We 
\index{unital}
call $V$ {\bf unital} if $V$ has an identity.
A {\it morphism} of unital vertex algebras 
is a morphism such that $1\mapsto 1$.
If $V$ is unital then a {\bf unital} vertex subalgebra of $V$ is assumed to
\index{unital!vertex subalgebra}
contain the identity of $V$.

If $1$ is an identity then $T1=0$ since $T$ is a derivation of $V\subast$
and hence $T1=T(1\cdot 1)=2T1$.
Thus $1$ is contained in the centre of $V\subslie$ by 
Proposition \ref{SS:finite confa}.
It follows that adjoining an identity $V\mapsto V\oplus\K 1$ is 
left adjoint to the forgetful functor 
from unital vertex algebras to vertex algebras.

An {\bf identity} of an {\it unbounded} vertex algebra $V$ is
\index{identity element}
an identity of $V\subast$ that lies in the centre of $V\subslie$.

A {\bf derivation} of a vertex algebra $V$ is a derivation of $V\subast$ 
\index{derivation!of a vertex algebra}
that is also a derivation of $V\subslie$.

A {\bf left ideal} of $V$ is
\index{ideal}
a left ideal of $V\subast$ that is also a left ideal of $V\subslie$.
Right ideals and (two-sided) ideals are defined in the same way.
If $I$ is an ideal then $V/I$ is 
\index{quotient}
the {\bf quotient} vertex algebra.

A {\bf gradation} of a vertex algebra $V$ is 
\index{gradation!of a vertex algebra}
a gradation of $V\subast$ that is also a gradation of $V\subslie$,
see section \ref{SS:confa}.
If $V$ is graded and $1$ an identity then $1\in V_0$ by Remark \ref{S:supalg}.

The intersection of vertex subalgebras is a vertex subalgebra.
If $S$ is a subset of a unital vertex algebra
\index{Saaa@$\sqbrack{S}$}
then $\sqbrack{S}$ denotes the unital vertex subalgebra generated by $S$.

A differential algebra $V$ is 
\index{CFT-type}
of {\bf CFT-type} if $V$ is unital and $\rho\N$-graded for some $\rho\in\Q_>$ 
such that $V_0=\K$.
A vertex algebra $V$
\index{CFT-type}
is of {\bf CFT-type} if
$V$ is graded such that $V\subast$ is of CFT-type.
Since $1\in V_0$ and $T1=0$, this is equivalent to 
$V\subslie$ being of CFT-type, see section \ref{SS:confa}.

A {\bf $U(1)$-vector} of a vertex algebra $V$ is 
\index{U1vector@$U(1)$-vector}
an even vector $J$ such that $J\subla J=k\la$ for some $k\in\K$. 
In other words, it is a $U(1)$-vector of $V\subslie$ such that
$\hk_J=J_1 J\in\K$, see section \ref{SS:chodos thorn}.

\subsection{Wick Formula}
\label{SS:wick}

We prove that the Wick formula implies that $a_t$ are differential operators.
One application of this is the BV-operator $G_1$, see section \ref{SS:bv algs}.

\bigskip

Recall that a conformal operator of a $\K[T]$-module $V$ is 
a power series $a\submu\in\fgl(V)\pau{\mu}$ such that 
$[T,a\submu]=-\mu a\submu$ and $a\submu b\in V[\mu]$ for any $b\in V$,
see section \ref{SS:conf op}.

Let $V$ be a vertex algebra.
The {\bf Wick formula} for a conformal operator $d\submu$ is
\index{Wick formula}
$$
d\submu(ab)
\; =\;
(d\submu a)b
\; +\;
\zeta^{da}\,a(d\submu b)
\; +\;
\int_0\upmu d\nu\, [(d\submu a)\subnu b].
$$
This is equivalent to
$$
d_t(ab)
\; =\;
(d_t a)b
\; +\;
\zeta^{da}\,
a(d_t b)
\; +\;
\sum_{i=0}^{t-1}\:
\binom{t}{i}\,
(d_i a)_{t-1-i}b
$$
since
\begin{align}
\notag
\int_0\upmu d\nu\, [(d\submu a)\subnu b]
\; &=\;
\sum_{t, i\geq 0} (d_i a)_t b\; \mu^{(t+i+1)}\binom{t+i+1}{i}
\\
\notag
\; &=\;
\sum_{t\geq 0}\:\sum_{i=0}^{t-1}
\binom{t}{i}\,(d_i a)_{t-1-i}b\; \mu^{(t)}.
\end{align}
The {\bf Wick formula} for a vertex algebra is 
\index{Wick formula!for a vertex algebra}
the special case $d\submu=[a\submu\;]$:
$$
[a\subla bc]
\; =\;
[a\subla b]c
\; +\;
\zeta^{ab}\,b[a\subla c]
\; +\;
\int_0\upla d\mu\, [[a\subla b]\submu c].
$$
We call this also 
\index{left!Wick formula}
the {\bf left Wick formula}.
There is a right Wick formula for $[ab\subla c]$, 
see section \ref{SS:right wick}.

Let $V$ be an algebra.
A {\bf homogeneous differential operator} of order $\leq n$ is 
\index{differential operator}
an operator $d$ such that $[d,a\cdot]-(da)\cdot$ is a 
\index{homogeneous differential operator}
homogeneous differential operator of order $\leq n-1$ for any $a\in V$.
Here $n\geq 1$ and, by definition, 
$d$ is a homogeneous differential operator of order $\leq 0$ if $d=0$.

An operator $d$ is a homogeneous differential operator of order $\leq 1$ iff 
$d$ is a derivation. 
Let $\cD_n$ be the space of homogeneous differential operators of 
order $\leq n$.
Then $\cD_n\subset\cD_{n+1}$, $\cD_n \cD_m\subset\cD_{n+m}$, and 
$[\cD_n,\cD_m]\subset\cD_{n+m-1}$.

\bigskip

{\bf Proposition.}\: {\it
Let $V$ be a vertex algebra that satisfies the Wick formula.
Then $a_t$ is a homogeneous differential operator of $V\subast$ of 
order $\leq t+1$ for any $a\in V$ and $t\geq 0$.
In particular, $a_0$ is a derivation of $V\subast$.
}

\bigskip

\begin{pf}
The Wick formula is 
$[a_t,b\cdot]-(a_t b)\cdot=\sum_{i=0}^{t-1}\binom{t}{i}(a_i b)_{t-1-i}$.
Thus the first claim follows by induction on $t$.
The operator $a_0$ is a derivation of the differential algebra $V\subast$
since $[T,a_0]=0$ by Remark \ref{SS:conf der}. 
\end{pf}

\subsection{Conformal Derivations of Vertex Algebras}
\label{SS:conf deriv va}

We prove that conformal derivations of a vertex algebra form
an unbounded vertex Lie algebra $\Der\subsc(V)$.

\bigskip

Recall that a conformal derivation of a conformal algebra is 
a conformal operator $d\submu$ such that 
$d\submu[a\subla b]=[(d\submu a)_{\la+\mu}b]+\zeta^{da}[a\subla(d\submu b)]$,
see section \ref{SS:conf der}.

A {\bf conformal derivation} of a vertex algebra $V$ is 
\index{conformal derivation}
a conformal derivation of $V\subslie$ that satisfies the Wick formula.
In Remark \ref{SS:conf j wick} we give an equivalent definition of 
conformal derivations.

From Remark \ref{SS:conf der} and by setting $\mu=0$ in 
the Wick formula, we obtain:

\bigskip

{\bf Remark.}\: {\it
Let $d\submu$ be a conformal derivation of a vertex algebra $V$.
Then $d_0$ is a derivation of $V$.
\hfill $\square$
}

\bigskip

{\bf Proposition.}\: {\it
The space $\Der\subsc(V)$ of conformal derivations of a vertex algebra $V$ is 
an unbounded vertex Lie subalgebra of $\Der\subsc(V\subslie)$.
}

\bigskip

\begin{pf}
That $\Der\subsc(V)$ is a $\K[T]$-submodule is clear since $T=-\mu\cdot$.
The Wick formula is equivalent to 
$[d\submu, a\cdot]=(d\submu a)\cdot+\int_0\upmu d\nu\, (d\submu a)\subnu$.
Thus the claim follows from 
\begin{align}
\notag
&[(d\submu)\subla d'\submu, a\cdot]
\; =\; 
[[d\subla, d'_{\mu-\la}], a\cdot]
\\
\notag
\; =\; 
&([d\subla, d'_{\mu-\la}]a)\cdot
+
\int_0\upla d\nu\, (d\subla d'_{\mu-\la}a)\subnu\;
+
\int_0^{\mu-\la} d\nu\, (d'_{\mu-\la}d\subla a)\subnu\; 
\\
\notag
&\qquad\qquad\qquad\,
+
\int_0^{\mu-\la} d\nu\, [d\subla,(d'_{\mu-\la}a)\subnu]
+
\int_0\upla d\nu\, [d'_{\mu-\la},(d\subla a)\subnu]
\\
\notag
=\;
&(((d\submu)\subla d'\submu)a)\cdot
+
\int_0\upmu d\nu\, (((d\submu)\subla d'\submu)a)\subnu.
\end{align}
In the last step we used that $[d\submu,a\subla]=(d\submu a)_{\la+\mu}$ and
applied the substitution formula from section \ref{SS:bracket vlie}.
\end{pf}

\subsection{Associative Vertex Algebras}
\label{SS:assoc va def}

We prove that if $V$ is an associative vertex algebra then 
$V\subast$ is almost commutative and $V\subast/C_2(V)$ is a Poisson algebra,
called
\index{Zhu Poisson algebra}
the {\bf Zhu Poisson algebra}. 
In section \ref{SS:zhu alg2}
we show that another quotient $A(V)=V/O(V)$ is a very interesting 
associative algebra.

\bigskip

The commutator $[\, ,]\subast$ of a pre-Lie algebra is a Lie bracket
by Remark \ref{SS:prelie}\,\itii.
On the other hand, any vertex Lie algebra has 
a natural Lie bracket $[\, ,]\subslie$ by Proposition \ref{SS:bracket vlie}.

\bigskip

{\bf Definition.}\: 
A 
\index{associative!vertex algebra}
vertex algebra $V$ is {\bf associative}
if $V\subast$ is a pre-Lie algebra, $V\subslie$ is a vertex Lie algebra, 
$[\, ,]\subast=[\, ,]\subslie$, and $[a\subla\;]\in\Der\subsc(V)$ for any $a$.

\bigskip

In other words, an associative vertex algebra is a differential pre-Lie algebra
with a Lie $\la$-bracket such that $[\, ,]\subast=[\, ,]\subslie$
and $[a\subla\;]\in\Der\subsc(V)$.

What we have defined are {\it left} associative vertex algebras:
differential {\it left} pre-Lie algebras with a {\it left} $\la$-bracket 
such that \dots\,.
We will have no occasion to consider right associative vertex algebras.

Let $V$ be an associative vertex algebra.
The {\bf Borcherds Lie algebra} of $V$ is 
\index{Borcherds Lie algebra}
$\fg(V):=\fg(V\subslie)$, the Borcherds Lie algebra of $V\subslie$, 
see section \ref{SS:vlie to loclie}.

We note that the unital vertex algebra $V\oplus\K 1$ is 
associative by Remark \ref{SS:prelie}\,\iti.
Any one-sided ideal $I$ of $V$ is a two-sided ideal because of 
conformal skew-symmetry and because 
$ab=ba+[a,b]$ and $[a,b]\in I$ for any $a\in V, b\in I$.

\medskip

Let $V$ be a differential algebra.
We call $V$ {\bf almost commutative} if
\index{almost commutative}
$[T^n V, T^m V]\subset T^{n+m+1}V$ for any $n, m\geq 0$.
In particular, $[V,V]\subset TV$.

Let $C_n(V)\subset V$ be the ideal generated by $T^{n-1}V$ where $n\geq 2$.
Of course, $C_{n+1}(V)\subset C_n(V)$.
If $V$ is associative and almost commutative then 
$C_n(V)=(T^{n-1}V)V$ since $a(T^m b)c=(T^m b)ac+[a,T^m b]c$.

Let $V$ be an associative vertex algebra.
Then $V\subast$ is almost commutative because
$[T^n a,T^m b]=\int_{-T}^0d\la\: (-\la)^n (T+\la)^m[a\subla b]$.

Define $C_n(V):=C_n(V\subast)$.
In Proposition \ref{SS:q assoc zhu poisson}\,\iti\ we prove that 
$C_n(V)=(T^{n-1}V)V$.
A vertex algebra $V$ is {\bf $C_n$-cofinite} if 
\index{Cn-cofinite@$C_n$-cofinite}
$\dim V/C_n(V)<\infty$.

Recall that a {\bf Poisson algebra} is 
\index{Poisson algebra}
a commutative algebra with a Lie bracket $\set{\, ,}$ such that 
$\set{a,bc}=\set{a,b}c+\paraab b\set{a,c}$.

\bigskip

{\bf Proposition.}\: {\it
Let $V$ be an associative unital vertex algebra.
The algebra $V\subast/C_2(V)$ is a Poisson algebra with Poisson bracket
induced by $a_0 b$.
}

\bigskip

\begin{pf}
The multiplication of $V\subast/C_2(V)$ is commutative since 
$V\subast$ is almost commutative and hence $[V,V]\subset TV$. 
By Remark \ref{SS:prelie}\,\itiv\ 
the pre-Lie algebra $V/C_2(V)$ is associative.

By Remark \ref{SS:vlie}, the $0$-th product induces a Lie bracket on $V/TV$.
The subspace $C_2(V)$ is a left ideal with respect to $a_0 b$ since
$a_0$ is a derivation of $V\subast$ by Proposition \ref{SS:wick}.
It is a right ideal since $a_0 b+b_0 a\in TV$.
Since $a_0$ is a derivation of $V\subast$, it follows that $V/C_2(V)$ is
a Poisson algebra. 
\end{pf}

\subsection{Commutative Vertex Algebras and Cosets}
\label{SS:comm va}

We prove that commutative vertex algebras and 
commutative differential algebras are equivalent notions.
We show that if $V$ is an $\N$-graded associative vertex algebra then
$V_0$ is a commutative algebra.

Examples of commutative vertex algebras are vertex Poisson algebras and 
cohomology vertex algebras, 
see sections \ref{SS:poisson va} and \ref{SS:topl va}. 
We treat commutativity again in section \ref{SS:commut centre}.

\bigskip

Elements $a, b$ of an algebra {\bf commute} if $ab=\paraab ba$.
Elements $a, b$ of a conformal algebra commute if 
$[a\subla b]=[b\subla a]=0$, see section \ref{SS:finite confa}.

Elements $a, b$ of 
\index{commuting!elements}
a vertex algebra $V$ {\bf commute} if $a, b\in V\subast$ and 
$a, b\in V\subslie$ commute.
If $V$ is associative then $a, b$ commute iff $a, b\in V\subslie$ commute 
since $[a\subla b]=0$ implies $[a,b]\subast=[a,b]\subslie=0$.
Note that $a$ need {\it not} commute with $a$, even if $V$ is associative.

\bigskip

{\bf Definition.}\: 
A 
\index{commutative!vertex algebra}
vertex algebra is {\bf commutative} if its elements commute.

\bigskip

After this subsection all {\bf commutative} vertex algebras are 
\index{commutative!vertex algebra}
assumed to be associative and unital.

\bigskip

{\bf Proposition.}\: {\it
The functor $V\mapsto V\subast$ is an isomorphism of categories
from commutative vertex algebras to differential algebras with
commutative multiplication.
Moreover, $V$ is associative iff $V\subast$ is. 
}

\bigskip

\begin{pf}
The first claim is clear.
A commutative vertex algebra $V$ is associative iff 
$V\subast$ is a pre-Lie algebra iff $V\subast$ is associative
by Remark \ref{SS:prelie}\,\itiv.
\end{pf}

\bigskip

We use the Proposition to identify the objects of these two categories.

The {\bf coset} of a subset $S$ of a vertex algebra $V$ is
\index{coset}
$$
C_V(S)
\; :=\;
\set{a\in V\mid \text{$a$ commutes with any $b\in S$}}.
$$
Two subsets $S$ and $S'$ {\bf commute} if 
\index{commuting!subsets}
$S\subset C_V(S')$.

The {\bf centre} of $V$ 
\index{centre!of a vertex algebra}
is $C(V):=C_V(V)$.
The elements of the centre are 
\index{central element!of a vertex algebra}
called {\bf central}.
The Wick formula and the conformal Jacobi identity imply:

\bigskip

{\bf Lemma.}\: {\it
Let $V$ be an associative vertex algebra.
If $a, b$ and $a, c$ commute then $a, bc$ and $a, b_t c$ and $a, Tb$ commute
for any $t\geq 0$.
\hfill $\square$
}

\bigskip

Assume that $V$ is associative.
The Lemma implies that $C_V(S)$ is a vertex subalgebra and 
$C_V(S)=C_V(\sqbrack{S})$.
In particular, $C(V)$ is a commutative vertex subalgebra.
If the elements of $S$ pairwise commute then 
$\sqbrack{S}$ is commutative and equal to the differential subalgebra
generated by $S$.

Proposition \ref{SS:finite confa} shows that
the torsion $\K[T]$-submodule is contained in $C(V)$. 
In particular, $\ker T\subset C(V)$ and 
any finite-dimensional associative vertex algebra is commutative.

If $V$ is $\N$-graded then any $a, b\in V_0$ commute
since $a_t b\in V_{-t-1}=0$ for $t\geq 0$. 
Thus $V_0$ is a commutative associative subalgebra of $V\subast$
by Remark \ref{SS:prelie}\,\itiv.
This follows also from Proposition \ref{SS:assoc va def} and $C_2(V)_0=0$.

\subsection{Vertex Poisson Algebras}
\label{SS:poisson va}

We prove that if $V$ is a vertex Poisson algebra then $V/C_2(V)$ is 
a Poisson algebra.

Many associative vertex algebras possess filtrations such that
the associated graded vertex algebra is a vertex Poisson algebras.
Any associative vertex algebra has a decreasing filtration, 
any $\Z$-graded associative vertex algebra has an increasing filtration,
and any enveloping vertex algebra has an increasing filtration of this type,
see sections \ref{S:filt span sets} and \ref{SS:vertex envelope}.

\bigskip

Recall again that a Poisson algebra is 
a commutative algebra with a Lie bracket $\set{\, ,}$ such that 
$\set{a,bc}=\set{a,b}c+\paraab b\set{a,c}$.

\bigskip

{\bf Definition.}\: 
A {\bf vertex Poisson algebra} is 
\index{vertex Poisson algebra}
a commutative vertex algebra with a Lie $\la$-bracket $\set{\subla}$ such that 
$\set{a\subla bc}=\set{a\subla b}c+\paraab\,b\set{a\subla c}$.

\bigskip

The $\la$-bracket of a vertex Poisson algebra is 
\index{Poisson!lambdabracket@$\la$-bracket}
called a {\bf Poisson $\la$-bracket}.

Note that the identity 
$\set{a\subla bc}=\set{a\subla b}c+\paraab\,b\set{a\subla c}$
is equivalent to $\set{a\subla\;}\in\Der\subsc(V)$, since $V$ is commutative.
If we write $\set{a\subla\;}=\sum a_t\,\la^{(t)}$
then $a_t$ is, in general, {\it not} a derivation of 
the commutative vertex algebra $V$ since, in general, $[T,a_t]\ne 0$ for $t>0$.

The identity $\set{a\subla bc}=\set{a\subla b}c+\paraab\,b\set{a\subla c}$ 
is equivalent to
$$
\set{ab\subla c}
\; =\;
(e^{T\del\subla}a)\set{b\subla c}
\; +\;
\paraab\, 
(e^{T\del\subla}b)\set{a\subla c}
$$
since $e^{T\del\subla}\set{c_{-\la}ab}=
e^{T\del\subla}(\set{c_{-\la}a}b+a\set{c_{-\la}b})$.

Recall that if $V$ is a commutative vertex algebra then $C_2(V)=(TV)V$ is 
the ideal generated by $TV$, see section \ref{SS:assoc va def}.

\bigskip

{\bf Proposition.}\: {\it
Let $V$ be a vertex Poisson algebra.
Then $V/C_2(V)$ is a Poisson algebra with Poisson bracket
induced by $\set{a_0 b}$.
}

\bigskip

\begin{pf}
By Remark \ref{SS:vlie}, the $0$-th product induces a Lie bracket on $V/TV$.
By Remark \ref{SS:conf der} the map $\set{a_0\;}$ is a $\K[T]$-module morphism.
It is clear that $\set{a_0\;}$ is a derivation of the commutative algebra $V$. 
Thus $C_2(V)=(TV)V$ is a left ideal with respect to $\set{a_0 b}$.
It is a right ideal since $\set{a_0 b}+\set{b_0 a}\in TV$.
It follows that $V/C_2(V)$ is a Poisson algebra. 
\end{pf}

\section{Modules}
\label{S:state f corres}

In sections \ref{SS:fields}--\ref{SS:id elem}
we discuss the unbounded endomorphism vertex algebra, modules, and the
adjoint module.
We show that vertex algebras and
bounded $\Z$-fold algebras with a translation operator are 
equivalent notions.

In sections \ref{SS:ratl fcts} and \ref{SS:z fold fields}
we show that $M$ is a module iff $M$ satisfies the associativity formula.

In sections \ref{SS:j id conf j id} and \ref{SS:conf j wick0}
we show that the conformal Jacobi identity and the Wick formula
are special cases of the commutator and associativity formula, 
which in turn are special cases of the Jacobi identity.
In section \ref{SS:skew sym}
we prove that conformal skew-symmetry and $[\, ,]\subast=[\, ,]\subslie$
are equivalent to skew-symmetry.

\subsection{The Unbounded Endomorphism Vertex Algebra}
\label{SS:fields}

The algebra $\End(E)\lau{z}$ does not contain 
non-commuting local distributions, 
while $\End(E)\pau{z\uppm}$ is not an algebra. 
We show that the space $\End\subsv(E)$ of fields is an unbounded vertex 
algebra.

In section \ref{S:vas of distr} we prove that 
$\End\subsv(E)$ is a pre-Lie algebra, that the conformal Jacobi identity
and the Wick formula hold, and that two fields satisfy conformal skew-symmetry 
and $[\,,]\subast=[\,,]\subslie$ iff they are local.

\bigskip

Let $E$ be a vector space.
In section \ref{SS:alg of distr} we remarked that
the space $\Endv(E)$ of fields is 
an unbounded conformal subalgebra of $\fgl(E)\pau{z\uppm}$.
Recall that a field is 
\index{field}
a distribution $a(z)\in\End(E)\pau{z\uppm}$ such that $a(z)b\in E\lau{z}$ 
for any $b\in E$.
We have $T=\del_z$ and $[a(w)\subla b(w)]=\res_z e^{(z-w)\la}[a(z),b(w)]$.

The space $\Endv(E)$ contains the algebra $\End(E)\lau{z}$.
These two spaces are equal iff $\dim E<\infty$.
For example, $\sum_{i\geq 0}\del_x^i\, z^{-i}$ is a field on $\K[x]$ 
that is not contained in $\End(\K[x])\lau{z}$.

The product $a(z)b(z)$ of two fields is in general not well-defined.
In other words, the operator product $a(z)b(w)$ is singular for $z=w$.
In order to explain this, let us introduce the following notation. 

Denote by $\cF[z\uppm](E)$ the space of distributions 
$a(z)\in\End(E)\pau{z\uppm}$ such that $a(z)b\in E[z\uppm]$ for any $b\in E$.
The spaces $\cF[z](E), \cF\lau{z,w}(E)$ etc. are defined in the same way.
Thus $\Endv(E)=\cF\lau{z}(E)$.

Distributions in $\cF\lau{z_1, \dots, z_r}(E)$ are 
called {\bf regular} at $z_1=\ldots=z_r$.
If $c(z_1, \dots, z_r)$ is regular then $c(z, \dots, z)$ is 
a well-defined field.
In particular, if $a(z)b(w)$ is regular at $z=w$ then 
the product $a(z)b(z)$ is well-defined. 

In general, the operator product $a(z)b(w)$ of two fields is {\it not}
contained in
$\cF\lau{z,w}(E)$, but only in the larger space $\cF\lau{z}\lau{w}(E)$.
An example for $E=\K[x]$ is
$$
\Big(\sum_{i\geq 0} \del_x^i z^{-i}\Big)
\Big(\sum_{j\geq 0} (x^{(j)}\cdot)w^j\Big)1
\; =\; 
\sum_{i, j\geq 0} x^{(j-i)}\: z^{-i}w^j.
$$ 

However, $a(z)b(w)$ is regular at $z=w$ if $a(z), b(z)\in\End(E)\lau{z}$.
More generally, this is true if $a(z)\in\End(E)\lau{z}$ and $b(z)\in\Endv(E)$.
It is also true if $a(z)\in\Endv(E)$ and $b(z)\in\cF[z\uppm](E)$.
We make use of this observation as follows.

There is a decomposition $E\pau{z\uppm}=z\inv E\pau{z\inv}\oplus E\pau{z}$ 
of $\K[\del_z]$-modules.
In particular, we have $\Endv(E)=\Endv(E)_+\oplus\End(E)\pau{z}$ where 
\index{EndvEplus@$\Endv(E)_+$}
$\Endv(E)_+:=z\inv\cF[z\inv](E)$.
Denote by $a(z)\mapsto a(z)\subpm$ the projections so that 
$a(z)_+=\sum_{t\geq 0}a_t z^{-t-1}$.

The {\bf normal ordered product} of 
\index{normal ordered product!of fields}
two fields $a(z)$ and $b(z)$ is
$$
\normord{a(z)b(w)}
\; :=\;
a(z)_-b(w)\, +\,\paraab b(w)a(z)_+.
$$
It is regular at $z=w$ and hence $\normord{a(z)b(z)}$ is a field, 
which we also call the {\bf normal ordered product}.
Thus $\Endv(E)$ is an unbounded vertex algebra.
We call $\Endv(E)$ the 
\index{unbounded!endomorphism vertex algebra}
{\bf unbounded endomorphism vertex algebra} of $E$.
The field $\id_E$ is an identity of $\Endv(E)$.

For example, if $E\ne 0$ then $\K\subset\End(E)$ and we can identify
$\K\lau{z}$ with a subspace of $\Endv(E)$.
Since $[a(z),b(w)]=0$ for $a(z)\in\K\lau{z}$ and $b(z)\in\Endv(E)$,
we have $\normord{a(z)b(z)}=a(z)b(z)=\;\normord{b(z)a(z)}$ and
$\K\lau{z}$ is contained in the centre of $\Endv(E)$.
Moreover, the commutative differential algebra $\K\lau{z}$ with $T=\del_z$
is a commutative vertex subalgebra of $\Endv(E)$.
If $\dim E=1$ then $\Endv(E)=\K\lau{z}$.

If $E$ is a $\K[T]$-module then we denote by 
$\Endv(E)_T:=\Endv(E)\cap\End(E)\pau{z\uppm}_T$ 
\index{EndvET@$\Endv(E)_T$}
the space of 
\index{translation covariant!field}
translation covariant fields.
This is an unbounded vertex subalgebra because
$[T,\, ]$ and $\del_z$ are derivations.

\subsection{Modules}
\label{SS:modules}

We reformulate modules in terms of a multiplication $a\cdot$ and a
$\la$-action $a\subla$ 
and consider modules over a commutative vertex algebra as an example.

\bigskip

{\bf Definition.}\: 
A {\bf module} over a vertex algebra $V$ is 
\index{module}
a vector space $M$ together with a morphism $Y_M: V\to\End\subsv(M)$.

\bigskip

If $V$ is unital then we shall always assume that $Y_M(1)=\id_M$.
We often just write $a(z)$ instead of $Y_M(a)$.
In section \ref{SS:z fold fields}
we show that $M$ is a module iff $M$ satisfies the associativity formula.
In Corollary \ref{SS:skew sym fields} we prove that if $V$ is associative
then any $V$-module is a $V\subslie$-module.

Let $V$ be a $\K[T]$-module and $M$ a vector space.
To give a $\K[T]$-module morphism $Y_M: V\to\End(M)\pau{z\uppm}$  
is equivalent to giving two $\K[T]$-module morphisms $V\to\End(M)\pau{z}$ and 
$V\to z\inv\End(M)\pau{z\inv}, a(z)\mapsto a(z)_+$. 

To give a $\K[T]$-module morphism $V\to\End(M)\pau{z}, a\mapsto a(z)$, 
is equivalent to giving an even linear map $V\to\End(M), 
a\mapsto a\cdot =a(0)$.
In fact, we have
$$
a(z)=e^{z\del_w}a(w)|_{w=0}=(e^{zT}a)(w)|_{w=0}=(e^{zT}a)\cdot.
$$
Conversely, the map $a\mapsto(e^{zT}a)\cdot$ is a $\K[T]$-module morphism
since $(e^{zT}Ta)\cdot\linebreak[0]=\del_z (e^{zT}a)\cdot$.

In section \ref{SS:unb confa power series} we showed that 
$\fgl(M)\pau{\la}$ is an unbounded vertex Lie algebra with $T=-\la\cdot$ and 
that there is a $\K[T]$-module isomorphism 
$$
z\inv\fgl(M)\pau{z\inv}\:\to\:\fgl(M)\pau{\la},\quad
a(z)\:\mapsto\:\res_z e^{z\la}a(z).
$$
We use it to identify these two spaces. 
Thus we write $a(z)_+=a\subla$.  

Summarizing, to give a $\K[T]$-module morphism $Y_M: V\to\End(M)\pau{z\uppm}$
is equivalent to giving an even linear map $V\to\End(M)$ and a morphism 
$V\to\fgl(M)\pau{\la}$. 
We have
$$
a(z)
\; =\;
(e^{zT}a)\cdot 
\; +\;
a\subla.
$$

In Proposition \ref{SS:q assoc assoc} we prove that 
a $\K[T]$-module morphism $V\to\Endv(M)$ defines a module iff
$a\cdot$ and $a\subla$ satisfy four identities.
The formulation of modules in terms of $\Endv(M)$ seems to be better.

As an example, let us consider a commutative vertex algebra $V$ and
$V$-modules $M$ such that $a\subla=0$.
In other words, $Y_M(V)\subset\End(M)\pau{z}$.
Such $M$ are the same as $V\subast$-modules. 
In fact, we have $[a(z)\subla b(z)]=0$ and 
$\normord{a(z)b(z)}=a(z)b(z)$ for $a(z), b(z)\in\End(M)\pau{z}$.
Thus the claim follows from $e^{zT}(ab)=(e^{zT}a)(e^{zT}b)$.

But not all $V$-modules satisfy $a\subla=0$.
For example, in section \ref{SS:fields} we have shown that 
$\Endv(\K)=\K\lau{z}$ as commutative vertex algebras.
Thus $\K$ is a module over the commutative vertex algebra $\K\lau{z}$.

\subsection{The Adjoint Module}
\label{SS:adj module}

We prove that vertex algebras and bounded $\Z$-fold algebras with 
a translation operator are equivalent notions, and
we discuss the relation between the properties that a vertex algebra $V$ is 
associative and that $V$ is a $V$-module.

\bigskip

To give an algebra is equivalent to giving a vector space $V$ together with 
an even linear map $V\to\End(V), a\mapsto a\cdot$.
The algebra is associative iff this map is an algebra morphism iff
$V$ is a $V$-module.

Recall that a $\Z$-fold module over a vector space $V$
is a vector space $M$ together with an even linear map 
$$
Y_M:\: V\;\to\;\End(M)\pau{z\uppm}, \quad a\;\mapsto\; a(z)\;=\;
\sum a_t\, z^{-t-1},
$$
see section \ref{SS:module vlie}.
By definition, $M$ is bounded iff $Y_M(V)\subset\Endv(M)$.

In the case $M=V$ we obtain the notion of a $\Z$-fold algebra.
The map $Y:=Y_V$ is
\index{state-field correspondence}
the {\bf state-field correspondence} of $V$.

Let $V$ be a $\Z$-fold algebra.
An even operator $T$ is 
\index{translation!endomorphism}
a {\bf translation endomorphism} if $(Ta)(z)=\del_z a(z)$;
equivalently, $(Ta)_t b=-t\, a_{t-1}b$.
It is 
\index{translation!generator}
a {\bf translation generator} if $[T,a(z)]=\del_z a(z)$; 
in other words, any $a(z)\in Y(V)$ is translation covariant.
It is a derivation of $V$ iff $[T,a(z)]=(Ta)(z)$.
Any two of these three properties imply the third one.
In this case $T$ is 
\index{translation!operator}
a {\bf translation operator}.

From section \ref{SS:modules} follows that
to give a $\Z$-fold algebra $V$ with a translation endomorphism is 
equivalent to giving an algebra $V\subast$ with an even operator $T$ and 
an even linear map $V\otimes V\to V\pau{\la}, a\otimes b\mapsto [a\subla b]$, 
such that $[Ta\subla b]=-\la[a\subla b]$.
For $t\geq 0$, we have 
$$
(T^{(t)}a)b
\,=\,
a_{-1-t}b, \qquad
[a\subla b]
\,=\,
\sum a_t b\,\la^{(t)}, 
\qquad
a(z)
\,=\,
(e^{zT}a)\cdot\;+\;[a\subla\;].
$$

\bigskip

{\bf Proposition.}\: {\it
The functor $V\mapsto (V\subast,[\subla])$ is an isomorphism of categories 
from $\Z$-fold algebras with a translation operator to 
unbounded vertex algebras.
}

\bigskip

\begin{pf}
The functor is a bijection on objects because the derivation property
$[T,a(z)]=(Ta)(z)$ is obviously compatible with the projections
$a(z)\mapsto a(z)\subpm$ and 
we have $[T,(e^{zT}a)\cdot]=(e^{zT}Ta)\cdot$ iff 
$T$ is a derivation of $V\subast$.
It is a bijection on morphisms since $a_{-1-t}b=(T^{(t)}a)b$.
\end{pf}

\bigskip

We use the Proposition to identify the objects of these two categories.

The above remarks show that to give a vertex algebra is equivalent to 
giving a $\K[T]$-module $V$ with a $\K[T]$-module morphism $Y: V\to\Endv(V)_T$.
By definition, $Y$ is a morphism of unbounded vertex algebras iff
$V=(V,Y)$ is a $V$-module. 
This is
\index{adjoint module}
the {\bf adjoint} module of $V$.

In Proposition \ref{SS:assoc va} we prove that if $V$ is associative then
$V$ is a $V$-module and the converse is true if conformal skew-symmetry and 
$[\, ,]\subast=[\, ,]\subslie$ are satisfied.
In Proposition \ref{SS:unital assoc va} we prove that if $V$ is unital then 
$V$ is associative iff $V$ is a $V\subslie$-module.
In Proposition \ref{SS:eq commut} we prove that
if $V$ is unital and has a Virasoro vector with $c\ne 0$ then
$V$ is associative iff $V$ is a $V$-module.

However, there exist vertex algebras $V$ such that $V$ is a $V$-module 
but $V$ is {\it not} associative.
For example, consider a vertex algebra $V$ with $[\subla]=0$.
Then $V$ is a $V$-module iff $V$ is a $V\subast$-module iff
$V\subast$ is associative, see section \ref{SS:modules}.
The vertex algebra $V$ is associative iff 
$V\subast$ is associative and commutative.

Let $V$ be a vertex algebra.
A vector space gradation $V=\bigoplus V_h$ is a vertex algebra gradation iff 
$TV_h\subset V_{h+1}$ and $(V_h)_t V_k\subset V_{h+k-t-1}$ for any $t\in\Z$.
If $V$ is unital then $TV_h\subset V_{h+1}$ follows from $Ta=a_{-2}1$ and 
$1\in V_0$.

\subsection{Identity Element}
\label{SS:id elem}

We compare the notions of an identity, morphism, and derivation
for vertex algebras with those for $\Z$-fold algebras.

\bigskip

Let $V$ be a $\Z$-fold algebra.
An even vector $1$ is a {\bf left identity} if
\index{left identity}
$1(z)=\id_V$.
It is a {\bf right identity} if 
\index{right identity}
$a(z)1=e^{zT}a$ for some even operator $T$.
Comparing coefficients of $z$ shows that $Ta=a_{-2}1$.
An {\bf identity} is 
\index{identity!of a Zfold algebra@of a $\Z$-fold algebra}
a left and right identity.
We call $V$ {\bf unital} if 
\index{unital!Zfold algebra@$\Z$-fold algebra}
$V$ has an identity.

Recall that an identity of an unbounded vertex algebra $V$
is an identity $1$ of $V\subast$ that lies in the centre of $V\subslie$,
see section \ref{SS:va}.
In this case $T1=0$ and $Ta=a_{-2}1$.
From $a(z)=(e^{zT}a)\cdot+[a\subla\;]$ follows that 
$1$ is an identity of $V$ iff $1$ is an identity of the $\Z$-fold algebra $V$.

\bigskip

{\bf Remark.}\: {\it
Let $V$ and $V'$ be unital unbounded vertex algebras.

\smallskip

\iti\:
A linear map $\phi: V\to V'$ is a morphism of unital unbounded vertex algebras
iff $\phi$ is a morphism of unital $\Z$-fold algebras.

\smallskip

\itii\:
An operator $d$ is a derivation of the unbounded vertex algebra $V$ iff 
$d$ is a derivation of the $\Z$-fold algebra $V$.
}

\bigskip

\begin{pf}
The implication `$\Rightarrow$' follows from $a_{-1-t}b=(T^{(t)}a)b$. 
The implication `$\Leftarrow$' follows from $Ta=a_{-2}1$ and $d1=0$.
\end{pf}

\bigskip

Part \iti\ shows that if $V$ is a unital vertex algebra then 
a bounded $\Z$-fold $V$-module $M$ is a $V$-module iff
$Y_M$ is a unital $\Z$-fold algebra morphism.

\subsection{Rational Functions}
\label{SS:ratl fcts}

We discuss the canonical morphism from rational functions to Laurent series.

\bigskip

The quotient field of $\K\pau{z}$ is $\K\lau{z}$.
The field $\K(z_1,\dots,z_r)$ of {\bf rational functions} is 
\index{rational function}
the free field over $\K$ generated by $z_1, \dots, z_r$.
Thus there exists a unique morphism of fields over $\K$ 
$$
T_{z_1,\dots,z_r}:
\;
\K(z_1,\dots,z_r)
\;\to\;
\K\lau{z_1}\dots\lau{z_r}
$$
such that $z_i\mapsto z_i$. 
The maps $T_{z_1,\dots,z_r}$ and $\del_{z_i}$ commute because
both of their compositions are derivations such that $z_j\mapsto\de_{ij}$.

\bigskip

{\bf Proposition.}\: {\it
Let $f(z)\in\K(z)$ be regular at $0$ and $n\in\Z$. 
Then
$$
T_z(z^n f(z))
\; =\;
z^n e^{z\del_w}f(w)|_{w=0}.
$$
}

\bigskip

\begin{pf}
We may assume that $n=0$.
Let $R$ be the ring of rational functions that are regular at $0$.
The map $T: R\to\K\pau{z}, f(z)\mapsto e^{z\del_w}f(w)|_{w=0}$,
is an algebra morphism because
$e^{z\del_w}(f(w)g(w))=(e^{z\del_w}f(w))e^{z\del_w}g(w)$
by the product formula. 
Since $Tz=z$, we get $T=T_z$.
\end{pf}

\bigskip

The morphism $T_{z_i}:\K(z_i)\to\K\lau{z_i}$ induces $2^{r-1}$ morphisms
$$
T_{z_i;z_1, \dots, z_r}:\:
\K\cB_1\dots (z_i)\dots\cB_r
\; \to\;
\K\cB_1\dots\lau{z_i}\dots\cB_r
$$
where $\cB_j$ stands either for $(z_j)$ or $\lau{z_j}$.
Due to its universal property, $T_{z_1,\dots,z_r}$
is equal to the composition of the morphisms
$T_{z_1;z_1, \dots, z_r}, \dots, T_{z_r;z_1, \dots, z_r}$ in any order.
For example, there is a commutative diagram 
$$
\xymatrix{
\K(z,w)=\K(z)(w)
 \; \ar[rr]^{T_{z;z,w}}
   \ar[d]_{T_{w;z,w}}   
   \ar[rrd]^{T_{z,w}}  & &
\;\K\lau{z}(w)
   \ar[d]^{T_{w;z,w}}
\\
\K(z)\lau{w}
\;\ar[rr]^{T_{z;z,w}} & &
\;\K\lau{z}\lau{w}.
}
$$

For $f(z)\in\K(z)$, the Proposition yields
$$
T_{z,w}f(z+w)
\; =\;
T_{z;z,w}T_{w;z,w}f(z+w)
\; =\;
T_{z;z,w}e^{w\del_z}f(z)
\; =\;
e^{w\del_z}T_z f(z).
$$
In section \ref{SS:delta dis} we defined $a(z+w):=e^{w\del_z}a(z)$ for 
$a(z)\in E\pau{z\uppm}$.
The last identity for $f(z)=z^t$ implies 
$$
a(z+w)
\; =\;
\sum_{t\in\Z}\: a_t\, T_{z,w}(z+w)^{-t-1}.
$$
This can be viewed as a Taylor formula for $a(z)$.

We will usually write $f(z+w)$ instead of $T_{z,w}f(z+w)$.
The order of the summands of $z+w$ determines whether we apply
$T_{z,w}$ or $T_{w,z}$.
It will be clear from the context whether $f(z+w)$ stands for the
rational function or the Laurent series.
For $r\in\Z$, we have
$$
(z-w)^r
\; =\;
e^{-w\del_z}z^r
\; =\;
\sum(-1)^i\binom{r}{i}z^{r-i}w^i
$$ 
and $T_{w,z}(z-w)^r=(-1)^r(w-z)^r=\sum(-1)^{r+i}\binom{r}{i}z^i w^{r-i}$.

\subsection{Associativity Formula}
\label{SS:z fold fields}

We give a uniform formula for the $r$-th products of the $\Z$-fold algebra
$\Endv(E)$
and remark that $M$ is a module iff $M$ satisfies the associativity formula.

Later we will show that the uniform formula gives the
coefficients of the OPE, see Propositions \ref{SS:tayl}
and \ref{SS:ope2} and section \ref{SS:ope heuris}.

\bigskip

Let $a(z), b(z)$ be fields on a vector space $E$. 
Then $a(z)b(w)$ is contained in the space $\cF\lau{z}\lau{w}(E)$ 
defined in section \ref{SS:fields}.
Since $\cF\lau{z}\lau{w}(E)$ is a vector space over $\K\lau{z}\lau{w}$,
the product
$$
f(z,w)a(z)b(w)
\; :=\;
T_{z,w}(f(z,w))\, a(z)b(w)
$$
is well-defined for any rational function $f(z,w)\in\K(z,w)$. 
Define
\begin{align}
\notag
f(z,w)[a(z),b(w)]
\, :=\,
&f(z,w)\, a(z)b(w)
\, -\,\paraab
f(z,w)\, b(w)a(z)
\\
\notag
=\,
&T_{z,w}(f(z,w))\, a(z)b(w)
\, -\,\paraab
T_{w,z}(f(z,w))\, b(w)a(z).
\end{align}

More generally, 
define the product of a rational function $f(z_1, \dots,z_r)$
with an expression involving operator products and Lie brackets of fields
$a_1(z_1), \linebreak[0]\dots, a_r(z_r)$
by first replacing the Lie brackets by commutators, 
so that we obtain a sum of operator products
$A_{\si}:=\pm a_{\si 1}(z_{\si 1})\dots a_{\si r}(z_{\si r}), \si\in\bbS_r$,
and then multiplying the summand $A_{\si}$ by 
$T_{z_{\si 1},\dots,z_{\si r}}(f)$ for any $\si$.

We will refer to the above notational 
\index{radial ordering}
convention as {\bf radial ordering}.

\bigskip

{\bf Proposition.}\: {\it
For any fields $a(z), b(z)$ and $r\in\Z$, we have 
\begin{align}
\notag
a(w)_r b(w)
\; &=\;
\res_z (z-w)^r [a(z),b(w)]
\\
\notag
&=\;
\sum_{s\in\Z}
\sum_{i\geq 0}
(-1)^i \binom{r}{i}
\big( a_{r-i}b_{s+i}
-
\paraab (-1)^r\,  b_{s+r-i}a_i\big)\; w^{-s-1}.
\end{align}
Equivalently, $a(w)(x)b(w)=\res_z \de(z-w,x)[a(z),b(w)]$
where $a(w)\mapsto a(w)(x)$ is the state-field correspondence of 
$\End\subsv(E)$.
}

\bigskip

\begin{pf}
The first identity holds for $r\geq 0$ by definition.
It holds for $r=-1$ since
$$
\res_z (z-w)\inv a(z)b(w)
\; =\;
\res_z \sum_{i<0}z^i w^{-i-1}a(z)b(w)
\; =\;
a(w)_- b(w)
$$
and $-\res_z (z-w)\inv b(w)a(z)=\res_z\sum_{i\geq 0}z^i w^{-i-1} b(w)a(z)
=b(w)a(w)_+$.
It holds for $r<0$ because by the integration-by-parts formula we have
$\res_z (z-w)^r[\del_z a(z),b(w)]=-r\, \res_z (z-w)^{r-1} [a(z),b(w)]$.
The second identity follows from the computations at the end of section
\ref{SS:ratl fcts}. 
\end{pf}

\bigskip

Let $V$ be a vertex algebra and $M$ a bounded $\Z$-fold $V$-module.
The {\bf associativity formula} is 
\index{associativity formula}
$(a_r b)(z)=a(z)_r b(z)$.
This is an identity in $\Endv(M)$ for any $r\in\Z$. 
By the Proposition it is equivalent to 
$$
(a_r b)_s c
\; =\;
\sum_{i\geq 0}
(-1)^i \binom{r}{i}
\big( a_{r-i}b_{s+i}c
-
\paraab (-1)^r\,  b_{s+r-i}a_i c\big).
$$
This is the associativity formula for {\bf indices} $r, s$.
It is well-defined for $\Z$-fold modules if $r\geq 0$ or 
if $a, c$ and $b,c$ are weakly local.

Proposition \ref{SS:adj module} implies that 
$M$ is a $V$-module iff $Y_M(Ta)=\del_z Y_M(a)$ and $M$ satisfies 
the associativity formula.

\subsection{Jacobi Identity}
\label{SS:j id conf j id}

The Jacobi identity unifies the commutator formula and 
the associativity formula into one symmetric identity.

\bigskip

Let $V$ be a vertex algebra and $M$ a bounded $\Z$-fold $V$-module.
The associativity formula is 
$$
(a(x)b)(w)
\; =\;
\res_z \de(z-w,x)[a(z),b(w)].
$$
The commutator formula can be written 
$$
\res_x \de(z,w+x)\linebreak[0](a(x)b)(w)
\; =\; 
[a(z),b(w)].
$$
The product $\de(z,w+x)(a(x)b)(w)$ is well-defined because
the coefficient of $z^{-t-1}$ is $(w+x)^t (a(x)b)(w)$ and this product 
is well-defined for any $t\in\Z$ since $(a(x)b)(w)c\in M\lau{w}\lau{x}$
for any $c$. 

The {\bf Jacobi identity} is 
\index{Jacobi identity}
$$
\de(z,w+x)\, (a(x)b)(w)
\; =\;
\de(z-w,x)\, [a(z),b(w)].
$$
Since $\res_z \de(z,w)=1$, applying $\res_z$ to the Jacobi identity
we get the associativity formula.
Applying $\res_x$ we obtain the commutator formula.

The Jacobi identity is equivalent to 
$$
\sum_{i\geq 0}\binom{t}{i}(a_{r+i}b)_{s+t-i}c
=
\sum_{i\geq 0}(-1)^i \binom{r}{i}
\big(a_{t+r-i}b_{s+i}c-\paraab (-1)^r b_{s+r-i}a_{t+i}c\big)
$$
since
\begin{align}
\notag
\de(z,w+x)\, &(a(x)b)(w)c
\; =\;
\sum_{t\in\Z}\: 
(a(x)b)(w)c\: (w+x)^t\, z^{-t-1}
\\
\notag
&=
\sum_{r,s, t\in\Z}
\bigg(
\sum_{i\geq 0}
\binom{t}{i} (a_{r+i}b)_{s+t-i}c
\bigg)
z^{-t-1} w^{-s-1} x^{-r-1}
\end{align}
and similarly for the right-hand side.
We call the above identity the Jacobi identity for {\bf indices} $r, s, t$.
It is well-defined for $\Z$-fold modules over $\Z$-fold algebras
if $r, t\geq 0$ or if $a, b$ and $a, c$ and $b, c$ are weakly local.

The Jacobi identity for indices $r, s, 0$ is 
the associativity formula for indices $r, s$.
The Jacobi identity for indices $0, s, t$ is 
the commutator formula for indices $t, s$.

The Jacobi identity for indices $r, s, t$ is easier to memorize by noting 
that the indices $r, s, t$ correspond to the pairs $ab, bc$, and $ac$.

\subsection{Conformal Jacobi Identity and Wick Formula I}
\label{SS:conf j wick0}

We show that the conformal Jacobi identity and the Wick formula are 
equivalent to certain special cases of the Jacobi identity.

\bigskip

{\bf Proposition.}\: {\it
Let $V$ be a $\Z$-fold algebra.
The conformal Jacobi identity,
the Jacobi identity, the associativity formula,
and the commutator formula for indices $r, s, t\geq 0$ are all equivalent.
}

\bigskip

\begin{pf}
The Jacobi identity for indices $\geq 0$ is
$
[[a_{\la+\ka}b]_{\mu+\ka}c]=[a_{\la+\ka}[b_{\mu-\la}c]]-
[b_{\mu-\la}[a_{\la+\ka}c]]
$
because from $z\de(z,w+x)=(w+x)\de(z,w+x)$ we get
$$
\res_{z,w,x}e^{z\ka}e^{w\mu}e^{x\la}\de(z,w+x)(a(x)b)(w)
=
\res_{w,x}e^{w(\mu+\ka)}e^{x(\la+\ka)}(a(x)b)(w)
$$
and similarly for the right-hand side.

Setting $t=\ka=0$ we obtain the associativity formula for indices $\geq 0$
and the conformal Jacobi identity. 
Conversely, replacing $\la$ by $\la+\ka$ and $\mu$ by $\mu+\ka$ we get back
the Jacobi identity for indices $\geq 0$.

Setting $r=\la=0$ we obtain the commutator formula for indices $\geq 0$.
Replacing $\mu$ by $\mu-\ka$ we obtain the conformal Jacobi identity. 
\end{pf}

\bigskip

We have already shown in section \ref{SS:conf jac} that 
the conformal Jacobi identity is equivalent to 
the commutator formula for indices $\geq 0$. 

The fact that the associativity formula for indices $\geq 0$ 
is equivalent to the conformal Jacobi identity also follows from
the following two remarks.
A conformal algebra $R$ satisfies the conformal Jacobi identity
iff the map $R\to\fgl(R)\pau{\la}, a\mapsto [a\subla\;]$, 
is an unbounded conformal algebra morphism, see section \ref{SS:conf op}.
The projection $\fg\pau{z\uppm}\to\fg\pau{\la}, a(z)\mapsto a(z)_+$, 
is an unbounded conformal algebra morphism, 
see section \ref{SS:unb confa power series}. 

Let $V$ be a vertex algebra.
The commutator formula for indices $t, s$ is 
$[a_t,b_s]=\sum\binom{t}{i}(a_i b)_{t+s-i}$, see section \ref{SS:module vlie}.
The Wick formula is $a_t(bc)=(a_t b)c+\paraab\, b(a_t c)+
\sum_{i=0}^{t-1}\binom{t}{i}(a_i b)_{t-1-i}c$ 
where $t\geq 0$, see section \ref{SS:wick}.
Thus the Wick formula is the commutator formula for indices $t\geq 0$ and 
$s=-1$.

Together with the Proposition we see that 
the Jacobi identity for indices $r\geq 0, s\in\Z, t\geq 0$ implies
the conformal Jacobi identity and the Wick formula.
In Proposition \ref{SS:conf j wick} we prove that also the converse is true.

\subsection{Skew-Symmetry}
\label{SS:skew sym}

We prove that skew-symmetry is equivalent to conformal skew-symmetry and 
$[\, ,]\subast=[\, ,]\subslie$.

\bigskip

{\bf Skew-symmetry} for a vertex algebra is
\index{skew-symmetry}
$$
a(z)b
\; =\;
\paraab\, e^{zT}b(-z)a.
$$
This is equivalent to
$$
a_r b
\; =\;
\paraab\sum_{i\geq 0}\: (-1)^{r+1+i}\, T^{(i)}(b_{r+i}a).
$$
We call this skew-symmetry for {\bf index} $r$.
It is well-defined for $\Z$-fold algebras with an even operator $T$
if $b, a$ are weakly local.

Skew-symmetry for indices $r\geq 0$ is conformal skew-symmetry because
$$
\res_z e^{z\la}e^{zT}b(-z)a
\; =\;
e^{T\del\subla}\res_z e^{z\la}b(-z)a
\; =\;
-[b_{-\la-T}a].
$$
Skew-symmetry for index $-1$ is equivalent to $[b,a]\subast=[b,a]\subslie$ 
because
$$
\res_z z\inv e^{zT}b(-z)a
\, =\,
\res_z\big( z\inv+\int_0^T d\la\, e^{z\la}\big) b(-z)a
\, =\,
ba-\int_0^T d\la\,[b_{-\la}a].
$$

\bigskip

{\bf Proposition.}\: {\it
Let $V$ be a vertex algebra.
Then skew-symmetry is equivalent to conformal skew-symmetry and 
$[\, ,]\subast=[\, ,]\subslie$.
}

\bigskip

\begin{pf}
We have seen that conformal skew-symmetry and $[\, ,]\subast=[\, ,]\subslie$
are equivalent to skew-symmetry for indices $r\geq -1$.
Hence it suffices to show that skew-symmetry for index $r\ne 0$ 
implies skew-symmetry for index $r-1$.
This follows from 
\begin{align}
\notag
\res_z (\del_z z^r)a(z)b
\; &=\;
-\res_z z^r (Ta)(z)b
\; =\;
-\res_z z^r e^{zT}b(-z)Ta
\\
\notag
\; &=\;
-\res_z z^r \del_z (e^{zT}b(-z)a)
\; =\;
\res_z (\del_z z^r) e^{zT}b(-z)a
\end{align}
using that
$\del_z (e^{zT}b(-z)a)=e^{zT}T(b(-z)a)-e^{zT}(Tb)(-z)a=e^{zT}b(-z)Ta$.
\end{pf}

\section{Associative Vertex Algebras of Distributions}
\label{S:vas of distr}

In sections \ref{SS:fz complete} and \ref{SS:a z uppm h}
we construct the FZ-completion of a graded algebra and 
show that if $\cA$ is FZ-complete then the space $\cA\pau{z\uppm}_H$ of 
homogeneous distributions is an unbounded vertex algebra.

In sections \ref{SS:norm ord prod pre lie} and \ref{SS:wick fields}
we prove that $\cA\pau{z\uppm}_H$ and the unbounded endomorphism vertex algebra
$\Endv(E)$ are pre-Lie algebras and satisfy the Wick formula.
In sections \ref{SS:ope sing}--\ref{SS:skew sym fields}
we use the fact that the $t$-th products of $\Endv(E)$ and $\cA\pau{z\uppm}_H$
are the coefficients of the OPE
to prove that distributions are local iff they satisfy skew-symmetry.
In sections \ref{SS:dong fields} and \ref{SS:vas of fields}
we prove that local distributions generate associative vertex subalgebras of 
$\Endv(E)$ and $\cA\pau{z\uppm}_H$.

\medskip

{\bf Conventions.}\:
All associative algebras and all unbounded vertex algebras are assumed to be 
unital.
In sections \ref{SS:norm ord prod pre lie}--\ref{SS:vas of fields}
we denote by $E$ a vector space and by $\cA$ an FZ-complete algebra.
Let $\rho\in\Q\uptimes$ such that $\rho\inv\in\N$.

\subsection{FZ-Complete Algebras}
\label{SS:fz complete}

We prove that any $\Z$-graded associative algebra $\cA$ has an 
{\bf FZ-completion} $\hcA$.
The FZ-completion is used in section \ref{SS:fz alg}
to construct the Frenkel-Zhu algebra.

\bigskip

A {\bf topological} vector space is 
\index{topological!vector space}
a vector space with a topology
such that addition and scalar multiplication are continuous. 
Here we assume that $\K$ is endowed with the discrete topology.
Thus scalar multiplication $\K\times V\to V$ is continuous iff
scalar multiplication $\ka\cdot: V\to V$ is continuous for any $\ka\in\K$.

A {\bf topological} $\K$-graded algebra is 
\index{topological!Zgraded algebra@$\Z$-graded algebra}
a $\K$-graded algebra $\cA=\bigoplus\cA_h$ together with
a topology on each $\cA_h$ such that $\cA_h$ are topological vector spaces 
and the multiplications $\cA_h\times\cA_{h'}\to\cA_{h+h'}$ are continuous.
 
Let $V$ be a vector space and $V_i\subset V, i\in I$, a family of subspaces
where $I$ is a partially ordered set such that $V_i\subset V_j$ for $i\leq j$
and $I$ is directed: 
for any $i, j$ we have $k\leq i$ and $k\leq j$ for some $k$.
Then there exists a unique topology on $V$ such that $(V_i)$ is a 
fundamental system of neighbourhoods of $0$.
This is 
\index{linear topology}
the {\bf linear topology} defined by $(V_i)$.
The subspaces $V_i$ are both open and closed with respect to the linear 
topology.
The linear topology is Hausdorff iff $\bigcap_i V_i=0$.

An {\bf FZ-topological} algebra is
\index{FZ-topological algebra}
a topological $\rho\Z$-graded associative algebra $\cA$ such that 
the topology of $\cA_h$ is linear and coarser than the linear topology 
defined by the closures of the subspaces $\sum_{l<k}\cA_{h-l}\cA_l$ for 
$k\in\Z$.

\bigskip

{\bf Lemma.}\: {\it
The forgetful functor $\cA\mapsto\cA$ from FZ-topological algebras to
$\rho\Z$-graded associative algebras has a left adjoint $\cA\mapsto\cA$.
}

\bigskip

\begin{pf}
Let $\cA$ be a $\rho\Z$-graded associative algebra.
Endow $\cA_h$ with the linear topology defined by 
$U_{h,k}:=\sum_{l<k}\cA_{h-l}\cA_l$ for $k\in\Z$.
The subspaces $U_{h,k}$ satisfy $\cA_h U_{h',k}\subset U_{h+h',k}$ and
$U_{h,k-h'}\cA_{h'}\subset U_{h+h',k}$.
It follows that $\cA$ is a topological $\Z$-graded algebra since
$(a+U_{h,k-h'})(b+U_{h',k})\subset ab+U_{h+h',k}$
for $a\in\cA_h$ and $b\in\cA_{h'}$.
Thus $\cA$ is an FZ-topological algebra.

Let $\cB$ be an FZ-topological algebra and $\phi: \cA\to\cB$ a morphism of
$\rho\Z$-graded algebras.
Since $U_{h,k}$ is mapped to $\sum_{l<k}\cB_{h-l}\cB_l$,
the morphism $\phi$ is degreewise continuous and hence a morphism of
FZ-topological algebras.
Thus the claim follows.
\end{pf}

\bigskip

Let $V$ be a topological vector space with a linear topology defined by 
$(V_i)$. 
Then the completion $\hV$ of $V$ is the inverse limit of $V/V_i$ for $i\in I$.
In particular, the topology of $\hV$ is the linear topology defined by
$(\bV_i)$ where $\bV_i$ is the closure of the image of $V_i$ in $\hV$.

A topological $\rho\Z$-graded associative algebra is {\bf complete}
\index{complete algebra}
if the topological vector spaces $\cA_h$ are complete and Hausdorff.
We call a complete FZ-topological algebra 
\index{FZ-complete!algebra}
just a {\bf FZ-complete} algebra.

\bigskip

{\bf Proposition.}\: {\it
\iti\:
The forgetful functor $\cA\mapsto\cA$ from FZ-complete algebras to
FZ-topological algebras has a left adjoint $\cA\mapsto\hcA$.

\smallskip

\itii\:
The forgetful functor $\cA\mapsto\cA$ from FZ-complete algebras to
$\rho\Z$-graded associative algebras has a left adjoint $\cA\mapsto\hcA$.
}

\bigskip

\begin{pf}
\iti\:
Let $\cA$ be an FZ-topological algebra and $\hcA_h$ 
the completion of the topological vector space $\cA_h$.
Then $\hcA:=\bigoplus\hcA_h$ is 
a topological $\rho\Z$-graded associative algebra.

Any continuous map $f$ satisfies $f(\overline{M})\subset\overline{f(M)}$ 
and hence $\overline{f(M)}=\overline{f(\overline{M})}$.

Define $U_{h,k}:=\sum_{l<k}\cA_{h-l}\cA_l$ and 
denote by $\ovlM$ the closure of a subset $M\subset\hcA_h$.
Taking for $f$ the canonical morphism $\cA\to\hcA$, 
we see that the topology of $\hcA_h$ is the linear topology 
defined by $\ovlU_{h,k}$ for $k\in\Z$, 
since $\cA$ is an FZ-topological algebra.

Let $U_{h,k}'$ be the closure of $\sum_{l<k}\hcA_{h-l}\hcA_l$.
In order to show that $\hcA$ is FZ-complete it suffices to show that
$U_{h,k}'=\ovlU_{h,k}$.
We obviously have $\ovlU_{h,k}\subset U_{h,k}'$ since $\cA_l\subset\hcA_l$.
Conversely, we have $\ovlM+\ovlN\subset\overline{M+N}$ and 
$\ovlM\cdot\ovlN\subset\overline{MN}$.
Moreover, $\sum_{l<k}M_l=\bigcup_{i<k}\sum_{l=i}^{k-1} M_l$ and
$\bigcup_i \ovlM_i\subset\overline{\bigcup_i M_i}$.
From $\hcA_h=\overline{\cA}_h$ we thus get $U_{h,k}'\subset\ovlU_{h,k}$.

The universal property of $\hcA_h$ implies that $\cA\mapsto\hcA$ is left
adjoint.

\smallskip

\itii\:
This follows from \iti\ and the Lemma since the composition of left adjoint
functors is left adjoint.
\end{pf}

\subsection{The Unbounded Vertex Algebra $\cA\pau{z\uppm}_H$}
\label{SS:a z uppm h}

We show that if $\cA$ is an FZ-complete algebra then
$\cA\pau{z\uppm}_H$ is an unbounded vertex algebra.

\bigskip

Let $\cA$ be an FZ-complete algebra.
Define
$$
\cA\pau{z_1\uppm, \dots, z_r\uppm}_H
\; :=\;
\bigoplus_{h\in\rho\Z}
\cA\pau{z_1\uppm, \dots, z_r\uppm}_h
$$ 
where $\cA\pau{z_1\uppm, \dots, z_r\uppm}_h$ consists of distributions
$\sum c_t z^{-t_1-1}\dots z_r^{-t_r-1}$ such that 
$c_t\in\cA_{h-t_1\ldots-t_r-r}$ for $t\in\Z^r$.
Considering $\cA$ as a graded Lie algebra,
$\cA\pau{z\uppm}_H$ is a graded unbounded conformal algebra by 
section \ref{SS:alg of distr}.

A sum $\sum_{i\in I} d_i$ in $\cA_h$ is {\bf summable} if 
\index{summable family}
there exists $d\in\cA_h$ such that for any neighbourhood $U$ of $0$ 
there exists a finite subset $J\subset I$ such that 
$\sum_{i\in J'}d_i\in d+U$ for any finite subset $J'\subset I$ containing $J$.
This $d$ is unique and we define $\sum_{i\in I} d_i:=d$.

The Cauchy criterion states that a sum $\sum_i d_i$ is summable iff 
for any $U$ there exists a finite subset $J$ such that 
$\sum_{i\in J'}d_i\in U$ for any finite subset $J'\subset I\setminus J$.
Since the topology is linear, we may assume that $U$ is a subspace.
Then this is equivalent to $d_i\in U$ for any $i\in I\setminus J$.
In particular, $\sum_i d_i$ is summable if for any $k\in\Z$ 
there exists a finite subset $J$ such that $d_i\in\sum_{l<k}\cA_{h-l}\cA_l$ 
for any $i\in I\setminus J$.

We remark that a sum of the form $\sum_{i\geq n}d_i$ is summable 
iff the sequence of partial sums $\sum_{i=n}^m d_i$ converges.

We call $c(z_1, \dots, z_r)\in\cA\pau{z_1\uppm, \dots, z_r\uppm}_h$ 
\index{regular at $z_1=\ldots=z_r$}
{\bf regular} at $z_1=\ldots=z_r$
if the sum $\sum_{t_1+\ldots+t_r=n} c_t$ in $\cA_{h-n-r}$ is summable for 
any $n\in\Z$.
In this case $c(z, \dots, z)$ is well-defined and contained in 
$\cA\pau{z\uppm}_h$.
We extend the notion of regularity linearly
to $\cA\pau{z_1\uppm, \dots, z_r\uppm}_H$.

For $a(z)\in\cA\pau{z\uppm}_h$ and $b(z)\in\cA\pau{z\uppm}_{h'}$, 
the {\bf normal ordered product} 
\index{normal ordered product}
\begin{align}
\notag
\normord{a(z)b(w)}
\; :=&\;
a(z)_-b(w)\, +\,\paraab\, b(z)a(w)_+
\\
\notag
\; =&\;
\sum_{m\in\Z}\:\Big(\sum_{n<0}a_n b_m+\paraab\sum_{n\geq 0} b_m a_n\Big)\, 
z^{-n-1}w^{-m-1}
\end{align}
is contained in $\cA\pau{z\uppm,w\uppm}_{h+h'}$ and regular at $z=w$.
It follows that $\cA\pau{z\uppm}_H$ with the {\bf normal ordered product} 
$\normord{a(z)b(z)}$ is a graded unital unbounded vertex algebra.

We consider now an analogue for $\cA\pau{z\uppm}_H$ of 
the space $\cF\lau{z}\lau{w}(E)$, 
see sections \ref{SS:fields} and \ref{SS:z fold fields}.

Let $\cA\pau{z_1\uppm, \dots, z_r\uppm}_{>,h}$ be the subspace of
$\cA\pau{z_1\uppm, \dots, z_r\uppm}_h$ consisting of $c(z_1, \dots, z_r)$
such that for any neighbourhood $U$ of $0\in\cA_{h'}$ and any 
$p\in\set{2, \dots, r}$
there exists $k$ such that if $h-\sum t_i-r=h'$ and 
$\sum_{i=p}^r t_i\geq k$ then $c_t\in U$.
Let $\cA\pau{z_1\uppm, \dots, z_r\uppm}_>$ be the direct sum of these spaces.

If $a_i(z)\in\cA\pau{z\uppm}_H$ then 
$a_1(z_1)\dots a_r(z_r)\in\cA\pau{z_1\uppm, \dots, z_r\uppm}_>$.

The space $\cA\pau{z_1\uppm, \dots, z_r\uppm}_>$ is a module over
$\K[z_1\uppm, \dots, z_r\uppm, (z_i-z_j)\inv\mid i<j]$.
The element $f:=(z_i-z_j)\inv$ acts by multiplication with 
$T_{z_i,z_j}(f)=\sum_{n\geq 0} z_i^{-1-n} z_j^n$, 
see section \ref{SS:ratl fcts}.
In fact, if $d:=fc$ then $d_t$ is a summable sum of
elements $c_{t'}$ with $\sum t'_l=\sum t_l$ and 
$\sum_{l=p}^r t'_l\geq\sum_{l=p}^r t_l$.
Hence if $c_{t'}\in U$ then $d_t\in U$, since we may assume that
$U$ is a vector subspace and closed.

Distributions in the intersection of the spaces 
$\cA\pau{z_{\si 1}\uppm, \dots, z_{\si r}\uppm}_>$, $\si\in\bbS_r$,
are regular at $z_1=\ldots=z_r$.
This follows from the fact that the set of $t$ such that $\sum t_i=n$
and $t_i\leq k$ is finite.
The number $p$ above is taken to be $r$.

As in Proposition \ref{SS:z fold fields} one has the residue formula
$$
a(w)_i b(w)
\; =\;
\res_z (z-w)^i [a(z),b(w)].
$$

\subsection{Pre-Lie Identity for Distributions}
\label{SS:norm ord prod pre lie}

We prove that $\Endv(E)$ and $\cA\pau{z\uppm}_H$ are pre-Lie algebras.
More generally, we prove that the normal ordered product of 
an associative $3/4$-algebra yields a pre-Lie algebra.

\bigskip

An {\bf associative $3/4$-algebra} is a pair $V=(V_+,V_-)$ of 
associative algebras together with a $V_-$-$V_+$-bimodule structure on 
$V:=V_+\oplus V_-$ such that  
\index{associative!$3/4$-algebra}
the inclusion $V_{\pm}\subset V$ is a morphism of 
right/left $V_{\pm}$-modules and $a(b)=(a)b$ for $a\in V_-, b\in V_+$ 
where $a(b), (a)b$ denote the action of $a$ on $b$, of $b$ on $a$.

In other words, an associative $3/4$-algebra is a pair $(V_+,V_-)$
of vector spaces together with even linear maps 
$V\subpm\otimes V\subpm\to V\subpm$ and $V_-\otimes V_+\to V$ 
such that associativity $(ab)c=a(bc)$ holds for 
any of the four out of eight combinations of homogeneous elements 
$a, b, c$ for which both sides of the identity are well-defined.

The elements 
\index{annihilators}
of $V_+, V_-$ are 
\index{creators}
called {\bf annihilators} and {\bf creators}.
Annihilators stand always on the right of creators.

For example, if $V$ is an associative algebra and $V\subpm\subset V$
are subalgebras such that $V=V_+\oplus V_-$ as vector spaces then
$V$ is an associative $3/4$-algebra.

Our main examples of associative $3/4$-algebras are $\Endv(E)$ and 
$\cA\pau{z\uppm}_H$.
The decomposition $V=V_+\oplus V_-$ is induced from the decomposition
$E\pau{z\uppm}=z\inv E\pau{z\inv}\oplus E\pau{z}$.
In particular, $\Endv(E)_+=z\inv\cF[z\inv](E)$ as in section \ref{SS:fields}
and $\Endv(E)_-=\End(E)\pau{z}$.
The multiplications are $a(z)\otimes b(z)\mapsto a(z)b(z)$.
They are well-defined by sections \ref{SS:fields} and \ref{SS:a z uppm h}.

For an associative $3/4$-algebra $V$, let $V^{::}$ denote 
the vector space $V$ endowed with the {\bf normal ordered product}
\index{normal ordered product}
$$
a\otimes b
\;\;\mapsto\;\;\;
\normord{ab}
\;\; :=\;\;
a_-b\, +\,\paraab\,  ba_+
$$
where $a\mapsto a_{\pm}$ are the projections $V\to V_{\pm}$.
We have
\begin{alignat}{2}
\notag
\normord{a_- b_-}
\; &=\;
a_- b_-,
\qquad &
\normord{a_- b_+}
\; &=\;
a_- b_+,
\\
\notag
\normord{a_+ b_-}
\; &=\;
\paraab\, b_- a_+,
\qquad &
\normord{a_+ b_+}
\; &=\;
\paraab\, b_+ a_+.
\end{alignat}
In particular, $V_-$ and the opposite algebra of $V_+$ are 
subalgebras of $V^{::}$.

Of course, in the case of $\Endv(E)$ and $\cA\pau{z\uppm}_H$ 
the above definition of the normal ordered product 
agrees with the one from sections \ref{SS:fields} and \ref{SS:a z uppm h}.

\bigskip

{\bf Proposition.}\: {\it
The algebra $V^{::}$ is a pre-Lie algebra.
}

\bigskip

\begin{pf}
We have $\normord{a\normord{bc}}=a_-b_-c+ a_- cb_+ + b_- ca_+ + cb_+a_+$ and
\begin{align}
\notag
&\normord{\normord{ab}c}
\; =\;
\normord{ab}_- c\: +\:  c\normord{ab}_+
\\
\notag
=\;
&a_-b_-c \, +\,
(a_-b_+)_- c \, +\,
(b_-a_+)_-c \, +\,
c(a_-b_+)_+ \, +\,
c(b_-a_+)_+ \, +\,
c\,b_+a_+
\end{align}
because $(b_+a_+)_-=(a_-b_-)_+=0$. 
Thus the associator $\normord{\normord{ab}c}-\normord{a\normord{bc}}$
is symmetric in $a, b$.
Hence $V^{::}$ is a pre-Lie algebra by Proposition \ref{SS:prelie}
\end{pf}

\bigskip

If $V$ is an associative $3/4$-algebra such that $V_+, V_-$ are commutative
then the normal ordered product is commutative and hence 
$V^{::}$ is commutative and associative by Remark \ref{SS:prelie}\,\itiv.

One can show that if one defines instead $\normord{a_+ b_+}=a_+ b_+$
then $V^{::}$ is in general not a pre-Lie algebra.

\subsection{Wick Formula for Distributions}
\label{SS:wick fields}

We prove that $\Endv(E)$ and $\cA\pau{z\uppm}_H$ satisfy 
the conformal Jacobi identity and the Wick formula.
By Proposition \ref{SS:conf jac}
we already know that the conformal Jacobi identity holds.
But as the proof shows, it is more natural to consider both identities
together.

\bigskip

For any distributions $a(z), b(z), c(z)$ in $\Endv(E)$ or $\cA\pau{z\uppm}_H$, 
\index{radially ordered Leibniz identity}
the {\bf radially ordered Leibniz identity} holds:
\begin{align}
\notag
&f(z,w,y)\, [[a(z),b(w)],c(y)] 
\\
\notag
=\;
&f(z,w,y)\,
\big( [a(z),[b(w),c(y)]]
\; -\;
\paraab[b(w),[a(z),c(y)]]\big)
\end{align}
where $f(z,w,y)\in\K[(z-w)\uppm,(w-y)\uppm,(z-y)\uppm]$.
This is proven in the same way as the
Leibniz identity for the commutator of an associative algebra.
The usual cancellation of terms still takes place 
because the power series expansions of $f(z,w,y)$
only depend on the order in which $a(z), b(w), c(y)$
appear in an operator product.

\bigskip

{\bf Proposition.}\: {\it
The unbounded vertex algebras $\Endv(E)$ and $\cA\pau{z\uppm}_H$ satisfy 
the conformal Jacobi identity and the Wick formula.
}

\bigskip

\begin{pf}
By section \ref{SS:conf j wick0} it suffices to prove
the Jacobi identity for indices $r\geq 0, s\in\Z, t\geq 0$.
Define $A:=\res_{z,w}(z-w)^r (w-y)^s (z-y)^t\, [[a(z),b(w)],c(y)]$.
Expanding the $t$-th power of $z-y=(z-w)+(w-y)$ we obtain
\begin{align}
\notag
A\; =\;
&\sum_{i\geq 0}\: 
\binom{t}{i}\: 
\res_w
(w-y)^{s+t-i}\,
\big[\res_z (z-w)^{r+i}[a(z),b(w)],\, c(y)\big] 
\\
\notag
=\;
&\sum_{i\geq 0}\: \binom{t}{i}\:
(a(y)_{r+i}b(y))_{s+t-i}c(y). 
\end{align}

Applying the radially ordered Leibniz identity and expanding 
the $r$-th power of $z-w=(z-y)-(w-y)$ we see that $A$ is also equal to
\begin{align}
\notag
&\res_{z,w}
(z-w)^r (w-y)^s (z-y)^t
\big(
[a(z),[b(w),c(y)]]
-
[b(w),[a(z),c(y)]]
\big)
\\
\notag
= 
&
\sum_{i\geq 0}\: 
(-1)^i\binom{r}{i}\:
\res_{z,w}
\Big(
(w-y)^{s+i}(z-y)^{t+r-i}\, [a(z),[b(w),c(y)]]
\\
\notag
&\qquad\qquad\qquad\qquad\;\;
-(-1)^r
(w-y)^{s+r-i}(z-y)^{t+i}\, [b(w),[a(z),c(y)]]
\Big)
\\
\notag
=
&\sum_{i\geq 0}\, 
(-1)^i 
\binom{r}{i}
\big(
a(y)_{t+r-i}b(y)_{s+i}c(y)-(-1)^r\,  
b(y)_{s+r-i}a(y)_{t+i}c(y)\big).
\end{align}
\end{pf}

\bigskip

The conformal Jacobi identity and the Wick formula are in fact
equivalent to the Jacobi identity for indices $r\geq 0, s\in\Z, t\geq 0$,
see Proposition \ref{SS:conf j wick}.

\subsection{Operator Product Expansion}
\label{SS:ope sing}

We prove that the $r$-th products of the $\Z$-fold algebras $\Endv(E)$ and 
$\cA\pau{z\uppm}_H$ are the coefficients of the OPE.

\bigskip

{\bf Lemma.}\: {\it
Let $\fg$ be an algebra and $a(z), b(z)\in\fg\pau{z\uppm}$. 
Then $c^i(z)=a(z)_i b(z)$ are the unique distributions such that
$$
[a_t,b_s]
\; =\;
\sum_{i\geq 0}\:
\binom{t}{i}\,
c^i_{t+s-i}
$$
for any $t\geq 0$ and $s\in\Z$.
Moreover, this identity is equivalent to 
$$
[a(z)_+,b(w)]
\; =\;
\sum_{i\geq 0}\:
\frac{c^i(w)}{(z-w)^{i+1}}.
$$
}

\bigskip

\begin{pf}
The inverse of the triangular matrix
$\big(\binom{t}{i}  w ^{t-i}\big)_{t,i\geq 0}$ is 
$\big(\binom{i}{s}(-w)^{i-s}\big)_{i,s}$ because
$1=((1-x)+x)^t=\sum_{i,s}\binom{t}{i}\binom{i}{s}(-1)^{i-s}x^{t-i+i-s}$.
Thus the identities $[a_t,b(z)]=\sum_i\binom{t}{i}w^{t-i}c^i(z)$ 
are equivalent to
$$
c^i(w)
\; =\; 
\sum_{s\geq 0}\:\binom{i}{s}(-w)^{i-s}\, [a_s,b(w)]
\; =\; 
\res_z (z-w)^i [a(z),b(w)].
$$
Applying $\del_w^{(i)}$ to $(z-w)^{-1}=\sum_{t\geq 0} w^t z^{-t-1}$ we obtain
$$
(z-w)^{-i-1}
\; =\;
\sum_{t\geq 0}\:
\binom{t}{i}\, w^{t-i} z^{-t-1}
$$
for $i\geq 0$.
This yields the second claim.
\end{pf}

\bigskip

{\bf Corollary.}\: {\it
Let $V$ be a $\Z$-fold algebra.
Then the commutator formula for indices $t\geq 0, s\in\Z$ 
is equivalent to $(a_r b)(z)=a(z)_r b(z)$ for any $r\geq 0$.
\hfill $\square$
}

\bigskip

The Corollary is part of Proposition \ref{SS:conf j wick}.

The following result is 
\index{operator product expansion}
the {\bf operator product expansion} or 
\index{OPE}
{\bf OPE}.

\bigskip

{\bf Proposition.}\: {\it
Let $a(z), b(z)$ be in $\Endv(E)$ or $\cA\pau{z\uppm}_H$.
Then 
$$
a(z)b(w)
\; =\;
\sum_{i\geq 0}\:
\frac{a(w)_i b(w)}{(z-w)^{i+1}}
\; +\;
\normord{a(z)b(w)},
$$
$\normord{a(z)b(w)}$ is regular at $z=w$, and 
$a(w)_i b(w)$ for $i<0$ is equal to the Taylor coefficient
$\del_z^{(-1-i)}\normord{a(z)b(w)}|_{z=w}$.
}

\bigskip

\begin{pf}
The first identity follows from the Lemma because
$$
a(z)b(w)
\; =\;
[a(z)_+,b(w)]
\; +\; 
a(z)_- b(w)
 + 
b(w)a(z)_+.
$$
The second claim was already shown 
in sections \ref{SS:fields} and \ref{SS:a z uppm h}.
Since $\del_z(a(z)\subpm)=(\del_z a(z))\subpm$ we get
$$
a(w)_{-1-i}b(w)
\; =\;\;
\normord{\del_w^{(i)}a(w)b(w)}
\; =\;
\del_z^{(i)}\normord{a(z)b(w)}|_{z=w}.
$$
\end{pf}

\subsection{Taylor Coefficients of the Operator Product}
\label{SS:tayl}

We prove that if $a(z), b(z)$ are weakly local then their $r$-th products
are equal to the Taylor coefficients of the operator product $a(z)b(w)$.

\bigskip

{\bf Lemma.} (Taylor's Formula)\: {\it
Let $c(z,w)$ be in $\End(E)\pau{z\uppm,w\uppm}$ or in 
$\cA\pau{z\uppm,w\uppm}_H$.
If $c(z,w)$ is regular at $z=w$ then $c^i(w)=\del_z^{(i)}c(z,w)|_{z=w}$ are 
\index{Taylor's formula}
the unique distributions such that for any $n\geq 0$ we have
$$
c(z,w)
\; =\;
\sum_{i=0}^{n-1}\:
c^i(w)\, (z-w)^i
\; +\;
(z-w)^n\, r(z,w)
$$
for some $r(z,w)$ that is regular at $z=w$.
If $c^i(z)=0$ for any $i\geq 0$ then $c(z,w)=0$.
}

\bigskip

\begin{pf}
We recall that regularity was defined in sections \ref{SS:fields} and 
\ref{SS:a z uppm h}.

The case $c(z,w)\in\End(E)\pau{z\uppm,w\uppm}$ follows directly from
Taylor's formula for $E\lau{z,w}$, see Lemma \ref{SS:loc diff op}.
Thus it suffices to consider the case
$c(z,w)\in\cA\pau{z\uppm,w\uppm}_{h'}$ for some $h'\in\Z$.

We do induction on $n$.
For $n=1$ we have to show that if $d(z,w)$ vanishes for $z=w$ then 
$d(z,w)=(z-w)r(z,w)$ for some $r(z,w)$ that is regular at $z=w$. 

Let $U_{h,k}$ be the closure of $\sum_{l\leq k}\cA_{h-l}\cA_l$.
Since $\cA_h$ is complete and Hausdorff and 
the $U_{h,k}$ form a neighbourhood basis of $0\in\cA_h$, 
the space $\cA_h$ is the inverse limit of $\cA_h/U_{h,k}$.
Let $V_{h,k}$ be a vector space complement of $U_{h,k}$ in $U_{h,k+1}$ 
for $k\leq -1$ and $U_{h,0}\oplus V_{h,0}=\cA_h$.
Then $\cA_h=\prod V_{h,k}$.

Let $\pi_k: \cA_h\to V_{h,k}$ be the projections.
Then $e(z,w)\in\cA\pau{z\uppm,w\uppm}_{h'}$ is regular at $z=w$ iff 
$\pi_k\sum e_{i,n-i} x^i\in V_{h'-n-2,k}[x\uppm]$ for any $k$ and $n\in\Z$.
Therefore the claim for $n=1$ follows from the fact that for 
$p(x)\in V_{h,k}[x\uppm]$ we have $p(1)=0$ iff 
$p(x)=(x-1)q(x)$ for some $q(x)\in V_{h,k}[x\uppm]$.
This fact is proven by reducing it to $\K[x\uppm]$ using a basis of $V_{h,k}$.

Using the above arguments, 
the lemma is now proven in the same way as Taylor's formula for $E\lau{z,w}$,
see Lemma \ref{SS:loc diff op}.
\end{pf}

\bigskip

{\bf Proposition.}\: {\it
Let $a(z), b(z)$ be in $\Endv(E)$ or $\cA\pau{z\uppm}_H$
such that $a(z), b(z)$ are weakly local of order $\leq n$ for some $n\in\Z$.
Then $(z-w)^n a(z)b(w)$ is regular at $z=w$ and for any $i\in\Z$ we have 
$$
a(w)_i b(w)
\; =\;
\del_z^{(n-1-i)}((z-w)^n a(z)b(w))|_{z=w}.
$$
}

\bigskip

\begin{pf}
By assumption, $a(z)_i b(z)=0$ for $i\geq n$. 
By Proposition \ref{SS:ope sing} we get
$$
(z-w)^n a(z)b(w)
=
\sum_{i=0}^{n-1}
a(w)_i b(w)\, (z-w)^{n-1-i}
+
(z-w)^n\normord{a(z)b(w)}.
$$
Thus the claim follows from Taylor's formula, Proposition \ref{SS:ope sing}, 
and from $\del_z^{(n+i)}((z-w)^n\normord{a(z)b(w)})|_{z=w}=
\del_z^{(i)}\normord{a(z)b(w)}|_{z=w}$.
\end{pf}

\subsection{Skew-Symmetry and Locality for Distributions}
\label{SS:skew sym fields}

We prove that a pair of distributions is local iff it satisfies skew-symmetry.
We use this result to show that any $V$-module is a $V\subslie$-module.

\bigskip

{\bf Proposition.}\: {\it
Let $a(z), b(z)$ be in $\Endv(E)$ or $\cA\pau{z\uppm}_H$.
Then $a(z), b(z)$ are local iff they are weakly local and
satisfy skew-symmetry.
}

\bigskip

\begin{pf}
Locality implies weak locality; skew-symmetry implies that 
$b(z), a(z)$ are also weakly local;
thus we may assume that $a(z), b(z)$ and $b(z), a(z)$ are weakly local
of order $\leq n$.
Then $c(z,w):=(z-w)^n a(z)b(w)$ and $c'(z,w):=(z-w)^n b(w)a(z)$ are
regular at $z=w$ by Proposition \ref{SS:tayl}.

In general, we have the chain rule $\del_w c(w,w)=(\del_z+\del_w)c(z,w)|_{z=w}$
since $\del_w(w^{n+m})=(\del_z+\del_w)(z^n w^m)|_{z=w}$.
Thus we get
\begin{align}
\notag
\del_z^{(r)}c(z,w)|_{z=w}
\;=\;
&
\sum_{i\geq 0}\:
(\del_z+\del_w)^{(i)}(-\del_w)^{(r-i)}c(z,w)|_{z=w}
\\
\notag
=\;
&
\sum_{i\geq 0}\:
(-1)^{r+i}\: \del_w^{(i)}(\del_w^{(r-i)}c(z,w)|_{z=w}).
\end{align}
Replace $r$ by $-r-1+n$ and apply Proposition \ref{SS:tayl} to 
the left-hand side.
Since $e^{zT}d(-z)=e^{zT}d'(-z)$ iff $d=d'$, we see that 
$a(z), b(z)$ satisfy skew-symmetry 
iff $(-1)^n b(w)_{r+i}a(w)=\del_w^{(-r-1-i+n)}c(z,w)|_{z=w}$ for any 
$r\in\Z$ and $i\geq 0$.
On the other hand, 
$$
b(w)_r a(w)
\;=\;
\del_w^{(-r-1+n)}((w-z)^n b(w)a(z))
\;=\;
(-1)^n\del_w^{(-r-1+n)}c'(z,w).
$$
Thus skew-symmetry holds iff $c$ and $c'$ have the same Taylor coefficients.
This is equivalent to $c=c'$ by Lemma \ref{SS:tayl}.
\end{pf}

\bigskip

{\bf Corollary.}\: {\it
Let $V$ be an associative vertex algebra.

\smallskip

\iti\:
Let $\phi$ be a morphism from $V$ to $\Endv(E)$ or $\cA\pau{z\uppm}_H$.
Then $\phi V$ is local.
In particular, any $V$-module is a $V\subslie$-module.

\smallskip

\itii\:
Let $M$ be a bounded $\Z$-fold $V$-module.
Then $M$ is a $V$-module iff 
$M$ is a $V\subslie$-module and $Y_M(ab)=\:\normord{Y_M(a)Y_M(b)}$.
}

\bigskip

\begin{pf}
\iti\:
Since $\phi$ is a morphism, skew-symmetry and weak locality for $V$ imply
skew-symmetry and weak locality for $\phi V$.
Thus $\phi V$ is local by the Proposition.
The second claim follows from the first since, by definition, 
a $V\subslie$-module $M$ is given by a morphism
$Y_M: V\subslie\to\End(M)\pau{z\uppm}$ such that $Y_M(V)$ is local.

\smallskip

\itii\:
This follows from \iti.
\end{pf}

\bigskip

Recall that $V\subslie$-modules are the same as $\fg(V)$-modules
by Remark \ref{SS:module vlie}.

\subsection{Dong's Lemma}
\label{SS:dong fields}

We prove that the unbounded vertex algebra structures of 
$\Endv(E)$ and $\cA\pau{z\uppm}_H$ preserve locality.
A refinement of this result and an alternative proof is given 
in section \ref{SS:duality of tensor alg}.

\bigskip

{\bf Remark.}\: {\it
Let $a(z), b(z)$ be in $\Endv(E)$ or $\cA\pau{z\uppm}_H$.
Then $a(z), b(z)$ are local of order $\leq n$ iff
there exists a distribution $c(z,w)$ that is regular at $z=w$ such that
$$
a(z)b(w)
\; =\;
\frac{c(z,w)}{(z-w)^n}, \qquad
\paraab\, b(w)a(z)
\; =\;
\frac{c(z,w)}{(z-w)_{w>z}^n}.
$$
}

\bigskip

\begin{pf}
`$\Rightarrow$'\: 
The distribution $c:=(z-w)^n a(z)b(w)=(z-w)^n b(w)a(z)$ is regular at $z=w$
since $E\lau{z}\lau{w}\cap E\lau{w}\lau{z}=E\lau{z,w}$ and since
distributions in $\cA\pau{z\uppm, w\uppm}_>\cap\cA\pau{w\uppm, z\uppm}_>$
are regular at $z=w$, see section \ref{SS:a z uppm h}.
We have $a(z)b(w)=c(z,w)/(z-w)^n$ since 
$\cF\lau{z}\lau{w}(E)$ is a vector space over $\K\lau{z}\lau{w}$
and $\cA\pau{z\uppm,w\uppm}_>$ is a module
over $\K[z\uppm,w\uppm,(z-w)\inv]$.
The second identity is proven in the same way.

\smallskip

`$\Leftarrow$'\:
We have $c(z,w)=(z-w)^n a(z)b(w)$
since $\cF\lau{z}\lau{w}(E)$ and $\cA\pau{z\uppm,w\uppm}_>$ are modules
over $\K[z\uppm,w\uppm,(z-w)\inv]$.
Similarly, $c(z,w)=(z-w)^n b(w)a(z)$. 
Thus $a(z), b(z)$ are local.
\end{pf}

\bigskip

{\bf Dong's Lemma.}\: {\it
Let $a(z), b(z), c(z)$ be in $\Endv(E)$ or $\cA\pau{z\uppm}_H$.
If they are 
\index{Dong's lemma}
pairwise local of orders $r, s, t$ then $a(z)_i b(z)$ and $c(z)$ are 
local of order $\leq r+s+t-i-1$ for any $i\in\Z$.
}

\bigskip

\begin{pf}
Let $r, s, t$ be the orders of locality of $ab, bc$, and $ac$.
We may assume that $i<r$.
Define $f:=(z-w)^r(w-y)^s(z-y)^t$.
As in the proof of the Remark we see that
$d:=fa(z)b(w)c(y)=fc(y)a(z)b(w)=fc(y)b(w)a(z)$ is regular at $z=w=y$
and
$$
a(z)b(w)c(y)
\; =\;
\frac{d(z,w,y)}{f(z,w,y)},
\qquad
c(y)a(z)b(w)
\; =\;
\frac{d(z,w,y)}{f(z,w,y)_{y>z>w}}.
$$
The claim now follows by applying $\del_z^{(-i-1+r)}((z-w)^r\cdot\; )|_{z=w}$ 
to these two equations and using 
Remark \ref{SS:skew sym fields} and the Remark.
\end{pf}

\subsection{Associative Vertex Algebras of Distributions}
\label{SS:vas of fields}

We prove that local subsets of $\Endv(E)$ and $\cA\pau{z\uppm}_H$ generate
associative vertex algebras.

\bigskip

{\bf Proposition.}\: {\it
Let $V$ be a vertex subalgebra of $\Endv(E)$ or $\cA\pau{z\uppm}_H$.
Then $V$ is associative iff $V$ is local. 
}

\bigskip

\begin{pf} 
The spaces $\Endv(E)$ and $\cA\pau{z\uppm}_H$ are pre-Lie algebras 
and satisfy the conformal Jacobi identity and the Wick formula by 
Propositions \ref{SS:norm ord prod pre lie} and \ref{SS:wick fields}.
Thus $V$ is associative iff $V$ satisfies skew-symmetry.
By Proposition \ref{SS:skew sym fields} this is equivalent to $V$ being local.
\end{pf}

\bigskip

By Proposition \ref{SS:assoc va} a unital vertex algebra is associative iff 
locality or the Jacobi identity hold.
Using one of these properties instead, one can prove more directly that 
$V$ is associative if $V$ is local.
Namely, one proves that local fields satisfy locality, 
see Proposition \ref{SS:holom loc implies itself} and 
Remark \ref{SS:holom loc implies itself},
or that local fields satisfy the Jacobi identity, 
see Proposition \ref{SS:jacobi fields}.

In section \ref{SS:dual skew loc}
we give an alternative proof of the fact that 
associative vertex subalgebras of $\Endv(E)$ are local.

Recall that if $S$ is a subset of a unital unbounded vertex algebra then 
$\sqbrack{S}$ denotes the unital unbounded vertex subalgebra generated by $S$.

\bigskip

{\bf Corollary.}\: {\it
Let $S$ be a local subset of $\Endv(E)$ or $\cA\pau{z\uppm}_H$.
Then $\sqbrack{S}$ is an associative unital vertex algebra.
}

\bigskip

\begin{pf}
The set $\sqbrack{S}$ is local because of Dong's lemma and because
$1$ is local to any $a(z)$.
Thus the claim follows from the Proposition.
\end{pf}

\section{The Frenkel-Zhu Algebra}
\label{S:fz alg}

We construct the Frenkel-Zhu algebra $\cA(V)$ and prove that $V\mapsto\cA(V)$
is a fully faithful functor into the category of local associative algebras.

\medskip

{\bf Conventions.}\:
All associative algebras and all unbounded vertex algebras are assumed to be 
unital.
Let $\rho\in\Q\uptimes$ such that $\rho\inv\in\N$.

\subsection{The Frenkel-Zhu Algebra}
\label{SS:fz alg}

We construct from any $\Z$-graded vertex algebra $V$ 
an FZ-complete algebra $\cA(V)$.
We 
\index{Frenkel-Zhu algebra}
call $\cA(V)$ the {\bf Frenkel-Zhu algebra}.

\bigskip

Let $V$ be a $\rho\Z$-graded associative vertex algebra.
There exists an FZ-complete algebra $\cA(V)$ and a morphism 
$Y: V\to\cA(V)\pau{z\uppm}_H$ of graded unbounded vertex algebras such that
the following universal property holds:
for any FZ-complete algebra $\cA$ and any morphism 
$\rho: V\to\cA\pau{z\uppm}_H$ 
there exists a unique topological $\rho\Z$-graded algebra morphism 
$\phi: \cA(V)\to\cA$ such that $\rho=\phi\circ Y$.
The pair $(\cA(V),Y)$ is unique up to a unique isomorphism.

The universal property says that, roughly speaking, the functor 
$V\mapsto\cA(V)$ is left adjoint to $\cA\mapsto\cA\pau{z\uppm}_H$.

We construct $\cA(V)$ as follows.
Let $\cA'$ be the tensor algebra of the vector space $V\otimes\K[x\uppm]$.
The algebra $\cA'$ is $\rho\Z$-graded where $a_t:=a\otimes x^t$
has weight $h-t-1$ for $a\in V_h$.
Let $\cA''$ be the FZ-completion of $\cA'$, 
see Proposition \ref{SS:fz complete}\,\itii.
Let $\cA(V)$ be the quotient of $\cA''$ by the degreewise closures of the ideal
generated by $1_t=\de_{t,-1}1$ and the {\bf associativity formula}
$$
(a_r b)_s 
\; =\;
\sum_{i\geq 0}(-1)^i\binom{r}{i}(a_{r-i}b_{s+i}-\paraab (-1)^r\, b_{s+r-i}a_i)
$$
for $a, b\in V$ and $t, r, s\in\Z$.
The quotient $\cA(V)$ is an FZ-complete algebra.
Define $Y: V\to\cA(V)\pau{z\uppm}_H$ by $a\mapsto\sum a_t z^{-t-1}$.
Then $Y$ is a morphism of unital $\Z$-fold algebras and hence of
unbounded vertex algebras by Remark \ref{SS:id elem}\iti.
It is straightforward to see that $\cA(V)$ and $Y$ satisfy the
universal property stated above.

The {\bf weak topology} of $\End(E)$ is the linear topology defined by
the subspaces $\set{a\mid a|_F=0}$ where $F$ runs through the
finite-dimensional subspaces of $E$. 
The weak topology is Hausdorff.
From $(a\circ b)|_F=a\circ(b|_F)$ follows that $\End(E)$ is a 
topological algebra. 
It is easy to see that any Cauchy filter in $\End(E)$ converges.
Thus $\End(E)$ is complete. 

If $E$ is a $\K$-graded vector spaces then $\End(E)$ contains a 
$\K$-graded subalgebra 
$$
\End(E)_H
\; :=\;
\bigoplus_{h\in\K}\End(E)_h
$$
where $\End(E)_h$ consists of endomorphisms $a$ such that 
$a: E_{h'}\to E_{h'+h}$; in other words, $[H,a]=ha$.
We endow $\End(E)_h$ with the weak topology.
Then $\End(E)_H$ is a complete topological $\K$-graded algebra. 

A $\K$-graded {\bf module} over a topological $\K$-graded algebra $\cA$ is 
\index{module!over a topological $\K$-graded algebra}
a $\K$-graded vector space $M$ together with a morphism
$\cA\to\End(M)_H$ of topological $\K$-graded algebras.

A $\Q$-graded vector space $E$ is {\bf truncated} if $E_h=0$ for 
$h\ll 0$.
A $\Q$-graded (vertex) algebra or a $\Q$-graded module are 
{\bf truncated} if the underlying $\Q$-graded vector space is truncated. 

If $E, E'$ are $\Q$-graded vector spaces that only differ by a shift,
$E_h=E'_{h+h_0}$ for some fixed $h_0$ and any $h$, 
then $\End(E)_H=\End(E')_H$.
Thus the $\Q$-graded module structures on $E$ and on $E'$ are in bijection.
In particular, if $M$ is a truncated $\Q$-graded module then, 
modulo a shift of the grading, $M$ is $\Q_{\geq}$-graded.

Suppose that $M$ is a $\rho\N$-graded vector space where $\rho\in\Q$
such that $\rho\inv\in\N$.
Then for any finite-dimensional subspaces $F\subset M$
there exists $k\in\Z$ such that $a|_F=0$ for any 
$a\in\sum_{l<k}\End(M)_{h-l}\End(M)_l$.
Thus $\End(M)_H$ is an FZ-complete algebra.
Moreover, $\End(M)_H\subset\Endv(M)$.
The universal property of $\cA(V)$ applied to $\cA=\End(M)_H$ yields:

\bigskip

{\bf Proposition.}\: {\it
Let $V$ be a $\rho\Z$-graded associative vertex algebra.
There is an isomorphism of categories $M\mapsto M$ from
$\rho\N$-graded $V$-modules to $\rho\N$-graded $\cA(V)$-modules.
\hfill $\square$
}

\subsection{Local Associative Algebras}
\label{SS:loc ass alg}

We define local associative algebras as FZ-complete algebras $\cA$
with a local subset $F\subcA\subset\cA\pau{z\uppm}_H$.
This notion has nothing to do with the usual notion of a local ring.

\bigskip

An {\bf FZ-complete differential} algebra is
\index{FZ-complete!differential algebra}
an FZ-complete algebra $\cA$ together with an even derivation $T$ of degree $1$
such that $T:\cA_h\to\cA_{h+1}$ is continuous.

For an FZ-complete differential algebra $\cA$, let
$$
\cA\pau{z\uppm}_{H,T}
\; :=\;
\cA\pau{z\uppm}_H\cap\cA\pau{z\uppm}_T
$$
be the subspace of translation covariant distributions: $Ta(z)=\del_z a(z)$.
This is a graded unbounded vertex subalgebra of $\cA\pau{z\uppm}_H$ because 
$T$ and $\del_z$ are derivations.

\bigskip

{\bf Definition.}\: 
A {\bf local associative algebra} is 
\index{local!associative algebra}
an FZ-complete differential algebra
$\cA$ together with a homogeneous local subset 
$F\subcA\subset\cA\pau{z\uppm}_{H,T}$ such that the subalgebra generated by 
$\set{a_n\mid a(z)\in F\subcA, n\in\Z}$ is degreewise dense in $\cA$.

\bigskip

A {\it morphism} $\phi: \cA\to\cA'$ of local associative algebras
is a morphism $\phi: \cA\to\cA'$ of FZ-complete algebras such that 
$\phi F\subcA\subset\sqbrack{F_{\cA'}}$.
Such a morphism satisfies $\phi\circ T=T\circ\phi$ since 
$\phi(Ta_n)=-n\phi(a_{n-1})=T\phi(a_n)$ and since the $a_n$ generate
a dense subalgebra. 

By Corollary \ref{SS:vas of fields},
if $\cA$ is a local associative algebra then $V(\cA):=\sqbrack{F\subcA}$
is a $\rho\Z$-graded associative vertex algebra. 
Thus we obtain a functor $\cA\mapsto V(\cA)$ from local associative algebras
to $\rho\Z$-graded associative vertex algebras.

Conversely, let $V$ be a $\rho\Z$-graded associative vertex algebra.
The Frenkel-Zhu algebra $\cA(V)$ has a unique structure of an 
FZ-complete {\it differential} algebra such that 
$Y(V)\subset\cA(V)\pau{z\uppm}_{H,T}$.
In fact, the endomorphism $T$ of $V[x\uppm]$ given by $T(a_n):=(Ta)_n$
induces a continuous derivation $T$ of $\cA(V)$ such that $Ta(z)=\del_z a(z)$
for any $a\in V$. 
Note that $(Ta)_n=-n\, a_{n-1}$ since $Y\circ T=\del_z\circ Y$.
By Corollary \ref{SS:skew sym fields}, the space $Y(V)$ is local.
It follows that $\cA(V)$ together with $Y(V)$ is a local associative algebra.

\bigskip

{\bf Proposition.}\: {\it
The functor $V\mapsto\cA(V)$ from $\rho\Z$-graded associative vertex algebras 
to local associative algebras is left adjoint to $\cA\mapsto V(\cA)$.
The morphism $Y: V\to V(\cA(V))$ is the unit of adjunction.
}

\bigskip

\begin{pf}
Let $V$ be a $\rho\Z$-graded associative vertex algebra and $\cA$ 
a local associative algebra.
By definition, a local associative algebra morphism $\cA(V)\to\cA$ is 
an FZ-complete algebra morphism $\psi:\cA(V)\to\cA$ such that 
$\psi Y(V)\subset V(\cA)$. 
By the universal property of $\cA(V)$, 
the map $\psi\mapsto\psi\circ Y$ is a bijection onto the set of
morphisms $\phi: V\to\cA\pau{z\uppm}_H$ such that 
$\phi(V)=\psi Y(V)\subset V(\cA)$.
These $\phi$s are just the $\rho\Z$-graded vertex algebra morphisms 
$V\to V(\cA)$.
\end{pf}

\subsection{The Embedding}
\label{SS:embedd assoc}

We prove that $V\mapsto\cA(V)$ is a fully faithful functor from 
truncated vertex algebras to local associative algebras.

\bigskip

Recall that a $\rho\Z$-graded vertex algebra $V$ is truncated iff $V_h=$
for $h\ll 0$.

\bigskip

{\bf Remark.}\: {\it
If $V$ is a truncated $\rho\Z$-graded associative vertex algebra then the
map $V\to\cA(V), a\mapsto a_{-1}$, is injective.
}

\bigskip

\begin{pf}
As pointed out in section \ref{SS:fz alg}, $V$ can be viewed as a 
$\rho\N$-graded $V$-module and hence, by Proposition \ref{SS:fz alg}, 
as a $\rho\N$-graded $\cA(V)$-module.
In other words, there is a morphism $\cA(V)\to\End(V)_H, a_n\mapsto a_n$.
The composition of the map $a\mapsto a_{-1}$ with the morphism 
$\cA(V)\to\End(V)_H\to V, a_n\mapsto a_n 1$, is the identity.
Thus the claim.
\end{pf}

\bigskip

{\bf Proposition.}\: {\it
The functor $V\mapsto\cA(V)$ from truncated $\rho\Z$-graded associative 
vertex algebras to local associative algebras is fully faithful.
}

\bigskip

\begin{pf}
Let $\cC$ be the category of truncated $\rho\Z$-graded associative 
vertex algebras and $\cC'$ the category of local associative algebras $\cA$
such that $V(\cA)$ is truncated.
If $V\in\cC$ then $Y: V\to V(\cA(V))$ is an isomorphism by the Remark.
In particular, $\cA(V)\in\cC'$.
By Proposition \ref{SS:fz alg} we thus have a pair of adjoint functors
between $\cC$ and $\cC'$ and the unit of adjunction $Y$ is an isomorphism.
Therefore $V\mapsto\cA(V)$ is fully faithful.
\end{pf}

\section{Supplements}
\label{S:suppl va}

In sections \ref{SS:ope heuris} and \ref{SS:ope2}
we show that the residue formula for the 
$t$-th products of $\Endv(E)$ computes 
the coefficients of the OPE.

In sections \ref{SS:holom loc implies itself} and \ref{SS:jacobi fields}
we prove that local fields satisfy locality and the Jacobi identity.

\subsection{Residue Formula and OPE: Heuristics}
\label{SS:ope heuris}

We explain {\it heuristically} that the residue formula for the 
$t$-th products of $\Endv(E)$ computes 
the coefficients of the operator product expansion (OPE).

\bigskip

Suppose that we can expand the operator product of two fields $a(z), b(z)$ 
as a singular power series in $z-w$:
$$
a(z)b(w)
\; =\;
\sum_{n\in\Z}\: c^n(w)\, (z-w)^{-n-1}.
$$
Let $C, C_{z>w}$, and $C_{w>z}^-$ be the following contours 
around $w$ and $0$ in the $z$-plane:
\bigskip
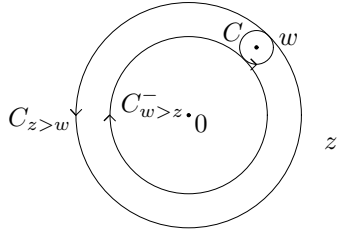
\begin{center}
   \setlength{\unitlength}{0.15cm}
   \begin{picture}(20,20)
     \put(10,10){\circle{20}}
     \put(10,10){\circle{14}}
     \put(10,10){\circle*{0.3}}
    \put(10,10){\circle*{0.2}}
     \put(10,10){\circle*{0.1}}

   \put(13.0,16.6){$C$}  
    \put(-6.0,9){$C_{z>w}$}  
   \put(4,10.3){$C_{w>z}^-$}  
  \put(22,7){$z$}  
   \put(18,16.2){$w$}  
    \put(10.5,8.5){$0$}
 
   \path(0,10)(-0.5,10.5)(0,10)(0.5,10.5)  
    \path(3,10)(2.5,9.5)(3,10)(3.5,9.5)  
  
     \put(16,16){\circle{3}}
      \put(16,16){\circle*{0.3}}
     \put(16,16){\circle*{0.2}}
     \put(16,16){\circle*{0.1}}
    \put(16,16){\circle*{0.25}}
     \put(16,16){\circle*{0.15}}
 \path(16,14.5)(15.5,15)(16,14.5)(15.5,14)  
    \end{picture}
\end{center}
\bigskip
Then the Cauchy residue theorem formally implies
\begin{align}
\notag
c^n(w)
=\;
&\res_{z=w}(z-w)^n a(z)b(w)
\\
\notag
=\;
&\oint_{C}(z-w)^n a(z)b(w)dz
\\
\notag
=\;
&\oint_{C_{z>w}}(z-w)^n a(z)b(w)dz
\; +\;
\oint_{C_{w>z}^-}(z-w)^n b(w)a(z)dz
\\
\notag
=\;
&\res_z (z-w)^n a(z)b(w)
\; -\;
\res_z (z-w)_{w>z}^n b(w)a(z)
\\
\notag
=\;
&a(w)_n b(w).
\end{align}
We have cheated twice. 
In the contour integral over $C_{w>z}^-$ we have changed the order of 
$a(z), b(w)$ because $a(z)b(w)$ is not well-defined in the domain $w>z$.
The residue theorem is not quite applicable to $C_{z>w}$ because of
the singularity at $z=w$.

\subsection{Residue Formula and OPE: A Proof}
\label{SS:ope2}

In Proposition \ref{SS:tayl} we proved that the $r$-th products of 
weakly local fields $a(z), b(z)$ are the Taylor coefficients of their 
operator product $a(z)b(w)$.
We give a new proof of this result in the case of local fields,
that is motivated by the heuristic argument of section \ref{SS:ope heuris}.

\bigskip

{\bf Remark.}\: {\it
For any $n\in\Z$, we have
$$
\del_w^{(n)}\de(z,w)
\; =\;
(z-w)^{-n-1}
\; -\;
(z-w)^{-n-1}_{w>z}.
$$
Equivalently, $\de(z,w+x)=\de(z-w,x)-\de(-w+z,x)$.
}

\bigskip

\begin{pf}
The identity holds true for $n<0$ because $(z-w)^m=(z-w)_{w>z}^m$ for 
$m\geq 0$.
It holds true for $n=0$ since 
$(z-w)^{-1}=\sum_{k\geq 0}\del_w^{(k)}(z-w)^{-1}|_{w=0} \, w^k=
\sum_{k\geq 0} w^k z^{-k-1}$ by Proposition \ref{SS:ratl fcts}
and $(z-w)_{w>z}^{-1}=-(w-z)^{-1}$.
Let $m\geq 0$.
Applying $\del_w$ to the above identity for $n=m$
we obtain the identity for $n=m+1$ times $m+1$.
Thus the claim follows by induction.
\end{pf}

\bigskip

{\bf Proposition.}\: {\it
If the fields $a(z), b(z)$ are local of order $\leq N$ then
for any $n\in\Z$ we have 
$$
a(w)_n b(w)
\; =\;
\del_z^{(-n-1+N)}((z-w)^N a(z)b(w))|_{z=w}.
$$
}

\bigskip

\begin{pf}
Proposition \ref{SS:ope sing} implies that $c(z,w):=(z-w)^N a(z)b(w)$ 
is in $\cF_2(E)$.
The Remark yields
\begin{align}
\notag
a(w)_n b(w)
\; &=\;
\res_z (z-w)^{n-N}c(z,w)-(z-w)_{w>z}^{n-N}c(z,w)
\\
\notag
\; &=\;
\res_z \del_w^{(-n-1+N)}\de(z,w)\, c(z,w)
\\
\notag
\; &=\;
\del_z^{(-n-1+N)}c(z,w)|_{z=w}
\end{align}
since $\res_z \del_w^k\de(z,w)c(z,w)=\res_z (-\del_z)^k\de(z,w)c(z,w)
=\res_z\de(z,w)\del_z^k c(z,w)\linebreak[0]=\del_z^k c(z,w)|_{z=w}$.
\end{pf}

\subsection{Locality for Local Fields}
\label{SS:holom loc implies itself}

Local fields generate an associative vertex algebra by 
Corollary \ref{SS:vas of fields}
and hence satisfy locality by Proposition \ref{SS:assoc va}.
We give two direct proofs of this fact.
Conversely, this fact provides a new proof of Corollary \ref{SS:vas of fields}
since a unital vertex algebra that satisfies locality is associative.

\bigskip

{\bf Proposition.}\: {\it
Let $a(z), b(z)$ be local fields.
Then the distributions $a(z)(Z)$ and $b(z)(Z)$ in 
$\End(\Endv(E))\pau{Z\uppm}$ are local.
}

\bigskip

\begin{pf}
Suppose that $a(z), b(z)$ are local of order $\leq n$ for some $n\geq 0$.
We have $a(w)(Z)c(w)=\res_z \de(z-w,Z)[a(z),c(w)]$ for $c(z)\in\Endv(E)$.
Thus
$$
a(x)(Z)\, b(x)(W)\, c(x)
\; =\;
\res_{z,w} \de(z-x,Z)\, \de(w-x,W)\, [a(z),[b(w),c(x)]].
$$
Since $z\de(z,w)=w\de(z,w)$ we have 
\begin{alignat}{2} 
\notag
&&(Z-W)^n\,\, &\de(z-x,Z)\, \de(w-x,W)
\\
\notag
=\;
&&((z-x)-(w-x))^n\,\, &\de(z-x,Z)\, \de(w-x,W)
\\
\notag
=\;
&&(z-w)^n\,\, &\de(z-x,Z)\, \de(w-x,W).
\end{alignat}
Together with the Leibniz identity we obtain
\begin{align}
\notag
&(Z-W)^n\, [a(x)(Z),b(x)(W)]c(x)
\\
\notag
=\;
&\res_{z,w} (z-w)^n\, \de(z-x,Z)\, \de(w-x,W)\, [[a(z),b(w)],c(x)]
\; =\;
0.
\end{align}
\end{pf}

\bigskip

{\bf Remark.}\: {\it
Let $V\subset\Endv(E)$ be a vertex subalgebra.
If $V$ is local then $Y(V)\subset\Endv(V)$ is local.
}

\bigskip

\begin{pf}
Let $a(z), b(z), c(z)\in V$ and $h_t\geq 0$ such that 
$a(z), b(z), c(z), a(z)_t c(z)$, and $b(z)_t c(z)$ are 
pairwise local of order $\leq h_t$.
Define $D_{h,k}:=(z-y)^h(w-y)^k a(z)b(w)c(y)$.
By Remark \ref{SS:skew sym fields} we have
$$
a(y)_t b(y)_s c(y)
\; =\;
\del_z^{(-t-1+h_s)}\del_w^{(-s-1+h_t)}D_{h_s,h_t}|_{z=w=y}.
$$
Using that $\del_z^{(n)}((z-w)a(z,w))|_{z=w}=\del_z^{(n-1)}a(z,w)|_{z=w}$
we see that the coefficient of $Z^{-t-1}W^{-s-1}$ of 
$(Z-W)^n a(y)(Z)b(y)(W)c(y)$ is 
\begin{align}
\notag
&\sum_{i\geq 0}(-1)^i\binom{n}{i}
\del_z^{(-t-n+i-1+h_s)}\del_w^{(-s-i-1+h_t)}
D_{h_s, h_t}|_{z=w=y} 
\\
\notag
=\;
&\sum_{i\geq 0}(-1)^i\binom{n}{i}
\del_z^{(-t-1+h_s)}\del_w^{(-s-1+h_t)}
D_{h_s+n-i, h_t+i}|_{z=w=y} 
\\
\notag
=\;
&\del_z^{(-t-1+h_s)}\del_w^{(-s-1+h_t)}
((z-w)^n D_{h_s, h_t})|_{z=w=y}.
\end{align}
Since the right-hand side is symmetric in $a(z), b(w)$ the claim follows.
\end{pf}

\subsection{Jacobi Identity for Local Fields}
\label{SS:jacobi fields}

Local fields generate a vertex algebra by Corollary \ref{SS:vas of fields}
and hence satisfy the Jacobi identity by Proposition \ref{SS:assoc va}.
We give a direct proof of this fact.

\bigskip

{\bf Lemma.}\: {\it
Let $t\in\Z$.

\smallskip

\iti\:
In $\K\lau{y}\pau{z,w}$ we have
$$
(z-y)^t_{y>z}
\; =\;
\sum_{i\geq 0}\:
\binom{t}{i}\,
(z-w)^i\,(w-y)_{y>w}^{t-i}.
$$

\smallskip

\itii\:
Let $a(z,w)\in\cF\lau{z}\lau{w}(E)+\cF\lau{w}\lau{z}(E)$ be local and
$b(z)\in\Endv(E)$ such that $(z-y)^n (w-y)^n [a(z,w),b(y)]=0$ for some 
$n\geq 0$.
Then
$$
(z-y)^t a(z,w)b(y)
\; =\;
\sum_{i\geq 0}\:
\binom{t}{i}\,
(z-w)^i\,(w-y)^{t-i}\, a(z,w)b(y).
$$
}

\bigskip

\begin{pf}
\iti\:
Define a morphism of fields over $\K$, $\io:\K(u,v)\to\K(z,w,y)$, by 
$u\mapsto z-w$ and $v\mapsto w-y$.
Let $(x_1,x_2,x_3)\in\set{(y,z,w), (y,w,z)}$.
Then there exists a unique continuous morphism $\tau$
such that we have a commutative diagram
$$
\xymatrix{
\K(u,v)\quad\ar[rr]^{\io} \ar[d]_{T_{v,u}} &&
\quad\K(z,w,y) \ar[d]^{T_{x_1,x_2,x_3}}
\\
\K\lau{v}\lau{u} \ar[rr]^{\tau\qquad} &&
\K\lau{x_1}\lau{x_2}\lau{x_3}.
}
$$
Indeed, uniqueness is clear. 
As for existence, let $f(u,v)\in\K\lau{v}\lau{u}$.
Then $f(z-w,-y+w')$ is in $\K\pau{z\uppm,w,w',y\uppm}$ and
$f(-w-z,-y+w')$ is in $\K\pau{y\uppm,w'}\lau{w}\pau{z}$.
Thus $f(z-w,-y+w)$ and $f(-w+z,-y+w)$ are well-defined.
Applying the diagram to $(u+v)^t$ we obtain the claim.

\smallskip

\itii\:
We may assume that $t<0$ and $(z-w)^n a(z,w)=0$. 
By \iti\ we have
\begin{align}
\notag
&(w-y)^{2n-t}(z-y)^n a(z,w)b(y)
\\
\notag
 =\,
&(w-y)^{2n-t}(z-y)^n b(y)a(z,w)
\\
\notag
 =\,
&\sum_{i\geq 0}\:
\binom{t}{i}\,
(z-w)^i\,(w-y)^{2n-i}\, (z-y)^{n-t}\, b(y)a(z,w)
\\
\notag
 =\,
&\sum_{i\geq 0}\:
\binom{t}{i}\,
(z-w)^i\,(w-y)^{2n-i}\, (z-y)^{n-t}\, a(z,w)b(y).
\end{align}
The main point is that the sum is finite with $i<n$ since
$a(z,w)$ is local.
Multiplying both sides by $(w-y)^{t-2n}(z-y)^{t-n}$
we obtain the assertion.
\end{pf}

\bigskip

{\bf Proposition.}\: {\it
Pairwise local fields satisfy the Jacobi identity.
}

\bigskip

\begin{pf}
The Lemma shows that the method of proof of the Wick formula for fields, 
see Proposition \ref{SS:wick fields}, also works in the present case.
More precisely, we apply the Lemma as it stands to 
$$
\res_{z,w}\,
(z-w)^r\, (w-y)^s\, (z-y)^t\, [[a(z),b(w)],c(y)]
$$
and obtain the left-hand side of the Jacobi identity.
(Part \iti\ of the Lemma is applied to the term involving $c(y)[a(z),b(w)]$
and part \itii\ is applied to $[a(z),b(w)]c(y)$.)
Then we apply the radially ordered Leibniz identity to the above expression
and apply the Lemma twice. 
Thus we obtain the right-hand side of the Jacobi identity. 
When the Lemma is applied to 
$$
\res_{z,w}\,
(z-w)^r\, (w-y)^s\, (z-y)^t\, [a(z),[b(w),c(y)]]
$$
the role of the variable $y$ in the Lemma is played by $z$
and we expand $(z-w)^r=((z-y)-(w-y))^r$ with $z>w$ and $z-y>w-y$. 
The point is that the field $a(z)$ corresponding to the variable $z$
always appears on the very left or very right 
of the four juxtapositions of $a(z), b(w), c(y)$. 
\end{pf}

\chapter{Basic Results}
\label{C:resul}

Sections \ref{S:va ids} and \ref{S:unital va ids}
are about implications between vertex algebra identities.
We prove that a vertex algebra is associative iff it satisfies the field
identities iff it satisfies skew-symmetry and one of the five three-element
field identities.
If there exists an identity element, 
both the Jacobi identity and locality alone are sufficient.
Under a mild further assumption, even the associativity formula suffices.
Many of these results are derived from two recursions and $\bbS_3$-symmetry
of the Jacobi identity.

\section{Field Identities I. Without Identity Element}
\label{S:va ids}

The vertex algebra identities are classified, first of all, according
to the number of elements for which they are defined.
For example, the Jacobi identity is 
\index{three-element identity}
a {\bf three-element} identity whereas skew-symmetry is 
\index{two-element identity}
a {\bf two-element} identity.
Moreover, the identities are divided into two groups as follows.

There are seven identities that are written in terms of the product and 
the $\la$-bracket: 
the pre-Lie identity and quasi-associativity, 
the conformal Jacobi identity, the left and right Wick formula, and
conformal skew-symmetry and $[\, ,]\subast=[\, ,]\subslie$.

The six {\bf field identities} are 
\index{field identity}
written in terms of fields:
the commutator and the associativity formula, locality and duality,
the Jacobi identity, and skew-symmetry.
The first four identities can be seen as the vertices of a square
whose vertical edges represent the fundamental recursion and
whose horizontal edges represent $\bbS_3$-symmetry.

Thus there are only three two-element identities.
They can be obtained from the three-element identities using an
identity element $1$, see section \ref{S:unital va ids}.

In section \ref{SS:loc and dual}
we discuss locality and duality.
In sections \ref{SS:fund recur} and \ref{SS:conseq fund rec}
we use the fundamental recursion to prove that 
the Jacobi identity is equivalent to locality and duality.
In section \ref{SS:s3 symm}
we use $\bbS_3$-symmetry to prove that if skew-symmetry holds then the other 
five field identities are equivalent to each other.
In sections \ref{SS:2nd recursion} and \ref{SS:assoc va}
we use a second recursion to prove that a vertex algebra is associative iff 
it satisfies the field identities.

In sections \ref{SS:q assoc zhu poisson} and \ref{SS:kind of assoc} 
we use quasi-associativity to show that $C_2(V)=(TV)V$ and that
$V_0/K$ is an associative algebra and 
we prove that duality implies a kind of associativity.

\subsection{Locality and Duality}
\label{SS:loc and dual}

We show that locality and duality are special cases of the Jacobi identity.

\bigskip

Let $V$ be a $\Z$-fold algebra and $M$ a $\Z$-fold $V$-module.
{\bf Locality} of order $\leq r$ is 
$$
(z-w)^r [a(z),b(w)]c
\; =\;
0.
$$
This is an identity in $M\pau{z\uppm,w\uppm}$.
If $r<0$, we assume that $M$ is bounded.
Then the identity involves radial ordering, see section \ref{SS:z fold fields}.
The coefficient of $z^{-t-1}w^{-s-1}$ is
$$
\sum_{i\geq 0}(-1)^i \binom{r}{i}
\big(a_{t+r-i}b_{s+i}c-\paraab (-1)^r\, b_{s+r-i}a_{t+i}c\big)
\; =\;
0.
$$
This is locality of order $\leq r$ for {\bf indices} $t, s$. 

We say that $M$ satisfies {\bf locality} if for 
\index{locality}
any $a, b, c$ there exists $r\geq 0$ such that 
$a, b, c$ satisfy locality of order $\leq r$. 
Of course, $Y_M(V)$ is local iff $M$ satisfies locality and $r$ depends
only on $a$ and $b$, but not on $c$.

For $a\in V$ and $b\in M$, define $o(a,b)\in\Z\cup\set{\pm\infty}$ to be the 
least $t\in\Z$ such that $a_s b=0$ for $s\geq t$.
If $V=M$ is bounded then $o:V^2\to\Z\cup\set{-\infty}$ is 
\index{locality function!of a vertex algebra}
the {\bf locality function} of $V$.
See section \ref{SS:loc fct} for properties of $o$.

If $r\geq o(a,b)$ then the Jacobi identity for indices $r, s, t$ is 
locality of order $\leq r$ for indices $t, s$.
Since $o(a,b)$ does not depend on $c$, the Jacobi identity implies that
$Y_M(V)$ is local.

Assume that $V$ and $M$ are bounded.
{\bf Duality} of order $\leq t$ is 
$$
(w+x)^t\, (a(x)b)(w)c
\; =\;
(x+w)^t\, a(x+w)b(w)c.
$$
Since $\binom{-r-1}{i}=(-1)^i\binom{r+i}{i}$,
the coefficient of $w^{-s-1}x^{-r-1}$ is
$$
\sum_{i\geq 0}\:
\binom{t}{i}\,
(a_{r+i}b)_{s+t-i}c
\; =\;
\sum_{i\geq 0}\:
(-1)^{i} \binom{r}{i}\,
a_{t+r-i}b_{s+i}c.
$$
This is duality of order $\leq t$ for {\bf indices} $r, s$.

We say that $M$ satisfies {\bf duality} if for 
\index{duality}
any $a, b, c$ there exists $t$ such that 
$a, b, c$ satisfy duality of order $\leq t$.
We call $M$ {\bf dual} if for any $a, c$ there exists $t$ such that
for any $b$ the elements $a, b, c$ satisfy duality of order $\leq t$.

If $t\geq o(a,c)$ then the Jacobi identity for indices $r, s, t$ is  
duality of order $\leq t$ for indices $r, s$.
Since $o(a,c)$ does not depend on $b$, the Jacobi identity implies that
$M$ is dual.

\subsection{The Fundamental Recursion}
\label{SS:fund recur}

We prove a recursion for the Jacobi identity that follows from $x=z-w$.

\bigskip

Let $E$ be a vector space and $e_i(z,w)\in E\pau{z\uppm,w\uppm}$ for 
$i=1, 2, 3$.
The {\bf Jacobi identity} for $e_1, e_2, e_3$ and indices $r, s, t$ is
\begin{align}
\notag
&\res_{w,x}(w+x)^t w^s x^r e_1(w,x)
\\
\notag
=\;
&\res_{z,w}z^t w^s ((z-w)^r e_2(z,w)-(-w+z)^r e_3(w,z)).
\end{align}
Here we assume either that $r, t\geq 0$ or that $e_i(z,w)\in E\lau{z}\lau{w}$.

Let $M$ be a $\Z$-fold module over a $\Z$-fold algebra $V$.
The Jacobi identity for $a, b\in V, c\in M$ is the Jacobi identity
for $(a(w)b)(z)c, a(z)b(w)c$, and $\paraab\, b(z)a(w)c$.
Likewise the conformal Jacobi identity and the Wick formula 
for a conformal operator $d\submu$ are the Jacobi identity for 
$(d\subnu a)(z)b, d\submu(a(w)b)$, and $a(z)d\subnu b$ and 
indices $r, s\geq 0, t=0$, resp., $r=0, s=-1, t\geq 0$.
Here we use the identification of $E\pau{\mu}$ with $z\inv E\pau{z\inv}$,
see section \ref{SS:unb confa power series}.

\bigskip

{\bf Lemma.}\: {\it
Let $E$ be a vector space and $e_i\in E\pau{z\uppm,w\uppm}$.
The Jacobi identity for two of the triples $(r+1,s,t), (r,s+1,t)$, 
and $(r,s,t+1)$ implies the Jacobi identity for the third triple.
}

\bigskip

\begin{pf}
Let $J_{r,s,t}$ denote the left-hand side of the Jacobi identity 
for indices $r, s, t$.
Then $J_{r,s,t+1}=J_{r+1,s,t}+J_{r,s+1,t}$.
The right-hand side satisfies the same identity.
This implies the claim.
\end{pf}

\bigskip

{\bf Proposition.}\: {\it
Let $E$ be a vector space and $e_i\in E\pau{z\uppm,w\uppm}$.

\smallskip

\iti\:
The Jacobi identity for indices $r=r_0, s\geq s_0, t\geq t_0$ implies 
the Jacobi identity for indices $r\geq r_0, s\geq s_0, t\geq t_0$.

\smallskip

\itii\:
The Jacobi identity for indices $r\geq r_1, s\in\Z, t\geq t_0$ and 
$r\geq r_0, s\in\Z, t\geq t_1$ implies
the Jacobi identity for indices $r\geq r_0, s\in\Z, t\geq t_0$.

\smallskip

These two results also hold true for any permutation of $r, s, t$.
}

\bigskip

\begin{pf}
\iti\:
By the Lemma the Jacobi identity for $r_0, s+1, t$ and for $r_0, s, t+1$ 
implies the Jacobi identity for $r_0+1, s, t$.

\smallskip
 
\itii\:
We may assume that $r_1> r_0$ and $t_1> t_0$.
By the Lemma the Jacobi identity for $r_1, s, t_1-1$ and $r_1-1, s, t_1$ 
implies the Jacobi identity for $r_1-1, s+1, t_1-1$.
Thus the claim follows by induction.
\end{pf}

\subsection{Consequences of the Recursion}
\label{SS:conseq fund rec}

We prove that the Jacobi identity is equivalent to locality and duality
and that the associativity formula is equivalent to duality and weak locality.

\bigskip

{\bf Proposition.}\: {\it
Let $M$ be a $\Z$-fold module over a $\Z$-fold algebra $V$.

\smallskip

\iti\:
The conformal Jacobi identity, the Jacobi identity, the associativity formula, 
and the commutator formula for indices $r, s, t\geq 0$ are all equivalent.

\smallskip

\itii\:
The commutator formula is equivalent to locality and 
$(a_r b)(z)c=a(z)_r b(z)c$ for $r\geq 0$.
In this case $Y_M(V)$ is local. 

\smallskip

\itiii\:
The associativity formula holds iff duality holds and 
$a(z), b(z)$ are weakly local on $c$.
In this case $M$ is dual.

\smallskip

\itiv\:
Duality and locality are equivalent to the Jacobi identity.
}

\bigskip

Parts \iti\ and \itii\ were already proven in Propositions 
\ref{SS:j id conf j id} and \ref{SS:module vlie}.

\bigskip

\begin{pf}
\iti\:
In section \ref{SS:modules} we have shown that the conformal Jacobi
identity is the associativity formula for indices $r, s\geq 0$.
Thus the claim follows from part \iti\ of Proposition \ref{SS:fund recur}
with $r_0=s_0=t_0=0$.

\smallskip

\itii\:
By part \iti\ of Proposition \ref{SS:fund recur} the commutator formula 
implies the Jacobi identity for indices $r\geq 0$ and $s, t\in\Z$.
For $r\geq o(a,b)$, we obtain that $Y_M(V)$ is local by 
section \ref{SS:loc and dual}.
For $t=0$, we obtain $(a_r b)(z)c=a(z)_r b(z)c$.
The converse follows from part \itii\ of Proposition \ref{SS:fund recur}.

\smallskip

\itiii\:
By part \iti\ of Proposition \ref{SS:fund recur} the associativity formula 
implies the Jacobi identity for indices $r, s\in\Z$ and $t\geq 0$.
For $t\geq o(a,c)$, we obtain that $V$ is dual by 
section \ref{SS:loc and dual}. 
Since $a_r b=0$ for $r\gg 0$, the associativity formula implies 
weak locality on $c$, that is, $a(z)_r b(z)c=0$ for $r\gg 0$. 
Conversely, weak locality implies that the associativity formula holds for 
indices $r\gg 0, s\in\Z$.
Thus the converse follows from part \itii\ of Proposition \ref{SS:fund recur}.

\smallskip

\itiv\:
This follows from part \itii\ of Proposition \ref{SS:fund recur}.
\end{pf}

\bigskip

{\bf Corollary.}\: {\it
Let $V$ be an associative vertex algebra and $M$ a $V$-module.
Then $M$ satisfies the Jacobi identity.
}

\bigskip

\begin{pf}
Corollary \ref{SS:skew sym fields} shows that $M$ satisfies locality.
By the Proposition 
locality and the associativity formula imply the Jacobi identity.
\end{pf}

\subsection{$\bbS_3$-Symmetry}
\label{SS:s3 symm}

We use $\bbS_3$-symmetry of the Jacobi identity
to prove that if skew-symmetry holds then the other five field identities are 
equivalent to each other.

\bigskip

{\bf Remark.} \: {\it
Let $E$ be a $\K[T]$-module.
Then $a(z)\in\End(E)\pau{z\uppm}$ is translation covariant iff
$$
a(z+w)
\; =\;
e^{wT}a(z)e^{-wT}.
$$
}

\bigskip

\begin{pf}
First an observation:
If $c\in E$ then $b(z)=e^{zT}c$ is the {\it unique}
power series such that $Tb(z)=\del_z b(z)$ and $b(0)=c$,
since $b(z)=e^{z\del_w}b(w)|_{w=0}=e^{zT}b(0)$.

Recall that $a(z+w)=e^{w\del_z}a(z)$.
Taking the coefficient of $w$ of the identity $a(z+w)=e^{wT}a(z)e^{-wT}$
we get $\del_z a(z)=[T,a(z)]$.

Conversely, both $a(z+w)$ and $e^{wT}a(z)e^{-wT}$ are power series in $w$, 
both have constant term $a(z)$, and 
both satisfy $\del_w b(w)=[T,b(w)]$.
Thus the converse follows from the observation.
\end{pf}

\bigskip

$\bbS_3$-symmetry of the Jacobi identity permutes $a, b, c$ and 
the indices $r, s, t$.
It is analogous to the $\bbS_3$-symmetry of the Leibniz identity 
$[[a,b],c]=[a,[b,c]]-[b,[a,c]]$.
Namely, if $[\, ,]$ is skew-symmetric then the Leibniz identity is 
equivalant to $[[b,a],c]=[b,[a,c]]-[a,[b,c]]$ and also to 
$[[a,c],b]=[a,[c,b]]-[c,[a,b]]$.

\bigskip

{\bf Lemma.}\: {\it
Let $V$ be a bounded $\Z$-fold algebra that satisfies skew-symmetry 
for some even operator $T$.

\smallskip

\iti\:
If $T$ is a translation generator then
the Jacobi identity for $a, b, c$ and indices $r$, any $s'\geq s$, and $t$ 
is equivalent to 
the Jacobi identity for $a, c, b$ and indices $t$, any $s'\geq s$, and $r$.

\smallskip

\itii\:
Let $M$ be a bounded $\Z$-fold $V$-module such that $Y_M(Ta)=\del_z Y_M(a)$.
Then the Jacobi identity holds for $a, b, c$ and indices $r, s, t$ 
iff it holds for $b, a, c$ and indices $r, t, s$.
}

\bigskip

\begin{pf}
Applying $\res_{z,w,x}z^{-t-1}w^{-s-1}x^{-r-1}$ to the Jacobi identity
$\de(z,w+x)\linebreak[5](a(x)b)(w)c=\de(z-w,x)[a(z),b(w)]c$, we obtain 
the Jacobi identity for indices $r, s, t$.
In the following we omit the fact that we actually consider this residue.
We use the Remark and that $\de(z,w)=\de(w,z)$ and $\de(z+x,w)=\de(z,w-x)$
by Proposition \ref{SS:delta dis}\,\itii,\itiii.

\smallskip

\iti\:
We show that skew-symmetry transforms the Jacobi identity for $a, b, c$ into 
the Jacobi identity for $a, c, b$ by interchanging the first and third term
and replacing $z, w, x$ by $x, -w, z$. 
In fact, we have
$$
\de(z,w+x)(a(x)b)(w)c
\; =\;
\de(w+x,z)e^{wT}c(-w)a(x)b,
$$
\begin{align}
\notag
\de(z-w,x)a(z)b(w)c
\; &=\;
\de(z-w,x)a(x+w)e^{wT}c(-w)b
\\
\notag
&=\;
\de(x+w,z)e^{wT}a(x)c(-w)b,
\end{align}
and 
$$
\de(-w+z,x)b(w)a(z)c
\; =\;
\de(x,-w+z)e^{wT}(a(z)c)(-w)b.
$$

\smallskip

\itii\:
We show that skew-symmetry transforms the Jacobi identity for $a, b, c$ into 
the Jacobi identity for $b, a, c$ by interchanging the second and third term
and replacing $z, w, x$ by $w, z, -x$. 
In fact, we have
\begin{align}
\notag
\de(z,w+x)(a(x)b)(w)c
\; &=\;
\de(z,w+x)(e^{xT}b(-x)a)(w)c
\\
\notag
&=\;
\de(z,w+x)e^{x\del_w}(b(-x)a)(w)c
\\
\notag
&=\;
\de(z,w+x)(b(-x)a)(w+x)c
\\
\notag
&=\;
\de(z,w+x)(b(-x)a)(z)c
\\
\notag
&=\;
\de(w,z-x)(b(-x)a)(z)c
\end{align}
and $\de(z-w,x)[a(z),b(w)]c=\de(w-z,-x)[b(w),a(z)]c$.
The last equation follows from $\de(z,w)=-\de(-z,-w)$ and from radial ordering.
\end{pf}

\bigskip

Recall that skew-symmetry is the only two-element field identity and 
there are five three-element field identities:
the commutator and the associativity formula, locality and duality, and
the Jacobi identity, see the introduction of section \ref{S:va ids}.

\bigskip

{\bf Proposition.}\: {\it
Let $V$ be a bounded $\Z$-fold algebra with a translation generator.
If skew-symmetry and one of the three-element field identities hold
then all field identities are satisfied.
}

\bigskip

\begin{pf}
By part \iti\ of the Lemma, locality and duality are equivalent.
Thus the claim follows from Proposition \ref{SS:conseq fund rec}.
\end{pf}

\subsection{A Second Recursion}
\label{SS:2nd recursion}

We prove a recursion for the associativity and commutator formula
that decreases the indices.
This recursion was proven for skew-symmetry in Proposition \ref{SS:skew sym}.

\bigskip

{\bf Proposition.}\: {\it
Let $V$ be an unbounded vertex algebra and $M$ a $\Z$-fold $V$-module with 
translation endomorphism.

\smallskip

\iti\:
The associativity formula for indices $r, s$ implies the associativity formula 
for indices $r-1, s$ if $r\ne 0$ and for indices $r, s-1$ if $s\ne 0$.

\smallskip

\itii\:
The commutator formula for indices $t, s$ implies the commutator formula 
for indices $t-1, s$ if $t\ne 0$ and for indices $t, s-1$ if $s\ne 0$.
}

\bigskip

\begin{pf}
\iti\:
The claim for $r\ne 0$ follows from 
$-r(a_{r-1}b)_s c=((Ta)_r b)_s c=((\del_z a(z))_r b(z))_s c
=-r(a(z)_{r-1}b(z))_s c$.

The claim for $s\ne 0$ follows from 
$-s(a_r b)_{s-1}c=(T(a_r b))_s c=((Ta)_r b+a_r Tb)_s c
=(\del_z(a(z)_rb(z)))_s c=-s(a(z)_r b(z))_{s-1}c$.

\smallskip

\itii\:
First, assume that $V$ is bounded.
We will use that 
$\del_x \de(z,w+x)=\del_x e^{x\del_w}\de(z,w)=
e^{x\del_w}\del_w \de(z,w)=\del_w \de(z,w+x)=-\del_z\de(z,w+x)$.
If $t\ne 0$ then the claim follows from $[(Ta)_t,b_s]c=-t[a_{t-1},b_s]c$ and
\begin{align}
\notag
&\res_{z,w,x}z^t w^s\,  \de(z,w+x)((Ta)(x)b)(w)c
\\
\notag
\; =\;
-&\res_{z,w,x}z^t w^s\,  \del_x \de(z,w+x)(a(x)b)(w)c
\\
\notag
\; =\;
-&\res_{z,w,x}\del_z(z^t)w^s\,  \de(z,w+x)(a(x)b)(w)c.
\end{align}
If $s\ne 0$ then the claim follows from $[a_t,(Tb)_s]c=-s[a_t,b_{s-1}]c$ and
\begin{align}
\notag
&\res_{z,w,x}z^t w^s\, \de(z,w+x) (a(x)Tb)(w)c
\\
\notag
\; =\;
&\res_{z,w,x}z^t w^s\, \de(z,w+x)\,(\del_w-\del_x) (a(x)b)(w)c
\\
\notag
\; =\;
&-\res_{z,w,x}z^t\del_w(w^s) \,\de(z,w+x) (a(x)b)(w)c.
\end{align}
If $V$ is not bounded then $t\geq 0$. 
In this case the commutator formula is 
$[a\subla,b(w)]c=\res_x e^{(w+x)\la}(a(x)b)(w)c$ and the claim
is proven similarly.
\end{pf}

\subsection{Axioms for Associative Vertex Algebras}
\label{SS:assoc va}

We prove that a vertex algebra is associative iff it satisfies 
the field identities.

\bigskip

Recall that there are six field identities: 
the commutator and the associativity formula, locality and duality, 
the Jacobi identity, and skew-symmetry, 
see the introduction of section \ref{S:va ids}.

\bigskip

{\bf Proposition.}\: {\it
A vertex algebra is associative iff the field identities hold.
}

\bigskip

\begin{pf}
By Proposition \ref{SS:skew sym} conformal skew-symmetry and 
$[\, ,]\subast=[\, ,]\subslie$ are equivalent to skew-symmetry.
Hence we may assume that $V$ satisfies skew-symmetry.
By Propositions \ref{SS:s3 symm} and \ref{SS:2nd recursion} it suffices to 
prove that $V$ is associative iff 
the commutator formula holds for indices $t, s\geq -1$.

The commutator formula for indices $t, s\geq 0$ is equivalent to
the conformal Jacobi identity by Proposition \ref{SS:j id conf j id}.
For indices $t\geq 0, s=-1$ it is the Wick formula, 
see section \ref{SS:conf j wick0}.
For indices $t=-1, s\geq 0$ it is equivalent to the Wick formula 
by $\bbS_3$-symmetry.
For indices $t=s=-1$ it is
$$
[a\cdot, b\cdot]
\; =\;
\sum\binom{-1}{i}\, (a_i b)_{-2-i}
\; =\;
\sum\: (-1)^i\, T^{(i+1)}(a_i b)\cdot.
$$
Since $[\, ,]\subast=[\, ,]\subslie$,
this identity is $[a\cdot, b\cdot]=[a, b]\cdot$.
\end{pf}

\subsection{Quasi-Associativity}
\label{SS:q assoc zhu poisson}

We use quasi-associativity to prove that $C_2(V)=(TV)V$ and that
if $V$ is graded then a certain quotient of $V_0$ is an associative algebra.

\bigskip

{\bf Quasi-associativity} for 
\index{quasi-associativity}
a vertex algebra $V$ is 
$$
(ab)c-abc
\; =\;
\Big(\int_0^T d\la \, a\Big) [b\subla c]
\; +\;
\paraab\Big(\int_0^T d\la \, b\Big) [a\subla c].
$$
This is equivalent to
$$
(ab)c-abc
\; =\;
\sum_{i\geq 0}\:
T^{(i+1)}(a)b_i c
\; +\;
\paraab\,
T^{(i+1)}(b)a_i c.
$$
Since the right-hand side is symmetric in $a$ and $b$,
quasi-associativity implies that $V\subast$ is a pre-Lie algebra.

Quasi-associativity is the associativity formula for indices $-1, -1$.
Hence it is satisfied for any associative vertex algebra.
In section \ref{SS:pre lie q assoc}
we prove quasi-associativity more directly and under weaker assumptions.

Here is an example.
Let $J\in V$ be a $U(1)$-vector such that $[J\subla J]=\la$, 
see section \ref{SS:va}.
Then quasi-associativity is $(JJ)J=J(JJ)+T^2 J$.

Recall that $C_n(V)$ is the ideal of $V\subast$ generated by $T^{n-1}V$, 
see section \ref{SS:assoc va def}.

\bigskip

{\bf Proposition.}\: {\it
Let $V$ be an associative unital vertex algebra.

\smallskip

\iti\:
We have $C_n(V)=(T^{n-1}V)V$ for $n\geq 2$.

\smallskip

\itii\:
Suppose that $V$ is $\K$-graded. 
Then $K:=\sum_{i\geq 1}(T^i V_0)V_{-i}$ is an ideal of $(V\subast)_0$ and 
$(V\subast)_0/K$ is an associative algebra.
}

\bigskip

\begin{pf}
\iti\:
Quasi-associativity implies that $(T^{n-1}V)V$ is a right ideal of $V\subast$
since $(T^m a)_i c=0$ for $i<m:=n-1$.
It is a left ideal since $aT^m(b)c=T^m(b)ac+[a,T^m b]c$ and 
$[a,T^m b]\in T^{m+1} V$ because $V\subast$ is almost commutative, 
see section \ref{SS:assoc va def}. 
Obviously, $T^{n-1}(V)V$ is a $\K[T]$-submodule.

\smallskip

\itii\:
Let $a, b\in V_0$ and $c\in V_{-i}$.
The subspace $K$ is a left ideal because $aT^i(b)c=T^i(b)ac+[a,T^i b]c$ and
$[a,T^i b]\in T^{i+1}V_{-1}\subset T^i V_0$ since 
$V\subast$ is almost commutative, see section \ref{SS:assoc va def}.  

Quasi-associativity implies that $K$ is a right ideal because
$T^{(j+1)}(T^i b)c_j a\linebreak[0]\in K$ and 
$T^{(j+1)}(c)(T^i b)_j a\in K$, using that $(T^i b)_j a=0$ for $i>j$.

Quasi-associativity implies that $V_0/K$ is associative.
\end{pf}

\bigskip

Part \iti\ is especially simple to prove for $n=2$. 
Quasi-associativity implies that $(TV)V$ is a right ideal.
It is a left ideal since $[V,V]\subset TV$. 

We mention that $C_n(V)=(T^{n-1}V)V$ is also an ideal with respect to $a_0 b$ 
because $a_0$ is a derivation and because of the right Wick formula, 
see section \ref{SS:right wick}.
Note that $K\subset C_2(V)_0$.

\subsection{A Kind of Associativity}
\label{SS:kind of assoc}

We prove that if duality is satisfied then any product $(a_r b)_s c$ is 
equal to a linear combination of products $a_p (b_q c)$.

\bigskip

{\bf Proposition.}\: {\it
Let $V$ be a bounded $\Z$-fold algebra, $M$ a bounded $\Z$-fold $V$-module, 
and $r, s\in\Z$.
Assume that $a, b\in V$ and $c\in M$ satisfy duality of order $\leq t$
for some $t\in\Z$. 

\smallskip

\iti\:
Let $N:=1-t-o(a,b)$. 
Then there exist $\la^i\in\Z$, 
depending only on $r, s, t, N$ but not on $a, b, c, V$, such that
$$
(a_r b)_s c
\; =\;
\sum_{i\geq N}\: \la^i\, a_{-i}(b_{s+r+i}c).
$$

\smallskip

\itii\:
Let $K:=1-t-o(b,c)$.
Then there exist $\ka^i\in\Z$,
depending only on $r, s, t, K$ but not on $a, b, c, V$, such that 
$$
a_r(b_s c)
\; =\;
\sum_{i\geq K}\: \ka^i\, (a_{r+s+i}b)_{-i}c.
$$
}

\bigskip

\begin{pf}
\iti\:
Duality of order $\leq t$ for indices $r, s-t$ is
$$
(a_r b)_s c
\; =\;
\sum_{i\geq 0} (-1)^{i} \binom{r}{i}\, a_{t+r-i}b_{s-t+i}c
\; -\;
\sum_{i>0} \binom{t}{i} (a_{r+i}b)_{s-i}c.
$$
Applying duality again for indices $r+j, s-t-j$
to the right-hand side, where $j=1, \dots, o(a,b)-r-1$,
one obtains the claim. 

\smallskip

\itii\:
This is proven in the same way.
\end{pf}

\bigskip

{\bf Corollary.}\: {\it
Let $V$ be a unital associative vertex algebra, $S\subset V$ a subset,
and $M$ a $V$-module.

\smallskip

\iti\:
We have 
$\sqbrack{S}=
\rspan\set{a^1_{n_1}\dots a^r_{n_r}1\mid a^i\in S, n_i\in\Z, r\geq 0}$.

\smallskip

\itii\:
Suppose that $M$ is generated by a subset $S'\subset M$.
Then $M=\rspan\set{a_t b\mid a\in V, b\in S', t\in\Z}$.

\smallskip

\itiii\:
If $V=\sqbrack{S}$ and $N\subset M$ is a subspace such that 
$a_t N\subset N$ for any $a\in S, t\in\Z$ then $N$ is a submodule.
}

\bigskip

\begin{pf}
This follows directly from the Proposition.
\end{pf}

\section{Field Identities II. With Identity Element}
\label{S:unital va ids}

In sections \ref{SS:loc skew sym} and \ref{SS:unital assoc va}
we prove that a unital vertex algebra is associative
iff it satisfies locality iff the Jacobi identity holds.

In sections \ref{SS:conf j wick}--\ref{SS:q assoc assoc}
we prove that the associativity formula is equivalent to 
four identities for $a\cdot$ and $a\subla$.
In sections \ref{SS:assoc alg bracket}--\ref{SS:eq commut}
we prove that the associativity formula implies associativity 
if $a_t b=1$ for some $a, b, t$.

\subsection{Locality implies Skew-Symmetry}
\label{SS:loc skew sym}

We prove that locality implies skew-symmetry.

\bigskip

Let $E$ be a $\K[T]$-module with an even vector $1$.
A distribution $a(z)\in\End(E)\pau{z\uppm}$ is 
\index{creative}
{\bf creative} if $a(z)1=e^{zT}a_{-1}1$.
We say that $1$ is 
\index{invariant vector}
{\bf invariant} if $T1=0$.

\bigskip

{\bf Lemma.} \: {\it
Suppose that $1$ is invariant.
If $a(z)\in\End(E)\pau{z\uppm}$ is translation covariant and 
$a(z)1\in E\lau{z}$ 
then $a(z)$ is creative. 
}

\bigskip

\begin{pf}
We have $a(z)1\in E\pau{z}$ because otherwise the pole order of $\del_z a(z)1$
is strictly smaller than the pole order of $[T,a(z)]1$.
By Remark \ref{SS:s3 symm} we have $a(z+w)1=e^{wT}a(z)e^{-wT}1$.
Setting $z=0$ we get $a(w)1=e^{wT}a_{-1}1$.
\end{pf}

\bigskip

Let $V$ be a $\Z$-fold algebra.
Recall that an even vector $1$ is a left identity if
$1(z)=\id_V$ and 
a right identity if $a(z)1=e^{zT}a$ for some even operator $T$,
see section \ref{SS:id elem}.

If $V$ is bounded and skew-symmetry holds for some even operator $T'$ 
then $1$ is a left identity iff $1$ is a right identity and $T'a=a_{-2}1$.
Indeed, from $1(z)a=a$ we get $a(z)1=e^{zT'}1(-z)a=e^{zT'}a$.
Conversely, we have $1(z)a=e^{zT'}a(-z)1=e^{zT'}e^{-zT}a=a$.

A {\bf weak right identity} is 
\index{weak right identity}
an even vector $1$ such that $a(z)1\in V\pau{z}$ and $a(0)1=a$ 
for any $a\in V$.

\bigskip

{\bf Proposition.}\: {\it
Let $V$ be a $\Z$-fold algebra with a translation generator $T$ and 
an invariant weak right identity. 
If $a, b, 1$ satisfy locality  
then $a, b$ are weakly local and satisfy skew-symmetry for $T$.
}

\bigskip

\begin{pf}
Since $a(z)1, b(z)1\in V\pau{z}$, 
we get $(z-w)^r a(z)b(w)1=(z-w)^r b(w)a(z)1\linebreak[0]\in V\pau{z,w}$
for $r\gg 0$.
Setting $w=0$ yields $z^r a(z)b\in V\pau{z}$.
Thus $a, b$ are weakly local.
Likewise $z^r b(z)a\in V\pau{z}$.
Together with the Lemma we obtain 
\begin{align}
\notag
(z-w)^r a(z)b(w)1
\, &=\, 
(z-w)^r b(w)a(z)1
\, =\, 
(z-w)^r b(w)e^{zT}a
\\
\notag
&=\, 
(z-w)^r e^{zT}b(w-z)a
\, =\, 
(z-w)^r e^{zT}b(-z+w)a.
\end{align}
Setting $w=0$ yields $a(z)b=e^{zT}b(-z)a$.
\end{pf}

\subsection{Axioms for Associative Unital Vertex Algebras}
\label{SS:unital assoc va}

We prove that a unital vertex algebra is associative iff locality holds iff
the Jacobi identity holds.
In fact, we prove a stronger claim about $\Z$-fold algebras.

\bigskip

{\bf Proposition.}\: {\it
Let $V$ be a $\Z$-fold algebra. 
The following is equivalent: 

\begin{enumerate}

\item[\iti]
$V$ is an associative unital vertex algebra;

\smallskip

\item[\itii]
$V$ is bounded, has a left identity, and satisfies duality and skew-symmetry;

\smallskip

\item[\itiii]
$V$ is bounded, has a weak right identity, and satisfies the Jacobi identity;

\smallskip

\item[\itiv]
there exist a translation generator and 
an invariant weak right identity and $V$ satisfies locality.
\end{enumerate}
}

\bigskip

\begin{pf}
A more precise statement of \iti\ is: 
$V$ is the $\Z$-fold algebra underlying an associative unital vertex algebra. 
This unital vertex algebra is unique since $1$ is unique and hence 
$T$ is unique, $Ta=(Ta)1=a_{-2}1$.

\smallskip

\iti$\Rightarrow$\itii\:
This follows from Proposition \ref{SS:assoc va}.

\smallskip

\itii$\Rightarrow$\iti, \itiii\:
Let $1$ be a left identity.
By section \ref{SS:loc skew sym} skew-symmetry for an even operator $T$ 
implies that $1$ is a right identity and $Ta=a_{-2}1$.
In particular, $T1=1_{-2}1=0$.

From $(x+w)^t a(x+w)b(w)1=(w+x)^t (a(x)b)(w)1$ for some $t\geq 0$ we get
$$
(x+w)^t a(x+w)e^{wT}b
\;=\;
(w+x)^t e^{wT}a(x)b.
$$
Because $a(x+w)e^{wT}b$ and $e^{wT}a(x)b$ are power series in $w$,
we can divide by $(x+w)^t$ and get $a(x+w)e^{wT}b=e^{wT}a(x)b$.
By Remark \ref{SS:s3 symm} this implies that $T$ is a translation generator.

By Proposition \ref{SS:s3 symm} the field identities hold.
Since $(Ta)(z)=(a_{-2}1)(z)=a(z)_{-2}\id_V=\del_z a(z)$,
we see that $V$ is a unital vertex algebra.
It is associative by Proposition \ref{SS:assoc va}.

\smallskip

\itiii$\Rightarrow$\itiv\:
Let $1$ be a weak right identity and define $Ta:=a_{-2}1$.
The Jacobi identity for $1, 1, 1$ and indices $-1, -1, -1$ is $T1=T1+T1$. 
Thus $T1=0$.
The coefficient of $w$ of the associativity formula
$$
(a(x)b)(w)1
\; =\;
\res_z\de(z-w,x)[a(z),b(w)]1
\; =\;
a(x+w)b(w)1
$$ 
is $T(a(x)b)=\del_x a(x)b+a(x)Tb$.
Thus $T$ is a translation generator.

\smallskip

\itiv$\Rightarrow$\itii\:
Proposition \ref{SS:loc skew sym} shows that $V$ is bounded and 
satisfies skew-symmetry for the translation generator $T$.
By Proposition \ref{SS:s3 symm} duality holds.
By Lemma \ref{SS:loc skew sym} the invariant weak right identity $1$ 
satisfies $a(z)1=e^{zT}a$.
Hence $1$ is a right identity and $Ta=a_{-2}1$.
By section \ref{SS:loc skew sym} skew-symmetry implies that 
$1$ is a left identity.
\end{pf}

\bigskip

There exist vertex algebras that satisfy the Jacobi identity but
which are {\it not} associative.
For example, take any vertex algebra $V$ with $[\subla]=0$
and such that $V\subast$ is non-commutative with $(ab)c=abc=0$ for any 
$a, b, c$.

\subsection{Conformal Jacobi Identity and Wick Formula II}
\label{SS:conf j wick}

We prove that the conformal Jacobi identity and the Wick formula
are equivalent to special cases of the commutator and associativity formula 
and the Jacobi identity.
We also characterize conformal derivations in terms of $Y$.

\bigskip

{\bf Proposition.}\: {\it
Let $V$ be an unbounded vertex algebra and $M$ a $\Z$-fold $V$-module 
such that $Y_M(Ta)=\del_z Y_M(a)$.
The following sets of identities are equivalent:

\smallskip

\iti\:
the conformal Jacobi identity and the Wick formula;

\smallskip

\itii\:
the Jacobi identity for indices $r\geq 0, s\in\Z, t\geq 0$;

\smallskip

\itiii\:
the commutator formula for indices $t\geq 0, s\in\Z$;

\smallskip

\itiv\:
the associativity formula for indices $r\geq 0, s\in\Z$;

\smallskip

\itv\:
the associativity formula for indices $r\geq 0, s\geq -1$.
}

\bigskip

\begin{pf}
The Wick formula is the Jacobi identity for indices $r=0, s=-1, t\geq 0$.
Together with the conformal Jacobi identity we get the Jacobi identity
for indices $r=0, s\geq -1, t\geq 0$.
By Proposition \ref{SS:2nd recursion} this yields 
the Jacobi identity for indices $r=0, s\in\Z, t\geq 0$. 
By Proposition \ref{SS:fund recur}\,\iti\ the Jacobi identity holds
for indices $r\geq 0, s\in\Z, t\geq 0$.

The other implications are trivial or follow likewise from
Propositions \ref{SS:fund recur}\,\iti\ and \ref{SS:2nd recursion}.
\end{pf}

\bigskip

The equivalence of \itiii\ and \itiv\ was already proven in 
Corollary \ref{SS:ope sing}.

An operator $d$ of an algebra is a derivation iff $[d,a\cdot]=(da)\cdot$.
The analogous result for conformal derivations is provided by 
the following remark.
Recall that we identify $E\pau{\mu}$ with 
$z\inv E\pau{z\inv}$, see section \ref{SS:unb confa power series}.

\bigskip

{\bf Remark.}\: {\it
A conformal operator $d\submu$ of a vertex algebra $V$ 
is a conformal derivation iff
$[(d\submu)\subla a(z)]=(d\subla a)(z)$ for any $a\in V$.
}

\bigskip

\begin{pf}
We have $d\submu\in\Der\subsc(V)$ iff 
the commutator formula is satisfied for $d\submu$, any $a, b\in V$ 
and indices $t\geq 0, s\geq -1$.

On the other hand, the identity $[(d\submu)\subla a(z)]=(d\subla a)(z)$ 
is the associativity formula for $d\submu$, any $a, b\in V$ and 
indices $r\geq 0, s\in\Z$.
It is equivalent to the commutator formula for indices $t\geq 0, s\in\Z$.
This fact is analogous to the equivalence of parts \itiii\ and \itiv\ of the
Proposition and is proven in the same way using 
Proposition \ref{SS:fund recur}\,\iti.

Hence the claim follows from an analogue of 
Proposition \ref{SS:2nd recursion}\,\itii.
The proposition is still valid in this case 
since its proof only uses that $(a(x)Tb)(w)c\linebreak[0]=
(\del_w-\del_x)(a(x)b)(w)c$ and this identity is still
satisfied since $[T,d\submu]=-\mu d\submu$ and $(Ta)(w)=\del_w a(w)$.
\end{pf}

\subsection{Right Wick Formula}
\label{SS:right wick}

We prove that the right Wick formula is a special case of 
the associativity formula.

\bigskip

Let $V$ be a vertex algebra and $M$ a bounded $\Z$-fold $V$-module such that
$Y_M(Ta)=\del_z Y_M(a)$.
The {\bf right Wick formula} is
\index{right Wick formula}
$$
[ab\subla c]
\; =\;
(e^{T\del\subla}a)[b\subla c]
\; +\;
\paraab(e^{T\del\subla}b)[a\subla c]
\; +\;
\paraab\int_0\upla d\mu\, [b\submu [a_{\la-\mu}c]].
$$
The substitution formula implies that
$\int_0\upla d\mu[b\submu [a_{\la-\mu}c]]=
\int_0\upla d\mu[b_{\la-\mu}[a\submu c]]$, see section \ref{SS:bracket vlie}.

\bigskip

{\bf Proposition.}\: {\it
The right Wick formula is the associativity formula for indices $r=-1$ and
$s\geq 0$:
$$
(ab)_s c
\; =\;
\sum_{i\geq 0}
T^{(i)}(a)b_{s+i}c
\; +\;
\paraab\, 
T^{(i)}(b)a_{s+i}c
\; +\;
\paraab
\sum_{i=0}^{s-1}
b_{s-1-i}a_i c.
$$
}

\bigskip

\begin{pf}
We need to prove that the right-hand side of the right Wick formula
is equal to $\res_z e^{z\la}(a(z)_-b(z)+b(z)a(z)_+)c$.
Since $a(z)_-=(e^{zT}a)\cdot$ and $e^{z\la}e^{zT}=e^{z(\la+T)}=
e^{T\del\subla}e^{z\la}$ we get the first term of the right Wick formula:
$$
\res_z e^{z\la}a(z)_-b(z)
\, =\,
\res_z (e^{T\del\subla}a)e^{z\la}b(z)
\, =\,
(e^{T\del\subla}a)b\subla.
$$
Since $a(z)_+=a_{-\del_z}z\inv$ and $[\del_z,e^{z\la}\cdot]=\la e^{z\la}$, 
we have
$$
e^{z\la}a(z)_+
\, =\,
e^{z\la}a_{-\del_z}z\inv
\, =\,
a_{\la-\del_z}e^{z\la}z\inv
\, =\,
a_{\la-\del_z}\Big(z\inv+\int_0\upla d\mu\, e^{z\mu}\Big)
$$
using that $\int_0\upla d\mu\, e^{z\mu}=\sum\la^{(n+1)}z^n$.
The $z\inv$ yields the second term:
\begin{align}
\notag
\res_z b(z)a_{\la-\del_z}z\inv
\, =\,
\res_z b(z)e^{-\del_z\del\subla}a\subla z\inv
\, =\,
&\res_z (e^{T\del\subla}b)(z)a\subla z\inv
\\
\notag
\, =\,
&(e^{T\del\subla}b)a\subla.
\end{align}
Applying again $[\del_z,e^{z\la}\cdot]=\la e^{z\la}$, we obtain the third term:
$$
\int_0\upla d\mu\,\res_z b(z)a_{\la-\del_z}e^{z\mu}c
=
\int_0\upla d\mu\,\res_z b(z)e^{z\mu}a_{\la-\del_z-\mu}c
=
\int_0\upla d\mu\, b\submu a_{\la-\mu}c.
$$
\end{pf}

\bigskip

One can show directly that $\int_0\upla d\mu\, b\submu a_{\la-\mu}c=
\sum_{t,s}k_{t,s}\,b_t a_s c\,\la^{(t+s+1)}$
where $k_{t,s}:=\sum_i (-1)^i\binom{t+i}{i}\binom{t+s+1}{s-i}$.
Thus the non-trivial part of the Proposition amounts to the binomial identity
$k_{t,s}=1$.
We give a direct proof of this identity in 
Proposition \ref{S:bin coeff}\,\itii.

The difference between the ``correction term'' of the left and
right Wick formula is the binomial coefficient and the bracketing
of the $t$-th products.

\subsection{Decomposition of the Associativity Formula}
\label{SS:q assoc assoc}

We prove that $M$ is a module iff the multiplication $a\cdot$ and the
$\la$-action $a\subla$ on $M$ satisfy four identities.

\bigskip

{\bf Remark.}\: {\it
Let $V$ be a vertex algebra that satisfies skew-symmetry.
Then the left and right Wick formulas are equivalent to each other.
}

\bigskip

\begin{pf}
The Wick formula for $a, b, c$ is the Jacobi identity for $a, b, c$
and indices $r=0, s=-1, t\geq 0$. 
By $\bbS_3$-symmetry it is equivalent to the Jacobi identity for $b, a, c$ 
and indices $r=0, s\geq 0, t=-1$ and hence to the Jacobi identity for 
$b, c, a$ and indices $r=-1, s\geq 0, t=0$. 
This is the right Wick formula for $b, c, a$ by 
Proposition \ref{SS:right wick}.
\end{pf}

\bigskip

In Remark \ref{SS:left right wick}
we give a direct proof of the Remark.

Let $V$ be a vertex algebra.
In section \ref{SS:modules} we showed that 
to give a $\Z$-fold $V$-module $M$ satisfying $Y_M(Ta)=\del_z Y_M(a)$
is equivalent to giving an even linear map $V\to\End(M), a\mapsto a\cdot$, and 
a $\K[T]$-module morphism $V\to\End(M)\pau{\la}, a\mapsto a\subla$.
Of course, $M$ is bounded iff $a\subla\in\Endv(M)_+$ for any $a$, 
see section \ref{SS:fields}.

A $\Z$-fold module $M$ is a module iff $M$ satisfies 
the associativity formula, see section \ref{SS:z fold fields}.
The next result shows that this is equivalent to four identities for
$a\cdot$ and $a\subla$.

\bigskip

{\bf Proposition.}\: {\it
Let $V$ be a vertex algebra and $M$ a bounded $\Z$-fold $V$-module such that
$Y_M(Ta)=\del_z Y_M(a)$.
Then $M$ is a $V$-module iff $M$ satisfies the conformal Jacobi identity, 
the left and right Wick formula, and quasi-associativity.
}

\bigskip

\begin{pf}
The associativity formula is equivalent to the associativity formula for 
indices $r, s\geq -1$ by Proposition \ref{SS:2nd recursion}.
For indices $r\geq 0, s\geq -1$ it is equivalent to 
the conformal Jacobi identity and the Wick formula by 
Proposition \ref{SS:2nd recursion}.
For indices $r=-1, s\geq 0$ it is the right Wick formula 
by Proposition \ref{SS:right wick}.
For indices $r=s=-1$ it is quasi-associativity.
\end{pf}

\subsection{Associative Algebras with a Bracket}
\label{SS:assoc alg bracket}

We prove two results about algebras with a bracket.
In sections \ref{SS:diff alg leibniz la brack} and \ref{SS:eq commut}
we prove analogous results for vertex algebras.

\bigskip

Let $V$ be an algebra together with a bracket.
The {\bf left Leibniz identity} is 
\index{left!Leibniz identity}
$[a,bc]=[a,b]c+\paraab\, b[a,c]$.
In other words, $[a,\;]$ is a derivation.
The {\bf right Leibniz identity} is 
\index{right!Leibniz identity}
$$
[ab,c]
\; =\;
a[b,c]
\; +\;
\paraab\,b[a,c]
\; -\;
\paraab\,[b,[a,c]].
$$
The left and right Leibniz identity are equivalent if $[\,,]=[\,,]\subast$.
This fact is the analogue of Remark \ref{SS:q assoc assoc}.

\bigskip

{\bf Proposition.}\: {\it
Let $V$ be a unital algebra together with a bracket such that 
the left and right Leibniz identity hold and 
there exist $a, b$ such that $[a,b]=1$.

\smallskip

\iti\:
If $[\, ,]$ satisfies the Leibniz identity then $[\, ,]$ is skew-symmetric.

\smallskip

\itii\:
If $V$ is associative then $[\,,]=[\,,]\subast$.
}

\bigskip

\begin{pf}
\iti\:
Take any $a, a', b, c\in V$.
We evaluate $[[aa', b],c]$ in two ways.
One way is first to apply the Leibniz identity for $[\, ,]$ and then
to apply the left and right Leibniz identity.
The other way is to do the same thing the other way round. 
The result is the identity
$$
([a,b]+[b,a])\,[a',c]
\; =\;
-(a\leftrightarrow a').
$$

Suppose that $[a',c]=1$. 
Taking $a=a'$ we get that $[a',b]=-[b,a']$ for any $b\in V$. 
Thus the right-hand side of the above identity vanishes 
for any $a, b\in V$ and $a', c$.
This shows that $[\,,]$ is skew-symmetric.

\smallskip

\itii\:
We evaluate $[aa', bc]$ in two ways.
One way is first to apply the left Leibniz identity
and then the right Leibniz identity.
The other way is to do the same thing the other way round. 
The result is the identity
$$
([a,b]\subast-[a,b])\,[a',c]
\; =\;
-(a\leftrightarrow a').
$$
The claim is now proven in the same way as \iti.
\end{pf}

\subsection{Wick and Vertex Leibniz imply Vertex Lie if $c\ne 0$}
\label{SS:diff alg leibniz la brack}

We prove that the Wick formulas and the conformal Jacobi identity imply 
conformal skew-symmetry if there exists a Virasoro vector with $c\ne 0$,
see section \ref{SS:conf va}.

\bigskip

{\bf Proposition.}\: {\it
Let $V$ be a unital vertex algebra such that $a_t b=1$ for some $a, b\in V$ 
and $t\geq 0$.
Then the left and right Wick formula and the conformal Jacobi identity 
imply conformal skew-symmetry.
}

\bigskip

\begin{pf}
Take any $a, a', b, c\in V$.
We evaluate $((aa')\subla b)_{\mu+\la}c$ in two ways.
One way is first to apply the conformal Jacobi identity and then
the Wick formulas.
The other way is to do the same thing the other way round. 
We will show that the result is the identity
$$
(b\submu a+a_{-\mu-T}b)\, a'\subla c
\; =\;
-(a\leftrightarrow a').
$$
The claim is now proven in the same way as Proposition 
\ref{SS:assoc alg bracket}.

On the left-hand side of the conformal Jacobi identity we have
\begin{align}
\notag
&((aa')\subla b)_{\mu+\la}c
\\
\notag
=
&\Big( 
(e^{T\del\subla}a)a'\subla b +
(e^{T\del\subla}a')a\subla b +
\int_0\upla d\nu\, a'\subnu a_{\la-\nu}b\Big)_{\la'} c|_{\la'=\mu+\la}
\\
\notag
=
&(e^{T\del\subla}a)(a'\subla b)_{\mu+\la}c+
(e^{T\del\subla}(a'_{-\mu-T}b))a_{\mu+\la}c
\\
\notag
&+
\int_0^{\mu+\la}d\nu\, 
(a'_{\nu-\mu}b)\subnu a_{\mu+\la-\nu}c+
(a\leftrightarrow a')+
\int_0\upla d\nu\, (a'\subnu a_{\la-\nu}b)_{\mu+\la} c
\end{align}
because
$(e^{T\del_{\la'}}e^{T\del\subla}a)(a'\subla b)_{\la'} c|_{\la'=\mu+\la}
=(e^{T\del\subla}a)(a'\subla b)_{\mu+\la} c$,
\begin{align}
\notag 
&(e^{T\del_{\la'}}(a'\subla b))(e^{T\del\subla}a)_{\la'}c|_{\la'=\mu+\la}
=(e^{T\del_{\la'}}(a'\subla b))e^{-\la'\del\subla}a_{\la'}c|_{\la'=\mu+\la}
\\
\notag
=&(e^{T\del_{\la'}}e^{-(\la''+T)\del\subla}(a'\subla b))a_{\la'}c
|_{\la'=\la''=\mu+\la}=
(e^{T\del\subla}(a'_{-\mu-T}b))a_{\mu+\la}c,
\end{align}
and 
$$
\int_0^{\la'}d\nu\, (a'\subla b)\subnu (e^{T\del\subla}a)_{\la'-\nu}c
\; =\;
\int_0^{\la'}d\nu\, (a'_{\la-\la'+\nu}b)\subnu a_{\la'-\nu}c.
$$

The two terms on the right-hand side of the conformal Jacobi identity are
$$
(aa')\subla b\submu c
=
(e^{T\del\subla}a)a'\subla b\submu c +
(e^{T\del\subla}a')a\subla b\submu c +
\int_0\upla d\nu\, a'\subnu a_{\la-\nu}b\submu c
$$
and
\begin{align}
\notag
b\submu (aa')\subla c 
=
&b\submu (e^{T\del\subla}a)a'\subla c +
b\submu (e^{T\del\subla}a')a\subla c +
\int_0\upla d\nu\, b\submu a'\subnu a_{\la-\nu}c
\\
\notag
=
&(e^{T\del\subla}(b\submu a))a'_{\la+\mu}c +
(e^{T\del\subla}a)b\submu a'\subla c +
\int_0\upmu d\nu\, (b\submu a)\subnu a'_{\la+\mu-\nu}c 
\\
\notag
&\qquad\qquad\qquad\qquad\qquad\qquad
+(a\leftrightarrow a')+\int_0\upla d\nu\, b\submu a'\subnu a_{\la-\nu}c
\end{align}
because $(b\submu e^{T\del\subla}a)a'\subla c
=e^{\mu\del\subla-\mu\del\subla}(b\submu e^{T\del\subla}a)a'\subla c
=((e^{T\del\subla}b)\submu (e^{T\del\subla}a))a'_{\la+\mu}c
=(e^{T\del\subla}(b\submu a))a'_{\la+\mu}c$
and 
$(b\submu e^{T\del\subla}a)\subnu a'\subla c=
(e^{T\del\subla}(b\submu a))\subnu a'_{\la+\mu}c=
e^{-\nu\del\subla}(b\submu a)\subnu\linebreak[0] a'_{\la+\mu}c=
(b\submu a)\subnu a'_{\la+\mu-\nu}c$.

\smallskip

Denote by $X_i, Y_i, Z_i$ the $i$-th summand of the first, second,
and third term of the conformal Jacobi identity. 
Thus $\sum_{i=1}^7 X_i=\sum_{i=1}^3 Y_i-\sum_{i=1}^7 Z_i$.
The conformal Jacobi identity for $a, b, c$ and $a', b, c$ shows that
the summands $X_4, Y_2, Z_5$ and $X_1, Y_1, Z_2$ cancel.
Applying the conformal Jacobi identity to 
$X_3, X_6, X_7, Y_3, Z_3, Z_6, Z_7$ we obtain $14$ summands of the
form $\int_0^{\de}f\subal g\subbe h_{\ga} c$ where
$\set{f, g, h}=\set{a, a', b}$ and $\de\in\set{\la, \mu, \la+\mu}$.
The substitution formula from section \ref{SS:bracket vlie}
implies that these $14$ terms cancel.
The remaining terms are $X_2, X_5, Z_1, Z_4$. 
We thus get
$$
(e^{T\del\subla}(b\submu a)\: +\: e^{T\del\subla}(a_{-\mu-T}b))a'_{\la+\mu}c
\; =\;
(a\leftrightarrow a').
$$
Applying $e^{-T\del\subla}$ to this equation and replacing $\la$ by $\la+T-\mu$
we obtain the identity at the beginning.
\end{pf}

\subsection{Associativity Formula implies Associativity if $c\ne 0$}
\label{SS:eq commut}

We prove that the associativity formula implies associativity if 
there exists a Virasoro vector with $c\ne 0$, see section \ref{SS:conf va}.

\bigskip

{\bf Proposition.}\: {\it
Let $V$ be a unital vertex algebra such that $a_t b=1$ for some $a, b\in V$ 
and $t\geq 0$.
Then the associativity formula implies that $V$ is associative.
}

\bigskip

\begin{pf}
By Proposition \ref{SS:q assoc assoc}
the conformal Jacobi identity, the left and right Wick formula,
and quasi-associativity are satisfied.
Because of Propositions \ref{SS:skew sym} and 
\ref{SS:diff alg leibniz la brack} 
it suffices to show that $[\,,]\subast=[\,,]\subslie$.

Take any $a, a', b, c\in V$.
We evaluate $(aa')\subla(bc)$ in two ways.
One way is first to apply the Wick formula and then the right Wick formula.
The other way is to do the same thing the other way round. 
We will show that the result is the identity
$$
([a,b]\subast-[a,b]\subslie)\, a'\subla c
\; =\;
-(a\leftrightarrow a').
$$
The claim is now proven in the same way as Proposition 
\ref{SS:assoc alg bracket}.

By first applying the Wick formula we get
\begin{align}
\notag
&(aa')\subla(bc)
\\
\notag
=\;
&
((aa')\subla b)c
\; +\;
b((aa')\subla c)
\; +\;
\int_0\upla d\mu\, ((aa')\subla b)\submu c
\\
\notag
=\;
&
\Big((e^{T\del\subla}a)a'\subla b
\; +\;
(e^{T\del\subla}a')a\subla b
\; +\;
\int_0\upla d\mu\, a'\submu a_{\la-\mu}b
\Big)
c
\\
\notag
&+\;
b\Big((e^{T\del\subla}a)a'\subla c
\; +\;
(e^{T\del\subla}a')a\subla c
\; +\;
\int_0\upla d\mu\, a'\submu a_{\la-\mu}c
\Big)
\int_0\upla d\mu
\\
\notag
&+\;
\int_0\upla d\mu
\Big((e^{T\del_{\la'}}a)a'_{\la'} b
\; +\;
(e^{T\del_{\la'}}a')a_{\la'} b
\; +\;
\int_0\upla d\nu\, a'\subnu a_{\la-\nu}b
\Big)\submu c|_{\la'=\la}
\\
\notag
=\;&
((e^{T\del\subla}a)a'\subla b)c
\; +\;
((e^{T\del\subla}a')a\subla b)c
\; +\;
\int_0\upla d\mu\, (a'\submu a_{\la-\mu}b)c
\\
\notag
&+\;
b(e^{T\del\subla}a)a'\subla c
\; +\;
b(e^{T\del\subla}a')a\subla c
\; +\;
\int_0\upla d\mu\, b\, a'\submu a_{\la-\mu}c
\\
\notag
&+\;
\int_0\upla d\mu
\Big(
(e^{T(\del_{\la'}+\del\submu)}a)(a'_{\la'} b)\submu c
\; +\;
(e^{T\del\submu}(a'_{\la'} b))(e^{T\del_{\la'}}a)\submu c
\\
\notag
&\qquad\qquad\qquad\qquad\qquad\qquad\qquad\quad+\;
\int_0\upmu d\nu\, (a'_{\la'} b)\subnu(e^{T\del_{\la'}}a)_{\mu-\nu}c
\Big)|_{\la'=\la}
\\
\notag
&+\;
(a\leftrightarrow a')
\; +\;
\int_0\upla d\mu\int_0\upla d\nu\, (a'\subnu a_{\la-\nu}b)\submu c.
\end{align}
Denote the thirteen terms on the right-hand side by
$A_r, A_r', D_r, G, G', F_r, H, \linebreak[0]I, \de, H',  I', \de', \ep$.

By first applying the right Wick formula we get
\begin{align}
\notag
&(aa')\subla(bc)
\\
\notag
=
&(e^{T\del\subla}a)a'\subla (bc)
\; +\;
(e^{T\del\subla}a')a\subla (bc)
\; +\;
\int_0\upla d\mu\, a'\submu a_{\la-\mu}(bc)
\\
\notag
=
&(e^{T\del\subla}a)
\Big((a'\subla b)c
\; +\;
b(a'\subla c)
\; +\;
\int_0\upla d\mu\, (a'\subla b)\submu c
\Big)
\; +\;
(a\leftrightarrow a')
\\
\notag
&+
\int_0\upla d\mu\, a'\submu 
\Big(
(a_{\la-\mu} b)c
\; +\;
b(a_{\la-\mu} c)
\; +\;
\int_0^{\la-\mu} d\nu\, (a_{\la-\mu} b)\subnu c
\Big)
\\
\notag
=
&(e^{T\del\subla}a)
(a'\subla b)c
\; +\;
(e^{T\del\subla}a)
b(a'\subla c)
\; +\;
(e^{T\del\subla}a)
\int_0\upla d\mu\, (a'\subla b)\submu c
+
(a\leftrightarrow a')
\\
\notag
&
+
\int_0\upla d\mu
\Big(
(a'\submu a_{\la-\mu} b)c
\; +\;
a_{\la-\mu}b(a'\submu  c)
\; +\;
\int_0\upmu d\nu\, (a'\submu a_{\la-\mu}b)\subnu c
\Big)
\\
\notag
&+
\int_0\upla d\mu
\Big(
(a'\submu b)a_{\la-\mu} c
\; +\;
b\, a'\submu a_{\la-\mu}c
\; +\;
\int_0\upmu d\nu\, (a'\submu b)\subnu a_{\la-\mu}c
\Big)
\\
\notag
&+
\int_0\upla d\mu\int_0^{\la-\mu} d\nu\, a'\submu 
(a_{\la-\mu} b)\subnu c.
\end{align}
Denote the thirteen terms on the right-hand side by
$A_l, B, C, A_l', B', C', D_l, \linebreak[0] E, \al, E', F_l, \be, \ga$.

We claim that by taking the difference of these two evaluations
of $(aa')\subla(bc)$ we obtain
$$
((e^{T\del\subla}a)b)a'\subla c
\; -\;
(b(e^{T\del\subla}a))a'\subla c
\; -\; 
\bigg(\int_{-T}^0 d\mu\, (e^{T\del\subla}a)\submu b\bigg) a'\subla c
\; +\;
(a\leftrightarrow a').
$$
Denote these six terms by $B_f, G_f, A_f, B_f', G_f', A_f'$.
From vanishing of this last expression we get the first identity of this proof
if we replace $a, a'$ by $e^{-T\del\subla}a$ and $e^{-T\del\subla}a'$.
In other words, we replace $a, a'$ by $T^{(i)}a, T^{(i)}a'$,
apply $(-\del\subla)^i$ to the resulting equation, and sum over $i\geq 0$.

The pre-Lie algebra identity shows that the terms $B, G$ yield $B_f, G_f$.
Likewise, we get $B_f', G_f'$ from $B', G'$.
The terms $D_l, D_r$ cancel as do $F_l, F_r$. 
Applying quasi-associativity to $A_l, A_r$ we get
$$
-\bigg(\int_0^T d\mu\, (a'\subla b)\bigg) (e^{T\del\subla}a)\submu c
\; -\;
\bigg(\int_0^T d\mu\, e^{T\del\subla}a\bigg) (a'\subla b)\submu c.
$$
Denote these two terms by $A_1, A_2$. 
In the same way we obtain two terms $A_1', A_2'$ from $A_l', A_r'$.
We have $-A_1=E'-A_f'-I$ because
\begin{align}
\notag
&\bigg(\int_0^T d\mu (a'_{\la-\mu}b)\bigg)\, a\submu c
\; =\;
-\bigg(\int\subla^{\la-T} d\mu (a'\submu b)\bigg)\, a_{\la-\mu} c
\\
\notag
=\; 
&\bigg(\int_0\upla d\mu (a'\submu b)\bigg)\, a_{\la-\mu} c
\; +\;
\bigg(\int_{-T}^0 d\mu (a'\submu b)\bigg)\, a_{\la-\mu} c
\\
\notag
&\qquad\qquad\qquad\qquad\qquad\qquad\qquad\qquad\qquad
-\;
\bigg(\int_{-T}^{\la-T} d\mu (a'\submu b)\bigg)\, a_{\la-\mu} c
\end{align}
and
\begin{align}
\notag
\bigg(\int_{-T}^{\la-T} d\mu (a'\submu b)\bigg)\, a_{\la-\mu} c
\; =\;
&\int_0\upla d\mu (e^{-T\del\submu}(a'\submu b))\, a_{\la-\mu} c
\\
\notag
=\;
&\int_0\upla d\mu (e^{T\del\submu}(a'_{\la-\mu}b))\, a\submu c
\end{align}
In the same way $A_1'$ together with $E$ and $I'$ yields $A_f$.

The terms $A_2, C, H$ cancel each other because
\begin{align}
\notag
-A_2
= 
\bigg(\int_0^T d\mu\, e^{T\del\subla}&a\bigg) (a'\subla b)\submu c
=
\sum_{n,m,i}(T^{(m+1)}T^{(i)}a)(a'_n b)_m c\, \la^{(n-i)}
\\
\notag
=\
&\sum_{n,m,i}\binom{m+i+1}{i}
T^{(m+i+1)}(a)(a'_{n+i}b)_m c\, \la^{(n)},
\end{align}
\begin{align}
\notag
C=(e^{T\del\subla}&a)\int_0\upla d\mu\, (a'\subla b)\submu c
=
\sum_{i,n,m} \binom{n+m+1}{n}T^{(i)}(a)(a'_n b)_m c\,\la^{(n+m+1-i)}
\\
\notag
=
&\sum_{n,m;\, i\geq -1-m}\binom{n+m+i+1}{n+i}
T^{(m+i+1)}(a)(a'_{n+i}b)_m c\,\la^{(n)},
\end{align}
and $H$ is equal to 
\begin{align}
\notag
&\int_0\upla d\mu
(e^{T(\del_{\la'}+\del\submu)}a)(a'_{\la'} b)\submu c|_{\la'=\la}
\\
\notag
=
&\sum_{i,n,m;j\leq m}
(T^{(i)}T^{(j)}a)(a'_n b)_m c\, \la^{(n-i)}\la^{(m-j+1)}
\\
\notag
=
&\sum_{i,n,m;j\leq m}
\binom{i+j}{j}\binom{n+m-j+1}{m-j+1}
(T^{(i+j)}a)(a'_{n+i}b)_m c\, \la^{(n+m-j+1)}
\\
\notag
=
&\sum_{n>0;\, j\leq m;\, i\geq -1-m}
\binom{i+m+1}{j}\binom{n}{m-j+1}
T^{(i+m+1)}(a)(a'_{n+i}b)_m c\, \la^{(n)}
\\
\notag
=
&\sum_{n>0;\, j\leq m;\, i\geq -1-m}
\bigg( \binom{n+m+i+1}{m+1}-\binom{m+i+1}{m+1}\bigg)
\\
\notag
&\qquad\qquad\qquad\qquad\qquad\qquad\qquad\qquad\qquad
T^{(m+i+1)}(a)(a'_{n+i}b)_m c\, \la^{(n)}
\end{align}
where in the third equation we replaced $n$ by $n-m+j-1$ and 
$i$ by $i+m-j+1$.
In fact, the summands of $A_2$ and $C$ for $n=0$ cancel,
the summands of $A_2$ for $n>0$ cancel with the second term of $H$, and
the summands of $C$ for $n>0$ cancel with the first term of $H$.
In the same way $A_2', C', H'$ cancel each other.

Applying the conformal Jacobi identity to 
$\al, \be, \ga, \de, \de', \ep$ we obtain $16$ summands of the
form $\int\int f\subal g\subbe h_{\ga} c$ where 
$(f, g, h)$ is a permutation of $(a, a', b)$.
These $16$ terms cancel. 
This can be checked using the substitution formula and the Fubini formula from 
section \ref{SS:bracket vlie}.
In particular, these two formulas imply
$$
\int_0\upla d\mu\int_0\upmu d\nu\: p(\mu,\nu)
\; =\;
\int_0\upla d\mu\int_0\upmu d\nu\: p(\la+\nu-\mu,\nu).
$$
\end{pf}

\section{Associative Unital Vertex Algebras}
\label{S:unital va}

In sections \ref{SS:commut centre}--\ref{SS:indecomp va}
we discuss commuting and idempotent elements and the block decomposition. 
In sections \ref{SS:tensor va} and \ref{SS:tensor va2}
we discuss the tensor product and affinization.
Section \ref{SS:conf va}
is about conformal vertex algebras.

\medskip

{\bf Convention.}\:
In this section all algebras and vertex algebras are assumed to be unital.

\subsection{Commuting Elements and the Centre}
\label{SS:commut centre}

We give a criterion for two elements to commute.

\bigskip

By definition, elements $a, b$ of an unbounded vertex algebra commute if
$ab=\paraab ba$ and $[a\subla b]=[b\subla a]=0$, 
see section \ref{SS:comm va}.

\bigskip

{\bf Remark.}\: {\it
Let $E$ be a vector space and $a(z), b(z)\in\Endv(E)$.
The fields $a(z), b(z)$ commute iff $[a(z),b(w)]=0$.
}

\bigskip

\begin{pf}
If $[a(z),b(w)]=0$ then $[a(z)\subla b(z)]=[b(z)\subla a(z)]=0$ and
$\normord{a(z)b(z)}=a(z)b(z)=b(z)a(z)=\normord{b(z)a(z)}$.

Conversely, if $a(z), b(z)$ commute then they satisfy 
conformal skew-sym\-metry and $[a(z),b(z)]\subslie=0=[a(z),b(z)]\subast$.
By Proposition \ref{SS:skew sym} they satisfy skew-symmetry.
By Proposition \ref{SS:skew sym fields} they are local.
Since $[a(z)\subla b(z)]=0$,
the weak commutator formula implies that $[a(z),b(w)]=0$.
\end{pf}

\bigskip

Let $V$ be an associative unital algebra, $C$ its centre, and 
$\End_V(V)$ the endomorphism algebra of the left $V$-module $V$.
Then $a\mapsto a\cdot$ is an algebra isomorphism $C\to\End_V(V)$.
Part \itii\ of the following result is the vertex analogue of this fact.

\bigskip

{\bf Proposition.}\: {\it
Let $V$ be an associative vertex algebra.

\smallskip

\iti\:
Elements $a, b\in V$ commute iff $[a(z),b(w)]=0$.

\smallskip

\itii\:
The centre of $V$ is isomorphic to $\End_V(V)$ via $a\mapsto a\cdot$.
}

\bigskip

\begin{pf}
\iti\:
Since $Y: V\to\Endv(V)$ is a monomorphism,
the claim follows from the Remark.

\smallskip

\itii\:
The map $C\to\End(V), a\mapsto a\cdot$, is an algebra morphism because 
$(ab)c=abc$ for $a, b\in C$ by Remark \ref{SS:prelie}\,\itiii.
Its image is in $\End_V(V)$ since $a(b_t c)=b_t (ac)$ for $a\in C$ by \iti.
It is injective because $a1=a$.
It is surjective because if $\phi\in\End_V(V)$ then $\phi 1\in C$ and 
$\phi=(\phi 1)\cdot$ because $[a\subla\phi 1]=\phi[a\subla 1]=0$ and
$\phi a=\phi(a1)=a\phi 1=(\phi 1)a$ for any $a\in V$.
\end{pf}

\subsection{Product Vertex Algebras and Central Idempotents}
\label{SS:prod of va idemp}

We explain the correspondence between product vertex algebras
and central idempotents.

\bigskip

The {\bf product} $\prod V_i$ of a family of vertex algebras $V_i$ is 
\index{product!of vertex algebras}
an unbounded vertex algebra with $(a^i)_t (b^i):=(a^i{}_t b^i)$ for $t\in\Z$.
If the family is finite then $\prod V_i$ is a vertex algebra and 
if any $V_i$ is associative then so is $\prod V_i$.

Let $V$ be an algebra.
An element $e\in V$ is 
\index{idempotent!}
an {\bf idempotent} if $e^2=e$.
The {\bf trivial} idempotents 
are $0$ and $1$.
Two idempotents $e, f$ are {\bf orthogonal} 
\index{idempotent!orthogonal}
if $ef=fe=0$.
Idempotents $e_1, \dots, e_r$ are {\bf complementary} if 
\index{idempotent!complementary}
they are pairwise orthogonal and $\sum e_i=1$.

If $e$ is an idempotent then $e, 1-e$ are complementary idempotents.
If $e, e'$ are orthogonal idempotents then $f:=e+e'$ is an idempotent and
$e=fe\in fV$.

The above notions for idempotents are also defined for a vertex algebra $V$
by applying them to the algebra $V\subast$.

We mention that if $V$ is associative and $e, f$ are complementary idempotents
then we have the {\bf Peirce decomposition}
\index{Peirce decomposition}
$V=eVe\oplus eVf\oplus fVe \oplus fVf$ where $eVe, fVf$ are subalgebras and
$eVf, fVe$ are bimodules.
In the following we only consider {\bf central} idempotents, that is, 
idempotents 
\index{idempotent!central}
that lie in the centre. 
In this case the Peirce decomposition reduces to a product decomposition,
as we explain next.

Let $V$ be a (vertex) algebra and $V_i\subset V$ subsets.
We write $V=V_1\times\ldots\times V_r$ if $V_i$ are ideals and 
the canonical linear map $V_1\times\ldots\times V_r\to V$ is 
a (vertex) algebra isomorphism.
In this case if $1=\sum e_i$ with $e_i\in V_i$ then 
$e_i$ are complementary central idempotents, 
$e_i$ is an identity of $V_i$, and $V_i=e_i V$.

Conversely, if $V$ is an associative algebra and $e_i$ are complementary 
central idempotents then $V=e_1 V\times\ldots\times e_r V$.
For vertex algebras, we have:

\bigskip

{\bf Proposition.}\: {\it
Let $V$ be an associative vertex algebra and $e_i$ complementary central 
idempotents of $V$.
Then $V=e_1 V\times\ldots\times e_r V$.
}

\bigskip

\begin{pf}
It suffices to consider the case $r=2$.
Let $e:=e_1$ and $f:=e_2$.
The subspace $eV$ is a $\K[T]$-submodule because $Te=T(e^2)=2eTe$ and hence
$T(ae)=(Ta)e+2(aTe)e\in eV$ by Remark \ref{SS:prelie}\,\itiii.
Thus $eV$ is an ideal since $a_n (eb)=e(a_n b)$ for any $a, b\in V, n\in\Z$
by Proposition \ref{SS:commut centre}\,\iti.

We have $V=eV+fV$ since $a=a(e+f)=ae+af$.
We have $eV\cap fV=0$ since $ea=fb$ implies $ea=e^2 a=eea=efb=(ef)b=0$ 
because of Remark \ref{SS:prelie}\,\itiii.
We have $a_n b=0$ for $a\in eV, b\in fV, n\in\Z$ because $eV$ and $fV$
are ideals and $eV\cap fV=0$.
Hence $V=eV\times fV$.
\end{pf}

\subsection{Block Decomposition}
\label{SS:indecomp va}

We prove that if the centre of a vertex algebra $V$ is noetherian
then $V$ is the unique finite product of indecomposable vertex algebras.
This is 
\index{block decomposition}
the {\bf block decomposition} of $V$.

\bigskip

Let $V$ be an associative (vertex) algebra.
A central idempotent $f\ne 0$ is {\bf centrally primitive} if there do not 
exist orthogonal central idempotents $e, e'\ne 0$ such that $f=e+e'$.
Proposition \ref{SS:prod of va idemp} shows that $fV$ is an ideal.
The element $f$ is centrally primitive iff 
all central idempotents of $fV$ are trivial.

We call $V$ {\bf indecomposable} if $V\ne 0$ and
there do not exist (vertex) subalgebras $W, W'\ne 0$ such that $V=W\times W'$.
Proposition \ref{SS:prod of va idemp} implies that $V$ is indecomposable 
iff $1$ is centrally primitive.

The following statement is a standard result from noncommutative ring theory.

\bigskip

{\bf Lemma.}\: {\it
Let $V$ be an associative algebra.

\smallskip

\iti\:
Suppose that there exist complementary, centrally primitive idempotents $e_i$.
Then $\set{\sum\de_i e_i\mid\de_i\in\set{0,1}}$ is the set of central
idempotents and $\set{e_1, \dots, e_r}$ is the set of centrally primitive 
idempotents of $V$.

\smallskip

\itii\:
If $V$ is left noetherian then there exist indecomposable subalgebras $V_i$, 
unique up to permutation, such that $V=V_1\times\ldots\times V_r$.
}

\bigskip

\begin{pf}
\iti\:
Let $e$ be a central idempotent.
For any $i$, the product $ee_i$ is a central idempotent in $e_i V$
and hence $ee_i=\de_i e_i$ for $\de_i\in\set{0,1}$.
We get $e=e\sum e_i=\sum\de_i e_i$.
Thus the two claims follow.

\smallskip

\itii\:
Assume that $V$ is not a finite product of indecomposable subalgebras.
By induction one constructs two infinite sequences of subalgebras 
$V_i, V'_i\ne 0$ such that $V=V_1\times\ldots\times V_r\times V'_r$ for any 
$r$ and $V'_i$ is not a finite product of indecomposable subalgebras.
Then $V_1\times\ldots\times V_r$ is a strictly ascending sequence of 
proper ideals. 
Contradiction.
Uniqueness follows from \iti.
\end{pf}

\bigskip

{\bf Proposition.}\: {\it
Let $V$ be an associative vertex algebra with a noetherian centre.
Then there exist indecomposable vertex subalgebras $V_i$, 
unique up to a permutation, such that $V=V_1\times\ldots\times V_r$.
}

\bigskip

\begin{pf}
The central idempotents of $V\subast$ are just the idempotents
of the centre.
Thus the claim follows from the Proposition \ref{SS:prod of va idemp} and 
the Lemma.
\end{pf}

\bigskip

From Remark \ref{SS:conf va} follows that 
the centre of a VOA is noetherian.

\subsection{Tensor Product I. Construction}
\label{SS:tensor va}

We construct the tensor product $V\otimes V'$ of associative vertex algebras.
In section \ref{SS:tensor va2} we prove 
the universal property of $V\otimes V'$ stated below.

\bigskip

Let us recall the notion of the tensor product of 
associative algebras $V$ and $V'$.
This is an associative algebra $P$ with morphisms
$\io: V\to P$ and $\io': V'\to P$ such that $\io V, \io'V'$ commute and
for any other such triple $(W,\phi,\phi')$, there is a unique
morphism $\al: P\to W$ such that $\phi=\al\io$ and $\phi'=\al\io'$.

The triple $(P,\io,\io')$ is unique up to a unique isomorphism.
The map $V\otimes V'\to P, a\otimes a'\mapsto (\io a)(\io'a')$, is 
a vector space isomorphism. 
The induced product on $V\otimes V'$ is
$(a\otimes a')(b\otimes b')=\zeta^{a'b}\,ab\otimes a'b'$.
Moreover, if $M$ is a $V$-module and $M'$ a $V'$-module then
$M\otimes M'$ is a $V\otimes V'$-module.

The tensor product of vertex algebras satisfies the same universal property:

\bigskip

{\bf Proposition.}\: {\it
For associative vertex algebras $V$ and $V'$,
there exists an associative vertex algebra $P$ and morphisms
$\io: V\to P$ and $\io': V'\to P$ such that $\io V, \io'V'$ commute and
for any other such triple $(W,\phi,\phi')$, there is a unique
morphism $\al: P\to W$ such that $\phi=\al\io$ and $\phi'=\al\io'$.
Furthermore, the map 
$V\otimes V'\to P, a\otimes a'\mapsto (\io a)(\io'a')$, is a vector space
isomorphism. 
}

\bigskip

The triple $(P,\io,\io')$ is unique up to a unique isomorphism.
It is 
\index{tensor product}
the {\bf tensor product} of $V$ and $V'$ and denoted $V\otimes V'$.
We now construct $V\otimes V'$.

For vector spaces $E$ and $E'$,
there is a natural map $E\lau{z}\otimes E'\lau{z}\to(E\otimes E')\lau{z}$,
$$
a(z)\otimes a'(z)
\;\mapsto\; 
a(z)\otimes a'(z)
\; =\;
\sum_{t, i\in\Z}\: a_i\otimes a'_{t-i-1}\, z^{-t-1}.
$$
It induces a map $\Endv(E)\otimes\Endv(E')\to\Endv(E\otimes E'),
a(z)\otimes a'(z)\mapsto a(z)\otimes a'(z)$,
given by $b\otimes b'\mapsto \zeta^{a' b}\,a(z)b\otimes a'(z)b'$.

\bigskip

{\bf Remark.}\: {\it
For associative vertex algebras $V$ and $V'$,
the vector space $V\otimes V'$
with $Y(a\otimes a'):=a(z)\otimes a'(z)$ is an associative vertex algebra. 
}

\bigskip

\begin{pf}
We use Proposition \ref{SS:unital assoc va} to prove the claim.
The operator $T=T\otimes 1+1\otimes T$ is a translation generator of 
$V\otimes V'$ since 
$\del_z(a(z)\otimes a'(z))=\del_z a(z)\otimes a'(z)+a(z)\otimes\del_z a'(z)$.
The vector $1=1\otimes 1$ is an invariant weak right identity since
$(a(z)\otimes a'(z))1=a(z)1\otimes a'(z)1\in (V\otimes V')\pau{z}$
and $a(0)1\otimes a'(0)1=a\otimes a'$.
Finally, $V\otimes V'$ satisfies locality since for $r, s\gg 0$ we have
\begin{align}
\notag
&(z-w)^{r+s}\,(a(z)\otimes a'(z))\,(b(w)\otimes b'(w))
\\
\notag
=\;
&(z-w)^r a(z)b(w)\,\otimes\, (z-w)^s a'(z)b'(w)
\\
\notag
=\;
&(z-w)^r b(w)a(z)\,\otimes\, (z-w)^s b'(w)a'(z)
\\
\notag
=\;
&(z-w)^{r+s}\,(b(w)\otimes b'(w))\,(a(z)\otimes a'(z)).
\end{align}
\end{pf}

\bigskip

The Remark shows that if $V$ and $V'$ are 
$\K$-graded associative vertex algebras, then so is $V\otimes V'$.

\subsection{Tensor Product II. Universality, Modules, Affinization}
\label{SS:tensor va2}

We prove Proposition \ref{SS:tensor va} and discuss the $V\otimes V'$-module
$M\otimes M'$ and the affinization $\hV$.

\bigskip

{\bf Lemma.}\: {\it
Let $V$ be an associative vertex algebra and $S, S'\subset V$ be two
commuting subsets.
For any $a, b\in S$ and $a', b'\in S'$, we have
$$
(aa')(z)bb'\; =\;(a(z)b)(a'(z)b').
$$
}

\bigskip

\begin{pf}
If $c, c'\in V$ commute then $[c(z),c'(w)]=0$
by Proposition \ref{SS:commut centre}\,\iti.\linebreak[5]
Moreover, $\sqbrack{S}\subset C_V(\sqbrack{S'})$ by section \ref{SS:comm va}.
Thus we get $(aa')(z)bb'=\linebreak[5]
\normord{a(z)a'(z)}bb'=a(z)a'(z)bb'=a(z)b\, a'(z)b'=
a(z)(a'(z)b')b=(a'(z)b')a(z)b\linebreak[0]=(a(z)b)a'(z)b'$.
\end{pf}

\bigskip

We now prove that the vertex algebra $V\otimes V'$ from 
Remark \ref{SS:tensor va} satisfies the universal property
from Proposition \ref{SS:tensor va}.

The maps $\io: V\to V\otimes V', a\mapsto a\otimes 1$, and 
$\io': V'\to V\otimes V', a'\mapsto 1\otimes a'$, are vertex algebra
morphisms since $a(z)\otimes 1(z)=a(z)$ and hence 
$(\io a)_n(b\otimes b')=(a_n b)\otimes b'$.
By Proposition \ref{SS:commut centre}\,\iti\ the vertex subalgebras
$\io V, \io'V'$ commute since
$(a\otimes 1)(z)(1\otimes b)(w)=a(z)\otimes b(w)=
(1\otimes b)(w)(a\otimes 1)(z)$.

Let $(W,\phi,\phi')$ be another such triple.
The Lemma shows that the linear map 
$\al: V\otimes V'\to W, a\otimes a'\mapsto (\phi a)(\phi'a')$,  
is a morphism of vertex algebras. 
In fact, this is a straightforward calculation: 
\begin{align}
\notag
&\al((a\otimes a')(z)b\otimes b')=\al(a(z)b\otimes a'(z)b')=
\phi(a(z)b)\phi'(a'(z)b')
\\
\notag
=&((\phi a)(z)\phi b)(\phi'a')(z)\phi'b'=
((\phi a)\phi'a')(z)(\phi b)\phi'b'=(\al(a\otimes a'))(z)\al(b\otimes b').
\end{align}
The map $\al$ obviously satisfies $\phi=\al\io$ and $\phi'=\al\io'$.
It is the unique such morphism because $a\otimes a'=(\io a)(\io'a')$.
This proves Proposition \ref{SS:tensor va}.

Let $M$ be a $V$-module and $M'$ a $V'$-module.
Then $N:=M\otimes M'$ is a $V$-module with $Y_N a:=a(z)\otimes 1$ and 
a $V'$-module with $Y_N a':=1\otimes a'(z)$.
The associative vertex subalgebras $Y_N V$ and $Y_N V'$ commute and 
generate an associative vertex subalgebra of $\Endv(N)$ by
Proposition \ref{SS:vas of fields} and Corollary \ref{SS:vas of fields}.
By Proposition \ref{SS:tensor va}
there is thus a morphism $Y_N: V\otimes V'\to\Endv(N)$.
In other words, $M\otimes M'$ is a $V\otimes V'$-module.
We have $Y_N(a\otimes a')=\linebreak[0]
\normord{(Y_N a)(Y_N a')}=(Y_N a)(Y_N a')=
a(z)\otimes a'(z)$.

We now consider $V\otimes V'$ in the special case where $V'$ is 
an even commutative vertex algebra $C=(C,\del)$.
We claim that $(C\otimes V)\subslie$ is equal to the 
vertex Lie algebra $C\otimes V\subslie$ defined in section \ref{SS:affin vlie}.
In fact, the $\la$-bracket of $(C\otimes V)\subslie$ is 
\begin{align}
\notag
[fa\subla gb]
\; =\;
&\res_z\, e^{z\la}(fa)(z)gb
\; =\;
\res_z\,  e^{z\la}\, (e^{z\del}f)g\; a(z)b
\\
\notag
\; =\;
&\res_z\, (e^{\del\del\subla}f)g\; e^{z\la}a(z)b
\; =\;
(e^{\del\del\subla}f)g\; [a\subla b].
\end{align}
Here we used that the ``Fourier transform'' of $z\cdot$ is $\del\subla$, 
see section \ref{SS:unb confa power series}.

The algebra $\K[x\uppm]$ is a $\Z$-graded commutative
differential algebra with $T=\del_x$ and $\K[x\uppm]_h=\K x^{-h}$. 
Viewed as a commutative vertex algebra,
we have $(p(x))(z)=(e^{z\del_x}p(x))\cdot=p(x+z)\cdot$ for $p(x)\in\K[x\uppm]$.

The {\bf affinization} of an associative vertex algebra $V$ is 
\index{affinization}
the tensor product vertex algebra $\hV:=V\otimes\K[x\uppm]$.
We have $(a_t)(z)=a(z)(x+z)^t$ where $a_t:=a\otimes x^t$ for $a\in V$ and
$t\in\Z$.
Therefore $(a_t)_r(b_s)=\sum_{i\geq 0}\binom{t}{i}(a_{r+i}b)_{t+s-i}$.
If $V$ is graded then $\hV$ is also graded with $a_t\in\hV_{h_a-t}$.
We have $(\hV)\subslie=\widehat{V\subslie}$.

\subsection{Conformal Vertex Algebras}
\label{SS:conf va}

We discuss Virasoro and conformal vectors of vertex algebras.

\bigskip

A {\bf Virasoro vector} of a vertex algebra $V$ is 
an even vector $L$ such that $L\subla L=(T+2\la)L+(c_L/2)\la^{(3)}$ 
for some $c_L\in\K$.
In other words, $L$ is a Virasoro vector of $V\subslie$ and $\hc_L\in\K$,
see section \ref{SS:vir vlie}.
The number $c_L$ is called {\bf central charge}.

A {\bf dilatation operator} of $V$ is 
\index{dilatation operator}
an even diagonalizable operator $H$ such that $[H,T]=T$ and
$H(a_t b)=(Ha)_t b+a_t(Hb)-(t+1)a_t b$ for any $t\in\Z$.
To give a gradation of $V$ is equivalent to giving a dilatation operator.

A {\bf conformal vector} of $V$ is a Virasoro vector $L$ such that 
$L_{(0)}=T$ and $H:=L_{(1)}$ is a dilatation operator.
In other words, $L$ is a conformal vector of $V\subslie$ and the gradation
defined by $L_{(1)}$ is an algebra gradation of $V\subast$.
In this case $L\in V_2$ and hence $T=L_{-1}, H=L_0$.

A {\bf conformal vertex algebra} is a vertex algebra together
with a conformal vector. 
A {\bf conformal vector} of a {\it graded} vertex algebra 
is a Virasoro vector such that $L_{(0)}=T$ and $L_{(1)}=H$.

A {\bf vertex operator algebra} or 
\index{vertex operator algebra}
{\bf VOA} is a conformal vertex algebra $V$ such that $\dim V_h<\infty$,
$V\even=\bigoplus_{h\in\Z}V_h, V\odd=\bigoplus_{h\in 1/2+\Z}V_h$, and
$V_h=0$ for $h\ll 0$.

Recall that the centre of a simple associative algebra is a field.

\bigskip

{\bf Remark.}\: {\it
Let $V$ be an associative conformal vertex algebra with centre $C$.
Then $C=\ker T$ and $C\subset V_0$.
If $V$ is simple then $C$ is a field.
If $V$ is a simple VOA and $\K=\C$ then $C=\C$.
}

\bigskip

\begin{pf}
By Proposition \ref{SS:finite confa} we have $\ker T\subset C$.
We have $C\subset\ker T$ since $T=L_{(0)}$.
We have $C\subset V_0$ since $H=L_{(1)}$.

If $V$ is simple then $V$ is a simple $V$-module since
$T=L_{-1}$ and hence any submodule is an ideal.
It follows that $\End_V(V)$ is a division algebra by Schur's lemma.
By Proposition \ref{SS:commut centre}\,\itii\ we have $C=\End_V(V)$.

If $V$ is a simple VOA and $\K=\C$
then $C\subset V_0$ is a finite field extension of $\C$.
Thus $C=\C$.
\end{pf}

\bigskip

{\bf Lemma.}\: {\it
Let $V$ be an associative vertex algebra and 
$S\subset V$ a generating subspace.
If $L_{(0)}|_S=T$ and $L_{(1)}|_S$ is diagonalizable
then $L_{(0)}=T$ and $L_{(1)}$ is a dilatation operator.
}

\bigskip

\begin{pf}
By Remark \ref{SS:conf deriv va}
the operator $L_{(0)}$ is a derivation of $V$.
Thus $\ker(T-L_{(0)})$ is a vertex subalgebra of $V$.
Since $S\subset\ker(T-L_{(0)})$ we get $L_{(0)}=T$.

The commutator formula implies
$$
[L_{(1)},a_{(t)}]b
\; =\;
(L_{(0)}a)_{(t+1)}b+(L_{(1)}a)_{(t)}b
\; =\; 
(L_{(1)}a)_{(t)}b-(t+1)a_{(t)}b.
$$
Moreover, $[L_{(1)},T]=-[T,L_{(1)}]=L_{(0)}=T$.
This shows that the span of the eigenvectors of $L_{(1)}$ is 
a vertex subalgebra.
Hence $L_{(1)}$ is diagonalizable.
The above two identities also show that $L_{(1)}$ is a dilatation operator.
\end{pf}

\bigskip

The Remark is the analogue for associative vertex algebras
of Proposition \ref{SS:vir vlie}.
As in Proposition \ref{SS:griess}, it implies:

\bigskip

{\bf Proposition.}\: {\it
Let $V$ be an associative vertex algebra of CFT-type and $S\subset V$ 
a generating subset.
If $L\in V_2$ is an even vector such that $L_{-1}|_S=T$ and $L_0|_S=H$ 
then $L$ is a conformal vector.
\hfill $\square$
}

\section{Filtrations and Generating Subspaces}
\label{S:filt span sets}

In sections \ref{SS:filtr alg}--\ref{SS:filt gen subset}
we discuss differential and invariant filtrations and 
give spanning sets for filtrations generated by a subset.

In sections \ref{SS:standard inv filt}--\ref{SS:no repeat va mod}
we use the standard invariant filtration to prove that 
generating subspaces of $V\subast$ satisfy the PBW-property
and are the complements of $C_1(V)$;
that the complements of $C_2(V)$ generate $V\subast$ without repeats;
and that the latter result generalizes to modules. 

In sections \ref{SS:canonical filt} and \ref{SS:no repeat va}
we use the canonical differential filtration
to extend the result about $C_2(V)$ to vertex algebras with negative weights.

\medskip

{\bf Convention.}\:
All algebras and vertex algebras in this section are assumed to be unital.

\subsection{Filtrations of Algebras}
\label{SS:filtr alg}

We prove a lemma about spanning sets for filtered vector spaces and 
explain the relation between filtered algebras and their associated
graded algebras.

\bigskip

An {\bf increasing filtration} of a vector space $E$ is 
\index{filtration!of an algebra}
a sequence $F$ of subspaces $F_n, n\in\Z$, such that $F_n\subset F_{n+1}$.
A {\bf decreasing} filtration satisfies $F_n\supset F_{n+1}$.
If $F$ is an increasing filtration then $(F_{-n})$ is a decreasing filtration.

Many statements about increasing filtrations, but not all, hold 
mutatis mutandis for decreasing filtrations. 
An exception is, for example, if $E$ is a $\K[T]$-module and 
$TF_n\subset F_{n+1}$.
In order to avoid repetitions, we often consider just one type of filtration.

A filtration $F$ is {\bf exhaustive} if $\bigcup F_n=E$, 
{\bf separated} if $\bigcap F_n=0$, and {\bf zero-stationary} if $F_n=0$
for some $n$.

If $F$ is an increasing filtration then the {\bf associated graded} space is 
$\rgr E=\rgr^F E:=\bigoplus\rgr_n E$ where 
\index{associated graded!space}
$\rgr_n E:=F_n/F_{n-1}$.
Let $\si_n: F_n\to F_n/F_{n-1}$ be the quotient map.

\bigskip

{\bf Lemma.}\: {\it
Let $E$ be a vector space with an exhaustive filtration $F$ and 
$S\subset E$ a subset 
such that $\si_n(S\cap F_n)\subset\rgr_n E$ is a spanning set for any $n$.
Then $S$ is a spanning set of $E$ if 

\begin{enumerate}
\item[\iti]
$F$ is zero-stationary or 

\smallskip

\item[\itii]
$E=\bigoplus E_h$ is $\Q$-graded such that $F_n$ and $S$ are graded subspaces 
and for any $h$ there exists $n$ such that $F_n\subset\bigoplus_{k>h}E_k$.
\end{enumerate}
}

\bigskip

\begin{pf}
Part \iti\ is clear.
Part \itii\ follows from \iti\ because
$F\cap E_h$ is an exhaustive and zero-stationary filtration of $E_h$
and $S\cap E_h$ is a subset such that $\si_n(S\cap E_h\cap F_n)$ is 
a spanning set of $\rgr_n E_h$. 
\end{pf}

\bigskip

A {\bf filtration} of an algebra $V$ is 
\index{filtration!of an algebra}
a vector space filtration $F$ such that $F_n F_m\subset F_{n+m}$ and 
$1\in F_0$.
Then $\rgr V$ with multiplication $(\si_n a)(\si_m b):=\si_{n+m}(ab)$ is 
a $\Z$-graded algebra.
This is 
\index{associated graded!algebra}
the {\bf associated graded} algebra.
If $V$ is associative then so is $\rgr V$.

Suppose that $V$ is associative and $F$ is increasing.
Then $\rgr V$ is commutative iff $[F_n,F_m]\subset F_{n+m-1}$ iff 
$F'_n:=F_{n+1}$ is an algebra filtration with respect to $[\, ,]$. 
In this case $\rgr^{F'}V$ is a Lie algebra.
Hence $\rgr^F V=\rgr^{F'}V$ has a commutative multiplication 
and a Lie bracket, $\set{\si_n a,\si_m b}:=\si_{n+m-1}[a,b]$.
From $[a,bc]=[a,b]c+\paraab b[a,c]$ we get
$\set{a,bc}=\set{a,b}c+\paraab b\set{a,c}$.
Thus $\rgr V$ is a Poisson algebra.

\subsection{Filtrations of Vertex Algebras}
\label{SS:filtr va}

We explain the relation between filtered vertex algebras and
their associated graded vertex algebras.

\bigskip

Let $V$ be a differential algebra.
A {\bf differential} filtration of $V$ is 
\index{differential filtration}
an algebra filtration $F$ such that $TF_n \subset F_{n+1}$.
In this case $\rgr V$ is a $\Z$-graded differential algebra with 
$T(\si_n a):=\si_{n+1} Ta$.

An {\bf invariant} filtration is 
\index{invariant filtration}
an algebra filtration $F$ such that $TF_n \subset F_n$.
In this case $\rgr V$ is a differential algebra with $T(\si_n a):=\si_n Ta$.

Let $R$ be a conformal algebra.
A {\bf differential} filtration of $R$ is 
\index{differential filtration!of a conformal algebra}
a filtration $F$ such that $TF_n\subset F_{n+1}$ and 
$(F_n)_k F_m\subset F_{n+m-k-1}$ for $k\geq 0$.
In this case $\rgr R$ is a conformal algebra with $T(\si_n a):=\si_{n+1}Ta$ and
$k$-th product $(\si_n a)_k(\si_m b):=\si_{n+m-k-1}(a_k b)$.
If $R$ is a vertex Lie algebra then so is $\rgr R$ because 
conformal skew-symmetry and the conformal Jacobi identity are homogeneous,
e.g.~we have
\begin{align}
\notag
(\si_n a)_r(\si_m b)
=
\si_{n+m-r-1}(a_r b)
=\,
&\si_{n+m-r-1}\sum (-1)^{r+1+i}T^{(i)}(b_{r+i}a)
\\
\notag
=\,
&\sum (-1)^{r+1+i}T^{(i)}((\si_m b)_{r+i}(\si_n a)).
\end{align}

An {\bf invariant} filtration is 
\index{invariant filtration!of a conformal algebra}
a filtration $F$ such that $TF_n\subset F_n$ and 
$(F_n)\subla F_m\subset F_{n+m}[\la]$.
In this case $\rgr R$ is a conformal algebra with $T(\si_n a):=\si_n Ta$ and
$\la$-bracket $(\si_n a)\subla(\si_m b):=\si_{n+m}(a\subla b)$.
If $R$ is a vertex Lie algebra then so is $\rgr R$.

Let $V$ be an associative vertex algebra.
A {\bf differential} filtration of $V$ is 
a differential filtration of $V\subast$ that is also 
a differential filtration of $V\subslie$.
{\bf Invariant} filtrations of $V$ are defined in the same way.
In both cases $\rgr V$ is an associative vertex algebra.

From $a_{-1-k}b=T^{(k)}(a)b$ and $Ta=a_{-2}1$ follows that 
a vector space filtration $F$ is a differential filtration iff 
$(F_n)_k F_m\subset F_{n+m-k-1}$ for any $k\in\Z$.
It is an invariant filtration iff
$(F_n)_k F_m\subset F_{n+m}$ for any $k\in\Z$.

Let $F$ be a decreasing differential filtration.
Then $\rgr V$ is commutative iff $(F_n)_k F_m\subset F_{n+m-k}$ for $k\geq 0$
iff $F'_n:=F_{n-1}$ is a differential filtration of $V\subslie$. 
In this case $\rgr^F V=\rgr^{F'}V$ is a vertex Lie algebra with 
$k$-th product $\set{(\si_n a)_k(\si_m b)}=\si_{n+m-k}(a_k b)$.
Since $[a\subla\;]$ satisfies the Wick formula in $V$ and the 
Wick formula is homogeneous, it follows that 
$\set{a\subla\;}$ satisfies the Wick formula in $\rgr V$.
Since $\rgr V$ is commutative,
the Wick formula for $\set{a\subla\;}$ is equivalent to 
$\set{a\subla bc}=\set{a\subla b}c+b\set{a\subla c}$.
Thus $\rgr V$ is a vertex Poisson algebra.

Let $F$ be an increasing invariant filtration.
Then $\rgr V$ is commutative iff $(F_n)\subla F_m\subset F_{n+m-1}[\la]$
iff $F'_n:=F_{n+1}$ is an invariant filtration of $V\subslie$. 
In the same way as above one shows that in this case 
$\rgr V$ is a vertex Poisson algebra with Poisson $\la$-bracket 
$\set{(\si_n a)\subla(\si_m b)}=\si_{n+m-1}[a\subla b]$.

We sometimes consider filtrations $F=(F_n)$ with $n\in\rho\Z$ where
$\rho\in\Q_>$ such that $\rho\inv\in\N$.
In this case the associated graded space $\rgr V$ is $\rho\Z$-graded
with $\rgr_n V:=V_n/V_{n-\rho}$ if the filtration is increasing.
Rescaling $F'_n:=F_{\rho n}$ yields a filtration $F'$ indexed by $\Z$.
Rescaling preserves the sets of algebra and invariant filtrations, 
but {\it not} the set of differential filtrations.

\subsection{Filtrations Generated by a Subset}
\label{SS:filt gen subset}

We give spanning sets for filtrations generated by a subset.

\bigskip

Let $V$ be an associative differential algebra.
A filtration $F$ is {\bf finer} than 
\index{finer filtration}
a filtration $F'$ if $F_n\subset F'_n$ for any $n$.
The intersection of a family of algebra filtrations 
is again an algebra filtration. 
It is the infimum of the family.
The same is true for invariant and differential filtrations.

Let $S\subset V$ be a subset together with a map $S\to\rho\Z, a\mapsto h_a$.
The increasing algebra filtration {\bf generated} by $S$ is the finest
increasing algebra filtration $F^S$ such that $a\in F^S_{h_a}$ for any 
$a\in S$. 
The filtration $F^S$ is indexed by $n\in\rho\Z$.
We have $F^S_n=\rspan\set{a_1\dots a_r\mid a_i\in S, \sum h_{a_i}\leq n}$.

The invariant and the differential filtration generated by $S$ 
are defined in the same way.
The increasing invariant filtration is given by 
$$
F^S_n
\; =\;
\rspan\big\{ T^{n_1}(a_1)\dots T^{n_r}(a_r)\mid a_i\in S, \sum h_{a_i}\leq n
\big\}.
$$
Replacing $\sum h_{a_i}$ by $\sum (h_{a_i}+n_i)$ we obtain 
a spanning set for the differential filtration.
Replacing $\leq n$ by $\geq n$ we obtain spanning sets for decreasing
filtrations. 
In all four cases the differential algebra $\rgr V$ is generated by
$\si_{h_a}a$ for $a\in S$.

Let $V$ be an associative vertex algebra and $S\subset V$ a subset
together with a map $S\to\rho\Z, a\mapsto h_a$.
The invariant filtration {\bf generated} by $S$ is 
\index{invariant filtration!generated by a subset}
the finest invariant filtration $F^S$ such that $a\in F^S_{h_a}$ for 
any $a\in S$.
The differential filtration generated by $S$ is defined in the same way.

Let $F^S$ be the invariant or the differential filtration generated by $S$.
Suppose that $S'\subset V$ is a subset that contains $S$ and $a\mapsto h_a$
is an extension to $S'$. 
Then $F^{S'}=F^S$ iff $a\in F^S_{h_a}$ for any $a\in S'\setminus S$.

\bigskip

{\bf Proposition.}\: {\it
Let $V$ be an associative vertex algebra and $S$ a subset with 
a map $S\to\rho\Z, a\mapsto h_a$.
The increasing invariant filtration generated by $S$ is 
$$
F^S_n
\; =\;
\rspan\big\{a^1_{n_1}\dots a^r_{n_r}1\:\big|\:
a^i\in S, n_i\in\Z, r\geq 0, \sum h_{a^i}\leq n\big\}.
$$ 
The elements $\si_{h_a}a, a\in S$,
generate $\rgr V$, and $\bigcup F^S_n=\sqbrack{S}$.

The same is true for the increasing differential filtration generated by $S$
if $\sum h_{a^i}$ is replaced by $\sum(h_{a^i}-n_i-1)$.
The same is also true for decreasing filtrations if $\leq n$ is replaced by
$\geq n$.
}

\bigskip

\begin{pf}
We only consider the increasing invariant filtration.
The other three cases are proven in exactly the same way.

Let $F_n$ be the above span of vectors.
From $a\in F^S_{h_a}$ and $1\in F^S_0$ we get 
$a^1_{n_1}\dots a^r_{n_r}1\in F^S_h\subset F^S_n$ if 
$h:=\sum h_{a^i}\leq n$. Hence $F_n\subset F^S_n$.
Since $a\in F_{h_a}$, it is enough to prove that 
$F$ is an invariant filtration to get $F=F^S$.

We prove by induction on $r$ that 
$(a^1_{n_1}\dots a^r_{n_r}1)_k b\linebreak[0]\in F_{h+m}$ if 
$b\in F_m$ and $k\in\Z$.
This is clear for $r=0$.
Define $a':=a^2_{n_2}\dots a^r_{n_r}1$.
By Proposition \ref{SS:kind of assoc} we have
$$
(a^1_{n_1} a')_k b
\; =\;
\sum_{i\geq N}\:\la_i\: a^1_{n_1+k-i}a'_i b
$$
for some $N, \la_i\in\Z$.
Hence by induction $(a^1_{n_1} a')_k b\in F_{h+m}$.

The elements $\si_{h_a}a, a\in S$, generate $\rgr V$
because $\rgr_n V$ is spanned by
$$
\si_n(a^1_{n_1}\dots a^r_{n_r}1)
\; =\;
(\si_{h_{a^1}}a^1)_{n_1}\dots(\si_{h_{a^r}}a^r)_{n_r}1
$$
for $a^i\in S$ and $n_i\in\Z$ such that $\sum h_{a^i}=n$.
The identity $\bigcup F^S_n=\sqbrack{S}$ 
follows from Corollary \ref{SS:kind of assoc}\,\iti.
\end{pf}

\subsection{The Standard Filtrations}
\label{SS:standard inv filt}

We prove that $\rgr V$ is commutative for the standard invariant filtration.

\bigskip

If $V$ is a $\rho\Z$-graded vector space then we endow 
any graded subset $S\subset V$ with the map $S\to\rho\Z$
associating any $a\in V_h$ its weight $h_a:=h$.

Let $V$ be a $\rho\Z$-graded associative differential algebra.
The {\bf standard differential} filtration is $F_n:=\bigoplus_{h\leq n}V_h$. 
It is equal to the increasing differential filtration generated by 
$\bigcup V_h$. 
We have $\rgr V=V$ as differential algebras.

The {\bf standard invariant} filtration $F\upssi$ is 
the finest increasing invariant filtration such that
$V_h\subset F\upssi_h$ for any $h$.
In other words, it is the increasing invariant filtration generated 
by $\bigcup V_h$.
We have $F\upssi_n=\rspan\set{T^{n_1}(a_1)\dots T^{n_r}(a_r)\mid 
a_i\in V, n_i\geq 0, \sum h_{a_i}\leq n}$.

Recall that a differential algebra $V$ is almost commutative if
$[T^n V, T^m V]\subset T^{n+m+1}V$, see section \ref{SS:assoc va def}.

\bigskip

{\bf Remark.}\: {\it
Let $V$ be a $\rho\Z$-graded almost commutative associative differential 
algebra.
Then $\rgr^{F\upssi}V$ is commutative.
}

\bigskip

\begin{pf}
The algebra $\rgr^{F\upssi}V$ is generated by $\si_{h_a}T^n a$ for 
$a\in V, n\geq 0$.
The $\si_{h_a}T^n a$ commute since
$[T^n a,T^m b]\in T^{n+m+1}(V_{h_a+h_b-1})\subset F\upssi_{h_a+h_b-\rho}$.
\end{pf}

\bigskip

Let $V$ be a $\rho\Z$-graded associative vertex algebra.
The {\bf standard differential} filtration is $F_n:=\bigoplus_{h\leq n}V_h$. 
It is equal to the increasing differential filtration generated by 
$\bigcup V_h$. 
We have $\rgr V=V$ as graded vertex algebras.

The {\bf standard invariant} filtration $F\upssi$ is 
the finest increasing invariant filtration such that 
$V_h\subset F\upssi_h$ for any $h$.
In other words, it is the increasing invariant filtration generated 
by $\bigcup V_h$.

\bigskip

{\bf Proposition.}\: {\it
Let $V$ be a $\rho\Z$-graded associative vertex algebra.
Then the vertex algebra $\rgr^{F\upssi}V$ is commutative.
}

\bigskip

\begin{pf}
The vertex algebra $\rgr V$ is generated by $\si_{h_a}a, a\in V$, 
by Proposition \ref{SS:filt gen subset}.
The $\si_{h_a}a$ commute since 
$a_k b\in V_{h_a+h_b-k-1}\subset F\upssi_{h_a+h_b-\rho}$ for $k\geq 0$.
\end{pf}

\subsection{Generators with the PBW-Property and $C_1(V)$}
\label{SS:pbw generat}

We prove that if $V$ is $\rho\N$-graded then 
a subspace $S$ generates $V\subast$ iff 
$S$ generates $V\subast$ with the PBW-property iff 
$S$ is a complement of $C_1(V)$.

\bigskip

Let $V$ be a differential algebra.
A subset $S$ generates $V$ with 
\index{PBW-property}
the {\bf PBW-property} if 
$V=\rspan\set{T^{n_1}(a_1)\dots T^{n_r}(a_r)\mid 
(a_1,n_1)\geq\ldots\geq(a_r,n_r)}$ for any total order on $S\times\N$. 
Of course, if $V$ is commutative and $S$ generates $V$ then 
$S$ generates $V$ with the PBW-property.

The following remark is proven in the same way as the proposition below.

\bigskip

{\bf Remark.}\: {\it
Let $V$ be an $\rho\N$-graded almost commutative associative differential 
algebra and $S$ a graded subset that generates $V$.
Then $S$ generates $V$ with the PBW-property.
\hfill $\square$
}

\bigskip

{\bf Proposition.}\: {\it
Let $V$ be an associative vertex algebra, 
$S$ a subset with a map $S\to\rho\N$, 
and $F^S$ the increasing invariant filtration generated by $S$.

\smallskip

\iti\:
Suppose that $S$ generates $V$ and $[a\subla b]\in F^S_{h_a+h_b-\rho}[\la]$ 
for any $a, b\in S$.
Then $S$ generates $V\subast$ with the PBW-property.

\smallskip

\itii\:
Suppose that $V$ is $\rho\N$-graded and $S$ generates $V\subast$.
Then $F^S=F\upssi$ and the assumptions of \iti\ are satisfied.
}

\bigskip

\begin{pf}
\iti\:
By Proposition \ref{SS:filt gen subset} 
the elements $\si_{h_a}a, a\in S$, generate $\rgr V$.
By assumption, they commute.
Thus $\rgr V$ is commutative and the elements $\si_{h_a}a, a\in S$, 
generate $\rgr V$ with the PBW-property.
The filtration $F^S$ is exhaustive because $S$ generates $V$.
Since $h_a\geq 0$, we have $F^S_n=0$ for $n<0$.
Thus Lemma \ref{SS:filtr alg}\iti\ implies that 
$S$ generates $V\subast$ with the PBW-property.

\smallskip

\itii\:
The space $V_h$ is the span of products of $T^{n_i}(a_i)$, with all possible
bracketings, where $a_i\in S$ and $n_i\geq 0$ such that $\sum(h_{a_i}+n_i)=h$.
Since $T^n a\in F^S_{h_a}$ and $\sum h_{a_i}\leq h$, we get $V_h\subset F^S_h$.
Hence $F^S=F\upssi$.
By Proposition \ref{SS:standard inv filt} the vertex algebra 
$\rgr V$ is commutative.
This implies $[a\subla b]\in F^S_{h_a+h_b-\rho}[\la]$.
\end{pf}

\bigskip

We will apply part \iti\ to the enveloping vertex algebra $U(R)$ of a vertex
Lie algebra $R$, see section \ref{SS:vertex envelope}.
In this case $S=R$ and $h_a=1$.

For a $\rho\N$-graded differential algebra $V$, 
define $C_1(V):=TV+V_> V_>$ where $V_>:=\bigoplus_{h>0}V_h$.

\bigskip

{\bf Lemma.}\: {\it
Let $V$ be a differential algebra of CFT-type and $S\subset V_>$ a 
graded subspace.
Then $S$ generates $V$ iff $S+C_1(V)=V_>$.
Moreover, $S$ is a minimal generating subspace iff $S\oplus C_1(V)=V_>$.
}

\bigskip

\begin{pf}
Let $V'$ be the differential subalgebra generated by $S$.
Any element of $V'$ is an iterated product of elements of $\K[T]S$.
Hence $V'_>\subset S+TV+V_>V_>=S+C_1(V)\subset V_>$.
Therefore $V=V'$ implies $S+C_1(V)=V_>$.

Suppose that $S+C_1(V)=V_>$.
We prove that if $V_k=V'_k$ for $k<h$ then $V_h=V'_h$.
We have 
$$
V_h
\; =\;
S_h+C_1(V)_h
\; =\;
S_h+(TV)_h+\bigoplus_{k+k'=h}V_k V_{k'}
\;\subset\; 
V'_h
$$
where $k, k'>0$, because $(TV)_h=T(V_{h-1})=T(V'_{h-1})\subset V'_h$
and because $V_k=V'_k$ since $k<h$.

The second claim follows from the first and from $V'_>\subset S+C_1(V)$.
\end{pf}

\subsection{Generators without Repeats and $C_2(V)$ I}
\label{SS:no repeat va si filtr}

We prove that the complements of $C_2(V)$ generate $V\subast$ without repeats.

\bigskip

Let $V$ be a differential algebra.
Recall that $C_2(V)$ is the ideal of $V$ generated by $TV$, 
see section \ref{SS:assoc va def}.
If $V$ is of CFT-type then $C_2(V)\subset C_1(V)$ since 
$C_1(V)=TV+V_> V_>$ is an ideal.

A subset $S$ generates $V$ {\bf without repeats}
\index{generating!without repeats}
if
$$
V
\; =\;
\rspan\set{T^{n_1}(a_1)\dots T^{n_r}(a_r)\mid 
a_i\in S, n_1>\ldots>n_r\geq 0, r\geq 1}.
$$
In this case $S+C_2(V)=V$.

\bigskip

{\bf Lemma.}\: {\it
Let $C$ be a commutative differential algebra that is generated 
by a subalgebra $C'$.
Then $C'$ generates $C$ without repeats. 
}

\bigskip

\begin{pf}
For any subset $S\subset\bigcup_{r\geq 1}\N^r$, define 
$$
C_S
\; :=\;
\rspan\set{T^{n_1}(a_1)\dots T^{n_r}(a_r) 
\mid a_i\in C', (n_i)\in S}.
$$
Fix $l, n\geq 1$.
Define $S:=
\set{(n_i)\mid n_1\geq\ldots\geq n_r\geq 0, r\leq l, \sum n_i=n}$ 
and $S':=\set{(n_i)\in S\mid n_1>\ldots>n_r}$.
Since $C$ is commutative and generated by $C'$, 
it suffices to prove that $C_S=C_{S'}$.

The set $S$ is finite and totally ordered: 
define $(n_1, \ldots, n_r)<(m_1, \ldots, m_s)$ if
$r<s$ or $r=s$ and there is $k$ such that $n_i=m_i$ for $i<k$ and $n_k>m_k$.

Let  $S_x:=\set{y\in S\mid y<x}$ for $x\in S$.
We prove that if $x\in S$ such that $C_{S_x}\subset C_{S'}$ 
then $C_{\set{x}}\subset C_{S'}$.
This implies that $C_S=C_{S'}$.

We may assume that $x\notin S'$.
Then there exists $j$ such that $n_j=n_{j+1}$ where $x=(n_i)$.
If $n_j=0$ then $x\in C_{S_y}$ for $y:=(\dots, n_j, n_{j+2}, \dots)<x$
since $C'$ is a subalgebra.
Suppose that $n_j>0$. 
Solving
$$
T^{(2n_j)}(a_j a_{j+1})
\; =\;
\sum
T^{(k)}(a_j)T^{(2n_j-k)}(a_{j+1})
$$
for $T^{(n_j)}(a_j)T^{(n_j)}(a_{j+1})$,
we see that $C_{\set{x}}\subset C_{S_x}\subset C_{S'}$.
\end{pf}

\bigskip

{\bf Remark.}\: {\it
Let $V$ be a $\rho\N$-graded almost commutative associative differential 
algebra and $S$ a graded subspace.
Then $S+C_2(V)=V$ iff $S$ generates $V$ without repeats. 
\hfill $\square$
}

\bigskip

This is proven in the same way as the following proposition.

\bigskip

{\bf Proposition.}\: {\it
Let $V$ be a $\rho\N$-graded associative vertex algebra and $S$ 
a graded subspace.
Then $S+C_2(V)=V$ iff $S$ generates $V\subast$ without repeats. 
}

\bigskip

\begin{pf}
The vertex algebra $\rgr^{F\upssi}V$ is commutative 
and generated by $C':=\set{\si_{h_a}a\mid a\in V}$
by Propositions \ref{SS:standard inv filt} and \ref{SS:filt gen subset}.
The subspace $C'$ is a subalgebra since 
$(\si_{h_a}a)(\si_{h_b}b)=\si_{h_a+h_b}ab$.

Recall that $C_2(V)=(TV)V$ by Proposition \ref{SS:q assoc zhu poisson}\,\iti.
If $a=(Tb)c\in C_2(V)$ then $a\in F\upssi_{h_a-1}$ and hence $\si_{h_a}a=0$.
Thus $C'=\set{\si_{h_a}a\mid a\in S}$.

Since $h_a\geq 0$, we have $F\upssi_n=0$ for $n<0$.
Therefore the claim follows from the Lemma and Lemma \ref{SS:filtr alg}\iti.
\end{pf}

\subsection{Generators without Repeats and $C_2(V)$ II. Modules}
\label{SS:no repeat va mod}

We extend Proposition \ref{SS:no repeat va si filtr} to weak modules.

\bigskip

Let $V$ be an associative vertex algebra with an invariant filtration $F$.
An {\bf invariant} filtration of a $V$-module $M$ is 
\index{invariant filtration}
a vector space filtration $F^M$ of $M$
such that $(F_n)_k F^M_m\subset F^M_{n+m}$ for any $k\in\Z$.
In this case $\rgr M$ is a module over $\rgr V$.

Assume that $V$ is $\rho\Z$-graded and endow $V$ with the 
standard invariant filtration.
The invariant filtration {\bf generated} by $c\in M$ is 
\index{filtration!generated by $c$}
the finest invariant filtration $F^M$ of $M$ such that $c\in F^M_0$.
As in Proposition \ref{SS:filt gen subset} one shows that
$F^M_n=\rspan\set{a^1_{n_1}\dots a^r_{n_r}c\mid 
a^i\in V, n_i\in\Z, r\geq 0, \sum h_{a^i}\leq n}$.

\bigskip

{\bf Proposition.}\: {\it
Let $V$ be a $C_2$-cofinite $\rho\N$-graded associative vertex algebra, 
$S$ a graded subspace such that $S+C_2(V)=V$ and $\dim S<\infty$, 
and $M$ a weak $V$-module that is generated by $c\in M$. 
Then
$$
M
\; =\;
\rspan\set{a^1_{n_1}\dots a^r_{n_r}c\mid 
a^i\in S, n_1 < \ldots < n_r < N, r\geq 0}
$$
where $N\in\Z$ is such that $a_n c=0$ for any $a\in S$ and $n\geq N$.
}

\bigskip

\begin{pf}
Endow $V$ with the standard invariant filtration and $M$ with the invariant
filtration $F^M$ generated by $c$.
Since $c$ generates $M$, $F^M$ is exhaustive.
We have $F^M_n=0$ for $n<0$ because $h_a\geq 0$ for any $a\in V$.
Therefore any spanning set of $\rgr M$ yields a spanning set of $M$ by
Lemma \ref{SS:filtr alg}\,\iti.

For any subset $I\subset\bigcup_{r\geq 0}\Z^r$, define 
$$
N_I
\; :=\;
\rspan\set{(\si_{h_{a^1}}a^1)_{n_1}\ldots (\si_{h_{a^r}}a^r)_{n_r}\si_0 c
\mid a^i\in V, (n_i)\in I}.
$$
Fix $l\geq 1$ and $n\in\Z$.
Define $I:=\set{(n_i)\mid n_1\leq\ldots\leq n_r<N, r\leq l, \sum n_i=n}$ 
and $I':=\set{(n_i)\in I\mid n_1<\ldots<n_r}$.

We have $\si_{h_a}a=0$ for $a\in C_2(V)$ as in the proof of 
Proposition \ref{SS:no repeat va si filtr}.
Therefore it suffices to take $a^i\in S$ in the definition of $N_I$.
Moreover, since $\rgr V$ is commutative, the operators $a_k$ of $\rgr M$
commute. 
Thus it suffices to prove that $N_I=N_{I'}$. 

The set $I$ is finite and totally ordered: 
define $(n_1, \ldots, n_r)<(m_1, \ldots, m_s)$ if
$r<s$ or $r=s$ and there is $k$ such that $n_i=m_i$ for $i<k$ and $n_k<m_k$.

Let  $I_x:=\set{y\in I\mid y<x}$ for $x\in I$.
We prove that if $x\in I$ such that $N_{I_x}\subset N_{I'}$ 
then $N_{\set{x}}\subset N_{I'}$.
This implies that $N_I=N_{I'}$.

We may assume that $x\notin I'$.
Then there exists $j$ such that $n:=n_j=n_{j+1}$ where $x=(n_i)$.
Let $a:=\si_{h_{a^j}}a^j$ and $b:=\si_{h_{a^{j+1}}}a^{j+1}$.
Solving the equation
$$
(ab)_{1+2n}
\; =\;
\sum_{i\geq 0}
(a_{-1-i}b_{1+2n+i}\, +\, b_{2n-i}a_i)
$$
for $a_n b_n=b_n a_n$,
we see that $N_{\set{x}}\subset N_{I_x}\subset N_{I'}$.
\end{pf}

\subsection{The Canonical Filtration}
\label{SS:canonical filt}

We give a spanning set for the canonical filtration $F\upsc$ and 
prove that $\rgr V$ is commutative and $F\upsc_1=C_2(V)$.

\bigskip

Let $V$ be an associative differential algebra.
The {\bf canonical} filtration $F\upsc$ is 
\index{canonical filtration}
the finest decreasing differential filtration such that $F\upsc_0=V$.
In other words, $F\upsc$ is the decreasing differential filtration 
generated by $V$ with $h_a=0$ for any $a$.

We have $F\upsc_n=\rspan\set{T^{n_1}(a_1)\dots T^{n_r}(a_r)\mid 
a_i\in V, n_i\geq 0, \sum n_i\geq n}$.
We have $C_n(V)\subset F\upsc_{n-1}$ for $n\geq 2$
because $T^{n-1}V\subset F\upsc_{n-1}$ 
and because $F\upsc_{n-1}$ is an ideal since $F\upsc_0=V$.
Moreover, $F\upsc_1=C_2(V)$.

The algebra $\rgr^{F\upsc}V$ is generated by 
$\si_n T^n a$ for $a\in V, n\geq 0$.
If $V$ is almost commutative then these elements commute since 
$[T^n a, T^m b]\in T^{n+m+1}V\subset F\upsc_{n+m+1}$.
Hence $\rgr^{F\upsc}V$ is commutative.

Let $V$ be a vertex algebra.
The {\bf canonical} filtration $F\upsc$ is 
\index{canonical filtration}
the finest decreasing differential filtration such that $F\upsc_0=V$.
In other words, $F\upsc$ is the decreasing differential filtration 
generated by $V$ with $h_a=0$ for any $a$.

\bigskip

{\bf Proposition.}\: {\it
Let $V$ be an associative vertex algebra.
Then $\rgr^{F\upsc} V$ is commutative and 
$F\upsc_n=\rspan\set{T^{n_1}(a_1)\dots T^{n_r}(a_r)\mid 
a_i\in V, n_i\geq 0, \sum n_i\geq n}$.
}

\bigskip

\begin{pf}
The vertex algebra $\rgr V$ is generated by $\si_0(a), a\in V$,
by Proposition \ref{SS:filt gen subset}.
These elements commute since $\si_0(a)_k\si_0(b)=\si_{-k-1}(a_k b)$ and
$\rgr_{-k-1} V=0$ for $k\geq 0$.
Thus $\rgr V$ is commutative.

Let $F_n$ be the above span of vectors.
By Proposition \ref{SS:filt gen subset} we need to show that 
$a^1_{n_1}\dots a^r_{n_r}1\in F_n$ for any $n_i\in\Z$ such that 
$\sum(-n_i-1)\geq n$.
Since $1\in F_0$, 
it suffices to show that $a_k F_n\subset F_{n-k-1}$ for $a\in V, k\in\Z$.
This is clear for $k<0$.
We prove by induction on $r$ that 
$a_k T^{n_1}(b^1)\dots T^{n_r}(b^r)\in F_{n-k-1}$ if 
$\sum n_i\geq n$ and $k\geq 0$.
For $r=0$, we have $a_k 1=0$.
Let $b':=T^{n_2}(b^2)\dots T^{n_r}(b^r)$.
Then 
$$
a_k T^{n_1}(b^1)b'
\; =\;
T^{n_1}(b^1)a_k b'
\; +\;
\sum\binom{k}{i}(a_i b^1)_{k-1-n_1-i}b'.
$$
By induction $T^{n_1}(b^1)a_k b'\in F_{n-k-1}$ and
$(a_i b^1)_{k-1-n_1-i}b'\in F_{n-k}$ because $(n-n_1)+(n_1+i-k)\geq n-k$.
\end{pf}

\bigskip

Let $V$ be an associative vertex algebra.
We have $C_n(V)\subset F\upsc_{n-1}$ for $n\geq 2$
as in the case of associative differential algebras.
The Proposition implies $F\upsc_1=C_2(V)$.

The Proposition shows that $\rgr^{F\upsc} V$ is a vertex Poisson algebra.
It is $\N$-graded as a commutative vertex algebra and 
$(\rgr_n V)_k(\rgr_m V)\subset\rgr_{n+m-k}V$ for $k\geq 0$.
In particular, $\rgr_0 V=V/C_2(V)$ has a commutative
multiplication and a bracket $\set{a,b}:=a_0 b$.
Thus we recover the Zhu Poisson algebra from section \ref{SS:assoc va def}.

\subsection{Generators without Repeats and $C_2(V)$ III}
\label{SS:no repeat va}

We extend Proposition \ref{SS:no repeat va si filtr} to vertex algebras
with elements of negative weight.

\bigskip

{\bf Remark.}\: {\it
Let $V$ be an almost commutative associative differential algebra.

\smallskip

\iti\:
We have $F\upsc_m\subset C_n(V)$ for $n\geq 3$ and $m\geq (n-2)2^{n-2}$.

\smallskip

\itii\:
Suppose that $V$ is $\Q$-graded such that $V_h=0$ for $h<h_0$ and 
let $S$ be a graded subspace.
Then $S+C_2(V)=V$ iff $S$ generates $V$ without repeats. 
${}_{}$${}_{}$${}_{}$${}_{}$${}_{}$${}_{}$${}_{}$${}_{}$
${}_{}$${}_{}$${}_{}$${}_{}$${}_{}$${}_{}$${}_{}$${}_{}$
\hfill $\square$
}

\bigskip

The Remark is proven in the same way as the following two results.

\bigskip

{\bf Lemma.}\: {\it
Let $V$ be an associative vertex algebra and $n\geq 3$.
Then $F\upsc_m\subset C_n(V)$ for $m\geq (n-2)2^{n-2}$.
}

\bigskip

\begin{pf}
Let $a\in F\upsc_m$.
By Proposition \ref{SS:canonical filt} we may assume that
$a=b^1\dots b^r$ where $b^i=T^{n_i}(a_i)$ such that $\sum n_i\geq m$.
If $n_i\geq n-1$ for some $i$ then $a\in C_n(V)$ since 
$C_n(V)$ is the ideal of $V\subast$ generated by $T^{n-1}V$.
Thus we may assume that $n_i\leq n-2$ for any $i$.

From $\sum n_i\geq m\geq (n-2)2^{n-2}$ we get 
$\#\set{i\mid n_i\ne 0}\geq 2^{n-2}$.
We prove by induction on $n$ that if $\#\set{i\mid n_i\ne 0}\geq 2^{n-2}$ 
then $a\in C_n(V)$.
It is true for $n=2$ since $C_2(V)$ is the ideal generated by $TV$.

Let $n\geq 2$ and $\#\set{i\mid n_i\ne 0}\geq 2^{n-1}$. 
Then by induction $a=b^1\dots b^s b$ for some $s<r$ and $b\in C_n(V)$ such 
that $\#\set{i\mid 1\leq i\leq s, n_i\ne 0}\geq 2^{n-2}$.
The commutator formula $[x_t,y_s]=\sum \binom{t}{i}(x_i y)_{t+s-i}$ implies
that
$$
b^1\dots b^s d_{-n}e
\;\equiv\;
d_{-n}b^1\dots b^s e\;\;\text{mod}\;\; C_{n+1}(V)
$$
for any $d, e\in V$ and by induction $b^1\dots b^s e\in C_n(V)$. 
Thus it suffices to prove that $a_{-n}b_{-n}c\in C_{n+1}(V)$
for any $a, b, c\in V$.
This follows from
$$
(ab)_{1-2n}c
\; =\;
\sum
(a_{-1-i}b_{1-2n+i}c+b_{-2n-i}a_i c)
$$
since all terms, except possibly $a_{-n}b_{-n}c$ for $i=n-1$, lie in 
$C_{n+1}(V)$.
\end{pf}

\bigskip

{\bf Proposition.}\: {\it
Let $V$ be a $\Q$-graded associative vertex algebra, 
such that $V_h=0$ for $h<h_0$, and $S$ be a graded subspace.
Then $S+C_2(V)=V$ iff $S$ generates $V\subast$ without repeats. 
}

\bigskip

\begin{pf}
By Proposition \ref{SS:canonical filt}
the vertex algebra $\rgr^{F\upsc}V$ is commutative and generated by
the subalgebra $\rgr_0 V=\si_0(S)$.
The Lemma shows that 
$F\upsc_n\subset C_m(V)\subset\bigoplus_{h\geq 2h_0+m-1}V_h$ 
for $n\geq (m-2)2^{m-2}$.
Therefore the claim follows from Lemma 
and Lemma \ref{SS:filtr alg}\itii.
\end{pf}

\bigskip

{\bf Corollary.}\: {\it
Let $V$ be a $C_2$-cofinite $\Q$-graded associative vertex algebra
such that $V_h=0$ for $h<h_0$.
Then $\dim V_h<\infty$ for any $h$ and $V$ is $C_n$-cofinite for any $n$.
}

\bigskip

\begin{pf}
Let $n_1>\ldots>n_r\geq 0$.
The vector $T^{n_1}(a_1)\dots T^{n_r}(a_r)$ has weight 
$\geq r(r-1)/2+rh_0$. 
This implies the first claim.
We have $T^{n_1}(a_1)\dots T^{n_r}(a_r)\in C_n(V)$ if $n_1\geq n-1$.
Thus $\dim V/C_n(V)<\infty$.
\end{pf}

\section{Supplements}
\label{S:suppl resul}

In sections \ref{SS:pre lie q assoc}---\ref{SS:left right wick}
we give a second proof of the fact that an associative vertex algebra 
satisfies the field identities,
we prove that duality and skew-symmetry imply locality for modules, and
we give a second proof of the fact that the left and right Wick
formulas are equivalent if skew-symmetry holds.

In section \ref{SS:lemma tensor prod}
we prove a lemma for tensor products.

In sections \ref{SS:bv algs}--\ref{SS:exten lie algebroi}
we discuss algebras and algebroids obtained from cohomology vertex algebras
and $\N$-graded vertex algebras.

In section \ref{SS:h Jac id res}
we rewrite the Jacobi identity in terms of residues.

\subsection{Pre-Lie Algebras and Quasi-Associativity}
\label{SS:pre lie q assoc}

We give a second proof of Proposition \ref{SS:assoc va}
stating that a vertex algebra is associative iff the field identities hold.

\bigskip

{\bf Proposition.}\: {\it
Let $V$ be a vertex algebra that satisfies the Wick formula and skew-symmetry. 
Then $V\subast$ is a pre-Lie algebra iff quasi-associativity holds.
}

\bigskip

\begin{pf}
By section \ref{SS:q assoc zhu poisson},
quasi-associativity implies the pre-Lie identity.
Conversely, suppose that $V\subast$ is a pre-Lie algebra.
Then 
\begin{align}
\notag
(ab)c-abc
\,&=\, 
cab-acb+[ab,c]-a[b,c]
\\
\notag
&=
-b[a,c]+[b,[a,c]]+[ab,c]-a[b,c].
\end{align}
By Remark \ref{SS:q assoc assoc} the right Wick formula holds.
Thus we get $[ab,c]=X(a,b)+X(b,a)+X$ where
\begin{align}
\notag
X(a,b)
\, &:=\,
\int_{-T}^0 d\la(e^{T\del\subla}a)b\subla c,
\quad
X
\, :=\,
\int_{-T}^0 d\la\int_0\upla d\mu\, b\submu a_{\la-\mu}c.
\end{align}
Using only properties of the translation operator we will show that 
$X(a,b)-a[b,c]= (\int_0^T d\la\, a)b\subla c$ and $X=-[b,[a,c]]$.
Thus we get quasi-associativity.

Since $[b,c]\subast=[b,c]\subslie$ is equivalent to skew-symmetry 
for $c, b$ and index $-1$ and since $T$ is a derivation, we have
\begin{align}
\notag
-a[b,c]
\; =\;
-a\int_{-T}^0 d\la\, b\subla c
\; =\;
&\res_z z\inv ae^{zT}b(-z)c
\; -\; 
a(bc)
\\
\notag
=\;
&\res_z z\inv e^{zT}(e^{-zT}a)b(-z)c
\; -\; 
a(bc).
\end{align}
Because $z\inv e^{zT}=z\inv +\int_0^T d\la\, e^{z\la}$ the first term
on the right-hand side is the sum of
\begin{align}
\notag
\res_z z\inv (e^{-zT}a)b(-z)c
\; =\;
\res_z (z\inv e^{zT}a)b(z)c
\; =\;
a(bc)
\; +\; 
\Big(\int_0^T d\la\, a\Big) b\subla c 
\end{align}
and, since $e^{z\la}e^{-zT}=e^{-T\del\subla}e^{z\la}$, of
\begin{align}
\notag
\res_z \int_0^T d\la\, e^{z\la}(e^{-zT}a)b(-z)c
\; =\;
&\int_0^T d\la\, (e^{-T\del\subla}a)\res_z e^{z\la}b(-z)c
\\
\notag
=\;
&\int_0^{-T} d\la\, (e^{T\del\subla}a)b\subla c
\; =\;
-X(a,b).
\end{align}
From $b\submu Td=(T+\mu)b\submu d$, the substitution formula, and 
the Fubini formula we obtain
\begin{align}
\notag
[b,[a,c]]\, =\,&\int_{-T}^0 d\mu\: b\submu\int_{-T}^0 d\la\, a\subla c
\, =\,
\int_{-T}^0 d\mu\int_{-T-\mu}^0 d\la\, b\submu a\subla c
\\
\notag
=\,
&\int_{-T}^0 d\mu\int_{-T}\upmu d\la\, b\submu a_{\la-\mu} c
\, =\, 
\int_{-T}^0 d\la\int\subla^0 d\mu\, b\submu a_{\la-\mu} c\, =\, -X.
\end{align}
\end{pf}

\bigskip

{\bf Corollary.}\: {\it
A vertex algebra is associative iff the field identities hold.
}

\bigskip

\begin{pf}
Suppose that $V$ is associative.
Then $V$ satisfies skew-symmetry by Proposition \ref{SS:skew sym}.
By Propositions \ref{SS:s3 symm} it suffice to prove the associativity formula.

By definition, the conformal Jacobi identity and the Wick formula hold.
By Remark \ref{SS:q assoc assoc} the right Wick formula holds.
By Proposition \ref{SS:pre lie q assoc} quasi-associativity holds.
Thus the associativity formula holds by Proposition \ref{SS:q assoc assoc}.

The converse is proven as in Proposition \ref{SS:assoc va}.
\end{pf}

\subsection{Duality and Locality for Modules}
\label{SS:dual skew loc}

We prove that duality and skew-symmetry imply locality for $\Z$-fold modules.
This result can be used in three proofs as a replacement of the fact 
that fields are local if they satisfy skew-symmetry,
see Proposition \ref{SS:skew sym fields}.
Namely, it can be used in the proof that $Y_M V$ is local for modules, 
that associative vertex subalgebras of $\Endv(E)$ are local, and 
that $U(R)$ is a quotient of $T\upast R$, 
see Corollary \ref{SS:skew sym fields} and 
Propositions \ref{SS:vas of fields} and \ref{SS:vertex envelope tensor alg}.

\bigskip

{\bf Proposition.}\: {\it
Let $V$ be a vertex algebra and $M$ a bounded $\Z$-fold $V$-module 
such that $Y_M(Ta)=\del_z Y_M(a)$.
Suppose that $a, b\in V$ satisfy skew-symmetry and 
$a, b, c$ and $b, a, c$ satisfy duality for any $c\in M$.
Then $Y_M(a)$ and $Y_M(b)$ are local of order $\leq o(a,b)$.
}

\bigskip

\begin{pf}
The idea of the proof is to argue that $a(bc)=(ab)c=(ba)c=b(ac)$.

Let $a, b, c$ and $b, a, c$ satisfy duality of order $\leq t$.
Let $r\geq o(a,b)$ and $s\geq o(b,c)$.
We have
\begin{align}
\notag
z^t w^s (z-w)^r a(z)b(w)c
\; &=\;
e^{-w\del_z}((z+w)^t w^s z^r a(z+w)b(w)c)
\\
\notag
&=\;
e^{-w\del_z}((w+z)^t w^s z^r (a(z)b)(w)c)
\\
\notag
&=\;
e^{-w\del_z}((w+z)^t w^s z^r (e^{zT}b(-z)a)(w)c)
\\
\notag
&=\;
e^{-w\del_z}e^{z\del_w}(w^t (w-z)^s z^r (b(-z)a)(w)c).
\end{align}
Define $p(z,w):=w^t (w-z)^s z^r (b(-z)a)(w)c$. The above
equations show that $p(z,w)\in V\pau{z,w}$. 
On the other hand, we have
\begin{align}
\notag
z^t w^s (z-w)^r b(w)a(z)c
\; &=\;
e^{-z\del_w}(z^t (w+z)^s (-w)^r b(w+z)a(z)c)
\\
\notag
&=\;
e^{-z\del_w}(z^t (z+w)^s (-w)^r (b(w)a)(z)c)
\\
\notag
&=\;
e^{-z\del_w}p(-w,z).
\end{align}
Thus the claim follows from 
$e^{-w\del_z}e^{z\del_w}p(z,w)=e^{-w\del_z}p(z,w+z)=
e^{-w\del_z}\linebreak[0]p(z,z+w)=p(z-w,z)$ and 
$e^{-z\del_w}p(-w,z)=p(-w+z,z)=p(z-w,z)$.
\end{pf}

\subsection{The Left and Right Wick Formulas}
\label{SS:left right wick}

We prove Remark \ref{SS:q assoc assoc} directly, 
without using $\bbS_3$-symmetry.

\bigskip

{\bf Remark.}\: {\it
Let $V$ be a vertex algebra that satisfies skew-symmetry.
Then the left Wick formula is equivalent to the right Wick formula:
$$
[ab\subla c]
\; =\;
(e^{T\del\subla}a)[b\subla c]
\; +\;
\paraab(e^{T\del\subla}b)[a\subla c]
\; +\;
\paraab\int_0\upla d\mu\, [b\submu [a_{\la-\mu}c]].
$$
}

\bigskip

\begin{pf}
The left Wick formula is equivalent to 
$$
e^{T\del\subla}[c_{-\la}ab]
\; =\;
e^{T\del\subla}\Big( [c_{-\la}a]b\; +\; a[c_{-\la}b]\; +\;
\int_0^{-\la}d\mu\,[[c_{-\la}a]\submu b]\Big).
$$
Note that $[c_{-\la}a]b=b[c_{-\la}a]+\int_{-T}^0 d\mu\, [[c_{-\la}a]\submu b]$.
Since 
$[g\submu b]=e^{\mu\del\subla}[e^{T\del\subla}g\submu b]$, 
we have
\begin{align}
\notag
e^{T\del\subla}
\int_{-T}^{-\la}d\mu\,[[c_{-\la}a]\submu b]
\; =\;
&\int_{-T}^{-\la-T}d\mu\,e^{(T+\mu)\del\subla}[[c_{-\la-T}a]\submu b]
\\
\notag
=\;
&-\int_0\upla d\mu\,[[c_{-(\la-\mu)-T}a]_{-\mu-T}b].
\end{align}
Using skew-symmetry, this shows that the above identity is also 
equivalent to the right Wick formula.
\end{pf}

\subsection{A Lemma for Tensor Products}
\label{SS:lemma tensor prod}

We give an alternative proof of Lemma \ref{SS:tensor va}.

\bigskip

{\bf Lemma.}\: {\it
Let $V$ be an associative vertex algebra and $S, S'\subset V$ be two
commuting subsets.
For any $a, b\in S, a', b'\in S'$, and $r\in\Z$, we have
$$
(aa')_r(bb')
\; =\;
\sum_{i\in\Z}\: \zeta^{a'b}\,(a_i b)(a'_{r-i-1}b').
$$
}

\bigskip

\begin{pf}
Since $Y: V\to\Endv(V)$ is a monomorphism, it is equivalent to prove the
identity for $YS, YS'$.
Let $a(z), b(z)\in YS$ and $a'(z), b'(z)\in YS'$ be local of order $\leq n$
and $\leq n'$. 
Then
$c(z,w):=(z-w)^n a(z)b(w)$ and $c'(z,w):=(z-w)^{n'}a'(z)b'(w)$
are fields by Remark \ref{SS:skew sym fields}.

By Proposition \ref{SS:commut centre}\,\iti\ we have $[a(z),a'(w)]=0$.
Thus $a(z)a'(w)=a'(w)a(z)$ is also a field and 
$\normord{a(z)a'(z)}=a(z)a'(z)$.
Moreover, $c(z,w)c'(z,w)\linebreak[0]=c'(z,w)c(z,w)$ is a field.
We have
$$
a(z)a'(z)\,b(w)b'(w)
\, =\, 
a(z)b(w)\,a'(z)b'(w)
\, =\,
\frac{c(z,w)c'(z,w)}{(z-w)^{n+n'}}
$$
and similarly 
$b(w)b'(w)a(z)a'(z)=c(z,w)c'(z,w)/(z-w)^{n+n'}_{w>z}$.
Moreover, we have the product formula
$$
\del_z^{(k)}(c(z,w)c'(z,w))
\; =\;
\sum\: \del_z^{(j)}c(z,w)\,\del_z^{(k-j)}c'(z,w).
$$
Thus the claim follows from Remarks \ref{SS:dong fields} and 
\ref{SS:skew sym fields} by taking $k=n+n'-r-1=(n-i-1)+(n'-(r-i-1)-1)$.
\end{pf}

\subsection{Batalin-Vilkovisky Algebras}
\label{SS:bv algs}

We prove that any odd, self-commuting element $G$ of a vertex algebra 
yields a BV-operator and thus an odd Leibniz bracket $[\, ,]_G$.

\bigskip

An {\bf odd bracket} on a vector space $\fg$ is 
\index{odd!bracket}
an odd linear map $\fg\otimes\fg\to\fg$.
Equivalently, it is a bracket on $\rPi\fg$,
where $\rPi$ is the parity-change functor.
An {\bf odd Leibniz bracket} on  $\fg$ is 
\index{odd!Leibniz bracket}
a Leibniz bracket on $\rPi\fg$. 
In other words, we have 
$$
[[a,b],c]
\; =\;
[a,[b,c]]
\; -\;
\zeta^{(a+1)(b+1)}\, [b,[a,c]].
$$
An {\bf odd Lie bracket} is 
\index{odd!Lie bracket}
a Lie bracket on $\rPi\fg$.
In other words, we have in addition $[a,b]=-\zeta^{(a+1)(b+1)}[b,a]$.
A {\bf derivation} of an odd bracket is 
\index{derivation!of an odd bracket}
a derivation of the algebra $\rPi\fg$.
Thus $d[a,b]=[da,b]+\zeta^{(a+1)d}[a,db]$.

A {\bf Gerstenhaber algebra} is a commutative algebra with an odd Lie
bracket such that $[a,\;]$ is a derivation of the commutative algebra
of parity $\ta+\oone$.
Thus $[a,bc]=[a,b]c+\zeta^{(a+1)b}b[a,c]$.
Roughly speaking, Gerstenhaber algebras are 
Poisson algebras for which the Poisson bracket is odd.

A {\bf BV-operator} on an algebra $V$ is an odd homogeneous 
differential operator $G$ of order $\leq 2$ such that $G^2=0$.
In other words, $G$ is an odd operator such that
$\del_a:=[G,a\cdot]-(Ga)\cdot$ is a derivation for any $a\in V$ and $[G,G]=0$.
Note that the parity of $\del_a$ is $\ta+1$.

For a BV-operator $G$, define an odd bracket by 
$$
[a,b]_G
\; :=\; 
\zeta^a\del_a b
\; =\;
\zeta^a\big( G(ab)-(Ga)b-\zeta^a a(Gb)\big).
$$
By definition, $\cdot$ and $[\, ,]_G$ satisfy the Leibniz identity.

A {\bf Batalin-Vilkovisky algebra} or {\bf BV-algebra} is 
a commutative algebra with a BV-operator.

For example, let $V$ be an associative vertex algebra and 
$G$ an odd element such that $G\subla G=0$.
Then $(G_{(t)})^2=[G_{(t)},G_{(t)}]/2=0$ for any $t\in\Z$.
In particular, $G_{(1)}$ is a BV-operator by Proposition \ref{SS:wick}.
The commutator formula implies $[a,b]_G=\zeta^a(G_{(0)}a)_{(0)}b$.

\bigskip

{\bf Remark.}\: {\it
Let $d$ be a differential of a vector space $E$, $a$ an operator, and 
$b:=[d,a]$.
Then $[b,d]=0$ and $b=0$ on $H_d(E)$.
}

\bigskip

\begin{pf}
By definition, $d$ is odd.
Thus $[d,d]=0$.
We have $[d,[d,a]]=[[d,d],a]-[d,[d,a]]$.
Hence $[b,d]=0$.
We have $b=0$ on $H_d(E)$ since $bc=dac$ for any $c\in\ker d$.
\end{pf}

\bigskip

{\bf Proposition.}\: {\it
Let $G$ be a BV-operator of an algebra $V$.
Then $[\, ,]_G$ is an odd Leibniz bracket and 
$G$ is a derivation of $[\, ,]_G$.
If the multiplication of $V$ is commutative then $[\, ,]_G$ is an 
odd Lie bracket.
In particular, any BV-algebra is a Gerstenhaber algebra.
}

\bigskip

\begin{pf}
The Remark and $G^2=0$ implies 
$[G,\del_a]=[G,[G,a\cdot]]-[G,(Ga)\cdot]=-\del_{Ga}$.
Thus $G$ is a derivation of $[\, ,]_G$.

The operator $G$ is also a BV-operator of the unital algebra $V\oplus\K 1$ 
where $G1:=0$.
Thus we may assume that $V$ is unital.
The claim that $[\, ,]_G$ is an odd Leibniz bracket is equivalent to
$[\del_a,\del_b]=\zeta^{a+1}\del_{\del_a b}$.
Both sides of this identity are derivations.
The space of derivations $d$ has zero intersection with the space $V\cdot$ of
left multiplications $a\cdot$ since $d1=0$ and $a1=a$.
Calculating modulo $V\cdot$, we have
$$[\del_a,\del_b]
\;\equiv\;
[\del_a,[G,b\cdot]]
\;\equiv\;
[[\del_a,G],b\cdot]+\zeta^{a+1}[G,(\del_a b)\cdot]
\;\equiv\;
\zeta^{a+1}\del_{\del_a b}
$$
using that $[\del_a,G]=-\zeta^{a+1}[G,\del_a]=\zeta^{a+1}\del_{Ga}$.
The remaining two claims are clear.
\end{pf}

\subsection{Topological Vertex Algebras and Chiral Rings}
\label{SS:topl va}

We prove that the cohomology of a topological vertex algebra is a BV-algebra.

\bigskip

Let $V$ be a vertex algebra. 
A {\bf differential} of $V$ is 
\index{differential!of a vertex algebra}
an odd derivation $d$ such that $d^2=0$.

If $d$ is a derivation of $V$ then $\ker d$ is a vertex subalgebra. 
If $d$ is a differential then $\im\, d$ is an ideal of $\ker d$,
since $d(ab)=(da)b$ and $d[a\subla b]=[(da)\subla b]$ for $b\in\ker d$.
Hence the {\bf cohomology} $H_d(V):=\ker d/\im\, d$ is
\index{cohomology!vertex algebra} 
a 
\index{vertex algebra!cohomology}
vertex algebra.
The 
\index{closed element}
elements of $\ker d$ and $\im\, d$ are 
\index{exact element}
called {\bf $d$-closed} and {\bf $d$-exact}, resp.

\bigskip

{\bf Remark.}\: {\it
Let $V$ be an associative conformal vertex algebra with a differential $Q$.
If $L$ is $Q$-exact then 
the cohomology vertex algebra $H_Q(V)$ is commutative with $T=0$.
}

\bigskip

\begin{pf}
If $L=QG$ then $[Q,G_{(0)}]=(QG)_{(0)}=T$. 
By Remark \ref{SS:bv algs} we get $T=0$ on $H_Q(V)$.
In particular, $H_Q(V)$ is commutative.
\end{pf}

\bigskip

Let $V$ be as in the Remark and $L=QG$ so that $[Q,G_{(0)}]=T$.
The fact that $H_Q(V)\subast$ is a commutative algebra can also be seen
by constructing so-called homotopies $\mu$ and $\nu$. 
Defining
$$
\mu(a,b)
\; :=\;
\sum (-1)^i G_{(0)}T^{(i)}(a_i b)/(i+1)
$$
for any $a, b\in V$, we have $Q\mu(a,b)=[a,b]-\mu(Qa,b)-\zeta^a\mu(a,Qb)$ since
$[\, ,]\subast=[\, ,]\subslie$.
Thus the multiplication of $H_Q(V)\subast$ is commutative.
Defining
$$
\nu(a,b,c)
\; :=\;
\sum\big( (G_{(0)}T^{(i)}a)b_i c\: +\: \paraab(G_{(0)}T^{(i)}b)a_i c\big)/(i+1)
$$
we have 
$Q\nu(a,b,c)=(ab)c-abc-\nu(Qa,b,c)-\zeta^a\nu(a,Qb,c)-\zeta^{a+b}\nu(a,b,Qc)$
because of quasi-associativity. 
Thus $H_Q(V)\subast$ is associative.

In section \ref{SS:n1 n2 top vir} we defined the 
topological Virasoro vertex Lie algebra $\tVir$.
It is generated by $L, J, Q, G$, and $\hd$.
A {\bf topological vertex algebra} is a conformal vertex algebra $V$
together with a morphism $\tVir\to V$ of conformal vertex Lie algebras
such that $J_0$ is diagonalizable on $V$.
Here and in the following we denote 
the image of $a\in\tVir$ in $V$ also by $a$.

Let $V$ be an associative topological vertex algebra.
The elements $Q, G$ of $V$ are odd, have weight $1$ and $2$, and 
satisfy $Q\subla G=L+J\la+\hd\la^{(2)}$ and $Q\subla Q=G\subla G=0$.
Thus $Q_0=Q_{(0)}$ is a differential, $Q_0G=L$, and
$G_0=G_{(1)}$ is a BV-operator.

The differential $Q_0$ is a derivation of $[\, ,]_{G_0}$ because
$[a,b]_{G_0}=\zeta^a(G_{(0)}a)_{(0)}b$ and 
\begin{align}
\notag
Q_0((G_{(0)}a)_{(0)}b)
\; &=\;
(Q_0G_{(0)}a)_{(0)}b\; +\; \zeta^{a+1}(G_{(0)}a)_{(0)}Q_0b
\\
\notag
\; &=\;
-(G_{(0)}Q_0a)_{(0)}b\; +\; \zeta^{a+1}(G_{(0)}a)_{(0)}Q_0b
\end{align}
where we have used that $[Q_0,G_{(0)}]=T$ and $(Tc)_{(0)}=0$.
Thus $[\, ,]_{G_0}$ induces an odd bracket on $H_{Q_0}(V)$.
The Remark and Proposition \ref{SS:bv algs} imply 
that $H_{Q_0}(V)$ is a BV-algebra.

Since $Q_0$ is homogeneous of degree $0$, the cohomology $H_{Q_0}(V)$
is graded with respect to $L_0$.
Since $[Q_0,G_0]=L_0$, we get $H_{Q_0}(V)=H_{Q_0}(V_0)$ 
by Remark \ref{SS:bv algs}.
We have $[Q_0,G_0]=0$ on $V_0$.

An $N$=2 {\bf superconformal vertex algebra} is 
\index{N2 superconformal vertex algebra@$N$=2 superconformal vertex algebra}
a conformal vertex algebra $V$ together with a morphism $\twoVir\to V$ 
such that $J_0$ is diagonalizable on $V$.
Of course, there is a bijection $V\mapsto V^A$ from $N$=2 superconformal
vertex algebras to topological vertex algebras obtained 
from the vertex Lie algebra isomorphism $\tVir\to\twoVir$,
see section \ref{SS:n1 n2 top vir}.
The mirror involution of $\twoVir$ yields another bijection $V\to V^B$,
see again section \ref{SS:n1 n2 top vir}.

The role of the differential $Q_0$ and the BV-operator $G_0$ is played
by $G^+_{-1/2}$ and $G^-_{1/2}$ or, alternatively, by $G^-_{-1/2}$ and
$G^+_{1/2}$.
The only difference is that $[G\suppm_{-1/2},G\supmp_{1/2}]=2L_0\mp J_0$ 
and hence $H_{G\suppm_{-1/2}}(V)=H_{G\suppm_{-1/2}}(V')$ where 
$V':=\set{a\mid 2L_0 a=\pm J_0 a}$.

Disregarding the fact that instead of $V$ one should consider the
full conformal field theory, 
the topological vertex algebras $V^A$ and $V^B$ are 
\index{topological twist}
the {\bf topological twists} or 
\index{Amodel@$A$-model}
the {\bf $A$-} and {\bf $B$-model} 
\index{Bmodel@$B$-model}
of $V$. 
The commutative algebras
$H_{G\suppm_{-1/2}}(V)$ are known as 
\index{chiral ring}
the {\bf chiral} and 
\index{anti-chiral ring}
{\bf anti-chiral ring}.

\subsection{1-Truncated Vertex Algebras and Courant Algebroids}
\label{SS:1trunc va}

We describe the algebraic structure on the subspaces $V_0$ and $V_1$
of an $\N$-graded vertex algebra and of an $\N$-graded vertex Poisson algebra.

We denote the $t$-th products by $a_t b$, not by $a_{(t)}b$.

\bigskip

Let $\cO$ be a commutative algebra and $\fa$ a 1-truncated vertex Lie algebra 
over the vector space $\cO$, see section \ref{SS:1trunc vlie}.
Denote the elements of $\cO$ by $f, g, h$ and the elements of $\fa$ by $x, y$.

Suppose that $\fa$ is endowed with an even linear map 
$\cO\to\End(\fa), f\mapsto f\cdot$.
Then we call $\fa$ a {\bf 1-truncated vertex algebra} over $\cO$ if 
\index{onetruncated vertex algebra@1-truncated vertex algebra}
$$
f(gx)
\; -\;
(fg)x
\; =\;
(xf)dg
\; +\;
(xg)df,
$$
and $1x=x$; 
the anchor $\fa\to\fgl(\cO)$ has image in $\Der(\cO)$ and 
satisfies $(fx)g=f(xg)$; 
we have $[x,fy]=(xf)y+f[x,y], (x,fy)=f(x,y)-y(xf)$,
and $d(fg)=fdg+gdf$.

For example, 1-truncated vertex algebras $\fa$ over $\K$ 
are exactly Lie algebras $\fa$ with 
an even invariant symmetric bilinear form since $d=0$ and $\Der(\K)=0$,
see section \ref{SS:1trunc vlie}.

If $V$ is an $\N$-graded associative unital vertex algebra then
$V_1$ is a 1-truncated vertex algebra over $V_0$.
Indeed, in section \ref{SS:comm va} we noted that 
$V_0$ is a commutative subalgebra of $V\subast$.
In section \ref{SS:1trunc vlie} we showed that $V_1$ is a 
1-truncated vertex Lie algebra over $V_0$.
The map $V_0\to\End(V_1), f\mapsto f\cdot$, is multiplication in $V\subast$.
The first axiom follows from quasi-associativity and from 
$f_{-2}g_0 x=-(xg)Tf+T((xg)_0 Tf)=-(xg)Tf$.
We have $(fx)g=\sum f_{-1-i}x_i g+x_{-1-i}f_i g=f(xg)$ and the 
commutator formula implies $(x,fy)=f(x,y)+(x_0 f)_0 y=f(x,y)-y(xf)$.
The remaining two non-trivial axioms follow from the fact that $a_0$
is a derivation.

A {\bf Courant algebroid} over $\cO$ is
\index{Courant algebroid}
a 1-truncated vertex Lie algebra $\fa$ over $\cO$ with an $\cO$-module
structure such that the anchor is an $\cO$-module morphism with image in
$\Der(\cO)$, $[x,fy]=(xf)y+f[x,y]$, the pairing $(x,y)$ is $\cO$-bilinear,
and $d:\cO\to\fa$ is a derivation.

Thus the difference between a Courant algebroid and a 1-truncated 
vertex algebra is that $\fa$ is a genuine $\cO$-module and the
pairing $(x,y)$ is $\cO$-bilinear. 

If $V$ is an $\N$-graded vertex Poisson algebra then
$V_1$ is a Courant algebroid over $V_0$.
This is similar to the case of associative vertex algebra.
We have $(fx)g=-g_0(fx)=f(xg)$, since $g_0 f=0$, and 
$(x,fy)=f(x,y)$, since $x_1 f=0$.

\subsection{Extended Lie Algebroids}
\label{SS:exten lie algebroi}

We construct extended Lie algebroids from 1-truncated vertex algebras
and from Courant algebroids.

\bigskip

Let $\cO$ be a commutative algebra.
Denote the elements of $\cO$ by $f, g, h$.

A {\bf Lie algebroid} over $\cO$ is
\index{Lie algebroid}
a Lie algebra $\cT$ together with an $\cO$-module structure and
a morphism $\cT\to\Der(\cO)$ of Lie algebras and $\cO$-modules such that
$[x,fy]=(xf)y+f[x,y]$ for any $x, y\in\cT$.

The morphism $\cT\to\Der(\cO)$ is 
\index{anchor}
the {\bf anchor}.
For example, $\Der(\cO)$ is a Lie algebroid over $\cO$.
Lie algebroids for which the anchor is zero are precisely Lie algebras over 
$\cO$.

An {\bf extended Lie algebroid} over $\cO$ is
\index{extended Lie algebroid}
a Lie algebroid $\cT$ together with a module $\Om$ over $\cO$ and 
over the Lie algebra $\cT$, an even $\cO$-bilinear map 
$\cT\otimes_{\cO}\Om\to\cO, x\otimes\om\mapsto (x,\om)$, and 
an even derivation and $\cT$-module morphism $d:\cO\to\Om$, such that
the $\cO$- and $\cT$-module structures on $\Om$ are related by 
$x(f\om)=x(f)\om+f(x\om)$ and $(fx)\om=f(x\om)+(x,\om)df$
and the pairing satisfies $(x,df)=xf$ and $x(y,\om)=([x,y],\om)+(y,x\om)$ 
for any $x, y\in\cT$ and $\om\in\Om$. 

For example, $\Der(\cO)$ with the module $\Om_{\cO/\K}^1$ of K\"ahler
differentials is an extended Lie algebroid.
We have $x(fdg)=(xf)dg+fd(xg)$ and $(x,fdg)=f(xg)$.

More generally, if $\cT$ is a Lie algebroid over $\cO$ then
$\cT$ with $\Om:=\Hom_{\cO}(\cT,\cO)$ is an extended Lie algebroid.
The $\cT$-module structure on $\Om$ is defined by 
$(y,x\om)=x(y,\om)-([x,y],\om)$,
the pairing $\cT\otimes_{\cO}\Om\to\cO$ is given by evaluation, and 
the derivation $d$ is defined by the formula $(x,df)=xf$.

\bigskip

{\bf Proposition.}\: {\it
Let $\fa$ be a 1-truncated vertex algebra or a Courant algebroid.
Define $\Om:=\rspan\set{fdg\mid f, g\in\cO}$ and $\cT:=\fa/\Om$.
Then $\cT$ with $\Om$ is an extended Lie algebroid.
}

\bigskip

\begin{pf}
We prove the claim only in the case that $\fa$ is a 1-truncated vertex algebra.
The case of a Courant algebroid is proven in the same way.

The morphism $\fa\to\Der(\cO)$ factorizes over $\cT$ because 
$(fdg)h=f((dg)h)\linebreak[0]=0$.
It follows that $\Om$ is invariant with respect to $f\cdot$ for any $f\in\cO$
and that $\Om$ and $\cT$ are $\cO$-modules.

The subspace $\Om$ is a left ideal and an abelian subalgebra of the
Leibniz algebra $\fa$ because $[x,fdg]=(xf)dg+fd(xg)$.
Since $[x,y]+[y,x]=T(x,y)$, we see that $\Om$ is a two-sided ideal, 
$\cT$ is a Lie algebra, and $\Om$ is a $\cT$-module. 
Thus $\cT$ with $\cT\to\Der(\cO)$ is a Lie algebroid over $\cO$.

From $(x,fdg)=f(x,dg)-(dg)(xf)=f(xg)$ we get $(\Om,\Om)=0$.
Thus the pairing $\fa\otimes\fa\to\cO$ induces a map 
$\cT\otimes\Om\to\cO, x\otimes\om\mapsto (x,\om)$.
This map is $\cO$-bilinear since $(x,fy)=f(x,y)-y(xf)$ for $x, y\in\fa$.

Finally, for $\om\in\Om$ we have $(fx)\om=f(x\om)+(x,\om)df$ because
$[fx,\om]=-[\om,fx]+d(fx,\om)=-(\om f)x-f[\om,x]+d(f(x,\om)-x(\om f))=
f[x,\om]+(x,\om)df$.
This proves the claim.
\end{pf}

\subsection{Jacobi Identity in terms of Residues}
\label{SS:h Jac id res}

We rewrite the Jacobi identity in terms of residues.

\bigskip

{\bf Proposition.}\: {\it
The following is equivalent:

\begin{enumerate}
\item[$\iti$]
$a, b, c$ satisfy the Jacobi identity;

\smallskip

\item[$\itii$]
for any $F(z,w,x)\in\K\lau{z,w,x}$ we have
$$
\res_x F(w+x,w,x)\, (a(x)b)(w)c
\; =\;
\res_z F(z,w,z-w)\, [a(z),b(w)]c;
$$

\smallskip

\item[\itiii]
the identity in \itii\ is satisfied
for any $F(z,w,x)\in\K[z\uppm,w\uppm,(z-w)\uppm]$;

\smallskip

\item[\itiv]
the identity in \itii\ is satisfied
for any $F(z,w,x)\in\K[z\uppm,(z-w)\uppm]$.
\end{enumerate}
}

\bigskip

\begin{pf}
\iti$\Rightarrow$\itii\:
Let $G(z,w,x)\in\K\lau{z,w,x}$.
The basic properties of the delta distribution yield
$$
\res_x G(w+x,w,x)\, (a(x)b)(w)c
\; =\;
\res_{z,x} G(z,w,x)\, \de(z,w+x)(a(x)b)(w)c
$$
and 
$$
\res_z G(z,w,z-w)\, [a(z),b(w)]c
\; =\;
\res_{z,x} G(z,w,x)\, \de(z-w,x)[a(z),b(w)]c.
$$
This shows that \iti\ implies \itii.

\smallskip

$\itii\Rightarrow\itiii$\:
Choose $G(z,w,x)\in \K[z\uppm,w\uppm,x\uppm]$
such that $G(z,w,z-w)=F(z,w,x)$.
Applying $\itii$ to $G(z,w,x)\in\K\lau{z,w,x}$
we obtain $\itiii$.

\smallskip

$\itiii\Rightarrow\itiv$\: This is trivial.

\smallskip

$\itiv\Rightarrow\iti$\:
Let $r, t\in\Z$ and $G(z,w,x):=z^t x^r$.
Applying $\itiv$ to $F(z,w,x)=G(z,w,z-w)$ for any $r, t\in\Z$
we obtain $\iti$ because of the two equations in the proof of 
\iti$\Rightarrow$\itii.
\end{pf}

\bigskip

{\bf Proposition.}\: {\it
Let $V$ be a bounded $\Z$-fold algebra.
Then the Jacobi identity is equivalent to locality and duality.
}

\bigskip

\begin{pf}
Suppose that $a, b, c$ satisfy locality of order $\leq r$ and duality
of order $\leq t$ and that $a(z)c\in z^{-t}V\pau{z}$.
Define $d(z,w,x):=x^{-r}(z-w)^r a(z)b(w)c\in z^{-t}x^{-r}V\pau{z}\lau{w}$.
Then $d(z,w,z-w)=a(z)b(w)c$ and $d(z,w,-w+z)=b(w)a(z)c$.
Moreover, we have
\begin{align}
\notag
(w+x)^t d(w+x,w,x)
\; &=\;
(x+w)^t d(x+w,w,x)
\\
\notag
&=\;
(x+w)^t a(x+w)b(w)c
\\
\notag
&=\;
(w+x)^t (a(x)b)(w)c.
\end{align}
Since $d(w+x,w,x), (a(x)b)(w)c\in V\lau{w}\lau{x}$ 
we get $d(w+x,w,x)=(a(x)b)(w)c$.
Proposition \ref{SS:delta dis}\,\iti\ and Remark \ref{SS:ope2} yield
\begin{align}
\notag
\de(z,w+x)(a(x)b)(w)c
\; &= \;
\de(z,w+x)d(z,w,x)
\\
\notag
\; &= \;
(\de(z-w,x)-\de(-w+z,x))d(z,w,x)
\\
\notag
\; &= \;
\de(x,z-w)[a(z),b(w)]c.
\end{align}
\end{pf}

\chapter{Enveloping Vertex Algebras}
\label{C:env}

There is a forgetful functor $V\mapsto V\subslie$ from 
associative vertex algebras to vertex Lie algebras.
This functor has a left adjoint $R\mapsto U(R)$.
The associative vertex algebra $U(R)$ is the enveloping vertex algebra 
of the vertex Lie algebra $R$. 
The category of bounded $R$-modules is isomorphic to the category
of $U(R)$-modules.

The purpose of this chapter is to explain three constructions of
$U(R)$ and to prove the Poincar\'e-Birkhoff-Witt theorem for $U(R)$.

In section \ref{S:env va} we construct $U(R)$ as a quotient of the
free vertex algebra generated by $R$ and as a quotient of the
tensor vertex algebra $T\upast R$.

In section \ref{S:PBW thm} we show that a certain Verma module $V(\fg)$
over a local Lie algebra $\fg$ has an associative vertex algebra structure 
and use the identity $U(R)=V(\fg(R))$ to prove the PBW-theorem.

\section{Free, Enveloping, and Tensor Vertex Algebras}
\label{S:env va}

In sections \ref{SS:loc fct}--\ref{SS:vertex envelope}
we construct enveloping vertex algebras as quotients of free vertex algebras.
In sections \ref{SS:tensor alg vleib}--\ref{SS:vertex envelope tensor alg}
we prove that the tensor algebra $T\upast R$ of a vertex Lie algebra $R$
has a vertex algebra structure and 
that $U(R)$ is the quotient of $T\upast R$ by the ideal
generated by $a\otimes b-b\otimes a-[a,b]\subslie$.

\subsection{Locality Function}
\label{SS:loc fct}

We prove some properties of the locality function of a vertex algebra.
They are used in the construction of free vertex algebras 
in section \ref{SS:free va}.

\bigskip

Recall that a locality function on a set $S$ is a map 
$o: S^2\to\Z\cup\set{-\infty}$.
The locality function $o$ of a vertex algebra $V$ is defined by 
$o(a,b):=\inf\set{n\in\Z\mid a_i b=0, i\geq n}$ for $a, b\in V$.

For fields $a(z), b(z)$, 
define $N(a(z),b(z))$ to be the least $n\in\Z$ such 
\index{local!pair of fields}
that $a(z), b(z)$ are {\bf local} of order $\leq n$: $(z-w)^n [a(z),b(w)]=0$.

\bigskip

{\bf Remark.}\: {\it
If $a(z), b(z)$ are local fields then $N(a(z),b(z))=o(a(z),b(z))$.
In particular, if $V$ is an associative unital vertex algebra 
then $o(a,b)=N(Ya,Yb)$.
}

\bigskip

\begin{pf}
It is clear that $o(a(z),b(z))\leq N(a(z),b(z))$.
Let $n, m\in\Z$ such that $n\geq N(a(z),b(z)), m\geq o(a(z),b(z))$, and 
$n\geq m$.
By Remark \ref{SS:skew sym fields} we have 
$d(z,w):=(z-w)^m a(z)b(w)\in\cF_2(E)$.
Since 
$$
(z-w)^{n-m}d(z,w)\; =\; (z-w)^n a(z)b(w)\; =\; (z-w)^n b(w)a(z)
$$
and both $d(z,w)$ and $b(w)a(z)$ are contained in $\cF\lau{w}\lau{z}(E)$ we get
$d(z,w)=(z-w)^m b(w)a(z)$.

If $V$ is an associative unital vertex algebra then $o(a,b)=o(Ya,Yb)$ 
since $Y$ is a monomorphism of unbounded vertex algebras.
Thus the second claim follows from the first.
\end{pf}

\bigskip

{\bf Proposition.}\: {\it
Let $o$ be the locality function of an associative vertex algebra.
Then:

\smallskip

\iti\: 
$o(a,b)=o(b,a)$;

\smallskip

\itii\:
$o(Ta,b)\leq o(a,b)+1$;

\smallskip

\itiii\:
$o(a_r b,c)\leq o(a,b)+o(a,c)+o(b,c)-r-1$ for any $r\in\Z$.
}

\bigskip

\begin{pf}
\iti, \itii\ follow from skew-symmetry and $(Ta)_r b=-r a_{r-1}b$.

\smallskip

\itiii\:
The Jacobi identity implies that $a, b, c$ satisfy duality of order 
$\leq o(a,c)$.
Let $s\geq o(a,b)+o(a,c)+o(b,c)-r-1$.
By Proposition \ref{SS:kind of assoc}, there are $\la^i\in\Z$ such that 
for $N:=1-o(a,c)-o(a,b)$ we have 
$$
(a_r b)_s c
\; =\;
\sum_{i\geq N}\: \la^i\, a_{-i}(b_{s+r+i}c).
$$
This is zero since $s+r+i\geq s+r-o(a,c)-o(a,b)+1\geq o(b,c)$.
\end{pf}

\subsection{Free Associative Vertex Algebras}
\label{SS:free va}

We construct the {\bf free} associative vertex algebra $V(S,o)$ generated 
\index{free!vertex algebra}
by a set $S$ with locality function $o$.

\bigskip

Recall that a morphism between sets with a locality function is a map $\phi$
such that $o(\phi a,\phi b)\leq o(a,b)$.

\bigskip

{\bf Proposition.}\: {\it
The functor $V\mapsto (V,o)$, from associative vertex algebras to 
sets with a locality function, has a left adjoint $(S,o)\mapsto V(S,o)$. 
}

\bigskip

\begin{pf}
Let $F$ be the free $\Z$-fold algebra with derivation generated by the set $S$.
A basis of $F$ is given by $B=\bigsqcup B_n$ where 
$B_1:=\set{T^k a\mid a\in S, k\geq 0}$ and 
$B_n:=\set{a_i b\mid a\in B_m, b\in B_{n-m}, 1\leq m\leq n-1, i\in\Z}$ for 
$n\geq 2$.

Define a map $o: B^2\to\Z\cup\set{-\infty}$ by induction on $n$: 
$o(T^k a,T^l b):=o(a,b)+k+l$ for $a, b\in S$, 
$o(T^k a,b_i c):=o(T^k a,b)+o(T^k a,c)+o(b,c)-i-1$ 
for $a\in S, b_i c\in B_n$, and 
$o(a_i b,c):=o(a,b)+o(a,c)+o(b,c)-i-1$ for $a_i b\in B_n, c\in B$.
Let $F'$ be the quotient of $F$ by the relations $a_i b=0$ and 
$(Ta)_n b=-n a_{n-1}b$ for $a, b\in B$ and $i\geq o(a,b), n\in\Z$.
Then $F'$ is a vertex algebra.
Let $V(S,o)$ be the quotient of $F'$ by the identities defining
an associative vertex algebra.

Let $V$ be an associative vertex algebra with a morphism 
$(S,o)\to (V,o)$.
Since $F$ is free, 
there exists a unique morphism $F\to V$ extending $S\to V$.
Proposition \ref{SS:loc fct} implies that $F\to V$ induces 
a morphism $F'\to V$.
Since $V$ is associative, the map $F'\to V$ induces a 
unique morphism $V(S,o)\to V$ extending $S\to V$.
\end{pf}

\bigskip

Note that $V:=V(S,o)\oplus\K 1$ is the free associative unital vertex algebra
generated by $S$ with locality function $o$.

Let $\cR\subset V$ be a subset and $I\subset V$ the ideal generated by $\cR$.
Then $V/I$ is the associative unital vertex algebra 
\index{relations!for a vertex algebra}
generated by $(S,o)$ with {\bf relations} $\cR$.

Suppose that $S$ is a $\K$-graded set, $S=\bigsqcup S_h$.
In the following we use the algebras $F, F', V(S,o)$ defined in
the proof of the Proposition. 
The gradation of $S$ induces a gradation on the free $\Z$-fold algebra $F$.
Thus $F'$ is also graded.
Since the associativity formula and skew-symmetry are homogeneous identities
it follows that $V(S,o)$ is a graded vertex algebra.

\subsection{Enveloping Vertex Algebras}
\label{SS:vertex envelope}

We construct the {\bf enveloping vertex algebra} $U(R)$ of a 
vertex Lie algebra.

\bigskip

Recall the definition of 
the enveloping algebra $U(\fg)$ of a Lie algebra $\fg$.
By definition, $\fg\mapsto U(\fg)$ is the left adjoint of 
the functor $V\mapsto (V,[\, ,]\subast)$ from associative unital algebras 
to Lie algebras.
The algebra $U(\fg)$ is constructed as the quotient of the 
free associative unital algebra $T\upast \fg$, generated by the
vector space $\fg$, by the ideal
generated by $[\io a,\io b]\subast-\io[a,b]$ for $a, b\in\fg$.
Here $\io:\fg\to T\upast \fg$ is the canonical map.

Recall also that a bounded module over a conformal algebra $R$ is 
a vector space $M$ together with 
a morphism $Y_M: R\to\Endv(M)\subslie$ such that $Y_M(R)$ is local,
see section \ref{SS:module vlie}.

\bigskip

{\bf Proposition.}\: {\it
The functor $V\mapsto V\subslie$, from associative unital vertex algebras to 
vertex Lie algebras, has a left adjoint $R\mapsto U(R)$.
Moreover,
the categories of bounded $R$-modules and of $U(R)$-modules are isomorphic.
}

\bigskip

\begin{pf}
Let $F$ be the free associative unital vertex algebra generated by
the $\K[T]$-module $R$ with locality function 
$o_+, o_+(a,b)=\inf\set{n\mid a_i b=0, i\geq n}$.
Let $U(R)$ be the quotient of $F$ by the ideal generated by
$[\io a\subla\io b]-\io[a\subla b]$ for $a, b\in R$.
Here $\io: R\to F$ is the canonical map.
Then the composition $R\to F\to U(R)$ is a morphism of vertex Lie algebras.

Let $V$ be an associative unital vertex algebra with a morphism 
$R\to V\subslie$.
Then there is a unique morphism $F\to V$ extending $R\to V$.
This morphism induces a unique morphism $U(R)\to V$ extending $R\to V$.

By Corollary \ref{SS:skew sym fields} 
any $U(R)$-module is a bounded $U(R)\subslie$-module.
Thus the canonical morphism $R\to U(R)\subslie$ yields a bounded $R$-module. 
Conversely, if $M$ is a bounded $R$-module then $\sqbrack{Y_M(R)}$
is an associative unital vertex algebra by Corollary \ref{SS:vas of fields}.
Thus there is a unique morphism $U(R)\to\sqbrack{Y_M(R)}$ extending $Y_M$.
Hence $M$ is an $U(R)$-module. 
\end{pf}

\bigskip

{\bf Remark.}\: {\it
Let $R$ be the free vertex Lie algebra generated by a set $S$ with 
a non-negative locality function $o$.
Let $V$ be the free associative unital vertex algebra generated by $S$
with locality function $o$.
Then $V=U(R)$.
}

\bigskip

\begin{pf}
This follows from the fact that both vertex algebras satisfy the
same universal property.
\end{pf}

\bigskip

Let $R, S, o, V$ be as in the Remark.
Let $R'$ be the vertex Lie algebra generated by $(S,o)$ with 
relations $\cR\subset R$.
Let $\cR'$ be the image of $\cR$ in $V$ and 
$V'$ the associative unital vertex algebra generated by $(S,o)$ with 
relations $\cR'$.
Then $V'=U(R')$ since both algebras satisfy the same universal property.

Let $R$ be a vertex Lie algebra and 
$F$ the increasing invariant filtration of $U(R)$ generated by $R$
with $h_a=1$ for any $a\in R$.
Since $R$ generates $U(R)$ and $[a\subla b]\in R[\la]\subset F_1[\la]$ 
for any $a, b\in R$,
Proposition \ref{SS:pbw generat}\,\iti\ and its proof
imply that $R$ generates $U(R)\subast$ with the PBW-property, 
$\rgr U(R)$ is a vertex Poisson algebra, 
and $F_n=\rspan\set{a_1\dots a_r\mid a_i\in R, 0\leq r\leq n}$.
In section \ref{SS:PBW thm} we prove that $\rgr U(R)=S\upast R$.

\subsection{Tensor Vertex Algebra}
\label{SS:tensor alg vleib}

We show that the tensor algebra $T\upast R$ of a vertex Leibniz algebra
has a natural vertex algebra structure.
We 
\index{tensor vertex algebra}
call $T\upast R$ the {\bf tensor vertex algebra} of $R$.

\bigskip

A {\bf Leibniz algebra} is
\index{Leibniz!algebra}
an algebra that satisfies the Leibniz identity,
see section \ref{SS:conf jac}.
A {\bf vertex Leibniz algebra} is
\index{vertex Leibniz algebra}
a conformal algebra that satisfies the conformal Jacobi identity.

The {\bf tensor algebra} of 
\index{tensor algebra}
a vector space $E$ is $T\upast E=\bigoplus_{n\geq 0} E^{\otimes n}$.
It is the free associative unital algebra generated by $E$.

\bigskip

{\bf Proposition.}\: {\it
Let $R$ be a vertex Leibniz algebra.
There exists a unique unital vertex algebra structure 
on $T\upast R$ such that 

\smallskip

\iti\:
the inclusion $R\to T\upast R$ is a conformal algebra morphism;

\smallskip

\itii\:
the identity is $1\in\K=T^0 R$ and 
$ab=a\otimes b$ for $a\in R, b\in T\upast R$;

\smallskip

\itiii\:
the associativity formula is satisfied for any $a\in R$ and 
$b, c\in T\upast R$.

\smallskip

\noindent
The increasing invariant filtration generated by $R$, 
with $h_a=1$ for any $a\in R$, is $F_n=\bigoplus_{i\leq n}R^{\otimes i}$.
Moreover, $(F_n)\subla F_m\subset F_{n+m-1}[\la]$.
}

\bigskip

\begin{pf}
The translation operator of $R$ induces a derivation $T$ 
of the tensor algebra $T\upast R$.
This must be the translation operator of $T\upast R$.

\smallskip

{\it Step 1.}\:
We first construct a linear map $Y:R\to\Endv(T\upast R)$ such that
$ab=a\otimes b$ for $a\in R, b\in T\upast R$ and such that  
$T$ is a translation operator.
As in section \ref{SS:modules},
it suffices to define $a\subla b$ for any $a\in R, b\in T\upast R$ 
such that $T$ is a translation operator for $a\subla b$.
We do this by induction as follows.
Define $a\subla 1:=0, a\subla b:=a\subla b$, and
$$
a\subla(bc)
\; :=\;
(a\subla b)c
\; +\;
b(a\subla c)
\; +\;
\int_0\upla d\mu\, (a\subla b)\submu c
$$
for $a, b\in R, c\in T\upast R$.
It is clear that $T$ is a derivation for $a\subla b$ and that
$(Ta)\subla b=-\la a\subla b$ for any $a\in R, b\in T\upast R$.
Moreover, $F_1 F_m\subset F_{1+m}$ and $(F_1)\subla F_m\subset F_m[\la]$.

We now prove the associativity formula for $a, b\in R, c\in T\upast R$
and indices $r\geq 0, s\in\Z$.
By Proposition \ref{SS:2nd recursion}
the associativity formula for indices $r\geq 0, s\in\Z$ is equivalent to
the conformal Jacobi identity and the Wick formula.
By construction, the Wick formula is satisfied for any
$a, b\in R, c\in T\upast R$. 
By induction on $n$, we prove the conformal Jacobi identity 
for any $a, b\in R$ and $c\in T^n R$.
This is clear for $n=0$.
For $n=1$ we have the conformal Jacobi identity for $R$. 
If $a, b, c\in R$ and $d\in T^n R$ then the induction hypothesis implies
\begin{align}
\notag
&
a\subla b\submu(cd)
\\
\notag
\; =\;
&a\subla((b\submu c) d)
\; +\;
a\subla(c(b\submu d))
\; +\;
\int_0\upmu d\nu\, a\subla((b\submu c)\subnu d)
\\
\notag
=\;
&(a\subla b\submu c) d
\, +\,
(b\submu c)(a\subla d)
\, +\,
\int_0\upla d\nu\, (a\subla b\submu c)\subnu d
\, +\,
(a\subla c)(b\submu d)
\, +\,
c(a\subla b\submu d)
\\
\notag
&+\,
\int_0\upla d\nu\, (a\subla c)\subnu (b\submu d)
\, +\,
\int_0\upmu d\nu\, (b\submu c)\subnu(a\subla d)
\, +\,
\int_0\upmu d\nu\, (a\subla b\submu c)_{\nu+\la} d.
\end{align}
Together with the substitution formula this yields
$[a\subla,b\submu](cd)=(a\subla b)_{\mu+\la}(cd)$.

\smallskip

{\it Step 2.}\:
Define $Y: T\upast R\to\Endv(T\upast R)$ by induction:
$1(z):=\id, a(z):=a(z)$ using step 1, and 
$(ab)(z):=\;\normord{a(z)b(z)}$ for $a\in R, b\in T\upast R$.
Because $\del_z$ is a derivation of $\Endv(T\upast R)$ and
the subspace of translation covariant fields is a $\Z$-fold subalgebra,
it follows that $T$ is a translation operator of $T\upast R$.
The left identity $1$ is also a right identity of $(T\upast R)\subast$
because $\normord{a(z)b(z)}1=a(z)_-b(z)1=(e^{zT}a)e^{zT}b$ for 
$a\in R, b\in T\upast R$.

We now prove the associativity formula for any $a\in R$ and 
$b, c\in T\upast R$.
By definition, it is satisfied for $a\in R, b, c\in T\upast R$
and indices $r=-1, s\in\Z$.
Hence by Proposition \ref{SS:2nd recursion}\,\iti\ it is satisfied for 
indices $r<0, s\in\Z$.

By induction on $n$, we prove that the associativity formula is 
satisfied for $a\in R, b\in T^n R, c\in T\upast R$ 
and indices $r\geq 0, s\in\Z$.
This is clear for $n=0$.
For $n=1$ this was proven in step 1.
Let $a, b\in R$ and $c\in T^n R$.
Then the induction hypothesis implies 
\begin{align}
\notag
&(a\subla(bc))(z)
\\
\notag
=\;
&((a\subla b)c)(z)
\; +\;
(b(a\subla c))(z)
\; +\;
\int_0\upla d\mu\, ((a\subla b)\submu c)(z)
\\
\notag
=\;
&\normord{b(z)(a(z)\subla c(z))}
\; +\;
\normord{(a(z)\subla b(z))c(z)}
\; +\;
\int_0\upla d\mu\, ((a(z)\subla b(z))\submu c(z)
\\
\notag
=\;
&a(z)\subla\normord{b(z)c(z)}
\end{align}
where we have used that $\Endv(T\upast R)$ satisfies the Wick formula.
Thus $Y$ satisfies \iti--\itiii.
The construction of $Y$ shows that $Y$ is unique.

By induction on $n$, we show that $F_n F_m\subset F_{n+m}$ and 
$(F_n)\subla F_m\subset F_{n+m-1}[\la]$.
The first claim follows from quasi-associativity and
the second claim follows from the right Wick formula for
$a\in R, b\in F_n$, and $c\in F_m$.
Clearly, $TF_n\subset F_n$.
Thus $F$ is an invariant filtration.
Since $R\subset F_1$, the invariant filtration $F'$ generated by $R$
is finer than $F$.
Conversely, $F_n\subset F'_n$.
Thus $F=F'$.
\end{pf}

\bigskip

If $R$ is abelian then the tensor vertex algebra $T\upast R$
is equal to the associative differential algebra $T\upast R$.
Indeed, the latter considered as a vertex algebra, 
with $(T\upast R)\subslie$ abelian,
satisfies properties \iti--\itiii\ of the Proposition.

In general, however, the algebras $(T\upast R)\subast$ and $T\upast R$ are 
very different. 
For example, quasi-associativity yields $(aa)c=a(ac)+2\sum(T^{(i+1)}a)a_i c$
for $a, c\in R$.

The vertex algebra $\rgr T\upast R$ is equal to the associative algebra 
$T\upast R$.
More precisely, $(\rgr T\upast R)\subast$ is equal to the associative
differential algebra $T\upast R$ because quasi-associativity for
$a\in R, b, c\in T\upast R$ implies associativity for $a\in F_1, 
b, c\in\rgr T\upast R$ and this implies associativity of 
$(\rgr T\upast R)\subast$ since then
$(a_1\dots a_r)(b_1\dots b_s)=a_1\dots a_r b_1\dots b_s$ for 
$a_i, b_j \in F_1$.
The Proposition implies that $(\rgr T\upast R)\subslie$ is abelian.

\subsection{Universality of the Tensor Vertex Algebra}
\label{SS:univ tensor alg}

We prove that $T\upast R$ is universal among all vertex algebras $V$
that satisfy the associativity formula for $a\in R, b, c\in V$.

\bigskip

{\bf Lemma.}\: {\it
Let $V$ and $W$ be unital vertex algebras and $S\subset V$ a subset 
such that $V=\rspan\set{a^1_{n_1}\dots a^r_{n_r}1\mid a^i\in S, n_i\in\Z}$.
Let $\phi: V\to W$ be an even linear map such that $\phi 1=1$ and
$\phi(a_n b)=(\phi a)_n(\phi b)$ for any $a\in S, b\in V$, and $n\in\Z$.
Then $\phi$ is a $\Z$-fold algebra morphism if either

\smallskip

\iti\:
$V$ and $W$ are associative or

\smallskip

\itii\:
duality holds for $a\in S$ and $b, c\in V$ and for $a\in \phi S$ and 
$b, c\in\phi V$.
}

\bigskip

\begin{pf}
This is proven by induction on $r$.
For $r=0$, we have $\phi 1=1$.
Let  $a\in S$ and $b, c\in V$ such that $\phi(b_n c)=(\phi b)_n(\phi c)$
for any $n\in\Z$.
By Proposition \ref{SS:kind of assoc}\,\iti\ for any $r, s\in\Z$ 
there exist $N, \la^i\in\Z$ such that
\begin{align}
\notag
\phi((a_r b)_s c)
\; =\;
\phi\sum_{i\geq N}\la^i\,a_{-i}b_{r+s-i}c
\; =\;
&\sum_{i\geq N}\la^i(\phi a)_{-i}(\phi b)_{r+s-i}(\phi c)
\\
\notag
\; =\;
&\phi(a_r b)_s\phi(c).
\end{align}
\end{pf}

\bigskip

{\bf Proposition.}\: {\it
Let $R$ be a vertex Leibniz algebra, $V$ a unital vertex algebra, and 
$\phi: R\to V$ a morphism such that 
$V$ satisfies the associativity formula for any $a\in\phi R, b, c\in V$.
Then $\phi$ has a unique extension $T\upast R\to V$.
In particular, $R\mapsto T\upast R$ is a functor.
}

\bigskip

\begin{pf}
Define $\phi: T\upast R\to V$ by induction: $\phi 1:=1, \phi a:=\phi a$, 
and $\phi(ab):=\phi(a)\phi(b)$ for $a\in R, b\in T\upast R$. 
Since $T$ is a derivation, it follows that $\phi: T\upast R\to V$
is a $\K[T]$-module morphism. 
Since $a_{-1-t}b=(T^{(t)}a)b$ for $t\geq 0$, 
this implies that $\phi(a_t b)=(\phi a)_t(\phi b)$ for any
$a\in R, b\in T\upast R$, and $t<0$.

By induction on $k$, we prove that $\phi(a\subla b)=(\phi a)\subla(\phi b)$ 
for any $a\in R$ and $b\in T^k R$.
This is clear for $k=0$.
For $k=1$ we have that $\phi: R\to V$ is a conformal algebra morphism.
Let $a, b\in R$ and $c\in T^k R$.
Then the Wick formula and the induction hypothesis imply
\begin{align}
\notag
&\phi(a\subla(bc))
=
\phi\Big((a\subla b)c
\, +\,
b(a\subla c)
\, +\,
\int_0\upla d\mu\, (a\subla b)\submu c\Big)
=
\phi(a)\subla\phi(bc).
\end{align}

Thus we have shown that $\phi(a_n b)=(\phi a)_n(\phi b)$ for any
$a\in R, b\in T\upast R, n\in\Z$.
By part \itii\ of the Lemma 
the map $\phi: T\upast R\to V$ is therefore a $\Z$-fold algebra morphism.
It is clear that $\phi$ is unique.
\end{pf}

\subsection{Tensor and Enveloping Vertex Algebras}
\label{SS:vertex envelope tensor alg}

We prove that $U(R)$ is a quotient of the tensor vertex algebra $T\upast R$.

\bigskip

Recall that the enveloping algebra $U(\fg)$ of a Lie algebra $\fg$ 
is the quotient of the tensor algebra $T\upast \fg$ by the ideal
generated by $a\otimes b-\paraab b\otimes a-[a,b]$ for $a, b\in\fg$.

In section \ref{SS:bracket vlie}
we showed that any vertex Lie algebra has a natural Lie bracket 
$[\, ,]\subslie$.

\bigskip

{\bf Proposition.}\: {\it
Let $R$ be a vertex Lie algebra.
Then $U(R)$ is the quotient of the tensor vertex algebra $T\upast R$ by 
the ideal generated by $a\otimes b-\paraab b\otimes a-[a,b]\subslie$ for 
$a, b\in R$.
}

\bigskip

\begin{pf}
Let $U$ be the above quotient of $T\upast R$.
Denote the image of $a\in R$ in $U$ also by $a$.
Then $[a,b]\subast=[a,b]\subslie$ in $U$ for any $a, b\in R$.
Conformal skew-symmetry is also satisfied for $a, b\in R$
since the canonical map $R\to U$ is a conformal algebra morphism.
Hence skew-symmetry holds for $a, b\in R$ by Proposition \ref{SS:skew sym}.

By Proposition \ref{SS:tensor alg vleib} any $a\in R, b, c\in U$
satisfy the associativity formula and $ab=a\otimes b$.
Thus $Y(U)\subset\sqbrack{Y(R)}$.
Doing Corollary \ref{SS:skew sym fields} elementwise for $a, b\in R$,
we see that $a(z), b(z)\in\Endv(U)$ are local for any $a, b\in R$.
Hence $\sqbrack{Y(R)}$ is local by Dong's lemma. 
This proves that $U$ is associative.

Let $V$ be an associative vertex algebra and $\phi: R\to V\subslie$ a 
morphism. 
By Proposition \ref{SS:univ tensor alg} the map $\phi$ has a unique 
extension $\phi: T\upast R\to V$.
Since $V$ satisfies skew-symmetry the morphism $\phi: T\upast R\to V$
induces a morphism $U\to V$.
Thus $U$ satisfies the universal property of the enveloping algebra.
\end{pf}

\bigskip

By Proposition \ref{SS:tensor alg vleib} 
the invariant filtration generated by $R\subset T\upast R$ is given by 
$F_n=\bigoplus_{i\leq n}R^{\otimes i}$.
It induces the invariant filtration generated by the image of $R$ in $U(R)$.
This filtration is discussed at the end of section \ref{SS:vertex envelope}.

\section{Poincar\'e-Birkhoff-Witt Theorem}
\label{S:PBW thm}

In sections \ref{SS:exis thm} and \ref{SS:va over loclie}
we use the existence theorem to construct 
a vertex algebra $V(\fg)$ from any local Lie algebra $\fg$.
In section \ref{SS:univ va over loclie}
we prove that $V(\fg)$ is the universal associative vertex algebra over $\fg$.
In section \ref{SS:PBW thm}
we prove for any vertex Lie algebra $R$ 
that the module categories for $R$, $\fg(R)$, and $U(R)$ are 
isomorphic and that $U(R)=V(\fg(R))$ and $\rgr U(R)=S\upast R$.

\subsection{Existence Theorem for Associative Vertex Algebras}
\label{SS:exis thm}

We prove that a vector space $V$ with a generating set $S\subset\Endv(V)$ 
of local fields has a vertex algebra structure such that $Y(V)=\sqbrack{S}$.
In section \ref{SS:univ va over loclie}
we prove that this vertex algebra structure is unique.

\bigskip

Let $E$ be a $\K[T]$-module with an even vector $1$.
Recall that $a(z)\in\End(E)\pau{z\uppm}$ is creative if 
$a(z)1=e^{zT}a_{-1}1$, see section \ref{SS:loc skew sym}.

\bigskip

{\bf Goddard's Uniqueness Theorem.}\: {\it
\index{Goddard's uniqueness theorem}
Let $\cF\subset\End(E)\pau{z\uppm}$ be a local subspace such that 
any $a(z)\in\cF$ is creative.
If $s_1:\cF\to E, a(z)\mapsto a_{-1}1$, is surjective 
then $s_1$ is an isomorphism.
}

\bigskip

\begin{pf}
Since $a(z)1=e^{zT}a_{-1}1$, it suffices to show that 
the linear map $\cF\to E\pau{z}, a(z)\mapsto a(z)1$, is injective.
Let $a(z)\in\cF$ such that $a(z)1=0$.
We prove that $a(z)b=0$ for any $b\in E$.
Let $b(z)\in\cF$ such that $b_{-1}1=b$.
For $r\gg 0$, we have $(z-w)^r a(z)b(w)1=(z-w)^r b(w)a(z)1=0$.
Setting $w=0$ yields $a(z)b=0$.
\end{pf}

\bigskip

{\bf Lemma.} \: {\it
If $a(z), b(z)\in\Endv(E)$ are creative then 
$(a(z)_n b(z))_{-1}1=a_n b_{-1}1$ for any $n\in\Z$.
}

\bigskip

\begin{pf}
Since $a_i 1=b_i 1=0$ for $i\geq 0$, we have 
$$
(a(z)_n b(z))_{-1}1=
\sum_{i\geq 0} (-1)^i\binom{n}{i}(a_{n-i}b_{-1+i}1-(-1)^n b_{-1+n-i}a_i 1)=
a_n b_{-1}1.
$$
\end{pf}

\bigskip

Recall that $\Endv(E)_T$ denotes the unbounded vertex algebra of 
translation covariant fields, see section \ref{SS:fields}.

\bigskip

{\bf Existence Theorem.} \: {\it
Let $V$ be a vector space with an even vector $1$ and an even operator $T$
such that $T1=0$.
Suppose that $S\subset\Endv(V)_T$ is a local subset such that 
$$
V
\; =\;
\rspan\set{a^1_{n_1}\dots a^r_{n_r}1\mid a^i(z)\in S, n_i\in\Z, r\geq 0}.
$$
Then there exists an associative vertex algebra structure on $V$ 
with identity $1$ such that $Y(V)=\sqbrack{S}$.
}

\bigskip

\begin{pf}
The unital vertex algebra $\sqbrack{S}$ is associative and local
by Proposition \ref{SS:vas of fields} and Corollary \ref{SS:vas of fields}.

Since $\sqbrack{S}\subset\Endv(V)_T$,
any $a(z)\in\sqbrack{S}$ is creative by Lemma \ref{SS:loc skew sym}.
The Lemma and the fact that $V$ is spanned by 
$a^1_{n_1}\dots a^r_{n_r}1$ implies that the map 
$s_1: \sqbrack{S}\to V, a(z)\mapsto a_{-1}1$, is surjective. 
By Goddard's uniqueness theorem $s_1$ is bijective. 

Define $Y:=s_1\inv$.
The Lemma shows that $s_1: \sqbrack{S}\to(V,Y)$ is an isomorphism of
unital $\Z$-fold algebras.
Thus $V$ is an associative vertex algebra with identity $1$ 
such that $Y(V)=\sqbrack{S}$.
\end{pf}

\subsection{Vertex Algebras over Local Lie Algebras}
\label{SS:va over loclie}

We construct an associative vertex algebra $V(\fg)$ from 
a local Lie algebra $\fg$.
In section \ref{SS:univ va over loclie}
we prove that the vertex algebra structure on $V(\fg)$ is unique and that
$V(\fg)$ is the universal vertex algebra over $\fg$.

\bigskip

{\bf Definition.}\: 
Let $\fg$ be a local Lie algebra.
A {\bf vertex algebra over} $\fg$ is 
\index{vertex algebra!over a local Lie algebra}
a vertex algebra $V$ together with a $\fg$-module structure $\rho$ 
such that $\rho F\subfg\subset Y(V)$.

\bigskip

A {\it morphism} of vertex algebras over $\fg$
is a vertex algebra morphism that is also a $\fg$-module morphism.

The condition $\rho F\subfg\subset Y(V)$ is equivalent to 
$Y(a_{-1}1,z)=\rho a(z)$ and hence to $(a_{-1}1)_t b=a_t b$ 
for any $a(z)\in F\subfg, b\in V$, and $t\in\Z$.

Recall that a $\fg$-module $M$ is 
bounded if $\rho F\subfg\subset\Endv(M)$.

\bigskip

{\bf Lemma.}\: {\it 
Let $\fg$ be a local Lie algebra, $M$ a $\fg$-module,
and $W\subset M$ a subspace such that $U(\fg)W=M$ and 
$a(z)b\in M\lau{z}$ for any $a(z)\in \rho F\subfg, b\in W$.
Then $M$ is bounded.
}

\bigskip

\begin{pf}
The canonical filtration of $U(\fg)$ induces an exhaustive increasing
filtration $(M_i)_{i\geq 0}$ of $M$ such that $M_0=W$.
We prove the claim by induction on $i$.
The induction beginning $i=0$ holds true by assumption.

Let $a(z), b(z)\in F\subfg$ and $c\in M_i$.
Since $a(z), b(z)$ are local and $a(z)c, b(z)c\in M\lau{z}$ 
there exists $n\geq 0$ such that 
$$
d(z,w)
\; :=\; 
(z-w)^n a(z)b(w)c
\; =\;
(z-w)^n b(w)a(z)c
\; \in\;
M\lau{z,w}.
$$ 
The distribution $(z-w)^{-n}d(z,w)-a(z)b(w)c$ is a Laurent series in $w$
and annihilated by $(z-w)^n$.
By Remark \ref{SS:loc diff op} we get $(z-w)^{-n}d(z,w)=a(z)b(w)c$.
Hence $a(z)b(w)c\in M\lau{z}\lau{w}$. 
In other words, $a(z)(b_k c)\in M\lau{z}$ for any $k\in\Z$.
\end{pf}

\bigskip

Let $\fg$ be a local Lie algebra.
Define $\bfg_+:=\rspan\set{a_t\mid a(z)\in\bF\subfg, t\geq 0}$.
The weak commutator formula shows that $\bfg_+$ is a subalgebra.
Since $a(w)_r b(w)=\res_z (z-w)^r [a(z),b(w)]$,
the algebra $\bfg_+$ is contained in the subalgebra generated by 
$\fg_+$ and hence is equal to it.

\smallskip

Let $\K$ be the trivial one-dimensional $\bfg_+$-module so that 
$\bfg_+\K=0$.
Define the induced $\fg$-module 
$$
V(\fg)
\; :=\;
U(\fg)\otimes_{U(\bfg_+)}\K
$$ 
and $1:=1\otimes 1\in V(\fg)$.
The map $U(\fg)/U(\fg)\bfg_+\to V(\fg), a\mapsto a1$, 
is a $\fg$-module isomorphism with inverse $a\otimes\la\mapsto\la a$.

\bigskip

{\bf Proposition.}\: {\it
If $\fg$ is a local Lie algebra then the $\fg$-module $V(\fg)$ has the 
structure of an associative vertex algebra over $\fg$ with identity $1$.
}

\bigskip

\begin{pf}
We apply the existence theorem from section \ref{SS:exis thm}
to $V=V(\fg)$ and $S=\rho F\subfg$.
It is clear that $\rho F\subfg$ is local and $V(\fg)$ is spanned by
$a^1_{n_1}\dots a^r_{n_r}1$.
The Lemma shows that $\rho F\subfg\subset\Endv(V(\fg))$.

The derivation $T$ of $\fg$ induces a derivation $T$ of $U(\fg)$ which
in turn induces an operator $T$ of $V(\fg), a\otimes 1\mapsto Ta\otimes 1$,
because $T\bfg_+\subset\bfg_+$.
Since $T$ is a derivation of $U(\fg)$ we have $T1=0$ and 
$[T,a\cdot]=(Ta)\cdot$ for any $a\in\fg$.
It follows that $[T,\rho a(z)]=\rho (Ta)(z)=\del_z \rho a(z)$ 
for any $a(z)\in F\subfg$.
\end{pf}

\subsection{Universal Vertex Algebra over a Local Lie Algebra}
\label{SS:univ va over loclie}

We prove that the vertex algebra $V(\fg)$ is 
the universal vertex algebra over $\fg$.

\bigskip

{\bf Remark.} \: {\it
Let $V$ be a vector space, $1\in V\even$, and $S\subset\Endv(V)$ a subset 
such that $V$ is the span of $a^1_{n_1}\dots a^r_{n_r}1$ for 
$a^i(z)\in S, n_i\in\Z, r\geq 0$.
Then there exists at most one associative vertex algebra structure on $V$ 
with identity $1$ and such that $S\subset Y(V)$.
}

\bigskip

\begin{pf}
Since $S\subset Y(V)$ and $Y: V\to Y(V)$ is a bijection with inverse 
$a(z)\mapsto a_{-1}1$ we have $Y(a_{-1}1,z)=a(z)$ for any $a(z)\in S$.
This is equivalent to $(a_{-1}1)_n b=a_n b$ for any 
$a(z)\in S, b\in V, n\in\Z$.
Thus the associativity formula implies
$Y(a^1_{n_1}\dots a^r_{n_r}1,z)=a^1(z)_{n_1}\dots a^r(z)_{n_r}\id_V$.
\end{pf}

\bigskip

Let $\fg$ be a local Lie algebra.
By Proposition \ref{SS:va over loclie} the $\fg$-module $V(\fg)$ has 
an associative vertex algebra structure over $\fg$ with identity $1$.
The Remark shows that this vertex algebra structure is unique.

An associative unital vertex algebra $V$ over $\fg$ is {\bf universal} if  
\index{universal vertex algebra of $\fg$}
$V$ is an initial object of the category of 
associative unital vertex algebras over $\fg$.

\bigskip

{\bf Proposition.}\: {\it
If $\fg$ is a local Lie algebra
then $V(\fg)$ is the universal associative unital vertex algebra over $\fg$.
}

\bigskip

\begin{pf}
Let $W$ be an associative unital vertex algebra over $\fg$.
Then $\bfg_+ 1=0$ in $W$ because $\rho\bF\subfg\subset Y(W)$.
The universal property of the induced module $V(\fg)$ implies that
there exists a unique $\fg$-module morphism $\phi:V(\fg)\to W$
such that $1\mapsto 1$. 

In order to show that $\phi$ is a vertex algebra morphism 
we apply Lemma \ref{SS:univ tensor alg}\,\iti\ with
$S:=\rspan\set{a_{-1}1\mid a(z)\in F\subfg}$.
From $(a_{-1}1)_t b=a_t b$ follows that 
$V(\fg)=\rspan\set{a^1_{n_1}\dots a^r_{n_r}1\mid a^i\in S, n_i\in\Z, r\geq 0}$.
We have $\phi(e_t b)=(\phi e)_t(\phi b)$ for $e\in S$ and $b\in V(\fg)$ since
$\phi((a_{-1}1)_t b)=\phi(a_t b)=a_t\phi(b)=(a_{-1}1')_t\phi(b)=
\phi(a_{-1}1)_t\phi(b)$ where $1':=1\in W$.
\end{pf}

\bigskip

If $\fg\to\fg'$ is a morphism of local Lie algebras
then any $\fg'$-module is a $\fg$-module.
Together with the Proposition this shows that 
$\fg\mapsto V(\fg)$ is a functor.

\subsection{Poincar\'e-Birkhoff-Witt Theorem}
\label{SS:PBW thm}

We prove the Poincar\'e-Birkhoff-Witt theorem for enveloping vertex algebras.

\bigskip

{\bf Proposition.}\: {\it
Let $R$ be a vertex Lie algebra.
Then $U(R)=V(\fg(R))$.
}

\bigskip

\begin{pf}
We show that $R\mapsto V(\fg(R))$ is left adjoint to $V\mapsto V\subslie$.
Let $\fg:=\fg(R)$.

Let $V$ be an associative unital vertex algebra and 
$\phi: R\to V\subslie$ a morphism. 
Then $V$ with $Y\circ\phi$ is an $R$-module and hence a $\fg$-module
by Remark \ref{SS:module vlie}.
Thus $V$ is an associative unital vertex algebra over $\fg$ and 
we get a unique morphism $V(\fg)\to V$ that is compatible with $R$.

Conversely, the map 
$\io: R\to V(\fg)\subslie, a\mapsto a_{-1}1$, is a morphism
because it is the composition of $Y_{V(\fg)}: R\to\Endv(V(\fg))$ and 
$Y(V(\fg))\to V(\fg), a(z)\mapsto a(0)1$.
This map $\io$ is the unit of adjunction.
\end{pf}

\bigskip

The proof shows that to give an associative unital vertex algebra over $\fg(R)$
is equivalent to giving an associative unital vertex algebra $V$ together with
a morphism $R\to V\subslie$.
 
We denote by $S\upast E=\bigoplus S^n E$ the symmetric algebra of 
a vector space $E$.
It is the free commutative algebra generated by $E$.
It is an $\N$-graded algebra.

Recall that if $R$ is a vertex Lie algebra then 
$R$ with $[\, ,]\subslie$ is a Lie algebra, see section \ref{SS:bracket vlie}.
This Lie algebra is denoted by $R\subsmlie$.

\bigskip

{\bf PBW Theorem.}\: {\it
Let $R$ be a vertex Lie algebra.
Then there exists a unique isomorphism of $R\subsmlie$-modules 
$U(R)\to U(R\subsmlie)$ such that $1\mapsto 1$.
Moreover, $\rgr U(R)=S\upast R$ 
as $\N$-graded commutative differential algebras.
}

\bigskip

\begin{pf}
Let $\fg:=\fg(R)$.
By Proposition \ref{SS:bracket vlie}
there is a differential Lie algebra isomorphism
$R\subsmlie\to\fg_-, a\mapsto a_{-1}$. 

If $W$ is an associative vertex algebra and $\io: R\to W\subslie$ a 
morphism then $W$ is an $R$-module and hence a $R\subsmlie$-module,
see section \ref{SS:bracket vlie}.
The action is just $a\otimes b\mapsto (\io a)b$.

By the Proposition there is an isomorphism $\phi: U(R)\to V(\fg)$ 
that is compatible with the canonical morphisms
$R\to U(R)$ and $R\to V(\fg), a\mapsto a_{-1}1$.
Thus $\phi$ is an $R\subsmlie$-module isomorphism.

By definition, $\bfg_+=\fg_+$ and $V(\fg)=U(\fg)/U(\fg)\fg_+$.
Since $\fg=\fg_+\oplus\fg_-$, the Poincar\'e-Birkhoff-Witt theorem implies that
$V(\fg)=U(\fg_-)$ as $\fg_-$-modules.
Thus $U(R)=V(\fg)=U(\fg_-)=U(R\subsmlie)$.
Since $U(R)=U(R\subsmlie)$ as filtered $\K[T]$-modules,
we obtain $\rgr U(R)=S\upast R$ from the Poincar\'e-Birkhoff-Witt theorem 
for $U(R\subsmlie)$.
\end{pf}

\subsection{Universal Affine Vertex Algebras}
\label{SS:affine va}

We describe universal affine vertex algebras using the PBW Theorem.

\bigskip

{\bf Proposition.}\: {\it
Let $R$ be a vertex Lie algebra such that $R=R'\oplus\K\hk$ where $R'$ is
a $\K[T]$-submodule and $T\hk=0$.
Then $R'\subset R\subsmlie$ is a Lie subalgebra, $R'\subsmlie$, and
for any $k\in\K$ 
there exists a unique isomorphism of $R'\subsmlie$-modules 
$U(R)/(\hk-k)\to U(R'\subsmlie)$ such that $1\mapsto 1$.
}

\bigskip

\begin{pf}
We have $[\hk,a]\subslie=0$ and $[a,b]\subslie\in TR=TR'\subset R'$ for any
$a, b\in R$.
Thus $R'$ is a Lie subalgebra of $R\subsmlie$ and 
$R\subsmlie=R'\subsmlie\times\K\hk$. 

Since $T(\hk-k)=0$, the element $e:=\hk-k\in U(R)$ is central 
by section \ref{SS:comm va}.
The subspace $U(R)e$ is a $\K[T]$-submodule since $Te=0$.
By Proposition \ref{SS:commut centre}\,\iti\ 
we have $a_t(be)=e(a_t b)$ for any $a, b\in U(R)$ and $t\in\Z$.
Thus $U(R)e$ is the ideal generated by $e$.

By the PBW Theorem from section \ref{SS:PBW thm} 
there is an isomorphism $U(R)\to U(R\subsmlie)$ 
of $R\subsmlie$-modules such that $1\mapsto 1$.

By Remark \ref{SS:prelie}\,\itiii\
we have $(a_1\ldots a_r)e=a_1\ldots a_r e$ for any $a_i\in R$.
Thus the isomorphism $U(R)\to U(R\subsmlie)$ maps $U(R)e$ to the left ideal 
$I\subset U(R\subsmlie)$ generated by $e':=\hk-k$.
Since $e'$ is central, $I$ is a two-sided ideal.

From $R\subsmlie=R'\subsmlie\times\K\hk$ we get
$U(R\subsmlie)=U(R'\subsmlie)\otimes U(\K\hk)$.
Of course, the map $U(R'\subsmlie)\otimes U(\K\hk)\to U(R'\subsmlie),
a\otimes\hk^n\mapsto k^n a$, 
induces an algebra isomorphism $U(R\subsmlie)/I\to U(R'\subsmlie)$.
The resulting isomorphism $U(R)/(\hk-k)\to U(R'\subsmlie)$ is unique
since the $R'\subsmlie$-module $U(R'\subsmlie)$ is generated by $1$.
\end{pf}

\bigskip

We apply the Proposition to an affine vertex Lie algebra $R=\fg[T]\oplus\K\hk$,
see section \ref{SS:aff vlie}.
Thus $R$ is an $\N$-graded vertex Lie algebra and 
$\fg=R_1$ is a Lie algebra with an even invariant symmetric bilinear form.
We have $[a\subla b]=[a,b]+\hk(a,b)\la$ for $a, b\in\fg$.

The {\bf universal affine vertex algebra} of level $k\in\K$ associated to $R$ 
is 
\index{universal affine vertex algebra}
$V:=U(R)/(\hk-k)$.

Remark \ref{SS:finite confa} and section \ref{SS:vertex envelope} imply that
$V$ is the associative unital vertex algebra generated by 
the vector space $\fg$ with relations $[a\subla b]=[a,b]+k(a,b)\la$ for 
$a, b\in\fg$.

Let $\hfg=\fg(R)=\fg[x\uppm]\oplus\K\hk$ be the affine Lie algebra 
associated to $R$.
By Proposition \ref{SS:bracket vlie}
there is a Lie algebra isomorphism $R\subsmlie\to\hfg_-, a\mapsto a_{(-1)}$,
where $\hfg_-=\rspan\set{a_{(t)}\mid a\in R, t<0}$.

The Lie algebra $\hfg$ is a $\Z$-graded Lie algebra.
We have $a_n=a_{(n)}\in\hfg_{-n}$ for $a\in\fg$ and $n\in\Z$.
Furthermore, $\hk=\hk_{(-1)}\in\hfg_0$ and $\hk_{(t)}=0$ for $t\ne -1$.
Since $[a_n,b_m]=[a,b]_{n+m}$ for $n, m<0$, we get 
$\hfg_-=\hfg_>\times\K\hk$ where $\hfg_>:=\bigoplus_{h>0}\hfg_h$.
This decomposition corresponds to the decomposition 
$R\subsmlie=R'\subsmlie\times\K\hk$ in the proof of the Proposition.
In particular, $R'\subsmlie$ is isomorphic to $\hfg_>$.

The Proposition implies that there exists a unique isomorphism of 
$\hfg_>$-modules from the universal affine vertex algebra to 
$U(\hfg_>)$ such that $1\mapsto 1$.

\section{Supplements}
\label{S:suppl env}

\subsection{Duality for the Tensor Vertex Algebra}
\label{SS:duality of tensor alg}

We prove that the tensor vertex algebra $T\upast R$ is dual.

\bigskip

Elements $a, c$ of a vertex algebra $V$ are {\bf dual} if 
\index{dual!elements}
there exists $t\in\Z$ such that for any $b$ the elements $a, b, c$ satisfy 
duality of order $\leq t$.
By definition, $V$ is dual iff any $a, c\in V$ are dual, 
see section \ref{SS:loc and dual}.

For $a\in V$,
define the field $a\upsT(z)$ by $a\upsT(z)b:=\paraab e^{zT}b(-z)a$.

The following lemma is a version of $\bbS_3$-symmetry of the Jacobi identity.

\bigskip

{\bf Lemma.}\: {\it
The elements $a, c$ are dual iff $a(z), c\upsT(z)$ are local.
}

\bigskip

\begin{pf}
This follows from $(a(x)b)(w)c=e^{wT}c\upsT(-w)a(x)b$ and
$$
a(x+w)b(w)c
\; =\;
a(x+w)e^{wT}c\upsT(-w)b
\; =\;
e^{wT}a(x)c\upsT(-w)b.
$$
\end{pf}

\bigskip

The following result is a refinement of Dong's lemma from 
section \ref{SS:dong fields}.
Its proof uses the residue formula for the $n$-th products of fields 
instead of Taylor's formula. 
Thus it provides also an alternative proof of the original Dong's lemma.

\bigskip

{\bf Dong's Lemma.}\: {\it
Let $a(z), b(z), c(z)$ be fields such that $a(z), b(z)$ are weakly local 
\index{Dong's lemma}
and $a(z), c(z)$ and $b(z), c(z)$ are local.
Then $a(z)_n b(z)$ and $c(z)$ are local for any $n\in\Z$.
}

\bigskip

\begin{pf}
Let $m\geq 0$ such that $n+m\geq 0$ and 
$a(z), b(z), c(z)$ are pairwise (weakly) local of order $\leq m$. 
Then
\begin{align}
\notag
&(x-w)^{4m}\,\res_z(z-w)^n[[a(z),b(w)],c(x)]
\\
\notag
=\;
&\sum_i\binom{3m}{i} \res_z(x-w)^m\, (z-w)^{n+i}\, (x-z)^{3m-i}
[[a(z),b(w)],c(x)].
\end{align}
This is equal to $0$.
In fact, the summands for $i\geq 2m$ are $0$ 
since $n+i\geq m$ and $a(z), b(z)$ are weakly local. 
Here we expand $(x-z)^{3m-i}=((x-w)-(z-w))^{3m-i}$ and
commute $\res_z$ with $x-w$.
The summands for $i\leq 2m$ are $0$ since $a(z), c(z)$ and $b(z), c(z)$ are
local and we may apply the radially ordered Leibniz identity.
It is easy to see that in the case $i\leq 2m$ it does not matter that 
$n+i$ may be negative.
\end{pf}

\bigskip

{\bf Proposition.}\: {\it
If $R$ is a vertex Leibniz algebra then $T\upast R$ is dual.
}

\bigskip

\begin{pf}
By induction on $n$, we prove that $a, c$ are dual for any $a\in T^n R$ and 
$c\in T\upast R$.
This is clear for $n=0$.
For $n=1$ it follows from the associativity formula that is satisfied for 
$a\in R$ and $b, c\in T\upast R$ by Proposition \ref{SS:tensor alg vleib}.

The Lemma shows that $a, c$ are dual iff $a(z), c\upsT(z)$ are local. 
Let $a\in R$ and $b, c\in T\upast R$ such that $b(z), c\upsT(z)$ are local.
By Proposition \ref{SS:conseq fund rec}\,\itiii\
the associativity formula for $a, b$ and any $d\in T\upast R$ 
implies that $a(z), b(z)$ are weakly local.
From Dong's lemma follows that $(ab)(z)=a(z)_{-1}b(z)$ and 
$c\upsT(z)$ are local.
That is the induction step.
\end{pf}

\chapter{Representation Theory}
\label{C:reprva}

\section{Zhu Correspondence}
\label{S:zhu corres}

In sections \ref{SS:n grad lie mod} and \ref{SS:zhu corres lie}
we construct the $\fg$-modules $M'(U)$ and $L(U)$ and 
prove that the maps $M\mapsto M_0$ and $U\mapsto L(U)$ 
are mutually inverse bijections between 
simple $\N$-graded $\fg$-modules and simple $\fg_0$-modules.

In sections \ref{SS:zhu prod} and \ref{SS:zhu alg2}
we prove that if $M$ is an $\N$-graded $V$-module then $M_0$ is 
a module over the Zhu algebra $A(V)$.
This is an associative algebra and a quotient of $\fg(V)_0$.

In sections \ref{SS:zhu alg fgv0} and \ref{SS:zhu corresp}
we prove that the maps $M\mapsto M_0$ and $U\mapsto L(U)$ 
are mutually inverse bijections between 
simple $\N$-graded $V$-modules and simple $A(V)$-modules.

In section \ref{SS:ratl vas}
we prove that if $V$ is rational then $A(V)$ is semisimple.

\bigskip

{\bf Conventions.}\:
All vertex algebras in this section are associative and unital.
We always denote by $\fg$ a $\Z$-graded Lie algebra and
by $V$ a $\Z$-graded vertex algebra.

Recall that if $a\in V_h$ then $h_a:=h$.
We denote by $a_{(t)}b$ the $t$-th product, so that 
$Y(a)=\sum a_{(t)}z^{-t-1}$, 
and $a_n:=a_{(h_a-1+n)}$, so that $a_n: V_h\to V_{h-n}$.

\subsection{$\N$-Graded Lie Algebra Modules}
\label{SS:n grad lie mod}

We show that if an $\N$-graded $\fg$-module $M$ is simple then 
the $\fg_0$-module $M_0$ is simple and we discuss when the converse is true.

\bigskip

Let $\fg$ be a $\Z$-graded Lie algebra.
Define $\fg_>:=\bigoplus_{h>0}\fg_h$ and 
define the subalgebras $\fg_{\geq}, \fg_<$, and $\fg_{\leq}$ in the same way.

The $\Z$-gradation of $\fg$ induces a $\Z$-gradation of 
the enveloping algebra $U(\fg)$.
The Poincar\'e-Birkhoff-Witt theorem implies that multiplication
$U(\fg_>)\otimes U(\fg_{\leq})\to U(\fg)$ is an isomorphism of 
left $\fg_>$-modules.
We will often use this fact. 

Let $M$ be an $\N$-graded $\fg$-module. 
Denote by 
$$
\Om(M)
\; :=\; 
\set{a\in M\mid \fg_< a=0}
$$
the space 
\index{singular vector}
of {\bf singular vectors} of $M$.
This is an $\N$-graded $\fg_0$-submodule of $M$ 
since $\fg_<\fg_0 a\subset \fg_0\fg_< a+[\fg_<,\fg_0]a$.
We have $\Om(M)_0=M_0$.

There exists a unique $\N$-graded submodule $J(M)\subset M$ that is
maximal among all $\N$-graded submodules $N$ with $N_0=0$. 
We have $J(M/J(M))=0$.

\bigskip

{\bf Remark.}\: {\it
Let $M$ be an $\N$-graded $\fg$-module.
We have $J(M)=0$ iff $\Om(M)=M_0$ iff for any non-zero $\N$-graded submodule 
$N\subset M$ we have $N_0\ne 0$.
}

\bigskip

\begin{pf}
It is clear that the first and third condition are equivalent.
If $J(M)\ne 0$ then $h:=\min\set{k\mid J(M)_k\ne 0}>0$ and 
$\Om(M)_h\supset J(M)_h\ne 0$. 

Conversely, if $\Om(M)_h\ne 0$ for some $h>0$ then the $\N$-graded submodule 
$N:=U(\fg)\Om(M)_h=U(\fg_>)\Om(M)_h$ is non-zero and $N_0=0$.
\end{pf}

\bigskip

By definition, a module $M$ over an associative algebra is simple
if $M\ne 0$ and $0$ and $M$ are the only submodules of $M$.

An $\N$-graded $\fg$-module $M$ is
\index{simple graded module} 
{\bf simple} if $M_0\ne 0$ and $0$ and $M$ are the only $\N$-graded 
submodules of $M$.

\bigskip

{\bf Proposition.}\: {\it
Let $M$ be an $\N$-graded $\fg$-module.
Then $M$ is simple iff $M_0$ is a simple $\fg_0$-module, $J(M)=0$, and 
$U(\fg)M_0=M$.
In this case $0$ and $M$ are the only submodules of $M$ and 
the gradation is unique.
}

\bigskip

\begin{pf}
`$\Rightarrow$'\:
Let $U\subset M_0$ be a non-zero $\fg_0$-submodule.
Then $U(\fg_>)U=U(\fg)U=M$.
Thus $U=M_0$.

Since $M_0\ne 0$, we have $J(M)=0$ and $U(\fg)M_0=M$.

\smallskip

`$\Leftarrow$'\:
Let $N\subset M$ be a non-zero $\N$-graded submodule. 
Then $N_0\ne 0$.
Thus $N_0=M_0$ and $M=U(\fg)N_0\subset N$.

\smallskip

Suppose that $M$ is simple.
We now prove that $U(\fg)a=M$ for any non-zero $a\in M$.
Let $a=\sum a_k$ where $a_k\in M_k$ and $h:=\max\set{k\mid a_k\ne 0}$.
Because $U(\fg)a_h=M$ and $M_0\ne 0$, 
there exists $b\in U(\fg)_{-h}$ such that $ba_h\in M_0$ is non-zero.
Then $ba=ba_h$ and hence $U(\fg)a\supset U(\fg)ba_h=M$.

From $M_0=\Om(M)$ and $M_h=U(\fg)_h M_0$ follows that the gradation is unique.
\end{pf}

\subsection{Zhu Correspondence for $\Z$-Graded Lie Algebras}
\label{SS:zhu corres lie}

We prove that the map $M\mapsto M_0$ is a bijection between 
simple $\N$-graded $\fg$-modules and simple $\fg_0$-modules.

\bigskip

Let $U$ be a $\fg_0$-module.
Then $U$ is a $\fg_{\leq}$-module via the quotient map $\fg_{\leq}\to\fg_0$.
It is graded with $U_0=U$.
Define the $\Z$-graded $\fg$-module
$$
M'(U)
\; :=\;
U(\fg)\otimes_{U(\fg_{\leq})}U.
$$
It is clear that $U\mapsto M'(U)$ is a functor.
Since $U(\fg)=U(\fg_>)\otimes U(\fg_{\leq})$ as left $\fg_>$-modules, we have 
$M'(U)=U(\fg_>)\otimes U$ as $\Z$-graded left $\fg_>$-modules.
Thus $M'(U)$ is actually $\N$-graded and $M'(U)_0=U$. 

Define the $\N$-graded $\fg$-module
$$
L(U)
\; :=\;
M'(U)/J(M'(U)).
$$
A morphism $U\to U'$ induces a morphism $\al: M'(U)\to M'(U')$
such that $\al J(M'(U))\subset J(M'(U'))$. 
Thus $U\mapsto L(U)$ is a functor from $\fg_0$-modules to 
$\N$-graded $\fg$-modules.
We have $L(U)_0=U, U(\fg)U=L(U)$, and $J(L(U))=0$.

The forgetful functor from $\N$-graded $\fg$-modules to 
$\N$-graded $\fg_{\leq}$-modules has the functor 
$N\mapsto U(\fg)\otimes_{U(\fg_{\leq})}N$ as its left adjoint.
It follows that the functor $M\mapsto M_0$ from $\N$-graded $\fg$-modules to
$\fg_0$-modules has the functor $U\mapsto M'(U)$ as its left adjoint,
since $\Hom_{\fg}(M'(U),M)=\Hom_{\fg_{\leq}}(U,M)=\Hom_{\fg_0}(U,M_0)$.

\bigskip

{\bf Proposition.}\: {\it
Let $\fg$ be a $\Z$-graded Lie algebra and $f$ the functor $M\mapsto M_0$ from
$\N$-graded $\fg$-modules $M$, that satisfy $U(\fg)M_0=M$, to $\fg_0$-modules.

\smallskip 

\iti\:
The functor $L$ is right adjoint to $f$ and fully faithful.

\smallskip 

\itii\:
The functors $L$ and $f$ are mutually inverse equivalences between
the categories of simple $A(V)$-modules and simple $\N$-graded $V$-modules. 
}

\bigskip

\begin{pf}
\iti\: 
Let $M$ be an $\N$-graded $\fg$-module such that $U(\fg)M_0=M$. 
The counit of adjunction $\al:M'(M_0)\to M$ is surjective since $U(\fg)M_0=M$.
The kernel $K$ of $\al$ satisfies $K_0=0$.
Hence $K\subset J(M'(M_0))$. 
Thus the quotient map $M'(M_0)\to L(M_0)$ factorizes over $M$ and
we obtain a morphism $\be: M\to L(M_0)$.

The morphism $\be$ and the inverse of the canonical isomorphism $U\to L(U)_0$
are the unit and counit of adjunction for $L$ and $f$.
Because the counit is an isomorphism,
the right adjoint functor $L$ is fully faithful.

\smallskip

\itii\:
Proposition \ref{SS:n grad lie mod} shows that $L$ and $f$ map
simple objects to simple objects.
The counit of adjunction $L(U)_0\to U$ is always an isomorphism and 
the unit $M\to L(M_0)$ is an isomorphism because it is a non-zero morphism 
between simple modules.
\end{pf}

\bigskip

{\bf Corollary.}\: {\it
Let $U$ be a $\fg_0$-module and $M$ an $\N$-graded $\fg$-module such that 
$M_0=U, U(\fg)U=M$, and $J(M)=0$.
Then there exists a unique isomorphism $M\to L(U)$ compatible with $U$.
}

\bigskip

\begin{pf}
Because $L$ is right adjoint to $M\mapsto M_0$,
the identity map $M_0\to U$ induces a $\fg$-module morphism $\al: M\to L(U)$. 
It is clear that $\al$ is surjective.
The kernel $K$ of $\al$ satisfies $K_0=0$.
Thus $K=0$.
\end{pf}

\subsection{Zhu Algebra}
\label{SS:zhu prod}

We prove that if $M$ is an $\N$-graded $V$-module then $M_0$ is 
a module over the Zhu algebra $A(V)$.
This algebra is a quotient of $\fg(V)_0$.

\bigskip

Let $V$ be a $\Z$-graded vertex algebra.
Then the Borcherds Lie algebra $\fg(V)$ is a $\Z$-graded Lie algebra 
by section \ref{SS:vlie to loclie}.

Recall that $\fg(V)$-modules are the same as $V\subslie$-modules 
by Remark \ref{SS:module vlie}.
They satisfy the commutator formula by Proposition \ref{SS:module vlie}. 
Any $V$-module is a $\fg(V)$-module and
a $\fg(V)$-module $M$ is a $V$-module iff $Y_M(ab)=\normord{Y_M(a)Y_M(b)}$
by Corollary \ref{SS:skew sym fields}\,\itii.
These statements are also true for $\N$-graded modules.

Let $M$ be an $\N$-graded $V$-module.
There is a linear map $V\to\fg(V)_0\to\End(\Om(M)), 
a\mapsto a_0=a_{(h_a-1)}$. 
We want to define a multiplication $\ast$ on $V$ such that this map
is an algebra morphism: $(a\ast b)_0=a_0 b_0$.

Let $c\in\Om(M)$. 
By Corollary \ref{SS:dual skew loc} the Jacobi identity holds for $M$.
In particular, if $a, b\in V$ and $t\geq o(a,c)$ then duality of order 
$\leq t$ holds:
$$
\sum_{i\geq 0}\binom{t}{i}(a_{(r+i)}b)_{(s+t-i)}c
\; =\;
\sum_{i\geq 0}(-1)^i \binom{r}{i}
\,a_{(t+r-i)}b_{(s+i)}c,
$$
see section \ref{SS:loc and dual}.
Since $h_a\geq o(a,c)$, we may take $t=h_a$.
For $s=h_b-1$, the right-hand side is $a_{(h_a+r)}b_{(h_b-1)}c$.
Taking $r=-1$, we get $(a\ast b)_0=a_0 b_0$ where
$$
a\ast b
\; :=\;
\sum_{i\geq 0}
\binom{h_a}{i}a_{(-1+i)}b
$$
is
\index{Zhu product}
the {\bf Zhu product}.
Note that $a\ast b\in\bigoplus_{i\geq 0} V_{h_a+h_b-i}$ and 
$1\ast a=a\ast 1=a$.

The surjective linear map $V\to\fg(V)_0, a\mapsto a_0$, has kernel $(T+H)V$ 
since $\fg(V)_0$ is the quotient of 
$\bigoplus_h V_h\otimes\K x^{h-1}$ by the span of 
$(Ta)\otimes x^{h_a}+h_a a\otimes x^{h_a-1}$.
Because $a_0=a_{(h_a-1)}$, we have  
$$
[a_0,b_0]
\; =\;
\sum\:\binom{h_a-1}{i}\,(a_{(i)}b)_0.
$$

Let $O(V)$ be the ideal of $(V,\ast)$ generated by $(T+H)V$.
The {\bf Zhu algebra} is 
\index{Zhu algebra}
$A(V):=V/O(V)$.
Thus $\Om(M)$ is an $A(V)$-module and there is a surjective linear map
$\fg(V)_0\to A(V), a_0\mapsto [a]$.

\subsection{Zhu Algebra and Affinization}
\label{SS:zhu alg2}

We prove that $A(V)$ is associative and that $O(V)=(T+H)(V)\ast V=V\ast_{-2}V$.
To prove this we use that $(V,\ast)\cong (\hV\subast)_0$,
where $\hV$ is the affinization of $V$.

\bigskip

Recall that the affinization $\hV$ of $V$ is 
the tensor product vertex algebra $V\otimes\K[x\uppm]$, 
see section \ref{SS:tensor va2}.
Define $a_n:=a\otimes x^n$ for $a\in V$ and $n\in\Z$.
The vertex algebra $\hV$ is $\Z$-graded with $a_n\in\hV_{h_a-n}$.
We have $T(a_n)=(Ta)_n+na_{n-1}$ and 
$(a_n)_{(r)}(b_m)=\sum_{i\geq 0}\binom{n}{i}(a_{(r+i)}b)_{n+m-i}$ for any 
$r\in\Z$.

In the following $a_n b$ for $a\in V$ and $b\in\hV$ denotes 
the product $a_n\cdot b$ in $\hV\subast$.
We will avoid the notation $a_n b=a_{(h_a-1+n)}b$ for $a, b\in\hV$.

By Proposition \ref{SS:tensor va}
we can identify $V$ and $\K[x\uppm]$ with vertex subalgebras of $\hV$
such that $V$ and $\K[x\uppm]$ commute and $a_n=ax^n$.
It follows that the commutative vertex algebra $\K[x\uppm]$ is contained in 
the centre of $\hV$.
By Remark \ref{SS:prelie}\,\itiii\ we have $(a_n b_m)x^k=a_n b_m x^k$.
Thus $x^k a_n=a_n x^k=a_{n+k}$.
Quasi-associativity implies that $(xT(a_n))b_m=T(a_n)b_{m+1}$ since
$x_{(i)}b_m=0$ for $i\geq 0$ and since 
$T(a_n)_{(0)}b_m=0$ and $T^{i+1}x=0$ for $i\geq 1$. 

There are linear isomorphisms $\phi: V\to\hV_0, a\mapsto a_{h_a}$, and
$\phi'=x\phi: V\to\hV_{-1}, a\mapsto a_{h_a+1}$.
From the explicit formula for the $t$-th products of $\hV$ we see that
$\phi: (V,\ast)\to(\hV\subast)_0$ is an algebra isomorphism. 
In particular, this shows that $V$ is a pre-Lie algebra with respect to 
the Zhu product.

We have $\phi((T+H)a)=(Ta)_{h_a+1}+h_a a_{h_a}=xT(a_{h_a})=xT\phi a$.
Together with $(xT(a_n))b_m=T(a_n)b_{m+1}$ we obtain
$\phi((T+H)(a)\ast b)=(T\phi a)\phi'b$.
This shows that the image of $(T+H)(V)\ast V$ is $(T\hV_0)\hV_{-1}$.
Moreover, the explicit formula for the $t$-th products of $\hV$ implies that
$(T+H)(a)\ast b=a\ast_{-2}b$ where 
$$
a\ast_n b
\; :=\;
\sum_{i\geq 0}\: \binom{h_a}{i}\, a_{(n+i)}b
$$
is 
\index{Zhu product}
the {\bf $n$-th Zhu product}.
Here $n\in\Z$ and, of course, $\ast_{-1}=\ast$.

By Proposition \ref{SS:q assoc zhu poisson}\,\itii\ 
the subspace $K:=\sum_{i\geq 1}(T^i \hV_0)\hV_{-i}$ is an ideal of 
$(\hV\subast)_0$ and $(\hV\subast)_0/K$ is an associative algebra.
We claim that $K=(T\hV_0)\hV_{-1}$. 
This implies that $O(V)=((T+H)V)\ast V$ and 
that the Zhu algebra $A(V)$ is associative.
To prove our claim, we note that
$$
T^{(i)}(a_{h_a+1})b_{h_b+i-1}
\; =\;
T^{(i)}(a_{h_a})b_{h_b+i}
\; +\;
T^{(i-1)}(a_{h_a})b_{h_b+i-1}
$$
for $i\geq 1$ since $a_{h_a+1}=x a_{h_a}$, $T$ is a derivation, 
$(xT(a_n))b_m=T(a_n)b_{m+1}$, and $T^j x=0$ for $j>1$. 
Therefore $T^{(i)}(a_{h_a})b_{h_b+i}\in(T^{i-1} \hV_0)\hV_{-i+1}$.

The fact that $\Om(M)$ is an $A(V)$-module also follows 
from the isomorphism $\phi: V\to\hV_0$.
Namely, $\K$ is a $\K[x\uppm]$-module by section \ref{SS:modules} and hence
$M=M\otimes\K$ is a module over $\hV=V\otimes\K[x\uppm]$ 
by section \ref{SS:tensor va2}.
We have $(a_n)_{(t)}=a_{(n+t)}$ and thus $(a_{h_a})\subla c=0$ for 
$c\in\Om(M)$. 
Quasi-associativity implies that the map $\hV_0\to\End(\Om(M)), 
a_{h_a}\mapsto (a_{h_a})\cdot=a_{(h_a-1)}$, is an algebra morphism and
$((T^i a_{h_a})b_{h_b+i})c=0$ for $i\geq 1$.

\subsection{Modules over $A(V)$ and $\fg(V)_0$}
\label{SS:zhu alg fgv0}

We show that the map $M\mapsto M_0$ from simple $\N$-graded $V$-modules to 
simple $A(V)$-modules is injective.

\bigskip

In section \ref{SS:zhu prod}
we showed that there is a surjective linear map $\io: \fg(V)_0\to A(V)$
such that $a_0\mapsto [a]$.

\bigskip

{\bf Proposition.}\: {\it
The map $\io: \fg(V)_0\to A(V)$ is a Lie algebra morphism.
}

\bigskip

\begin{pf}
Since there is an algebra isomorphism $A(V)\to\hV_0/K, [a]\mapsto [a_{h_a}]$, 
it suffices to show that 
$[a_{h_a},b_{h_b}]\equiv\sum\binom{h_a-1}{i}(a_{(i)}b)_{h_a+h_b-i-1}$
modulo $K$.

Let $a_n\in\hV_{-i}$.
From $a_n=x^i a_{n-i}$ we get 
$T^{(i)}a_n=\sum_j T^{(j)}(a_{n-i})\binom{i}{j}x^j\equiv a_{n-i}$ 
modulo $K$.
By section \ref{SS:tensor va2} the vertex Lie algebra $(\hV)\subslie$ is 
equal to the affinization of $V\subslie$.
Thus we have 
\begin{align}
\notag
[a_{h_a},b_{h_b}]=
&\int_{-T}^0 d\la\, [(a_{h_a})\subla b_{h_b}]
=\int_{-T}^0 d\la\, (e^{\del_x\del\subla}x^{h_a})\, x^{h_b}\, [a\subla b]
\\
\notag
=
&\sum_{i, j\geq 0}
\int_{-T}^0 d\la\binom{h_a}{i}x^{h_a+h_b-i}\, a_{(j+i)}b\:\la^{(j)}
\\
\notag
\equiv
&\sum_{i, j\geq 0}
(-1)^j\binom{h_a}{i}\, (a_{(j+i)}b)_{h_a+h_b-i-j-1}.
\end{align}
The claim now follows from $\sum_{i=0}^n(-1)^{n+i}\binom{h}{i}=\binom{h-1}{n}$,
which is proven by induction using that 
$-\binom{h-1}{n}+\binom{h}{n+1}=\binom{h-1}{n+1}$.
\end{pf}

\bigskip

{\bf Corollary.}\: {\it
Let $L$ be a conformal vector of $V$. 
Then $L_0$ is contained in the centre of $\fg(V)_0$
and $[L]$ is contained in the centre of $A(V)$.
}

\bigskip

\begin{pf}
We have 
$[L_0,a_0]=[L_{(1)},a_{(h_a-1)}]=(L_{(0)}a)_{(h_a)}+(L_{(1)}a)_{(h_a-1)}=0$, 
since $L_{(0)}=T$ and $L_{(1)}=H$.
The Lemma implies the second claim.
\end{pf}

\bigskip

An $\N$-graded $V$-module $M$ is
\index{simple graded module} 
{\bf simple} if $M_0\ne 0$ and $0$ and $M$ are the only $\N$-graded 
submodules of $M$.

It is clear that $\N$-graded $V$-modules form a full subcategory of
the category of $\N$-graded $\fg(V)$-modules and that these two categories
have the same simple objects. 
The Proposition implies that $A(V)$-modules form a full subcategory of
the category of $\fg(V)_0$-modules and that these two categories
have the same simple objects.
In particular, for any $A(V)$-module $U$ the $\N$-graded $\fg(V)$-module $L(U)$
is defined.

By section \ref{SS:zhu prod}
the functor $M\mapsto M_0$ sends $\N$-graded $V$-modules to $A(V)$-modules.
Proposition \ref{SS:zhu corres lie}\,\itii\ implies that 
this functor induces an injective map from isomorphism classes of
simple $\N$-graded $V$-modules to isomorphism classes of
simple $A(V)$-modules.

\subsection{Zhu Correspondence}
\label{SS:zhu corresp}

We prove that the maps $M\mapsto M_0$ and $U\mapsto L(U)$
are mutually inverse bijections between 
simple $\N$-graded $V$-modules and simple $A(V)$-modules.

\bigskip

{\bf Remark.}\: {\it
Let $M$ be an $\N$-graded $\fg$-module and $p: M\to M_0$ the projection. 
Then
$$
J(M)
\; =\;
\set{a\in M\mid  p(U(\fg)a)=0}.
$$
}

\bigskip

\begin{pf}
Let $J$ be the right-hand side.
It is clear that $J(M)\subset J$ and that $J$ is a submodule. 
Because $M$ is an $\N$-graded module, 
it follows that $J$ is an $\N$-graded submodule.
Since $J_0=0$, we get $J\subset J(M)$.
\end{pf}

\bigskip

{\bf Proposition.}\: {\it
Let $U$ be an $\fg(V)_0$-module.
Then $U$ is an $A(V)$-module iff $L(U)$ is a $V$-module.
}

\bigskip

\begin{pf}
If $L(U)$ is a $V$-module then $U=L(U)_0$ is an $A(V)$-module by 
section \ref{SS:zhu prod}.
Suppose that $U$ is an $A(V)$-module.

Let $M:=M'(U)$ and $\fg:=\fg(V)$.
By the Remark we have to prove that 
$$
p(g(ab)(z)c)
\; =\;
p(g\hskip -2pt \normord{a(z)b(z)}c)
$$
for any $a, b\in V, c\in M$, and $g\in U(\fg)$. 
First, we consider the claim for $g=1$.

By Proposition \ref{SS:module vlie} 
the $V\subslie$-module $M$ satisfies the commutator formula.
Below we will prove in two steps that $p$ applied to duality holds:
\begin{equation}
\tag{$\ast$}
(w+x)^t\, p((a(x)b)(w)c)
\; =\;
(x+w)^t\, p(a(x+w)b(w)c)
\end{equation}
for any $a, b\in V, c\in M$, and $t\gg 0$.
The commutator formula and ($\ast$) imply that 
$p$ applied to the Jacobi identity is satisfied.
Indeed, this follows from Proposition \ref{SS:fund recur}
in the same way as in the proof of Proposition \ref{SS:conseq fund rec},
the only difference is that now we apply $p$.
Thus the claim is proven for $g=1$, since it is a special case of 
$p$ applied to the Jacobi identity.

Since $\fg$ generates $U(\fg)$,
it suffices to prove that if the claim holds for $g$ then it
holds for $ge_{(n)}$ where $e\in V$ and $n\in\Z$.

Let $r\in\Z$ such that $(ab)(z), \normord{a(z)b(z)}$, and $e(z)$
are pairwise local of order $\leq r$.
We have 
\begin{align}
\notag
(z-w)^r p(g e(z)(ab)(w)c)
\; &=\;
(z-w)^r p(g(ab)(w)e(z)c)
\\\
\notag
&=\;
(z-w)^r p(g\hskip -2pt \normord{a(w)b(w)}e(z)c)
\\
\notag
&=\;
(z-w)^r p(g e(z)\hskip -2pt \normord{a(w)b(w)}c).
\end{align}
Since $(ab)(z)c$ and $\normord{a(z)b(z)}c$ are in $M\lau{z}$, 
we may divide by $(z-w)^r$.
Thus the claim holds for $ge_{(n)}$.
To finish the proof, we need to show ($\ast$).

\smallskip

{\it 1st Step.}\:
We prove ($\ast$) in the special case $c\in U$.
Equation ($\ast$) holds iff for any $r\in\Z$ we have
\begin{equation}
\tag{$\ast_r$}
\res_x x^r (w+x)^t p((a(x)b)(w)c)
\; =\;
\res_x x^r (x+w)^t p(a(x+w)b(w)c).
\end{equation}

Because $U$ is an $A(V)$-module, we have
\begin{equation}
\notag
\sum_{i\geq 0}
\binom{h_a}{i}\,
p((a_{(-1+i)}b)_{(h_a+s-i)}c)
\; =\;
\sum_{i\geq 0}
(-1)^{i} \binom{-1}{i}\,
p(a_{(h_a-1-i)}b_{(s+i)}c)
\end{equation}
for $s=h_b-1$, and for degree reasons this also holds for any $s\ne h_b-1$.
This equation is $p$ applied to duality of order $\leq h_a$ 
for indices $r=-1$ and any $s\in\Z$. 
Thus we obtain ($\ast_r$) for $r=-1$ and $t=h_a$.

Since $M$ satisfies the commutator formula,
the Jacobi identity holds for indices $r\geq 0$ and $s, t\in\Z$. 
In particular, duality of order $\leq h_a$ holds for indices 
$r\geq 0$ and $s\in\Z$.
Hence ($\ast_r$) is satisfied for $r\geq -1$ and $t=h_a$.
By Lemma \ref{SS:fund recur} we obtain ($\ast_r$) for $r\geq -1$ and 
$t\geq h_a$.

The integration-by-parts formula yields
\begin{align}
\notag
&\res_x x^r (w+x)^{t+1} p(((Ta)(x)b)(w)c)
\\
\notag
=\;
-&\res_x (r x^{r-1}w(w+x)^t
+ rx^r(w+x)^t 
+(t+1)x^r(w+x)^t)
\\
\notag
&\qquad\qquad\qquad\qquad\qquad\qquad\qquad\qquad\qquad\quad
p((a(x)b)(w)c).
\end{align}
Since $(Ta)(x+w)= e^{w\del_x}(Ta)(x)=e^{w\del_x}\del_x a(x)=\del_x a(x+w)$,
we also have
\begin{align}
\notag
&\res_x x^r (x+w)^{t+1} p((Ta)(x+w)b(w)c)
\\
\notag
=\;
-&\res_x (rx^{r-1}w(x+w)^t + rx^r (x+w)^t + (t+1)x^r (x+w)^t)
\\
\notag
&\qquad\qquad\qquad\qquad\qquad\qquad\qquad\qquad\qquad
p(a(x+w)b(w)c).
\end{align}
By induction, we see that ($\ast_r$) is satisfied for any $r\in\Z$ and 
$t\geq h_a$.

\smallskip

{\it 2nd Step.}\:
We prove ($\ast$) in general.
By the first step, ($\ast$) holds for $c\in U$.
Since $M=U(\fg_>)U$, it thus suffices to prove that if ($\ast$) holds for
$c\in M$ then ($\ast$) also holds for $e_{(n)}c$ where $e\in V$ and $n<h_e-1$.

In the proof we apply the commutator formula four times and use ($\ast$) for 
$e_{(j)}a, b, c$ and for $a, e_{(i)}b, c$.
Since $p(e_{(n)}M)=0$, we obtain
\begin{alignat}{2}
\notag
&&&(w+x)^t p((a(x)b)(w)e_{(n)}c)
\\
\notag
=&
-&&(w+x)^t p\sum_{i\geq 0}\binom{n}{i}w^{n-i}(e_{(i)}a(x)b)(w)c
\\
\notag
=&
-&&(w+x)^t p\sum_{i\geq 0}\binom{n}{i}w^{n-i} \Big(
(a(x)e_{(i)}b)(w)c
\\
\notag
&&&\qquad\qquad\qquad\qquad\qquad\;\;\;
+\sum_{j\geq 0}\binom{i}{j}x^{i-j}((e_{(j)}a)(x)b)(w)c\Big)
\\
\notag
=&
-&&(w+x)^t p\Big( 
\sum_{i\geq 0}\binom{n}{i}w^{n-i} (a(x)e_{(i)}b)(w)c
\\
\notag
&&&\qquad\qquad\qquad
+\sum_{j\geq 0}\binom{n}{j}(w+x)^{n-j}((e_{(j)}a)(x)b)(w)c\Big)
\\
\notag
=&
-&&(x+w)^t p\Big( 
\sum_{i\geq 0}\binom{n}{i}w^{n-i} a(x+w)(e_{(i)}b)(w)c
\\
\notag
&&&\qquad\qquad\qquad
+\sum_{j\geq 0}\binom{n}{j}(x+w)^{n-j}(e_{(j)}a)(x+w)b(w)c\Big)
\\
\notag
=&&&
(x+w)^t p(a(x+w)b(w)e_{(n)}c)
\end{alignat}
where in the third step we used
\begin{align}
\notag
&\sum_{i,j\geq 0}\binom{n}{i}\binom{i}{j}w^{n-i}x^{i-j}
=
\sum_{i,j\geq 0}\binom{n}{i+j}\binom{i+j}{j}w^{n-i-j}x^i
\\
\notag
=
&\sum_{i,j\geq 0}\binom{n}{j}\binom{n-j}{i}w^{n-i-j}x^i
=
\sum_{j\geq 0}\binom{n}{j}(w+x)^{n-j}.
\end{align}
\end{pf}

\bigskip

{\bf Theorem.}\: {\it
Let $V$ be a $\Z$-graded vertex algebra.
The functors $M\mapsto M_0$ and $L$ are mutually inverse equivalences between
the categories of simple $\N$-graded $V$-modules and simple $A(V)$-modules. 
}

\bigskip

\begin{pf}
Because of the remarks at the end of section \ref{SS:zhu alg fgv0},
this follows from Proposition \ref{SS:zhu corres lie}\,\itii\
and the Proposition.
\end{pf}

\subsection{Rational Vertex Algebras}
\label{SS:ratl vas}

We prove that if $V$ is rational then $A(V)$ is semisimple.

\bigskip

A $\Z$-graded vertex algebra $V$ is 
\index{rational vertex algebra}
{\bf rational}
if for any $\N$-graded $V$-module $M$ there are simple $\N$-graded $V$-modules
$M_i$ and an isomorphism $M\to\bigoplus M_i$ of $V$-modules.

\bigskip

{\bf Proposition.}\: {\it
Let $V$ be a rational $\Z$-graded vertex algebra.
Then:

\smallskip

\iti\:
The Zhu algebra $A(V)$ is semisimple.
If $V$ is even then $V$ has only finitely many simple $\N$-graded modules
up to isomorphism.

\smallskip

\itii\:
If $\K=\C$, $V$ is even, and $\dim V_h<\infty$ for any $h$ then
$A(V)\cong\End(\C^{n_1})\times\ldots\times\End(\C^{n_r})$ for some unique
$n_i\geq 1$.
}

\bigskip

\begin{pf}
\iti\:
Let $U$ be an $A(V)$-module.
Then $L(U)\cong\bigoplus M_i$ for some simple $\N$-graded $V$-modules $M_i$.
The Remark implies $U=\Om(L(U))\cong\Om(\bigoplus M_i)=\bigoplus\Om(M_i)$.
We have $\Om(M_i)=(M_i)_0$ by Remark \ref{SS:zhu corresp}.
These are simple $A(V)$-modules by Theorem \ref{SS:zhu corresp}\,\itii.
Hence $A(V)$ is semisimple.

It is well-known that an even semisimple associative algebra has only
finitely many simple modules up to isomorphism. 
Thus the second claim follows from Theorem \ref{SS:zhu corresp}\,\itii.

\smallskip

\itii\:
The Wedderburn-Artin theorem implies that $A(V)$ is a product of
matrix algebras of division algebras over $\C$. 
Any division algebra $\ne\C$ contains the rational function field $\C(x)$
and hence has uncountable dimension, since the rational functions
$1/(x-a)$ for $a\in\C$ are linearly independent.
Thus the claim follows from the fact that $A(V)=V/O(V)$ has 
at most countable dimension since $\dim V_h<\infty$.
\end{pf}

\bigskip

The proof shows that if $\K=\C$ and $A$ is an even 
semisimple associative algebra then $\dim A<\infty$ iff 
$A(V)\cong\End(\C^{n_1})\times\ldots\times\End(\C^{n_r})$.

\section{Supplements}
\label{S:suppl reprva}

In section \ref{SS:gen verma vertex}
we construct a functor from $A(V)$-modules to $\N$-graded $V$-modules
that is left adjoint to $M\mapsto M_0$.

In sections \ref{SS:assoc zhu prod}--\ref{SS:commut zhu revisi}
we give a second proof of the facts that $O(V)$ is an ideal,
that $A(V)$ is associative, and that the map $\fg(V)_0\to A(V), a_0\mapsto a$, 
is a Lie algebra morphism.

In section \ref{SS:zhu affine lie} we compute the Zhu algebra
of universal affine vertex algebras.

\subsection{Verma $V$-Modules}
\label{SS:gen verma vertex}

We construct a functor from $A(V)$-modules to $\N$-graded $V$-modules
that is left adjoint to $M\mapsto M_0$.

\bigskip

Let $M$ be a $\fg(V)$-module.
Denote by $Q(M)$ the quotient of $M$ by the $\fg(V)$-submodule generated by 
the coefficients of $(ab)(z)c-\normord{a(z)b(z)}c$ for any $a, b\in V$ and 
$c\in M$. 
Corollary \ref{SS:skew sym fields}\,\itii\ implies that
$Q(M)$ is a $V$-module and the functor $M\mapsto Q(M)$ 
from $\fg(V)$-modules to $V$-modules is left adjoint to the
forgetful functor. 
The same is true for $\N$-graded modules over $\fg(V)$ and $V$.

By section \ref{SS:zhu corres lie} the functor $U\mapsto M'(U)$ 
from $\fg(V)_0$-modules to $\N$-graded $\fg(V)$-modules is left adjoint to 
the functor $M\mapsto M_0$.
By section \ref{SS:zhu alg fgv0} $A(V)$-modules form a full subcategory of 
the category of $\fg(V)_0$-modules.
Since the composition of left adjoint functors is left adjoint,
the functor $U\mapsto Q(M'(U))$ from $A(V)$-modules to $\N$-graded $V$-modules
is left adjoint to the functor $M\mapsto M_0$.

\subsection{Associativity Formula for the Zhu Products}
\label{SS:assoc zhu prod}

We prove an associativity formula for the Zhu products.
We use it in sections \ref{SS:zhu alg} and \ref{SS:zhu affine lie}
to prove results about $A(V)$.

\bigskip

Recall that the $n$-th Zhu product is 
$$
a\ast_n b
\; =\;
\sum_{i\geq 0}\: \binom{h_a}{i}\, a_{(n+i)}b
\; =\;
\res_z z^n(1+z)^{h_a}a(z)b.
$$

\bigskip

{\bf Proposition.}\: {\it
For $a, b, c\in V$ and $r, s\in\Z$, we have
\begin{align}
\notag
&(a\ast_r b)\ast_s c
\,=\,
\sum_{i, j\geq 0}
(-1)^i \binom{-r-1}{j}\binom{r}{i}
\\
\notag
&\qquad\qquad\qquad\qquad
(a\ast_{r-i}(b\ast_{s+i+j}c)-\paraab (-1)^r\, b\ast_{s+r-i+j}(a\ast_{i}c)).
\end{align}
}

\bigskip

\begin{pf}
For $a(z)\in V\pau{z\uppm}$, we have
\begin{align}
\notag
a(z+(w+x))
=
e^{(w+x)\del_z}a(z)
=
e^{w\del_z}e^{x\del_z}a(z)
&=
e^{w\del_z}a(z+x)
\\
\notag
&=
a((z+w)+x).
\end{align}
Thus we get
$$
(1+z)^h
\; =\;
(1+(w+z-w))^h
\; =\;
((1+w)+(z-w))^h
$$
for any $h\in\Z$.
The Jacobi identity in terms of residues with 
$F(z,w,x)=(1+w)^{h_a+h_b-r-i-1}w^s x^{r+i}$,
see Proposition \ref{SS:h Jac id res}, yields
\begin{align}
\notag
&(a\ast_r b)\ast_s c
\\
\notag
=\;
&\sum_{i\geq 0}
\binom{h_a}{i}
\res_{w,x}
(1+w)^{h_a+h_b-r-i-1} w^s x^{r+i} (a(x)b)(w)c
\\
\notag
=\;
&\sum_{i\geq 0}
\binom{h_a}{i}
\res_{z,w}
(1+w)^{h_a-i+h_b-r-1} w^s (z-w)^{r+i} [a(z),b(w)]c
\\
\notag
=\;
&\sum_{j\geq 0}
\binom{-r-1}{j}
\res_{z, w}
(1+z)^{h_a}(1+w)^{h_b} w^{s+j}(z-w)^r[a(z),b(w)]c
\\
\notag
=\;
&\sum_{i, j\geq 0}
(-1)^i \binom{-r-1}{j}\binom{r}{i}
\\
\notag
&\qquad\qquad\qquad\qquad
(a\ast_{r-i}(b\ast_{s+j+i}c)-(-1)^r\, b\ast_{s+j+r-i}(a\ast_i c)).
\end{align}
\end{pf}

\bigskip

We note that only the associativity formula of $V$ is used in the proof
because $F(z,w,x)\in\K\lau{w,x}$.

\subsection{Zhu Algebra. Revisited}
\label{SS:zhu alg}

We use the associativity formula for $\ast_n$ to prove 
that $O(V)=V\ast_{-2}V$ is an ideal and $A(V)$ is associative.
This was proven in section \ref{SS:zhu alg2} using $\hV$.

\bigskip

{\bf Remark.}\: {\it
For $a, b\in V$ and $n\in\Z$, we have
$$
(T+H+n+1)(a)\ast_n b
\; =\;
-n\, a\ast_{n-1}b.
$$
In particular, $V\ast_{n-1}V\subset V\ast_n V$ if $n\ne 0$ and 
$(T+H)(a)\ast b=a\ast_{-2} b$.
}

\bigskip

\begin{pf}
We have
\begin{align}
\notag
Ta\ast_n b 
\; &= \;
\res_z (1+z)^{h_a+1}z^n (Ta)(z)b
\\
\notag
&=\;
-\res_z\del_z((1+z)^{h_a+1}z^n)a(z)b
\\
\notag
&=\;
-\res_z((h_a+1)(1+z)^{h_a}z^n+n(1+z)(1+z)^{h_a}z^{n-1})a(z)b
\\
\notag
&=\;
-(h_a+1) a\ast_n b
\; -\;
n\, a\ast_{n-1}b
\; -\;
n\, a\ast_n b.
\end{align}
\end{pf}

\bigskip

{\bf Proposition.}\: {\it
The subspace $V\ast_{-2}V$ is a two-sided ideal of $(V,\ast)$
and the quotient $A(V)=V/(V\ast_{-2}V)$ is associative. 
}

\bigskip

\begin{pf}
By the Remark we have $V\ast_n V\subset V\ast_{-2}V$ for $n\leq -2$.
Thus Proposition \ref{SS:assoc zhu prod} implies
$$
(a\ast_r b)\ast_s c
\; -\;
\sum_{j=0}^{-r-1}
\binom{-r-1}{j} a\ast_r (b\ast_{s+j}c)
\;\in\;
V\ast_{-2}V
$$
for any $r, s<0$, since $s+r+j-i\leq s-1-i\leq -2$.
Taking $(r,s)=(-2,-1)$ and $(r,s)=(-1,-2)$ 
we see that $V\ast_{-2}V$ is a right and a left ideal.
Taking $(r,s)=(-1,-1)$ we see that $A(V)$ is associative.
\end{pf}

\subsection{Commutator of the Zhu Algebra}
\label{SS:commut zhu revisi}

We give a direct proof that the map $\fg(V)_0\to A(V), a_0\mapsto a$, is 
a Lie algebra morphism.
This was proven in Lemma \ref{SS:gen verma vertex} using the affinization 
$\hV$.

\bigskip

Let $E$ be a vector space, $a(z)\in E\lau{z}$, and 
$f(z)\in z\K\pau{z}, f(z)\ne 0$.
Then the sum 
$a(f(z)):=\sum_{n\in\Z} a_n f(z)^{-n-1} \in E\lau{z}$
is summable.

\bigskip

{\bf Remark.}\: {\it
Let $a(z)\in E\lau{z}$ and
$f(z)\in z\K\pau{z}$ such that $\del_z f(0)\ne 0$.
Then
$$
\res_z a(z)
\; =\;
\res_z a(f(z))\del_z f(z).
$$
}

\bigskip

\begin{pf}
For $n\in\Z, n\ne -1$, we have
$$
\res_z f(z)^n \del_z f(z)
\; =\;
\frac{\res_z \del_z(f(z)^{n+1})}{n+1}
\; =\;
0.
$$
Because
$f(z)^{-1}\in \del_z f(0)^{-1}z^{-1}+\K\pau{z}$ we get
$\res_z a_0 f(z)^{-1} \del_z f(z)=a_0$.
\end{pf}

\bigskip

{\bf Proposition.}\: {\it
The following identity holds in $A(V)$:
$$
a\ast b-\paraab\, b\ast a
\; =\;
\sum_{i\geq 0}\binom{h_a-1}{i}a_{(i)}b.
$$
}

\bigskip

\begin{pf}
For $c\in V$, we have $e^{zT}c=(1+z)^{-h_c}c$ in $A(V)\pau{z}$ because 
using $(T+H)V\subset V\ast_{-2}V$ we have by induction
$$
T^{(n)}c
=
-H T^{(n-1)}c/n
\equiv
-(h_c+n-1)\binom{-h_c}{n-1}c/n
=
\binom{-h_c}{n}c.
$$
Thus we obtain
\begin{align}
\notag
b(z)a
\; &=\;
\sum_{n\in\Z}
e^{zT}(a_{(n)}b)\: (-z)^{-n-1}
\\
\notag
\; &=\;
\sum_{n\in\Z}
(1+z)^{-h_a-h_b+n+1}(a_{(n)}b)\: (-z)^{-n-1}
\\
\notag
\; &=\;
(1+z)^{-h_a-h_b}\, a(-z(1+z)^{-1})b.
\end{align}
Let $w:=-z(1+z)^{-1}$. 
Then $1+w=(1+z)^{-1}$ and thus $z=-w(1+w)^{-1}$.
Moreover,
$\del_z(-z(1+z)^{-1})=(-(1+z)-(-z))(1+z)^{-2}=
-(1+z)^{-2}$.
The Remark yields
\begin{align}
\notag
b\ast a
\; &=\;
\res_z (1+z)^{h_b}z^{-1}b(z)a
\\
\notag
&=\;
\res_z (1+z)^{-h_a}z^{-1}a(-z(1+z)^{-1})b
\\
\notag
&=\;
\res_w (1+w)^{h_a-1}w^{-1}a(w)b.
\end{align}
Thus we get
\begin{align}
\notag
a\ast b-b\ast a
\; &=\;
\res_z((1+z)^{h_a}-(1+z)^{h_a-1})z^{-1}a(z)b
\\
\notag
&=\;
\res_z(1+z)^{h_a-1}a(z)b.
\end{align}
\end{pf}

\subsection{Zhu Algebra of Universal Affine Vertex Algebras}
\label{SS:zhu affine lie}

We compute the Zhu algebra of universal affine vertex algebras.

\bigskip

Let $R=\fg[T]\oplus\K\hk$ be an affine vertex Lie algebra.
The universal affine vertex algebra of level $k\in\K$ associated to $R$
is the quotient $U(R)/(\hk-k)$, see section \ref{SS:affine va}.

The affine Lie algebra $\hfg=\fg(R)=\fg[x\uppm]\oplus\K\hk$ 
is a $\Z$-graded Lie algebra.
Note that $a_n\in\hfg_{-n}$.
As usual, we define $\hfg_>:=\bigoplus_{h>0}\hfg_h$ and 
$\hfg_{\geq}:=\hfg_>\oplus\hfg_0$.

\bigskip

{\bf Proposition.}\: {\it
Let $V$ be the universal affine vertex algebra of level $k\in\K$
associated to an affine vertex Lie algebra $R=\fg[T]\oplus\K\hk$.
Then the map $\fg\to V, a\mapsto a_{-1}1$, induces an algebra isomorphism 
$U(\fg)\to A(V)$. 
}

\bigskip

\begin{pf}
Let $a, a', a^i\in\fg$ and $b\in V$.
Since $h_a=1$, we have
\begin{equation}
\tag{$\ast$}
a\ast_n b
\; =\;
\sum_{i\geq 0}\binom{h_a}{i}\: a_{(n+i)}b
\; =\;
a_n b
\; +\;
a_{n+1}b
\end{equation}
for $n\in\Z$.
By Remark \ref{SS:zhu alg} we have $V\ast_n V\subset V\ast_{-2}V$ for 
$n\leq -2$.
Thus $N:=\rspan\set{a_n b+a_{n+1}b\mid a\in\fg, b\in V, n\leq -2}\subset 
V\ast_{-2}V$.
The subspace $N$ is a $\hfg_{\geq}$-submodule because
$$
a_n(a'_m+a'_{m+1})b
\; =\;
(a'_m+a'_{m+1})a_n b
\; +\;  
([a,a']_{n+m}+[a,a']_{n+m+1})b.
$$

We now prove that $N=V\ast_{-2}V$.
We show that $b\ast_n c\in N$ for any $n\leq -2$ and $c\in V$
by induction on $h=h(b)$.
The case $h=0$ follows from $1\ast_n c=0$.
Suppose the claim is true for $h$.
Let $a\in\fg, b\in V$, and $m\leq -1$ such that $h(a_m b)=h+1$.
By ($\ast$) we have
$$
(a_m b)\ast_n c
\; =\;
(a\ast_m b)\ast_n c
\; -\;
(a_{m+1}b)\ast_n c.
$$
By induction, $(a_{m+1}b)\ast_n c\in N$.
Proposition \ref{SS:assoc zhu prod} yields 
\begin{align}
\notag
&(a\ast_m b)\ast_n c
\,=\,
\sum_{i, j\geq 0}
(-1)^i \binom{-m-1}{j}\binom{m}{i}
\\
\notag
&\qquad\qquad\qquad\qquad
(a\ast_{m-i}(b\ast_{n+j+i}c)-\paraab (-1)^m\, b\ast_{n+j+m-i}(a\ast_{i}c)).
\end{align}
By induction, $b\ast_{n+j+m-i}(a\ast_{i}c)\in N$ since $j\leq -m-1$.
Equation ($\ast$) shows that 
$a\ast_{m-i}(b\ast_{n+j+i}c)\in N$ if $m-i\leq -2$,
i.e.~if $m\leq -2$ or $i>0$. 
If $m=-1$ and $i=0$ then $j=0$.
By induction, $b\ast_n c\in N$.
Thus $a\ast_{-1}(b\ast_n c)\in N$ follows from ($\ast$) 
and the fact that $N$ is a $\hfg_{\geq}$-submodule.
The proof that $N=V\ast_{-2}V$ is complete.

Proposition \ref{SS:commut zhu revisi} and ($\ast$) imply that
$$
b\ast a
\; =\; 
a\ast b-\sum_{i\geq 0}\binom{h_a-1}{i}a_{(i)}b
\; =\; 
a\ast b-a_0 b
\; =\; 
a_{-1}b
\; \in\;  
A(V).
$$
We get
\begin{equation}
\tag{$\dagger$}
a^1\ast\dots\ast a^r 
\; =\;
a^r_{-1}\dots a^1_{-1}1.
\end{equation}

The map $\fg\to A(V), a\mapsto a$, induces an algebra morphism $\al$
from the tensor algebra of $\fg$ to $A(V)$.
Equation ($\dagger$) implies that
$$
\al:\quad
a^1\dots a^r
\;\mapsto\;
a^r_{-1}\dots a^1_{-1}1.
$$
From ($\dagger$) and $N\subset V\ast_{-2}V$ we obtain
$$
\al(aa'-a'a)=
[a'_{-1},a_{-1}]1=
[a',a]_{-2}1=
[a,a']_{-1}1=
\al([a,a']).
$$
Thus $\al$ induces a morphism $\al:U(\fg)\to A(V)$.

By section \ref{SS:affine va} we have an isomorphism $V\cong U(\hfg_>)$
of $\hfg_>$-modules.
The map $\be: \hfg_>\to\fg, a_n\mapsto (-1)^{n+1}a$, is an anti-morphism
of Lie algebras since 
$\be[a_n,a'_m]=(-1)^{n+m+1}[a,a']=-[\be a_n,\be a'_m]$.
Thus $\be$ induces an algebra anti-morphism $U(\hfg_>)\to U(\fg)$.
Hence we get a linear map
$$
\be:\:
V\to U(\fg),
\quad
a^1_{n_1}\dots a^r_{n_r}1
\; \mapsto \;
(-1)^{r+\sum_i n_i}\, a^r \dots a^1.
$$
It is clear that $\be$ induces a linear map 
$\be: A(V)=V/N\to U(\fg)$ such that $\be\circ\al=\id$. 
Since $N$ is a $\hfg_<$-submodule, we have 
$$
a^1_{n_1}\dots a^r_{n_r}1
\; = \;
(-1)^{r+\sum_i n_i}\, a^1_{-1}\dots a^r_{-1}1
\;\in\; 
V/N.
$$
This shows that $\al\circ\be=\id_{A(V)}$.
\end{pf}

\appendix
\chapter{Superalgebra}
\label{C:supalg}

\section{Superalgebra}
\label{S:supalg}

All {\bf vector spaces} are
\index{vector space}
assumed to be {\it super} vector spaces, that is, they are $\Z/2$-graded.
The $\Z/2$-gradation is written $E=E\even\oplus E\odd$.

The space $\Hom(E,F)$ of linear maps $E\to F$ is a super vector space
with $\Hom(E,F)_p:=\bigoplus_{q\in\Z/2}\Hom(E_q,F_{q+p})$ for $p\in\Z/2$.

A {\bf super vector subspace} of a super vector space $E$ is
\index{super vector subspace}
a super vector space $F$ together with an even monomorphism $\io: E\to F$.
We shall often identify a super vector subspace $(F,\io)$ with the
super vector space $\io F$ endowed with the inclusion $\io F\subset E$.

The {\bf parity-change functor} $\rPi$ is the endofunctor of 
\index{parity-change functor}
the supercategory of vector spaces defined by $(\rPi E)_p:=E_{p+\bar{1}}$.

Since $\Z\subset\K$, we can identify $\Z$-gradations with certain
$\K$-gradations by defining $V_h:=0$ for $h\notin\Z$.

A subset $S\subset V$ is {\bf homogeneous} if $S\subset\bigoplus (S\cap V_h)$.
It is {\bf graded} if $S\subset\bigcup (S\cap V_h)$.

\bigskip

{\bf Proposition.}\: {\it
Let $V$ be a $\K$-graded vector space with Hamiltonian $H$.
A subspace $E\subset V$ is homogeneous iff $HE\subset E$.
}

\bigskip

\begin{pf}
Suppose that $HE\subset E$.
Let $a=\sum a_h\in E$ with $a_h\in V_h$.
From $H^i a=\sum h^i a_h, i=1, \dots, \#\set{h\ne 0\mid a_h\ne 0}$,
and nonvanishing of the Vandermonde determinant follows that
$a_h$ for $h\ne 0$ is a linear combination of $H^i a\in E$.
From $a\in E$ we also get $a_0\in E$.
The converse is clear.
\end{pf}

\bigskip

A {\bf multiplication} on a vector space $V$ is
\index{multiplication}
an even linear map $V\otimes V\to V$.
The image of $a\otimes b$ is denoted by $ab$ or $[a,b]$.
In the latter case the multiplication is 
\index{bracket}
called a {\bf bracket}.
An {\bf algebra} is 
\index{algebra}
a vector space with a multiplication.
The {\bf opposite} algebra $V\upop$ of 
\index{opposite!algebra}
an algebra $V$ 
is the vector space $V$ with 
\index{opposite!multiplication}
the {\bf opposite} multiplication $a\otimes b\mapsto\paraab ba$.

A multiplication is {\bf commutative} if 
\index{commutative!multiplication}
$ab=\paraab ba$ and {\bf skew-symmetric} if 
\index{skew-symmetric!multiplication}
$ab=-\paraab ba$.
An algebra is {\bf unital} if 
\index{unital!algebra}
there exists an {\bf identity} $1$  such that $1a=a1=a$.
A {\bf commutative algebra} is
\index{commutative!algebra}
a unital algebra with commutative and associative multiplication.

For example, if $S$ is a set then the polynomial ring $\K[S]$
is the associative algebra generated by $S$ with relations $[a,b]=0$. 
In other words, it is the free commutative algebra generated by $S$.
In particular, 
the polynomial ring $\K[z_1, \dots, z_n, \zeta_1, \dots, \zeta_m]$
in the even 
\index{even!variable}
variables $z_i$ and the odd 
\index{odd!variable}
variables $\zeta_j$ is $\K[S]$
with $S\even=\set{z_1, \dots, z_n}$ and $S\odd=\set{\zeta_1, \dots, \zeta_m}$.

\bigskip

{\bf Remark.}\: {\it
Let $V$ be a $\K$-graded algebra and $1$ an identity.
Then $1\in V_0$.
}

\bigskip

\begin{pf}
Let $1=\sum 1_h$ with $1_h\in V_h$. 
Then $1_h=1_h 1=\sum_k 1_h 1_k$ and hence $1_h 1_k=0$ for $k\ne 0$.
Thus $1_h=1 1_h=\sum_k 1_k 1_h=0$ for $h\ne 0$.
\end{pf}

\bigskip

A {\bf bilinear form} on a vector space $E$ is 
\index{bilinear form}
an {\it even} linear map $E\otimes E\to\K$.
Equivalently, a bilinear form on $E$ is a pair of bilinear forms 
$E\even\otimes E\even\to\K$ and $E\odd\otimes E\odd\to\K$.

A bilinear form is {\bf symmetric} if $(a,b)=\paraab(b,a)$ and 
{\bf skew-symmetric} if $(a,b)=-\paraab(b,a)$.
A {\bf symplectic vector space} is 
\index{symplectic vector space}
an even vector space with a non-degenerate skew-symmetric bilinear form.

\section{Binomial Coefficients}
\label{S:bin coeff}

Let $a$ be an element of an associative unital algebra and $n\in\K$. 
The 
\index{divided power}
divided powers and the 
\index{binomial coefficient}
binomial coefficients are
$$
a^{(n)}:=
\begin{cases}
a^n/n! \\
\;\; 0 
\end{cases}
\qquad
\binom{a}{n}:=
\begin{cases}
(n!)\inv\prod_{i=0}^{n-1}(a-i) &\quad\text{if $n\in\N$}\\
\qquad\quad 0 &\quad\text{otherwise}.
\end{cases}
$$
We have $\binom{a+1}{n}=\binom{a}{n-1}+\binom{a}{n}$ and
$\binom{a}{n}=(-1)^n\binom{-a-1+n}{n}$.
If $a$ and $b$ are commuting even elements of $A$ then 
$\binom{a+b}{n}=\sum_{i\geq 0}\binom{a}{i}\binom{b}{n-i}$
and
$(a+b)^{(n)}=\sum_{i\geq 0}a^{(i)}b^{(n-i)}$.

For $n, m, i\in\N$, we have $\sum_{j=n}^m\binom{j}{n}=\binom{m+1}{n+1}$ and
$\binom{n}{i}\binom{i}{m}=\binom{n}{m}\binom{n-m}{i-m}$.

\bigskip

{\bf Proposition.}\: {\it
\iti\:
For $n, m\in\N$ with $m>0$, we have
$$
\sum_{i=0}^n\: (-1)^i\frac{\binom{n}{i}}{m+i}
\; =\;
\frac{n!}{\prod_{i=0}^n(m+i)}.
$$
\itii\:
For $n, m\in\N$, we have
$$
\sum_{i=0}^n\: (-1)^i\binom{m+i}{i}\binom{m+n+1}{n-i}
\; =\;
1.
$$
}

\bigskip

\begin{pf}
\iti\:
We do induction on $m$. For $m=1$, we have
\begin{align}
\notag
&\sum_{i=0}^n (-1)^i\frac{\binom{n}{i}}{1+i}
=
\sum_{i=0}^n (-1)^i \binom{n+1}{i+1}(n+1)\inv
\\
\notag
=
&(n+1)\inv -(1-1)^{n+1}(n+1)\inv
=
(n+1)\inv.
\end{align}
For $m\geq 1$, we have
\begin{align}
\notag
&\sum_{i=0}^n (-1)^i\frac{\binom{n}{i}}{m+1+i}
=
\sum_{i=0}^n \frac{(-1)^i\binom{n+1}{i+1}+(-1)^{i+1}\binom{n}{i+1}}{m+1+i}
\\
\notag
=
&\frac{1}{m}
-
\sum_{i=0}^{n+1} (-1)^i\frac{\binom{n+1}{i}}{m+i}
+
\sum_{i=0}^n (-1)^i\frac{\binom{n}{i}}{m+i}
-
\frac{1}{m}
\\
\notag
=
&-
\frac{(n+1)!}{\prod_{i=0}^{n+1}(m+i)}
+
\frac{n!}{\prod_{i=0}^n(m+i)}
=
\frac{n!}{\prod_{i=0}^n(m+1+i)}.
\end{align}

\itii\:
Using \iti, we have 
\begin{align}
\notag
&\sum_{i=0}^n\: (-1)^i\binom{m+i}{m}\binom{m+n+1}{m+i+1}
\\
\notag
\; =\;
&\sum_{i=0}^n\: (-1)^i\binom{m+i}{m}\binom{m+n}{m+i}\frac{m+n+1}{m+1+i}
\\
\notag
\; =\;
&\sum_{i=0}^n\: (-1)^i\binom{n}{i}\binom{m+n}{m}\frac{m+n+1}{m+1+i}
\\
\notag
\; =\;
&\frac{n!}{\prod_{i=0}^n(m+1+i)}\binom{m+n}{m}(m+n+1)
\; =\;
1.
\end{align}
\end{pf}

\chapter{Bibliographical Notes}
\label{C:notes}

\section{Chapter \ref{C:intro}}

The original definition of an associative vertex algebra is due to Borcherds 
\cite{borcherds.voa}, section 4.

Bakalov and Kac proved that to give an associative vertex algebra is 
equivalent to giving a vertex Lie algebra with a pre-Lie algebra structure 
such that the Wick formula holds and $[\,,]\subast=[\,,]\subslie$
\cite{bakalov.kac.field.algebras}, Theorem 7.9.

Introductary books about vertex algebras have been written 
by Frenkel, Lepowsky, and Meurman, 
by Frenkel, Huang, and Lepowsky, 
by Kac, 
by Xu, 
by Matsuo and Nagatomo, 
by Frenkel and Ben-Zvi, and
by Lepowsky and Li
\cite{frenkel.lepowsky.meurman.book,
frenkel.huang.lepowsky.axiom,
kac.beginners.first.ed,
kac.beginners,
xu.vertex.book,
matsuo.nagatomo.locality,
frenkel.benzvi.book,
frenkel.benzvi.book.2nd,
lepowsky.li.voa.book}.
These books also treat more advanced topics.

More specialized books about vertex algebras have been written 
by Feingold, Frenkel, and Ries,
by Tsukada, 
by Prevost, 
by Husu, 
by Dong and Lepowsky, 
by Huang, 
by Weiner, 
by Meurman and Primc, 
by Tamanoi, and 
by Barron
\cite{feingold.frenkel.ries.spinor.vertex.triality.e(1)8,
tsukada.string.path.integral.realization.voas,
prevost.vas.integral.bases.enveloping.algs.affine.lie,
husu.jacobi.extensions1,
dong.lepowsky.book,
huang.book,
weiner.bosonic.constr.voparaalgs.sympl.affine.algs,
meurman.primc.annihilating.fields.sl2.combinat.ids,
tamanoi.elliptic.genera.voas.book,
barron.n=1.supertubes}.

There is also a proceedings volume dedicated to vertex algebras
\cite{voas.math.phys.toronto.00}.

\section{Chapter \ref{C:vlie}}

\indent
{\bf Section \ref{S:loclie}.}\:
\ref{SS:loclie}.\:
Kac introduced the notion of a local Lie algebra.
He defined them in terms of locality instead of the weak commutator formula.
He called a Lie algebra that is spanned by the coefficients of
mutually local distributions a ``Lie algebra of formal distributions" and
also a ``formal distribution Lie algebra"
\cite{kac.beginners.first.ed}, Definition 2.7a, 
\cite{kac.beginners}, Definition 2.6b.
He called such an algebra ``regular'' if there exists 
a derivation $T$ such that $Ta_t=-ta_{t-1}$
\cite{kac.beginners.first.ed,kac.beginners}, equation (4.7.1).
Thus local Lie algebras are ``regular formal distribution Lie algebras''.

Dong, Li, and Mason defined a notion that is essentially equivalent to
a regular local Lie algebra
\cite{dong.li.mason.vertex.lie.poisson.algebras}, Definition 3.1.
They call their notion a ``vertex Lie algebra''.

Primc used the name ``local Lie algebra'' for the Lie algebras $\fg(R)$
constructed from vertex Lie algebras $R$ 
\cite{primc.lie}, abstract and after Theorem 4.1.
The Lie algebras $\fg(R)$ are in fact regular local Lie algebras, 
see section \ref{S:vlie to loclie}.

The name ``local Lie algebra'' is also used in two different contexts.
Kac used it in his work on the classification of transitive, irreducible 
graded Lie algebras of finite growth \cite{kac.kacmoody}.
There a local Lie algebra is 
a Lie algebra $\fg_0$ together with $\fg_0$-modules $\fg_{-1}, \fg_1$ 
and a $\fg_0$-module morphism $\fg_{-1}\otimes\fg_1\to\fg_0$.
Kac constructed $\Z$-graded Lie algebras from these data.

Kirillov used the name ``local Lie algebra'' for continuous Lie algebra 
structures on the space of smooth functions on a manifold satisfying 
a certain condition on the support \cite{kirillov.local.lie}. 
Local Lie algebras in this sense are a generalization of Poisson manifolds. 

It is not true that an operator $T$ of a formal distribution Lie algebra
satisfying $Ta_t=-ta_{t-1}$ must be a derivation.
This is erroneously claimed in 
\cite{kac.beginners.first.ed,kac.beginners}, after equation (4.7.1).
A counterexample is the formal distribution Lie algebra $\Vir$
with $F_{\Vir}:=\set{L(z),\hc e^z}$ and $TL_n=-(n+1)L_{n-1}$ and $T\hc=\hc$.

\smallskip

\ref{SS:loop affi lie}, \ref{SS:heis cliff lie}.\:
Affine Lie algebras and more general Kac-Moody algebras
are studied in the books 
\cite{kac.kmbook,kass.moody.patera.slansky.affine.lie,
moody.pianzola.triangular.book,wakimoto.infini.dim.lie.book,
wakimoto.lec.infini.dim.lie.book}.
See also \cite{pressley.segal.loopgrps,kumar.kac.moody.flag.repr.book} 
for books about the corresponding infinite-dimensional Lie groups.

Garland proved Proposition \ref{SS:loop affi lie} 
\cite{garland.arithmetic.loop.groups}, Theorem (2.36).
See also \cite{wilson.euclid.lie.univ.central.ext}, Theorem, 
and \cite{weibel.homological}, Theorem 7.9.11.

\smallskip

\ref{SS:witt alg}, \ref{SS:vir alg}.\:
The Virasoro algebra is studied in a book by Kac and Raina 
\cite{kac.raina.bombay}.
The discovery of the Virasoro algebra is usually attributed 
to Gelfand and Fuchs who computed the cohomology of the Lie algebra of 
smooth vector fields on the circle
\cite{gelfand.fuks.cohomologies.lie.vector.fields.circle}. 

The name of Virasoro is associated with this algebra 
because he introduced operators $O_n$ acting on the Fock space
of the bosonic string \cite{virasoro.virasoro}. 
Fubini and Veneziano expressed the operators $O_n$ in terms of operators $L_n$
and stated that $[L_n,L_m]=(m-n)L_{n+m}$ \cite{fubini.veneziano.virasoro}.
In a note added in proof they pointed out that 
the correct formula actually is 
$[L_n,L_m]=(m-n)L_{n+m}+c_n\de_{n+m}\id$ 
for some $c_n\in\C$.
In fact, the operators $L_n$ generate the Virasoro algebra
and play a fundamental role in string theory.

The fact that a Lie algebra has a universal central extension iff it
is perfect, is proven for example in \cite{weibel.homological}, Theorem 7.9.2.
The correspondence between extensions and $H^2(\fg,M)$,
used in the proof of Proposition \ref{SS:vir alg}, 
can be found in \cite{weibel.homological}, Theorem 7.6.3.

The fact that the Virasoro algebra is the unique non-trivial
$1$-dimensional central extension of $\Witt$ 
is proven for example in \cite{kac.raina.bombay}, Proposition 1.3,
and \cite{frenkel.lepowsky.meurman.book}, Proposition 1.9.4.

\smallskip

\ref{SS:super witt}, \ref{SS:super vir}.\:
Neveu and Schwarz discovered the contact Lie algebra $K_1$ 
in superstring theory \cite{neveu.schwarz.model}, equations (4.3) and (4.4).
Kac defined the Lie algebras $W_N, K_N$, and $S_N$ 
with $\K[z\uppm]$ replaced by $\K[z_1, \dots, z_m]$
\cite{kac.lie.super.algebras}, section 5.4.2.

The large $N$=4 super Virasoro algebra was discovered by 
Sevrin, Troost, and van Proeyen 
\cite{sevrin.troost.vanproeyen.2d.n=4.superconformal}.

\medskip

{\bf Section \ref{S:loc weak comm f}.}\:
\ref{SS:distr}.\:
Borcherds introduced the commutator formula for vertex algebras
\cite{borcherds.voa}, section 8.
The weak commutator formula is a weak version of it,
see section \ref{SS:module vlie}.

Kac observed that one can rewrite $\sum p_j(t)d^j_{t+s}$ as
$\sum\binom{t}{i}c^i_{t+s-i}$ and conversely, using that 
the polynomials $\binom{t}{i}$ form a basis of $\K[t]$
\cite{kac.beginners.first.ed,kac.beginners}, Theorem 2.3\,\itiv, \itvi.
But he missed to say that the polynomials $p_j(t)$ must have degree $\leq N-1$ 
and the sum over $j$ may be arbitrary.

Borcherds used the space $V\pau{z\uppm}$ of $V$-valued distributions in the
construction of the lattice vertex algebra $V$
\cite{borcherds.voa}, section 3.
He called distributions formal Laurent series.
Kac introduced the name ``formal distribution" 
\cite{kac.beginners.first.ed,kac.beginners}, equation (2.1.1).

Frenkel, Lepowsky, and Meurman discussed the algebra of distributions in detail
\cite{frenkel.lepowsky.meurman.book}, sections 2.1--2.2 and 8.1--8.3.
They also considered more general expressions of the form 
$\sum_{t\in\K}a_t z^t$. 
They just called them formal sums.

Frenkel, Lepowsky, and Meurman introduced the one-variable version 
$\de(x)=\sum_{t\in\Z}x^t$ of the delta distribution $\de(z,w)$ 
\cite{frenkel.lepowsky.meurman.book}, equation (2.1.22).
Kac defined $\de(z,w)$ 
\cite{kac.beginners.first.ed,kac.beginners}, equation (2.1.3).
He denoted it by $\de(z-w)$.

\smallskip

\ref{SS:diff op}
Frenkel, Lepowsky, and Meurman observed that distributions can be viewed
as linear maps defined on the space of Laurent polynomials
\cite{frenkel.lepowsky.meurman.book}, equations (A.3.1)--(A.3.5).

\smallskip

\ref{SS:delta dis diff op}, \ref{SS:loc diff op}.\:
Kac proved that locality and the weak commutator formula are equivalent
\cite{kac.beginners.first.ed,kac.beginners}, Theorem 2.3 \iti.

Chambert-Loir found a beautiful proof 
of the fact that locality implies the weak commutator formula
\cite{chambertloir.distributions.champs}, Proposition 2.2.
We present his proof.

\medskip

{\bf Section \ref{S:loclie to vlie}.}\:
\ref{SS:confa}.\:
Kac introduced the name ``conformal algebra'' for the structure underlying
a vertex Lie algebra \cite{kac.conf.algebras}, section 3.
Previously, he used the name ``conformal algebra'' for vertex Lie algebras.
In his new terminology, vertex Lie algebras are ``Lie conformal algebras''.
He also defined associative and commutative conformal algebras, see
section \ref{SS:assoc confa}.

Lian and Zuckerman defined the notion of weak locality for fields 
\cite{lian.zuckerman.classicalq}, Definition 3.5.
\cite{lian.zuckerman.quantum}, Definition 2.4.
They called it ``locality''.
Kac introduced the name ``weak locality'' 
\cite{kac.beginners.first.ed}, beginning of section 4.11,
\cite{kac.beginners}, equation (4.11.3).

\smallskip

\ref{SS:alg of distr}, \ref{SS:weak comm f ex}.\:
Kac defined the unbounded conformal algebra $\fg\pau{z\uppm}$
\cite{kac.beginners.first.ed,kac.beginners}, equation (2.3.8).
It generalizes the unbounded conformal algebra of fields, 
see section \ref{SS:fields}.

\smallskip

\ref{SS:conf jac}.\:
Kac proved the conformal Jacobi identity for distributions
\cite{kac.beginners.first.ed}, Proposition 2.3 (c),
\cite{kac.beginners}, Proposition 2.3 (d).

\smallskip

\ref{SS:delta dis}, \ref{SS:conf skew sym}.\:
Kac proved conformal skew-symmetry for local distributions 
\cite{kac.beginners.first.ed}, Proposition 2.3 (b),
\cite{kac.beginners} Proposition 2.3 (c).

He remarked that conformal skew-symmetry implies that the identities 
$[(Ta)\subla b]=-\la[a\subla b]$ and $[a\subla Tb]=(T+\la)[a\subla b]$ are 
equivalent \cite{kac.beginners.first.ed}, equation (2.7.2), 
\cite{kac.beginners}, before  Corollary 2.7.

\smallskip

\ref{SS:vlie}.\:
Getzler, Kac, and Primc defined independently from each other 
the notion of a vertex Lie algebra 
\cite{getzler.maninpairs}, after Proposition A.2, 
\cite{kac.beginners.first.ed}, Definition 2.7b, 
\cite{kac.beginners}, Definition 2.7, 
\cite{primc.lie}, equation (3.3).
They defined vertex Lie algebras as a part of the structure of an
associative vertex algebra. 

Kac formulated vertex Lie algebras in terms of a $\la$-bracket
\cite{kac.beginners}, equations (2.7.2) and (C1)$\subla$--(C3)$\subla$.

The name ``vertex Lie algebra'' is due to Primc. 
Kac called them conformal algebras in 
\cite{kac.beginners.first.ed,kac.beginners}
and Lie conformal algebras in \cite{kac.conf.algebras}, section 3.
Dong, Li, and Mason used the name ``vertex Lie algebra'' 
for a related notion
\cite{dong.li.mason.vertex.lie.poisson.algebras}, Definition 3.1.

Borcherds noted that if $V$ is an associative vertex algebra then $V/TV$ is 
a Lie algebra \cite{borcherds.voa}, section 4.
Kac and Primc remarked that this result is true for any vertex Lie algebra
\cite{kac.beginners.first.ed,kac.beginners}, Remark 2.7a, 
\cite{primc.lie}, Lemma 3.1.

Li discovered Dong's lemma, see section \ref{SS:dong fields}.

Kac constructed from a local Lie algebra $\fg$ the vertex Lie algebra
$R(\fg)$ \cite{kac.beginners.first.ed}, after Definition 2.7b, 
\cite{kac.beginners}, after Corollary 2.7.

\medskip

{\bf Section \ref{S:vlie to loclie}.}\:
\ref{SS:vlie to loclie}.\:
Borcherds constructed from any associative vertex algebra $V$ 
the Lie algebra $\fg(V)=\hV/T\hV$ \cite{borcherds.voa}, section 8.
Therefore we call $\fg(R)$ the Borcherds Lie algebra of $R$. 
It is analogous to the Frenkel-Zhu associative algebra $\cA(V)$,
see section \ref{SS:fz alg}.

Kac called $\fg(R)$ the maximal formal distribution Lie algebra
associated to $R$ \cite{kac.beginners}, before Theorem 2.7.

Roitman found the universal property of the Lie algebra $\fg(R)$ that we use 
as the definition of $\fg(R)$ \cite{roitman.michael.phd}, Theorem 1.1 (e).
The observation that the functors $\fg\mapsto R(\fg)$ and $R\mapsto\fg(R)$ 
are adjoint is new.

\smallskip

\ref{SS:module vlie}.\:
Borcherds discovered the commutator formula
\cite{borcherds.voa}, section 8.
He wrote it in terms of $t$-th products. 
He showed that any associative vertex algebra satisfies it. 
The name is due to Frenkel, Huang, and Lepowsky 
\cite{frenkel.huang.lepowsky.axiom}, equation (2.3.13).
Kac called it the Borcherds commutator formula 
\cite{kac.beginners}, equation (4.6.3).

Kac defined modules over a vertex Lie algebra $R$ in terms of a map
$R\to\fgl(M)\pau{\la}$, see section \ref{SS:conf op}.

Our definition of an $R$-module is new, and hence so is the 
correspondence between $R$-modules and $\fg(R)$-modules and 
the equivalence with the commutator formula. 

Getzler {\it defined} modules over a vertex Lie algebra $R$
to be modules over $\fg(R)$ \cite{getzler.maninpairs}, after Lemma A.3.

\smallskip

\ref{SS:affin vlie}.\:
Kac defined the affinization of a vertex Lie algebra 
\cite{kac.beginners}, equation (2.7.4).
It generalizes the affinization of associative vertex algebras,
see section \ref{SS:tensor va2}.
Kac also mentioned $C\otimes R$ \cite{kac.beginners}, Remark 2.7d.

The notion of a vertex Lie algebra over $C$ is new, and hence so is 
the universal property of $C\otimes R$.

\smallskip

\ref{SS:loclie of vlie}.\:
Borcherds' construction of the Lie algebra $\fg(V)=\hV/T\hV$
was generalized by Kac and by Primc to vertex Lie algebras
\cite{kac.beginners.first.ed}, after Definition 2.7b, 
\cite{kac.beginners}, after equation (2.7.7) and Remark 2.7b,
\cite{primc.lie}, Theorem 4.1.
Kac observed in addition that $\fg(R)$ is a local Lie algebra 
and $R\mapsto\fg(R)$ is a functor.
Primc observed that $\fg(R)$ is a differential Lie algebra.

Kac remarked that the map $R\to\fg_R, a\mapsto a_{-1}$, is injective
\cite{kac.beginners}, Lemma 2.7.
Kac defined the local Lie algebras $\fg$ and $\fg(R(\fg))$ to be equivalent 
and stated that the ``category of equivalence classes'' of local
Lie algebras is equivalent to the category of vertex Lie algebras
\cite{kac.beginners}, Theorem 2.7.

The observation, that the functor $R\mapsto\fg(R)$ is fully faithful, is new.

\smallskip

\ref{SS:vlie into loclie}.\:
Fattori and Kac showed that a local Lie algebra $\fg$ is regular if
$\bF\subfg=\K[\del_z]F\subfg$ and the $a_t$ form a basis of $\fg$
\cite{fattori.kac.classifi.finite.simple.lie.conf.superalgebras}, 
Proposition 3.2.

\smallskip

\ref{SS:bracket vlie}.\:
Bakalov and Kac and Gorbounov, Malikov, and Schechtman 
defined the bracket $[\, ,]\subslie$ and proved that it is a Lie bracket
\cite{bakalov.kac.field.algebras}, equation (7.24) and Proposition 7.10, 
\cite{gorbounov.malikov.schechtman.gerbes2}, equation (10.6.1) and 
Theorem 10.7,
\cite{gorbounov.malikov.schechtman.gerbes2.v2},
equation (8.21.1) and Theorem 8.22.

Bakalov and Kac observed that the map $R\subsmlie\to\fg(R)_-, a\mapsto a_{-1}$,
is an isomorphism of differential algebras
\cite{bakalov.kac.field.algebras}, before Theorem 7.12.

\medskip

{\bf Section \ref{S:ex vlie}.}\:
\ref{SS:unb confa power series}\:
The results about the unbounded vertex Lie algebra $\fg\pau{\mu}$
and the morphism $\fg\pau{z\uppm}\to\fg\pau{\mu}$ are new.

\smallskip

\ref{SS:conf op}\:
Kac defined conformal operators and, more generally, conformal linear maps 
$a\submu: R\to R'[\mu]$ between $\K[T]$-modules $R$ and $R'$ 
\cite{kac.beginners}, Definition 2.10.
He shows that conformal operators of a finite $\K[T]$-module form
a vertex Lie algebra \cite{kac.beginners}, after Remark 2.10d.

Kac defined conformal modules over a vertex Lie algebra in terms of
$t$-th products and noted that to give a conformal module is equivalent to
giving a morphism $R\to\fgl\subsc(M)$ 
\cite{kac.beginners}, Definition 2.8b and Remark 2.10e.

D'Andrea and Kac proved that any torsion conformal operator is zero and
that any conformal operator annihilates the torsion submodule 
\cite{dandrea.kac.finite.conformal.algebras}, Remark 3.1 and Proposition 3.2.

\smallskip

\ref{SS:conf der}\:
D'Andrea and Kac defined conformal derivations and inner conformal derivations 
\cite{dandrea.kac.finite.conformal.algebras}, Definition 3.2.

\smallskip

\ref{SS:free vlie}.\:
Roitman constructed free vertex Lie algebras 
\cite{roitman.free.conformal.vertex.algebras}, Proposition 3.1.
Our construction of them is different from Roitman's.
We explain his method in section \ref{SS:free vlie rev}.

\smallskip

\ref{SS:finite confa}.\:
Kac defined finite vertex Lie algebras and central elements of
vertex Lie algebras \cite{kac.beginners}, after Remark 2.7a.
He showed that vertex Lie algebras $R=E[\mu]$ can be defined
in terms of $E$ and a map $E\otimes E\to R[\la]$
\cite{kac.beginners.first.ed}, Theorem 2.7. 
This result is not contained in the second edition \cite{kac.beginners}.

D'Andrea and Kac proved that the torsion submodule of a vertex Lie algebra
is contained in the centre 
\cite{dandrea.kac.finite.conformal.algebras}, Proposition 3.1.

\smallskip

\ref{SS:aff vlie}.\:
Lian proved that for any associative vertex algebra $V$ of CFT-type 
the subspace $V_1$ is a Lie algebra with an invariant symmetric
bilinear form \cite{lian.simplevertex}, Theorem 3.7.
Previously, Frenkel, Lepowsky, and Meurman proved this result for
lattice vertex algebras \cite{frenkel.lepowsky.meurman.book}, Remark 8.9.1.

\smallskip

\ref{SS:witt vlie}.\:
D'Andrea and Kac proved that the Witt vertex Lie algebra is the only 
non-abelian even vertex Lie algebra of rank $1$
\cite{dandrea.kac.finite.conformal.algebras}, Proposition 3.3.
They acknowledged that the result is joint work with Wakimoto.
Kac explains this result in his book \cite{kac.beginners}, equation (2.7.10).

D'Andrea and Kac proved that any even finite simple vertex Lie algebra
is isomorphic to a loop or the Witt vertex Lie algebra  
\cite{dandrea.kac.finite.conformal.algebras}, Theorem 5.1.
Kac had conjectured this result \cite{kac.beginners.first.ed}, Conjecture 2.7.

\smallskip

\ref{SS:vir vlie}.\:
Frenkel and Zhu proved that if $V$ is a vertex operator algebra of CFT-type 
and $V'$ a conformal vertex subalgebra then $C_V(V')=\ker L'_{-1}$ 
\cite{frenkel.zhu.vertex}, Theorem 5.2.
In particular, the centre of $V$ is equal to the kernel of $T$.
Li also proved this last statement
\cite{li.invariant.forms.voas}, Proposition 3.3 (a).

\smallskip

\ref{SS:griess}.\:
Griess constructed the largest sporadic finite simple group $\M$, 
called the monster \cite{griess.friendly.giant}.
He constructed the monster as a group of automorphisms of a
commutative nonassociative algebra $B$ with an invariant, symmetric
bilinear form. 
As was shown later, $\M$ is in fact the automorphism group of $B$.
The dimension of $B$ is 196884.
The $\M$-module $B$ decomposes as $\K L\oplus B'$ where $L/2$ is an identity 
of $B$ and $B'$ is the smallest non-trivial simple $\M$-module.
Assuming the existence of $\M$ and the $\M$-module $B'$, 
Norton had shown earlier that $B$ has the structure of a commutative
nonassociative algebra with an invariant, symmetric bilinear form. 

The work of Griess played an important role in the construction of
the moonshine module $V\upnatural$ by Frenkel, Lepowsky, and Meurman
\cite{frenkel.lepowsky.meurman.moonshine.module,frenkel.lepowsky.meurman.book}.
The moonshine module is a conformal vertex algebra of CFT-type 
on which $\M$ acts such that $V\upnatural_1=0$ and $V\upnatural_2=B$.

Borcherds remarked that if $V$ is a conformal vertex algebra with 
a positive definite inner product and $V_1=0$ then 
the multiplication of the Griess algebra $V_2$ is commutative 
\cite{borcherds.voa}, section 9.

Frenkel, Lepowsky, and Meurman showed that if $V$ is a lattice vertex 
operator algebra 
associated to a positive definite lattice such that $V_1=0$ then 
a certain fixed point subspace of $V_2$ has a commutative product and 
an invariant, symmetric bilinear form 
\cite{frenkel.lepowsky.meurman.book}, Theorem 8.9.5.

Lam showed that for any vertex operator algebra $V$ of CFT-type with $V_1=0$
the above structure exists on $V_2$
\cite{lam.vertex.from.comm.assoc.algebras}, before Lemma 2.

Zagier observed that if $V$ is of CFT-type and $L\in V_2$ such that
$L_{-1}=T$ and $L_0=H$ then $L$ is a conformal vector
\cite{zagier.vertex.sl2}, Proposition 2.

\smallskip

\ref{SS:gko vector}.\:
Belavin, Polyakov, and Zamolodchikov defined the notions of a 
primary field, primary state, and quasi-primary field
in two-dimensional conformal field theory
\cite{belavin.polyakov.zamolodchikov.minimal}, 
equations (1.16), (3.8), (3.9), (3.22), (A.4).

The condition $L_n a=0$ for $n>0$ also arises in the quantization of
the bosonic string as a necessary condition for states to be physical
\cite{zwiebach.first.course.string.th}, equation (21.32).
Borcherds used this condition to construct Lie algebras from vertex algebras
\cite{borcherds.voa}, section 4.
He called the space of primary vectors the physical subspace.

Goddard, Kent, and Olive discovered the coset construction 
\cite{goddard.kent.olive.coset,goddard.kent.olive.unitary.coset}.

Frenkel and Zhu explained the coset construction in the context of 
vertex algebras \cite{frenkel.zhu.vertex}, section 5.
They proved that if $V$ is a vertex operator algebra of CFT-type and
$V'$ a subalgebra then $C_V(V')=\ker L'_{-1}$ and 
if $L'$ is quasi-primary then $L-L'$ is a conformal vector of $C_V(V')$
\cite{frenkel.zhu.vertex}, Theorems 5.2 and 5.1.

\smallskip

\ref{SS:chodos thorn}.\:
The original reference to the Chodos-Thorn construction is
the string theory paper \cite{chodos.thorn.linear.dilaton}.
The construction is most often applied to the Heisenberg vertex algebra. 
The resulting theory is sometimes called the linear dilaton theory
\cite{polchinski.book1}.
The mathematical reference for the general result is 
\cite{dong.lin.mason.voa.sl2.modules}, Proposition 4.1.
The Chodos-Thorn construction is important, for example,
because it allows to study the representation theory of 
the Virasoro algebra using representations of the simpler
Heisenberg algebra.

\smallskip

\ref{SS:n1 n2 top vir}.\:
The Chodos-Thorn construction applied to $\twoVir$ 
is a manifestation of the topological twist in $N$=2 supersymmetric
field theory. 
Witten discovered the topological twist for the $N$=2 sigma model
of a K\"ahler manifold \cite{witten.sigma}, after equation (2.19).
This is the $A$-model of mirror symmetry and is formulated 
mathematically in terms of Gromov-Witten theory. 
For Calabi-Yau manifolds, the sigma model has an $N$=2 super Virasoro symmetry 
and a second twist is well-defined on the quantum level \cite{witten.mirror}. 
This is the $B$-model.
The mirror involution identifies the two twists. 
Eguchi and Yang discussed the topological twist for general $N$=2 
superconformal field theories \cite{eguchi.yang.twisting}.

\medskip

{\bf Section \ref{S:vlie suppl}.}\:
\ref{SS:s3 symm conf J}.\: 
Li proved $\bbS_3$-symmetry of the conformal Jacobi identity 
\cite{li.vertex.alg.vertex.poisson.alg}, Proposition 2.8.
D'Andrea and Kac observed this property independently 
\cite{dandrea.kac.finite.conformal.algebras}, Remark 2.4.

\smallskip

\ref{SS:semidir prod}.\:
D'Andrea and Kac defined the semi-direct product of the Witt and
a loop vertex Lie algebra
\cite{dandrea.kac.finite.conformal.algebras}, Example 3.4.

\smallskip

\ref{SS:assoc confa}.\:
Kac defined associative and commutative conformal algebras
\cite{kac.beginners}, 
(A1)$\subla$, (A2)$\subla$ of section 2.10 and after Remark 2.10b.

\smallskip

\ref{SS:lambda commut}.\:
Kac defined the $\la$-commutator and claimed that it is a Lie
$\la$-bracket \cite{kac.beginners}, Remarks 2.10a and 2.10b.

\smallskip

\ref{SS:frob alg vlie}.\:
Lam constructed from a commutative Frobenius algebra $C$
an associative vertex algebra of CFT-type with $V_1=0$ and Griess algebra $C$
\cite{lam.vertex.from.comm.assoc.algebras}, end of section 5.2.

Rebout and Schechtman generalized Lam's construction to 
Virasoro algebroids 
\cite{rebout.schechtman.alg.frobenius.virasoro}, section 6.
Any commutative Frobenius algebra is a Virasoro algebroid over $\K$ and
any Virasoro algebroid $C$ over a commutative algebra $\cO$ 
yields an $\N$-graded associative vertex algebra $V$ with $V_0=\cO$, 
$V_1=0$, and $V_2=C$.

\smallskip

\ref{SS:free vlie rev}.\:
Roitman constructed free vertex Lie algebras using free Lie algebras 
\cite{roitman.free.conformal.vertex.algebras}, Proposition 3.1.

\smallskip

\ref{SS:1trunc conf}, \ref{SS:1trunc vlie}.\:
Gorbounov, Malikov, and Schechtman 
defined 1-truncated vertex Lie algebras and showed that the forgetful
functor from $\N$-graded vertex Lie algebras to 1-truncated vertex Lie algebras
has a left adjoint
\cite{gorbounov.malikov.schechtman.gerbes2.v2}, section 9.6 and Theorem 9.8. 

Li and Yamskulna simplified the definition of a 1-truncated vertex Lie algebra
\cite{li.yamskulna.vertex.algebroids}, Lemma 2.5 and Proposition 2.6.

Loday defined Leibniz algebras
\cite{loday.cyclic.homology,
loday.pirashvili.enveloping.leibniz.algebras.cohomology}.

\section{Chapter \ref{C:va}}

${}_{}$
\indent
{\bf Section \ref{S:va def}.}\:
\ref{SS:prelie}.\:
Since pre-Lie algebras are little known, we give a brief overview of the role 
of pre-Lie algebras in mathematics.
We also give a summary of some basic results about other classes of 
non-associative algebras.

There are two main examples of pre-Lie algebras:
an operadic construction of pre-Lie algebras that leads to 
important Lie algebras and 
a correspondence between pre-Lie algebras and affine structures.

Some authors claim that already Cayley studied pre-Lie algebras 
\cite{cayley.prelie.trees}.
Vinberg and Gerstenhaber defined pre-Lie algebras at about the same time 
\cite{vinberg.homogeneous.convex.cones}, 
chapter 2, Definition 2, \cite{gerstenhaber.cohomology}, equation (6).
Vinberg observed that the associator is symmetric in the first two arguments
and hence he called pre-Lie algebras 
\index{left-symmetric!algebra}
left-symmetric algebras
\cite{vinberg.homogeneous.convex.cones}, chapter 2, equation (5).
He showed that there is a correspondence between homogeneous convex domains
and certain pre-Lie algebras.

Gerstenhaber discovered a graded pre-Lie algebra structure on 
the Hochschild cochain complex of an associative algebra $A$
with coefficients in $A$. 
He proved that the commutator of a pre-Lie algebra is a Lie bracket
\cite{gerstenhaber.cohomology}, Theorem 1, and deduced that
the Hochschild cohomology has an odd Lie bracket. 
Together with the cup product the Hochschild cohomology is
thus a so-called Gerstenhaber algebra. 
Algebras for which the commutator is a Lie bracket
\index{Lie-admissible algebra}
are called Lie-admissible.

Gerstenhaber's construction of the pre-Lie multiplication uses
the operadic structure of the Hochschild complex.
He called this a pre-Lie system and proved that any pre-Lie system
yields a pre-Lie algebra \cite{gerstenhaber.cohomology}, Theorem 2.
Kapranov and Manin observed a variation of Gerstenhaber's result,
namely, they showed that for any operad $\cP$ the sum 
$\bigoplus\cP(n)$ is a pre-Lie algebra with $a\circ b:=\sum a\circ_i b$
\cite{kapranov.manin.morita.operads}, Theorem 1.7.3.
Gerstenhaber and Voronov pointed out earlier the underlying Lie algebra 
structure \cite{gerstenhaber.voronov.gerstenhaber}, equation (3).
Markl, Shnider, and Stasheff defined in their book
the pre-Lie algebra structure on 
the cochain complex computing the cohomology of algebras over an operad  
\cite{markl.shnider.stasheff.book}, Proposition 3.111.

Matsushima \cite{matsushima.affine.str.cplx.mfds} showed that
there exists a close relation between affine structures on manifolds
and pre-Lie algebra structures on the tangent bundle.
An affine structure is a connection $\nabla$ on the
tangent bundle that is torsion-free, $\nabla_X Y-\nabla_Y X=[X,Y]$,
and flat, $[\nabla_X,\nabla_Y]=\nabla_{[X,Y]}$.
In this case $XY:=\nabla_X Y$ defines a pre-Lie algebra structure on
the tangent bundle with $[X,Y]\subast=[X,Y]$.
Affine structures correspond to atlases with affine coordinate transformations.
Affine structures, in particular on solvable Lie groups, have
been studied from this point of view in a number of articles.
Matsushima proved that commutative pre-Lie algebras are associative
\cite{matsushima.affine.str.cplx.mfds}, Lemma 3.

There is the following algebraic version of the correspondence between 
pre-Lie algebras and affine structures.
Let $\cO$ be a commutative algebra such that $\Der(\cO)$ is a free
$\cO$-module with a basis of commuting derivations $\del_1, \dots, \del_n$.
Then $\nabla_X\del_i:=0$, for $X\in\Der(\cO)$ and $i=1, \dots, n$,
defines an affine structure and hence $\Der(\cO)$ is a pre-Lie algebra. 
For $\cO=\K[x_1, \dots, x_n]$, 
this pre-Lie algebra is 
\index{left-symmetric!Witt algebra}
called the left-symmetric Witt algebra.

More generally, if $\del_1, \dots, \del_n$ are commuting derivations
of a commutative algebra $\cO$ then the free $\cO$-module spanned by $\del_i$ 
is a pre-Lie algebra with $(f\del_i)(g\del_j):=f(\del_i g)\del_j$.
In the special case $n=1$, this algebra also satisfies $(ab)c=(ac)b$.
Such pre-Lie algebras are called 
\index{Novikov algebra}
Novikov algebras.
They appear in the theory of Hamiltonian systems 
\cite{gelfand.dorfman.hamilton.algebraic}, equation (6.3).

Connes and Kreimer defined a pre-Lie algebra product for graphs
and used it to study renormalization of quantum field theories,
see \cite{kreimer.feynman.graphs.hopf.algs.review}, Proposition 8.

Kupershmidt showed that if $V$ is a pre-Lie algebra then so is $V\oplus V\dual$
\cite{kupershmidt.non.abelian.phase.space.prelie}, Proposition 2.

Burde gave an overview of pre-Lie algebras \cite{burde.prelie.overview}.

Besides Lie algebras, the best known non-associative algebra
is probably Cayley's algebra of octonions. 
The algebra of octonions is 
\index{alternative algebra}
an {\it alternative} algebra which means that 
associativity is satisfied for $a, b, c$ if any two of the elements $a, b, c$
are equal or, equivalently, the associator is skew-symmetric.
Zorn proved that the only alternative finite-dimensional real division
algebras are $\R, \C$, Hamilton's quaternions $\bbH$, and the Cayley algebra.
Zorn's result is an extension of Frobenius' theorem stating that 
the only associative finite-dimensional real division algebras are 
$\R, \C$, and $\bbH$.

Emil Artin discovered the beautiful theorem stating that an algebra is 
alternative iff any subalgebra generated by two elements is associative,
see \cite{schafer.intro.nonassoc.algs}, Theorem 3.1.
The proof uses the equally remarkable 
\index{Moufang identity}
Moufang identity 
$\fa(a,bc,b)=-\fa(a,b,c)b$ for the associator of an alternative algebra.
Artin's theorem shows that any alternative algebra is {\it power-associative}
\index{power-associative algebra}
which means that any subalgebra generated by one element is associative,
or, equivalently, $a^n a^m=a^{n+m}$.

Another well-known class of non-associative algebras are {\it Jordan algebras}.
Physical observables in quantum mechanics are hermitian operators.
The product of hermitian operators $a, b$ is not hermitian, but
their anti-commutator $ab+ba$ is. 
The resulting algebra was proposed by Jordan as a foundation for
quantum mechanics. 
It is commutative and satisfies $(ab)a^2=a(ba^2)$. 
These are the two axioms for 
\index{Jordan algebra}
Jordan algebras.
Any associative algebra endowed with the anti-commutator is 
a Jordan algebra. 
Finite-dimensional simple Jordan algebras have been classified.
They are of interest because of their relation to exceptional simple 
Lie algebras.

Jordan algebras are power-associative.
This can be seen as a corollary of a theorem of Albert which states
that an algebra is power-associative iff $a a^2=a^2 a$ and $a^2 a^2=a(a a^2)$.

Other examples of non-associative algebras are Griess algebras, 
see section \ref{SS:griess}, and the 24-dimensional Chevalley algebra, 
that is used 
\nocite{chevalley.collected.works.spinors}
in triality \cite{chevalley.spinors}, section 4.2.

\smallskip

\ref{SS:va}.\: 
Bakalov and Kac introduced the notion of a (non-associative) unital 
vertex algebra \cite{bakalov.kac.field.algebras}, Definition 2.1.
They defined the notion in terms of the state-field correspondence
and called unital vertex algebras ``state-field correspondences''.

Lepowsky and Li used the notion of a (non-associative) unital 
unbounded vertex algebra \cite{lepowsky.li.voa.book}, Definition 5.3.1.
They called it a ``weak vertex algebra''.

\smallskip

\ref{SS:wick}.\:
Bais, Bouwknegt, Surridge, and Schoutens wrote down 
the Wick formula for ``operators'' $a(z)$
\cite{bais.bouwknegt.surridge.schoutens.casimir}, equation (A.6).
In fact, they wrote it in a way that only applies to fields.
See section \ref{SS:wick fields} for more references dealing with
the Wick formula for fields.

The Wick formula for fields is a generalization to {\it non-free} 
holomorphic fields of Wick's theorem \cite{wick.theorem.normal.ordering}.
This theorem is part of the standard material of any book about
quantum field theory. 
It deals with the normal ordered product of a finite family of free fields.

Getzler gave a direct proof of the fact that the commutator formula for 
indices $t\geq 0, s\in\Z$ follows from Borcherds' axioms for an associative 
vertex algebra \cite{getzler.maninpairs}, Proposition A.2.

Kac, following Getzler, introduced the Wick formula into the 
mathematical literature \cite{kac.beginners.first.ed}, equation (3.3.10), 
\cite{kac.beginners}, equation (3.3.11).
He called it the ```non-commutative' Wick formula''.
He showed how to write it in terms of the $\la$-bracket 
\cite{kac.beginners}, equation (3.3.12).

Akman proved that the $t$-th products $a_t$ of an associative vertex algebra
are differential operators \cite{akman.bv}, Theorem 2.2.

\smallskip

\ref{SS:conf deriv va}.\:
The notion of a conformal derivation of a vertex algebra is new.
Hence the result that conformal derivations form an 
unbounded vertex Lie algebra is new as well.

\smallskip

\ref{SS:assoc va def}.\:
Borcherds defined associative vertex algebras \cite{borcherds.voa}, section 4.
He defined them as unital bounded $\Z$-fold algebras with a
translation operator satisfying skew-symmetry and the associativity formula. 

Bakalov and Kac proved that to give an associative vertex algebra is 
equivalent to giving a vertex Lie algebra with a pre-Lie algebra structure 
such that the Wick formula holds and $[\,,]\subast=[\,,]\subslie$
\cite{bakalov.kac.field.algebras}, Theorem 7.9.

Before Bakalov and Kac proved their theorem, 
Dong, Li, and Mason observed that the commutator of the 
$(-1)$-st product of an associative vertex algebra is a Lie bracket 
\cite{dong.li.mason.vertex.lie.poisson.algebras}, Lemma 5.3.
Then Matsuo and Nagatomo and Akman pointed out that 
the $(-1)$-st product satisfies the pre-Lie algebra identity
\cite{matsuo.nagatomo.locality}, Remark 8.3.3, 
\cite{akman.htpgerstenhaber}, Proposition 6 in section 4.3.

Previously, the physicists Sasaki and Yamanaka and 
Bais, Bouwknegt, Surridge, and Schoutens
noted that the normal ordered product of fields satisfies 
the pre-Lie algebra identity, see section \ref{SS:norm ord prod pre lie}. 
This implies the result for the $(-1)$-st product since $Y$ is 
a monomorphism.

Zhu defined the subspace $C_2(V):=\rspan\set{a_{-2}b\mid a, b\in V}$ and 
remarked that $V/C_2(V)$ is a commutative algebra \cite{zhu.phd}, section 4.3.
He also noticed that $V/C_2(V)$ is in fact a Poisson algebra 
\cite{zhu.modular}, section 4.4.
Li defined the subspaces $C_n(V)$ \cite{li.finiteness.regular}, equation (3.1).

Zhu required that $C_2(V)$ has finite codimension in his proof of
the modularity of the space of characters of $V$. 
He called this the finiteness condition $C_2$ \cite{zhu.phd}, section 4.3.
He conjectured that any rational vertex operator algebra is $C_2$-cofinite
\cite{zhu.phd}, section 4.3, \cite{zhu.modular}, section 4.4.
This conjecture is still open.
Li proved that regular vertex operator algebras are $C_2$-cofinite
\cite{li.finiteness.regular}, Theorem 3.8.

\smallskip

\ref{SS:comm va}.\:
Borcherds observed that there is an equivalence between commutative vertex 
algebras and commutative differential algebras 
\cite{borcherds.voa}, section 4.

\smallskip

\ref{SS:poisson va}.\:
Enriquez and Frenkel defined vertex Poisson algebras 
in the framework of integrable systems
\cite{enriquez.frenkel.poisson}, section 1.3.
Frenkel and Ben-Zvi gave a definition in the setting of vertex algebras 
\cite{frenkel.benzvi.book}, Definition 15.2.1.

Dong, Li, Mason gave a slightly different formulation of 
vertex Poisson algebras 
\cite{dong.li.mason.vertex.lie.poisson.algebras}, Definition 3.7.
Beilinson and Drinfeld introduced the closely related notion of a 
coisson algebra \cite{beilinson.drinfeld.chiral}, section 1.4.18.

Dong, Li, Mason proved that $V/C_2(V)$ is a Poisson algebra 
\cite{dong.li.mason.vertex.lie.poisson.algebras}, Proposition 3.12.

\medskip

{\bf Section \ref{S:state f corres}.}\:
\ref{SS:fields}.\:
Li defined the notion of a field on a vector space
\cite{li.localsystems}, Definition 3.1.1. 
He called them ``weak vertex operators''.
Lepowsky and Li used the same terminology in their book
\cite{lepowsky.li.voa.book}, Definition 5.1.1.
Kac introduced the name ``field'' 
\cite{kac.beginners.first.ed}, equation (1.3.1). 

Wick invented the normal ordered product in quantum field theory
\cite{wick.theorem.normal.ordering}.
Therefore it is also called the Wick product.
Normal ordering is a fundamental concept in quantum field theory
and can be found in any textbook on the subject.

Normal ordering was used in the representation theory of affine algebras.
Igor Frenkel defined a normal ordering $\normord{a_n b_m}$ for
elements of an affine algebra as part of the Sugawara construction
\cite{frenkel.affine.lie.boson.fermion.corr}, equation (I.1.9).
Igor Frenkel considered a normally ordered product of a differential polynomial
of $U(1)$-currents \cite{frenkel.kac.moody.dual.resonance}, equation (3.21).
Borcherds used Frenkel's idea in his construction of lattice vertex
algebras \cite{borcherds.voa}, section 3.

Feigin and Edward Frenkel defined the normal ordered product 
$\normord{a(z)b(z)}\linebreak[0] =a(z)_-b(z)+b(z)a(z)_+$
for certain homogeneous distributions with coefficients in a 
completed enveloping algebra of an affine algebra 
\cite{feigin.frenkel.gelfanddickey}, before Lemma 1.

Lian and Zuckerman defined the normal ordered product $\normord{a(z)b(w)}$
for fields and noted that $a(z)_{-1-t}b(z)=\;\normord{\del_z^{(t)}a(z)b(z)}$ 
for any $t\geq 0$ \cite{lian.zuckerman.classicalq}, equations (3.4) and (3.10),
\cite{lian.zuckerman.quantum}, equations (2.2) and (2.4).

Previously, Lian defined $\normord{a(z)b(w)}$ for elements $a, b$ of
a vertex algebra \cite{lian.simplevertex}, Definition 2.4.

Li, Lian and Zuckerman, and Meurman and Primc 
introduced $\Endv(E)$ as a $\Z$-fold algebra 
\cite{li.localsystems}, Lemma 3.1.4,
\cite{lian.zuckerman.classicalq}, Definition 3.1,
\cite{lian.zuckerman.quantum}, Definition 2.1, 
\cite{meurman.primc.annihilating.fields.sl2.combinat.ids}, Proposition 2.3. 
Actually, only Li considered fields on a general vector space $E$ without 
a gradation.

Previously, Frenkel and Zhu gave the formula for the $t$-th products 
in a special situation, see section \ref{SS:a z uppm h}.

Li's work about $\Endv(E)$ was very influential and greatly simplified the
theory of vertex algebras.

\smallskip

\ref{SS:modules}.\:
Borcherds defined the notion of a module over a vertex algebra
\cite{borcherds.voa}, section 4. 
He defined modules in terms of the associativity formula. 
Borcherds also mentioned that $\K$ is a module over the
commutative vertex algebra $\K[z\uppm]$ \cite{borcherds.voa}, section 8.

\smallskip

\ref{SS:adj module}, \ref{SS:id elem}.\:
Borcherds pointed out that $V$ is a $V$-module \cite{borcherds.voa}, section 4.
Since he defined associative vertex algebras and modules in terms of the
associativity formula, this statement is trivial.
Frenkel, Huang, and Lepowsky called $V$ the ``adjoint module'' 
\cite{frenkel.huang.lepowsky.axiom}, Remark 4.1.4.

Bakalov and Kac showed that (non-associative) unital vertex algebras and
unital bounded $\Z$-fold algebras with a translation operator
are equivalent notions \cite{bakalov.kac.field.algebras}, 
Lemma 5.1 and Proposition 2.6.

\smallskip

\ref{SS:ratl fcts}.\:
Frenkel, Lepowsky, and Meurman discussed Taylor series expansions 
of rational functions in one and several variables 
\cite{frenkel.lepowsky.meurman.book}, equations (8.1.5), (8.1.6), and
(8.10.37).

\smallskip

\ref{SS:z fold fields}.\:
Lian and Zuckerman pointed out the relation between the 
normal ordered product and the $\Z$-fold algebra $\Endv(E)$,
see section \ref{SS:fields}.

Borcherds discovered the associativity formula \cite{borcherds.voa}, 
section 4.
He wrote it in terms of $t$-th products. 
It is an axiom in his definition of an associative vertex algebra.

The name ``associativity formula'' is due to Matsuo and Nagatomo
\cite{matsuo.nagatomo.locality}, equation (4.3.3).
Li and Lepowsky and Li called it the ``iterate formula'' 
\cite{li.localsystems,lepowsky.li.voa.book}.
Kac and Bakalov and Kac called it the ``$n$-th product axiom''
\cite{kac.beginners}, section 4.11, 
\cite{bakalov.kac.field.algebras}, equation (4.1).

\smallskip

\ref{SS:j id conf j id}.\:
Frenkel, Lepowsky, and Meurman introduced the Jacobi identity 
\cite{frenkel.lepowsky.meurman.book}, section 8.10.
It is the main axiom in their definition of a vertex operator algebra.
They proved the Jacobi identity for vertex operators
and for the moonshine module.

Kac called the Jacobi identity 
\index{Borcherds identity}
the {\it Borcherds identity} 
\cite{kac.beginners.first.ed,kac.beginners}, Theorem 4.8.

\smallskip

\ref{SS:conf j wick0}.\:
Primc proved that the Jacobi identity, the associativity formula, and 
the commutator formula for indices $r, s, t\geq 0$ are all equivalent
\cite{primc.lie}, Lemmas 6.1 and 6.2.

\smallskip

\ref{SS:skew sym}.\:
Borcherds discovered skew-symmetry \cite{borcherds.voa}, 
section 4.
He wrote it in terms of $t$-th products. 
It is an axiom in his definition of an associative vertex algebra.

Frenkel, Huang, and Lepowsky introduced the name 
\cite{frenkel.huang.lepowsky.axiom}, equation (2.3.19).
Kac called it first ``quasi-symmetry'' 
and later called it also skew-symmetry 
\cite{kac.beginners.first.ed,kac.beginners}, equation (4.2.3).

Bakalov and Kac proved that skew-symmetry is equivalent to 
conformal skew-symmetry and $[\, ,]\subast=[\, ,]\subslie$
\cite{bakalov.kac.field.algebras}, Theorem 7.9.

\medskip

{\bf Section \ref{S:vas of distr}.}\:
The main result of this section is that local distributions in
$\Endv(E)$ and $\cA\pau{z\uppm}$ satisfy the axioms of Bakalov and Kac for
an associative vertex algebra.
This was essentially proven by three groups of physicists around 1989.

The main mathematical contribution to this problem is due to Li. 
He proved that a set $S$ of local fields generates a vertex subalgebra of 
$\Endv(E)$ that is associative in the sense of Borcherds 
\cite{li.localsystems}, Corollary 3.2.11.

\smallskip

\ref{SS:fz complete}.\:
Matsuo, Nagatomo, and Tsuchiya defined FZ-topological algebras
\cite{matsuo.nagatomo.tsuchiya.quasi.finite.algs.hamilton.voas}, 
Definition 1.2.1.
They used a different name for it.
They required that the topology of $\cA_h$ is {\it equal} to the topology 
defined by the closures of $\sum_{l<k}\cA_{h-l}\cA_l$. 

They showed that any $\Z$-graded algebra yields an FZ-topological algebra,
defined FZ-complete algebras, and showed that any FZ-topological algebra
has a completion that is an FZ-complete algebra
\cite{matsuo.nagatomo.tsuchiya.quasi.finite.algs.hamilton.voas}, 
section 1.3, Definition 1.2.1, and Proposition 1.3.1.

\smallskip

\ref{SS:a z uppm h}.\:
Igor Frenkel and Zhu defined $t$-th products for $t\in\Z$ of 
an affine current and an arbitrary homogeneous distribution 
with coefficients in a completed enveloping algebra
of an affine algebra \cite{frenkel.zhu.vertex}, Definition 2.2.1.

Feigin and Edward Frenkel defined a normal ordered product 
in a similar context \cite{feigin.frenkel.gelfanddickey}, before Lemma 1.
The unbounded vertex algebra $\cA\pau{z\uppm}_H$ is a natural generalization
of the products of these authors.

\smallskip

\ref{SS:norm ord prod pre lie}.\:
The physicists Sasaki and Yamanaka and 
Bais, Bouwknegt, Surridge, and Schoutens
observed that the normal ordered product of chiral fields satisfies 
the pre-Lie algebra identity 
\cite{sasaki.yamanaka.virsinegordon}, equation (5.10),
\cite{bais.bouwknegt.surridge.schoutens.casimir}, equation (A.11).

Matsuo and Nagatomo stated this result as a mathematical proposition and
referred to the two physics papers
\cite{matsuo.nagatomo.locality}, Proposition 1.4.3.

The physicists Sevrin, Troost, Van Proeyen, and Spindel stated that 
primary fields satisfy a family of four-term identities indexed by 
$r, s\geq -1$
of which the identity for $r=s=-1$ is the pre-Lie algebra identity
\cite{sevrin.troost.vanproeyen.spindel.ext.susy.sigma.mod.grp.mfds2},
equation (A.5). 

Kac mentioned that the commutator of the normal ordered product of fields
is a Lie bracket 
\cite{kac.beginners.first.ed,kac.beginners}, before equation (3.1.5).
He wrote that he learned this fact from Radul.

\smallskip

\ref{SS:wick fields}.\:
Bais, Bouwknegt, Surridge, and Schoutens 
remarked that chiral fields satisfy the Wick formula
\cite{bais.bouwknegt.surridge.schoutens.casimir}, equation (A.6).

Sevrin, Troost, Van Proeyen, and Spindel stated that primary chiral fields 
satisfy the commutator formula for indices $t\geq 0, s\in\Z$ and noted that
the Wick formula is a special case of it
\cite{sevrin.troost.vanproeyen.spindel.ext.susy.sigma.mod.grp.mfds2},
equations (A.6) and (A.8). 

Lian and Zuckerman proved formulas for $\normord{a(z)b(z)}_t c(z)$ and for
$a(z)_t \normord{b(z)c(z)}$ where $a(z), b(z), c(z)$ are pairwise local
fields and $t\in\Z$
\cite{lian.zuckerman.quantum}, Lemma 3.2.
In particular, they proved the right Wick formula.

Kac proved that fields satisfy the Wick formula and, more generally, 
the commutator formula for indices $t\geq 0, s\in\Z$ 
\cite{kac.beginners.first.ed,kac.beginners}, Proposition 3.3\,\itc. 
Matsuo and Nagatomo proved that fields satisfy the Jacobi identity for 
indices $r\geq 0, s\in\Z, t\geq 0$  
\cite{matsuo.nagatomo.locality}, Theorem 3.2.1 and Corollary 3.2.2.
We present Matsuo and Nagatomo's proof.

Matsuo and Nagatomo showed how these results can be proven using
contour integrals \cite{matsuo.nagatomo.locality}, end of section B.2.

\smallskip

\ref{SS:ope sing}.\:
Wilson and Kadanoff showed the importance of the 
operator product expansion in quantum field theory
\cite{wilson.ope1,kadanoff.ope}.
It is now a fundamental concept treated in every textbook on 
quantum field theory.

Polyakov and Kadanoff put forward the idea to consider the coefficients of 
the OPE as products of an algebra 
\cite{polyakov.bootstrap,kadanoff.ope,kadanoff.ceva.operator.algebra.2d.ising}.

Lian and Zuckerman proved that the $t$-th products for $t\geq 0$ are
the coefficients of the singular part of the OPE 
\cite{lian.zuckerman.classicalq}, Proposition 3.2, 
\cite{lian.zuckerman.quantum}, Proposition 2.3.

Previously, Lian proved that $(a_t b)(z)$ for $t\geq 0$ are
the coefficients of the singular part of the OPE 
where $a, b$ are elements of an associative vertex algebra
\cite{lian.simplevertex}, Lemma 2.5\,\itii.

Roitman observed that the commutator formula for indices $t\geq 0, s\in\Z$
is equivalent to the associativity formula for indices $r\geq 0, s\in\Z$
\cite{roitman.free.conformal.vertex.algebras}, equation (1.3).

\smallskip

\ref{SS:tayl}.\:
The author observed that the $r$-th products of weakly local fields are 
equal to the Taylor coefficients of the operator product
\cite{rosellen.ope.algs}, equation (3.4).

Kac proved Taylor's formula for a field $a(z,w)$ 
\cite{kac.beginners.first.ed}, Lemma 3.1,
\cite{kac.beginners}, Proposition 3.1.

\smallskip

\ref{SS:skew sym fields}.\:
Bais, Bouwknegt, Surridge, and Schoutens showed that any
``operators'' $a(z), b(z)$ satisfy $[a(z),b(z)]\subast=[a(z),b(z)]\subslie$
\cite{bais.bouwknegt.surridge.schoutens.casimir}, equation (A.10).
They proved this by comparing the constant terms of the OPEs of 
$a(z)b(w)$ and $b(w)a(z)$. 
Similarly, comparing the singular part yields conformal skew-symmetry.

Sevrin, Troost, Van Proeyen, and Spindel stated that primary fields 
$a, b, c$ satisfy the identity 
$(a_r b)_s c=\sum (-1)^{r+i+1}(T^{(i)}(b_{r+i}a))_s c$ for $r\geq -1$ and
$s\geq 0$ \cite{sevrin.troost.vanproeyen.spindel.ext.susy.sigma.mod.grp.mfds2},
equation (A.7). 
They used double contour integrals to prove this.

Li proved that a set $S$ of local fields generates an associative 
vertex algebra \cite{li.localsystems}, Corollary 3.2.11.
In particular, the fields in $S$ satisfy skew-symmetry.

Lian and Zuckerman proved that weakly local fields are local iff
they satisfy skew-symmetry
\cite{lian.zuckerman.quantum}, Proposition 2.10.

Li proved that if $M$ is a $\Z$-fold module over an associative vertex algebra
then duality for $M$ implies the Jacobi identity for $M$ and hence, 
in particular, locality for $M$ \cite{li.localsystems}, Proposition 2.3.3.
Our proof of this fact, using that skew-symmetry implies locality, is new.
In Proposition \ref{SS:dual skew loc} we give yet another proof. 

Frenkel and Ben-Zvi proved that $V$-modules are equivalent to 
$\fg(V)$-modules such that $1(z)=\id_M$ and $(ab)(z)=\normord{a(z)b(z)}$
\cite{frenkel.benzvi.book}, Proposition 4.1.5, 
\cite{frenkel.benzvi.book.2nd}, Theorem 5.1.6.

\smallskip

\ref{SS:dong fields}.\:
Li and Meurman and Primc discovered Dong's lemma 
\cite{li.localsystems}, Proposition 3.2.7,
\cite{meurman.primc.annihilating.fields.sl2.combinat.ids}, Proposition 2.5.
Li acknowledged Dong for providing him with a proof of it.
Kac introduced the name ``Dong's lemma''
\cite{kac.beginners.first.ed,kac.beginners}, Lemma 3.2.

We present the author's proof of Dong's lemma that uses 
a reformulation of locality and the fact
that the $r$-th products of local fields are the Taylor coefficients of 
the operator product \cite{rosellen.ope.algs}, Lemma 4.1\,\iti.

Kac claimed that locality is equivalent to two identities of the form
$a(z)b(w)=c(z,w)/(z-w)^n$
\cite{kac.beginners.first.ed,kac.beginners}, Theorem 2.3\,\itvii.
But he missed that one needs to require that $c(z,w)$ is regular at $z=w$ 
in order that associativity holds when multiplying with $(z-w)^n$.

\smallskip

\ref{SS:vas of fields}.\:
Li proved that any maximal local subspace of $\Endv(E)$ is an
associative vertex subalgebra of $\Endv(E)$
\cite{li.localsystems}, Theorem 3.2.10.
Because of Zorn's and Dong's lemma, this is equivalent to the statement that
any local vertex subalgebra of $\Endv(E)$ is associative.

The observation, that any associative vertex subalgebra of $\Endv(E)$
is local, is new.

Li proved that a local subset of $\Endv(E)$ generates an 
associative vertex subalgebra of $\Endv(E)$
\cite{li.localsystems}, Corollary 3.2.11.

\medskip

{\bf Section \ref{S:fz alg}.}\:
\ref{SS:fz alg}.\:
Frenkel and Zhu introduced the Frenkel-Zhu algebra $\cA(V)$
\cite{frenkel.zhu.vertex}, section 1.3.
They called it the universal enveloping algebra.
They observed that there is an equivalence between $\N$-graded $V$-modules
and $\N$-graded $\cA(V)$-modules.

\smallskip

\ref{SS:loc ass alg}.\:
Our definition of the notion of a local associative algebra is new.
The statement that the functor $V\mapsto\cA(V)$ is left adjoint to
the functor $\cA\mapsto V(\cA)$ is new as well.

\smallskip

\ref{SS:embedd assoc}.\:
The remark that the map $V\to\cA(V), a\mapsto a_{-1}$, is injective and
the proposition that the functor $V\to\cA(V)$ is fully faithful are new.

\medskip

{\bf Section \ref{S:suppl va}.}\:
\ref{SS:ope heuris}.\:
Frenkel and Zhu remarked that the formula for the $t$-th products of
$\Endv(E)$ can be viewed as a rigorous formulation of the 
OPE-coefficient $\res_{z=w}(z-w)^n a(z)b(w)$
\cite{frenkel.zhu.vertex}, Remark after Definition 2.2.1.

\smallskip

\ref{SS:ope2}.\:
Frenkel, Lepowsky, and Meurman proved that
$\del_z^{(n)}\de(z)=(1-z)^{-n-1}-(1-z)^{-n-1}_{w>z}$ for any $n\geq 0$ where
$\de(z)=\sum_{i\in\Z}z^i$ \cite{frenkel.lepowsky.meurman.book}, 
Proposition 8.1.2. 

The proof that the $r$-th products of local fields are the Taylor 
coefficients of the operator product is new.

\smallskip

\ref{SS:holom loc implies itself}.\:
Li proved that local fields satisfy locality
\cite{li.localsystems}, Proposition 3.2.9.

\smallskip

\ref{SS:jacobi fields}.\:
Matsuo and Nagatomo proved that local fields satisfy the Jacobi identity 
\cite{matsuo.nagatomo.locality}, Theorem 3.4.1 and Corollary 3.4.2.

\section{Chapter \ref{C:resul}}

${}_{}$
\indent
{\bf Section \ref{S:va ids}.}\:
\ref{SS:loc and dual}.\:
Li called duality ``associativity'' \cite{li.localsystems}, (2.2.9).
Li and Lepowsky called it the ``weak associativity relation''
\cite{lepowsky.li.voa.book}, equation (3.1.16).

\smallskip

\ref{SS:fund recur}.\:
Matsuo and Nagatomo found the fundamental recursion of the Jacobi identity
\cite{matsuo.nagatomo.locality}, equation (4.3.1).

\smallskip

\ref{SS:conseq fund rec}.\:
Roitman observed that the commutator formula for indices $t\geq 0, s\in\Z$
is equivalent to the associativity formula for indices $r\geq 0, s\in\Z$
\cite{roitman.free.conformal.vertex.algebras}, equation (1.3).

Bakalov and Kac proved that the associativity formula is equivalent to
duality and weak locality for $a(z), b(z)$
\cite{bakalov.kac.field.algebras}, Theorem 4.7\,\itb.

\smallskip

\ref{SS:s3 symm}.\:
Frenkel, Huang, and Lepowsky proved that if skew-symmetry is satisfied then
the Jacobi identity holds for elements $a, b, c$ iff it holds for
any permutation of $a, b, c$ 
\cite{frenkel.huang.lepowsky.axiom}, Proposition 2.7.1.
They called this $\bbS_3$-symmetry of the Jacobi identity.

They observed that their argument also shows that if skew-symmetry holds then
the commutator formula and the associativity formula are equivalent
\cite{frenkel.huang.lepowsky.axiom}, Remark 2.7.2.

\smallskip

\ref{SS:2nd recursion}.\:
Bakalov and Kac proved the second recursion for the associativity formula
\cite{bakalov.kac.field.algebras}, proof of Theorem 5.4.

\smallskip

\ref{SS:assoc va}.\:
Bakalov and Kac proved that a vertex algebra is associative iff
it satisfies skew-symmetry and the associativity formula
\cite{bakalov.kac.field.algebras}, Theorem 7.9.

We present a new proof of their result. 
Instead of the associativity formula we prove the commutator formula,
which is equivalent to the associativity formula
because of skew-symmetry and $\bbS_3$-symmetry.
See Corollary \ref{SS:pre lie q assoc} for the original proof.

\smallskip

\ref{SS:q assoc zhu poisson}.\:
Kac introduced quasi-associativity and noted that it is a special
case of the associativity formula
\cite{kac.beginners}, equation (4.11.4).

Matsuo discovered that the subspace $K=\sum_{i\geq 1}(T^i V_0)V_{-i}$
is an ideal of $V_0$ and $V_0/K$ is associative (unpublished).
Miyamoto discovered this independently in the case that $V$ is the
affinization of a vertex algebra $W$, see section \ref{SS:zhu alg2}.
In this case $V_0/K$ is isomorphic to the Zhu algebra of $W$.

\smallskip

\ref{SS:kind of assoc}.\:
Li proved that $a_r b_s c$ is an explicit linear combination of $(a_i b)_j c$
\cite{li.zhu}, Lemma 3.12.
He wrote that this result is a reformulation of a result of
Dong, Li, and Mason \cite{dong.li.mason.higherzhu}.
Li and Lepowsky presented Li's result in their book
\cite{lepowsky.li.voa.book}, Proposition 4.5.7.

Dong and Mason proved that the submodule generated by an element $c$
is the span of $a_n c$ for $a\in V$ and $n\in\Z$
\cite{dong.mason.on.q.galois}, Proposition 4.1.
They acknowledged that Li proved this result previously
\cite{li.haisheng.phd}.

\medskip

{\bf Section \ref{S:unital va ids}.}\:
\ref{SS:loc skew sym}.\:
Li proved that if $a(z)$ is translation covariant then $a(z)$ is creative
\cite{li.invariant.forms.voas}, Proposition 3.3\,\itb.

Li proved that locality implies skew-symmetry 
\cite{li.localsystems}, proof of Proposition 2.2.4.

\smallskip

\ref{SS:unital assoc va}.\:
Li proved that a unital vertex algebra $V$ satisfies duality and
skew-symmetry iff $V$ satisfies locality iff $V$ satisfies the
Jacobi identity \cite{li.localsystems}, Propositions 2.2.4 and 2.2.6.

\smallskip

\ref{SS:conf j wick}.\:
Roitman observed that the commutator formula for indices $t\geq 0, s\in\Z$
is equivalent to the associativity formula for indices $r\geq 0, s\in\Z$
\cite{roitman.free.conformal.vertex.algebras}, equation (1.3).

The observation that a conformal operator $d\submu$ of a vertex algebra 
is a conformal derivation iff $[(d\submu)\subla a(z)]=(d\subla a)(z)$
is new.

\smallskip

\ref{SS:right wick}.\:
Bakalov and Kac introduced the right Wick formula
\cite{bakalov.kac.field.algebras}, equation (5.9).
de Sole gave it the name \cite{desole.thesis}, Remark 1.1.4. 

Bakalov and Kac proved that the right Wick formula
is the associativity formula for indices $r=-1$ and $s\geq 0$
\cite{bakalov.kac.field.algebras}, Lemma 5.2.

\smallskip

\ref{SS:q assoc assoc}.\:
Bakalov and Kac proved that if skew-symmetry holds then 
the left and right Wick formula are equivalent 
\cite{bakalov.kac.field.algebras}, proof of Theorem 7.9.
We give a new proof of their result using $\bbS_3$-symmetry.
See Remark \ref{SS:left right wick} for the original proof. 

Bakalov and Kac proved that a vertex algebra $V$ satisfies 
the associativity formula iff $V$ satisfies conformal Jacobi identity, 
the left and right Wick formula, and quasi-associativity
\cite{bakalov.kac.field.algebras}, Theorem 5.4.

\smallskip

\ref{SS:assoc alg bracket}.\:
The results about associative algebras with a bracket are new.
They are motivated by the work of Bakalov and Kac, see sections
\ref{SS:diff alg leibniz la brack} and \ref{SS:eq commut}.

\smallskip

\ref{SS:diff alg leibniz la brack}.\:
Bakalov and Kac proved that if a unital vertex algebra $V$
contains elements $a, b$ such that $a_t b=1$ for some $t\geq 0$ then the
Wick formulas and the conformal Jacobi identity imply conformal skew-symmetry
\cite{bakalov.kac.field.algebras}, Lemma 8.1 and proof of Theorem 8.4.

\smallskip

\ref{SS:eq commut}.\:
Bakalov and Kac proved that if a unital vertex algebra $V$
contains elements $a, b$ such that $a_t b=1$ for some $t\geq 0$ then the
associativity formula implies that $V$ is associative
\cite{bakalov.kac.field.algebras}, Theorem 8.4.

\medskip

{\bf Section \ref{S:unital va}.}\:
\ref{SS:commut centre}.\:
Dong and Mason proved that the centre $C$ of a vertex operator algebra $V$
is isomorphic to $\End_V(V)$
\cite{dong.mason.local.semilocal.voas}, Proposition 2.5.

\smallskip

\ref{SS:prod of va idemp}.\:
Dong and Mason showed that complementary central idempotents
of a vertex operator algebra $V$ yield a product decomposition of $V$
\cite{dong.mason.local.semilocal.voas}, equation (2.6).

\smallskip

\ref{SS:indecomp va}.\:
Dong and Mason proved the existence and uniqueness of the block decomposition
for a vertex operator algebra 
\cite{dong.mason.local.semilocal.voas}, Theorem 2.7.

The Lemma from noncommutative ring theory 
can be found in the book by Lam \cite{lam.noncomm.rings},
Propositions 22.1 and 22.2.

\smallskip

\ref{SS:tensor va}, \ref{SS:tensor va2}.\:
Borcherds defined the affinization of a vertex algebra
\cite{borcherds.voa}, section 8.
He also mentioned tensor products of vertex algebras and of modules
and noticed that affinization is the tensor product vertex algebra
$V\otimes\K[x\uppm]$.

Frenkel, Huang, Lepowsky constructed the tensor product of vertex algebras
and of modules \cite{frenkel.huang.lepowsky.axiom}, 
section 2.5 and Propositions 3.7.1 and 4.6.1.

Matsuo and Nagatomo proved the universal property of the tensor product 
of vertex algebras \cite{matsuo.nagatomo.locality}, Proposition 4.4.1.

\smallskip

\ref{SS:conf va}.\:
Frenkel and Zhu proved that if $V$ is a vertex operator algebra of CFT-type 
and $V'$ a conformal vertex subalgebra then $C_V(V')=\ker L'_{-1}$ 
\cite{frenkel.zhu.vertex}, Theorem 5.2.
In particular, the centre of $V$ is equal to the kernel of $T$.

Li proved that an element $c$ of an ordinary module $M$ satisfies
$Y_M(V)c\subset M\pau{z}$ iff $Tc=0$ 
\cite{li.invariant.forms.voas}, Proposition 3.3\,\ita.
Li also observed that the centre $C$ of a vertex operator algebra
is contained in $V_0$
\cite{li.invariant.forms.voas}, Corollary 4.2.

Li proved implicitly that if $\K=\C$ and $V$ is a simple 
vertex operator algebra then $C=\C$ 
\cite{li.invariant.forms.voas}, proof of Lemma 4.3 or of Proposition 4.8.

\medskip

{\bf Section \ref{S:filt span sets}.}\:
This section is based on two papers by Li and on results of 
Gaberdiel and Neitzke and of Buhl
\cite{li.vertex.alg.vertex.poisson.alg}, section 4, 
\cite{li.poisson.vas2,gaberdiel.neitzke.rationality.finite.w,
buhl.spanning.set.voa.modules}.
We have included analogous results for almost commutative associative 
differential algebras in order to demonstrate their similarity to
associative vertex algebras.

\smallskip

\ref{SS:filtr alg}.\:
Filtrations are, of course, a basic tool in the study of associative algebras,
in particular, of algebras that are ``almost'' commutative like
enveloping algebras and Weyl algebras. 
But most standard textbooks about ring theory do not discuss filtrations.

\smallskip

\ref{SS:filtr va}.\:
Karel and Li had the idea to use filtrations in the theory of vertex algebras
\cite{karel.li.generating.subspaces}, before Lemma 3.8.
They considered filtrations on the Frenkel-Zhu algebra.
They acknowledged the physicist Watts for the idea
to use filtrations \cite{watts.coset}, equation (11).
Karel and Li's and Watts' filtrations are analogous to invariant filtrations.

Li essentially defined invariant filtrations of vertex algebras, 
except for the fact that
he requires in addition that $(F_n)_k F_m\subset F_{n+m-1}$ for $k\geq 0$
\cite{li.vertex.alg.vertex.poisson.alg}, Definition 4.1.
Thus $\rgr V$ is always commutative.
Li calls 
\index{good filtration}
such filtrations {\it good}.
Arakawa gave the general definition of invariant filtrations a little later
\cite{arakawa.repr.th.w.algs}, Definition A.6.1.

Li introduced differential filtrations of vertex algebras 
\cite{li.poisson.vas2}, equation (2.21).
He does not give them a name.

Li proves for invariant and differential filtrations for which $\rgr V$ is 
commutative that $\rgr V$ is a vertex Poisson algebra 
\cite{li.vertex.alg.vertex.poisson.alg}, Proposition 4.2,
\cite{li.poisson.vas2}, Proposition 2.6.
Arakawa notes that $\rgr V$ is a vertex algebra for 
a general invariant filtration \cite{arakawa.repr.th.w.algs}, Appendix A.6.

\smallskip

\ref{SS:filt gen subset}.\: 
The notion of invariant and differential filtrations $F^S$ generated by 
a subset $S$ is new.
Hence the spanning sets for $F^S$ are new as well.

Li associates to a subset $S$ with a map $S\to\Z_>$ 
a vector space filtration $F'$ and proves that if $S$ is compatible with $F'$
then $F'$ is an invariant filtration such that $\rgr V$ is commutative
\cite{li.vertex.alg.vertex.poisson.alg}, Theorems 4.6.
In this case $F'=F^S$, but in general $F'$ is {\it not} 
an invariant filtration, in particular, $F'\ne F^S$.

\smallskip

\ref{SS:standard inv filt}.\:
Li defined the standard invariant filtration $F\upssi$ of a 
vertex algebra of CFT-type
\cite{li.vertex.alg.vertex.poisson.alg}, after Theorem 4.14.
He shows that $F\upssi$ is the finest good filtration such that 
$V_h\subset F\upssi_h$ \cite{li.vertex.alg.vertex.poisson.alg}, Theorem 4.14.
We define $F\upssi$ for $\rho\Z$-graded vertex algebras
and give a simpler proof that $\rgr V$ is commutative.

\smallskip

\ref{SS:pbw generat}.\:
Karel and Li considered generating subsets with the PBW-property.
They call them generating subsets of weak PBW-type 
\cite{karel.li.generating.subspaces}, Introduction.
Li calls them generating subsets with the PBW spanning property
\cite{li.vertex.alg.vertex.poisson.alg}, before Theorem 4.8.

Li introduced the subspace $C_1(V)$ and proved that if $S+C_1(V)=V$
then $S$ generates $V$ as a vertex algebra 
\cite{li.finiteness.regular}, equation (3.3) and Proposition 3.3.

Karel and Li proved that $S+C_1(V)=V$ iff $S$ {\it strongly generates} $V$,
that is, 
\index{strongly generating}
$V$ is the span of $a^1_{n_1}\dots a^r_{n_r}1$ for $a^i\in S, n_i<0$
\cite{karel.li.generating.subspaces}, Theorem 3.5. 
Moreover, they proved that if $V$ is of CFT-type and
$S$ strongly generates $V$ then $S$ generates $V$ with the PBW-property 
\cite{karel.li.generating.subspaces}, Corollary 3.13.
A result of this type appeared previously in 
the physics literature in a paper by Watts 
\cite{watts.coset}, after equation (7).

Li proved that $S$ strongly generates $V$ iff $S$ generates $V$ and
$S$ satisfies a compatibility with $F^S$
\cite{li.vertex.alg.vertex.poisson.alg}, Theorems 4.6 and Theorem 4.11.
In this case $S$ generates $V$ with the PBW-property and $F^S=F\upssi$ 
\cite{li.vertex.alg.vertex.poisson.alg}, Theorems 4.8 and  4.14.

The proof we give is much simpler.
The main point is that instead of strongly generating 
it suffices to assume that $S$ generates $V\subast$.
In this way we obtain a new proof of the result of Karel and Li.

\smallskip

\ref{SS:no repeat va si filtr}.\:
Gaberdiel and Neitzke introduced the notion of generators without repeats and
proved that complements of $C_2(V)$ generate $V$ without repeats
\cite{gaberdiel.neitzke.rationality.finite.w}, Proposition 8.
They obtain a series of consequences from this result.

Gaberdiel and Neitzke's proof has three main arguments: 
$\rgr V$ is commutative; $\si_{h_a}a=0$ for $a\in C_2(V)$; and
one gets rid of repetitions by using 
a special case of the associativity formula.
But they actually do not use the vertex algebra $\rgr V$
and their inductive argument is quite involved. 
In our proof their inductive argument is replaced by
a lemma about commutative differential algebras which 
is due to Li \cite{li.poisson.vas2}, Lemma 4.1.

\smallskip

\ref{SS:no repeat va mod}.\:
Buhl generalized the result of Gaberdiel and Neitzke to modules
\cite{buhl.spanning.set.voa.modules}, Theorem 1.
His proof is very, very complicated. 
Miyamoto later gave a much simpler proof by extending the three main 
arguments of Gaberdiel and Neitzke more directly to modules
\cite{miyamoto.mod.inv.c2}, Lemma 2.4.
Moreover, Miyamoto's result is stronger than Buhl's 
because Buhl has to allow finitely many repetitions.
However, for applications this difference does not matter.
Yamauchi extended the result to twisted modules 
\cite{yamauchi.modularity.voas.semisimple.primary.vectors}, Lemma 3.3. 

Independently of Buhl, Nagatomo and Tsuchiya proved a generalization
of Gaberdiel and Neitzke's result for certain modules over
the Frenkel-Zhu algebra $\cA(V)$
\cite{nagatomo.tsuchiya.cfts.regular.chiral.voas.1.genus0}, Theorem 3.2.7.

\smallskip

\ref{SS:canonical filt}.\:
Li defined the canonical filtration $F\upsc$ of a vertex algebra, 
and more generally of a module, in terms of spanning sets 
\cite{li.poisson.vas2}, Definition 2.7.
He proved that $F\upsc$ is a differential filtration and 
$\rgr V$ is commutative \cite{li.poisson.vas2}, Proposition 2.14.
The proof we give of these facts is a bit simpler and 
the description of $F\upsc$ as the differential filtration
generated by $V$ is new.

\smallskip

\ref{SS:no repeat va}.\:
Li generalized the result of Gaberdiel and Neitzke to vertex algebras
with elements of negative weight \cite{li.poisson.vas2}, Theorem 4.7.

The Lemma is due to Li \cite{li.poisson.vas2}, Theorem 3.5.
We follow his proof, with the only exception that we do not use his
smaller spanning sets for $F\upsc$ \cite{li.poisson.vas2}, Lemma 2.9.

The Corollary is due to Gaberdiel and Neitzke
\cite{gaberdiel.neitzke.rationality.finite.w}, Theorem 11 and remark
after Corollary 9.

\medskip

{\bf Section \ref{S:suppl resul}.}\:
\ref{SS:pre lie q assoc}.\:
Bakalov and Kac proved that if a vertex algebra $V$ satisfies 
the Wick formula and skew-symmetry then 
the pre-Lie identity is equivalent to quasi-associativity.
\cite{bakalov.kac.field.algebras}, proof of Theorem 7.9.

Bakalov and Kac proved that a vertex algebra is associative iff
it satisfies skew-symmetry and the associativity formula
\cite{bakalov.kac.field.algebras}, Theorem 7.9.
We present their proof.

\smallskip

\ref{SS:dual skew loc}.\:
Li proved that if $M$ is a $\Z$-fold module over an associative vertex algebra
then duality for $M$ implies the Jacobi identity for $M$ and hence, 
in particular, locality for $M$ \cite{li.localsystems}, Proposition 2.3.3.

We present a proof due to Li that duality for $M$ implies locality for $M$
\cite{li.higher.dim.analogues.vertex}, Proposition 3.13.
Li applied his proof, more generally, to modules over $G_n$-vertex algebras.

\smallskip

\ref{SS:left right wick}.\:
Bakalov and Kac proved that if skew-symmetry holds then 
the left and right Wick formula are equivalent 
\cite{bakalov.kac.field.algebras}, proof of Theorem 7.9.
We present their proof.

\smallskip

\ref{SS:lemma tensor prod}.\:
Matsuo and Nagatomo proved the universal property of the tensor product 
of vertex algebras \cite{matsuo.nagatomo.locality}, Proposition 4.4.1.
We give a new proof of the key identity used in the proof.

\smallskip

\ref{SS:bv algs}.\:
The fact that any BV-algebra yields a Gerstenhaber algebra is standard.
See for example the book by Huybrechts \cite{huybrechts.cplx.geom}, 
Proposition 6.A.2.

\smallskip

\ref{SS:topl va}.\:
Lian and Zuckerman proved that if $V$ is a topological vertex algebra
then $H_{Q_0}(V)$ is a BV-algebra 
\cite{lian.zuckerman.conjecture}, Theorem 2.2 and Lemma 2.1.

\smallskip

\ref{SS:1trunc va}.\:
Gorbounov, Malikov, and Schechtman introduced the notion of a
$1$-truncated vertex algebra and noted that there is a forgetful
functor from $\N$-graded vertex algebras to $1$-truncated vertex algebras
\cite{gorbounov.malikov.schechtman.gerbes2}, section 3.1 and equation (3.1.1).

Bressler defined the notion of a Courant algebroid
\cite{bressler.vertex.algebroids1}, section 4.1.

Yamskulna showed that any vertex Poisson algebra yields a Courant algebroid
\cite{yamskulna.v.poisson.courant.algebroids.deform1}, Theorem 3.8.

\smallskip

\ref{SS:exten lie algebroi}.\:
Pradines introduced Lie algebroids in 1967, see
\cite{cannasdasilva.weinstein.geom.models.noncomm.algs}, section 16.2.
This is a vector bundle on a smooth manifold with a Lie bracket and
additional structure.

Rinehart defined Lie algebroids in an algebraic setting, constructed
their enveloping algebras, and proved a Poincar\'e-Birkhoff-Witt
theorem for them \cite{rinehart.diff.forms.on.comm.algs}.
Therefore Lie algebroids are 
\index{Lie!Rinehart algebras}
also called {\it Lie-Rinehart algebras}
\cite{huebschmann.lie.rinehart.gerstenhaber.bv}.

Gorbounov, Malikov, and Schechtman introduced the notion of 
an extended Lie algebroid and showed that any $1$-truncated vertex algebra
yields an extended Lie algebroid
\cite{gorbounov.malikov.schechtman.gerbes2}, sections 3.2 and 2.1.

Previously, Malikov, Schechtman, and Vaintrob gave an outline of
this result \cite{malikov.schechtman.vaintrob.derham}, section 6.
Malikov and Schechtman provided more details
\cite{malikov.schechtman.derham}, section 2.2.

The analogous construction of an extended Lie algebroid from a 
Courant algebroid is new.

\smallskip

\ref{SS:h Jac id res}.\:
Frenkel, Lepowsky, and Meurman wrote the Jacobi identity in terms of
residues \cite{frenkel.lepowsky.meurman.book}, Proposition A.2.8.

\section{Chapter \ref{C:env}}

${}_{}$
\indent
{\bf Section \ref{S:env va}.}\:
\ref{SS:loc fct}, \ref{SS:free va}.\:
Borcherds stated that for any locality function there exists a 
free vertex algebra \cite{borcherds.voa}, section 4.

Roitman constructed free vertex algebras as enveloping algebras of
free vertex Lie algebras 
\cite{roitman.free.conformal.vertex.algebras}, end of section 3.1.

\smallskip

\ref{SS:vertex envelope}.\:
Primc proved that the functor $V\mapsto V\subslie$ has a left adjoint
$R\mapsto U(R)$ and that any bounded $\fg(R)$-module is a $U(R)$-module 
\cite{primc.lie}, Theorems 5.5 and 5.8.

Dong, Li, and Mason proved that any bounded module over a regular local
Lie algebra $\fg$ is a module over an associative vertex algebra contructed 
from $\fg$
\cite{dong.li.mason.vertex.lie.poisson.algebras}, Theorem 4.8.

The construction of $U(R)$ as a quotient of a free vertex algebra is new.

\smallskip

\ref{SS:tensor alg vleib}.\:
Bakalov and Kac proved that $T\upast R$ has a unique vertex algebra
structure \cite{bakalov.kac.field.algebras}, Theorem 6.1.

The observations about the invariant filtration of $T\upast R$ are new.

Loday defined Leibniz algebras
\cite{loday.cyclic.homology,
loday.pirashvili.enveloping.leibniz.algebras.cohomology}.

\smallskip

\ref{SS:univ tensor alg}.\:
Bakalov and Kac proved that $T\upast R$ satisfies a universal property.
\cite{bakalov.kac.field.algebras}, before Remark 6.5.

We prove a slightly stronger universal property that implies that
$R\mapsto T\upast R$ is a functor. 
This last observation is new.

\smallskip

\ref{SS:vertex envelope tensor alg}.\:
Gorbounov, Malikov, and Schechtman proved that if $I$ is the ideal
of the ordinary tensor algebra $T\upast R$ that is 
generated by $a\otimes b-\paraab b\otimes a-[a,b]\subslie$
then $T\upast R/I$ has a unique associative vertex algebra structure 
and $U(R)=T\upast R/I$ 
\cite{gorbounov.malikov.schechtman.gerbes2}, Theorem 10.3, 
\cite{gorbounov.malikov.schechtman.gerbes2.v2}, Theorem 8.3.

Bakalov and Kac stated that $U(R)$ is the quotient of the tensor 
vertex algebra $T\upast R$ by the ideal generated by 
$a\otimes b-\paraab b\otimes a-[a,b]\subslie$
\cite{bakalov.kac.field.algebras}, Remark 7.13.

\medskip

{\bf Section \ref{S:PBW thm}.}\:
\ref{SS:exis thm}.\:
Goddard proved the uniqueness theorem in the framework of 
meromorphic conformal field theory
\cite{goddard.meromorphic}, Theorem 1.

Frenkel, Kac, Radul, and Wang and Meurman and Primc proved 
the existence theorem 
\cite{frenkel.kac.radul.wang.w}, Proposition 3.1.
\cite{meurman.primc.annihilating.fields.sl2.combinat.ids}, Theorem 2.6.

Xu proved a much weaker version of the existence theorem
\cite{xu.spinor.voas.modules}, Theorem 2.4.

\smallskip

\ref{SS:va over loclie}.\:
Kac proved that if $\fg$ is a local Lie algebra then
the Verma module associated to a Lie algebra morphism
$\set{a\in\fg\mid T^n a=0,\; n\gg 0}\to\C$ has a unique
associative vertex algebra structure
Theorem 4.7.

Primc and Dong, Li, and Mason proved this result for any
{\it regular} local Lie algebra $\fg$ and for the Verma module associated to
the decomposition $\fg=\fg_-\oplus\fg_+$
\cite{primc.lie}, Theorems 5.3,
\cite{dong.li.mason.vertex.lie.poisson.algebras}, Theorem 4.8.

The notion of a vertex algebra over a local Lie algebra is new.

\smallskip

\ref{SS:univ va over loclie}.\:
The result that $V(\fg)$ is the universal vertex algebra over $\fg$ is new.

\smallskip

\ref{SS:PBW thm}.\:
Primc proved that $U(R)=V(\fg(R))$ \cite{primc.lie}, Theorems 5.5.

Gorbounov, Malikov, and Schechtman and Bakalov and Kac proved 
the PBW theorem $U(R)=U(R\subsmlie)$
\cite{gorbounov.malikov.schechtman.gerbes2}, Corollary 10.8, 
\cite{gorbounov.malikov.schechtman.gerbes2.v2}, Corollary 8.27,
\cite{bakalov.kac.field.algebras}, Theorem 7.12.

Li proved that $S\upast R$ is a vertex Poisson algebra and
$\rgr U(R)=S\upast R$ as vertex Poisson algebras
\cite{li.vertex.alg.vertex.poisson.alg}, Propositions 3.7 and 4.7.

Frenkel and Ben-Zvi stated the fact that $S\upast R$ is 
a vertex Poisson algebra
\cite{frenkel.benzvi.book}, Remark 15.2.2,
\cite{frenkel.benzvi.book.2nd}, Remark 16.2.2.

\smallskip

\ref{SS:affine va}.\:
Frenkel and Zhu constructed the universal affine vertex algebras $V$
in the case that $\fg$ is a finite-dimensional simple Lie algebra
and proved that $V$ has a conformal vector if $k$ is not equal to the
dual Coxeter number \cite{frenkel.zhu.vertex}, Theorem 2.4.1 and Remark 2.4.1.

Lian generalized the work of Frenkel and Zhu to arbitrary Lie algebras
$\fg$ with an invariant, symmetric bilinear form
\cite{lian.simplevertex}, Theorem 4.8.

Frenkel and Zhu and Lian constructed a vertex algebra structure on a 
certain Verma module $V$, using completed topological algebras and 
correlation functions, but they did not prove any universal property of 
the vertex algebra $V$.

\medskip

{\bf Section \ref{S:suppl env}.}\:
\ref{SS:duality of tensor alg}.\:
Bakalov and Kac proved that $T\upast R$ is dual
\cite{bakalov.kac.field.algebras}, Theorem 6.1.

They also proved the refinement of Dong's lemma that states that
if $a(z), b(z)$ are weakly local and $a(z), c(z)$ and $b(z), c(z)$ are local
then $a(z)_t b(z)$ and $c(z)$ are local for any $t\in\Z$
\cite{bakalov.kac.field.algebras}, Lemma 3.5\,\ita.

They introduced the field $a\upsT(z)$ and proved that 
$a, c$ are dual iff $a(z), c\upsT(z)$ are local
\cite{bakalov.kac.field.algebras}, Proposition 4.1\,\itb.

\section{Chapter \ref{C:reprva}}

${}_{}$
\indent
{\bf Section \ref{S:zhu corres}.}\:
Zhu discovered the relationship between $\N$-graded $V$-modules and 
$A(V)$-modules \cite{zhu.phd,zhu.modular}. 
Essentially all of the results of this section are due to him.
He used correlation functions to prove his main result.

Li found another approach to the Zhu correspondence \cite{li.haisheng.phd}. 
Instead of correlation functions he used the Borcherds Lie algebra $\fg(V)$
and standard constructions of Lie modules. 
We follow his approach. 
Dong, Li, and Mason gave the first published account of Li's work 
\cite{dong.li.mason.twisted.zhus.AV}.
Moreover, they generalized the Zhu correspondence to twisted modules.

\smallskip

\ref{SS:n grad lie mod}, \ref{SS:zhu corres lie}.\:
Zhu introduced the notion of an $\N$-graded module over a 
vertex operator algebra 
\cite{zhu.phd}, Definition 1.2.2, \cite{zhu.modular}, Definition 1.2.3.
Previously, modules were always assumed to be graded by the eigenvalues of 
$L_0$ \cite{frenkel.lepowsky.meurman.book}, equations (8.10.23) and (8.10.32),
\cite{frenkel.huang.lepowsky.axiom}, Definition 4.1.1.

Li introduced the standard constructions of $M'(U)$ and $L(U)$ and the
$\fg_0$-module $\Om(M)$ into vertex algebra theory \cite{li.haisheng.phd}, 
after Remark 2.2.9, Introduction to Chapter 3, and equation (2.2.12).

Li pointed out that the Zhu correspondence is analogous to a
correspondence between simple $\N$-graded $\fg$-modules and simple
$\fg_0$-modules \cite{li.haisheng.phd}, Introduction to Chapter 3.
He mentioned also that such a correspondence exists for 
$\Z$-graded associative algebras $A$.
Instead of $A_0$ one has to take the quotient algebra 
$A_0/\sum_{h\leq -1}A_{-h}A_h$.

The results of sections \ref{SS:n grad lie mod} and \ref{SS:zhu corres lie}
are essentially well-known and elementary.
Dong, Li, and Mason observed that $\Om(L(U))=U$ for any vertex operator
algebra $V$ and any $A(V)$-module $U$ 
\cite{dong.li.mason.twisted.zhus.AV}, Theorem 6.3.
Li proved this result under the assumption that $L_0$ acts semisimply on $U$
\cite{li.haisheng.phd}, Proposition 3.2.9.

\smallskip

\ref{SS:zhu prod}.\:
Zhu defined the Zhu product and the subspace $O(V)=V\ast_{-2}V$ and 
proved that if $M$ is an $\N$-graded $V$-module then $M_0$ is an
$A(V)$-module \cite{zhu.phd,zhu.modular}, Definition 2.1.2 and Theorem 2.1.2.

Li observed that $\fg(V)_0\cong V/(T+H)V$ \cite{li.haisheng.phd},
Proposition 2.2.7.

\smallskip

\ref{SS:zhu alg2}.\:
Zhu proved that $O(V)$ is an ideal and $A(V)=V/O(V)$ is 
an associative algebra \cite{zhu.phd,zhu.modular}, Theorem 2.1.1 (1).

Miyamoto gave a proof of Zhu's result using the affinization $\hV$ of $V$ 
(Proceedings of 17th Algebraic Combinatoric Symposium 
at Tsukuba, 19--21 June 2000), Lemmas 3 and 4. 
Yamauchi explained it to me and Matsuo found it independently.

Akman suggested that the Zhu product might be a pre-Lie product
\cite{akman.htpgerstenhaber}, at the end.
Since $(V,\ast)\cong\hV_0$, this is in fact true.

\smallskip

\ref{SS:zhu alg fgv0}.\:
Zhu gave a formula for the commutator of the Zhu algebra 
\cite{zhu.phd}, Lemma 2.1.7, \cite{zhu.modular}, Lemma 2.1.3.
He deduced from it that a conformal vector is contained in the centre of 
$A(V)$ \cite{zhu.phd,zhu.modular}, Theorem 2.1.1 (3).

Li observed that Zhu's formula is equivalent to the
fact that the surjection $\fg(V)_0\to A(V), a_0\mapsto [a]$,
is a Lie algebra morphism \cite{li.haisheng.phd}, Lemma 3.2.1.
We give a new proof of this fact using $\hV$.

\smallskip

\ref{SS:zhu corresp}.\:
Zhu proved that for any $A(V)$-module $U$ there exists 
an $\N$-graded $V$-module $M$ such that $M_0=U$ and $J(M)=0$ 
\cite{zhu.phd,zhu.modular}, Theorem 2.2.1.
Li proved that if $U$ is an $A(V)$-module then the $\fg(V)$-module $L(U)$ is 
a $V$-module \cite{li.haisheng.phd}, Theorem 3.2.6.

Zhu deduces from his result that there is a bijection between simple
$\N$-graded $V$-modules and simple $A(V)$-modules 
\cite{zhu.phd,zhu.modular}, Theorem 2.2.2.

\smallskip

\ref{SS:ratl vas}.\:
Zhu defined the notion of a rational vertex algebra
\cite{zhu.phd}, Definition 1.2.3, \cite{zhu.modular}, Definition 1.2.4.
He proved that if $V$ is rational then $A(V)$ is semisimple
\cite{zhu.phd,zhu.modular}, Theorem 2.2.3.

Zhu required in his definition of a rational vertex algebra two more 
conditions.
He required in addition that $V$ has only finitely many simple $\N$-graded 
modules up to isomorphism and that $\dim M_h<\infty$ for any such module $M$.

Dong, Li, and Mason showed that these two conditions are redundant
\cite{dong.li.mason.twisted.zhus.AV}, Theorem 8.1 (b), (c).
The subspaces $M_h$ are eigenspaces of $L_0$ with eigenvalue $h_0+h$
for some fixed $h_0$.
Dong, Li, and Mason proved that if $V$ is an even rational vertex operator 
algebra and $\K=\C$ then $A(V)$ is finite-dimensional
\cite{dong.li.mason.twisted.zhus.AV}, Theorem 8.1 (a).
They proved all these results more generally for twisted modules.

\medskip

{\bf Section \ref{S:suppl reprva}.}\:
\ref{SS:gen verma vertex}.\:
Li gave a construction of the $V$-module $Q(M'(U))$ and 
remarked that the functor $U\mapsto Q(M'(U))$ from $A(V)$-modules
to $V$-modules is left adjoint to the functor $M\mapsto\Om(U)$
\cite{li.haisheng.phd}, Remark 3.2.7 and Corollary 3.2.8.

Dong, Li, and Mason extended this result to twisted modules
\cite{dong.li.mason.twisted.zhus.AV}, Theorem 6.2.

\smallskip

\ref{SS:assoc zhu prod}, \ref{SS:zhu alg}.\:
Matsuo and Nagatomo defined the $n$-th Zhu product $\ast_n$ for any $n\in\Z$,
proved the relation between $\ast_n$ and $\ast_{n-1}$ and 
the associativity formula for the Zhu products, 
and used these results to give a new proof that $O(V)$ is an ideal and 
$A(V)$ is associative \cite{matsuo.nagatomo.locality},
Lemmas 8.4.1 and 8.4.2 and Proposition 8.4.3.

\smallskip

\ref{SS:commut zhu revisi}.\:
We present Zhu's computation of the commutator of the Zhu algebra 
\cite{zhu.phd}, Lemma 2.1.7, \cite{zhu.modular}, Lemma 2.1.3.

\smallskip

\ref{SS:zhu affine lie}.\:
Frenkel and Zhu calculated the Zhu algebra of the universal affine
vertex algebra \cite{frenkel.zhu.vertex}, Theorem 3.1.1.
We follow their proof and 
use the results of Matsuo and Nagatomo from sections
\ref{SS:assoc zhu prod} and \ref{SS:zhu alg} to show that $N=O(V)$.

\backmatter

\bibliographystyle{alpha}
\bibliography{lit/bib,lit/collections}


\include{referenc}
\printindex


\end{document}